%% file: ANNcalculus022.tex
\documentclass[a4paper,12pt]{article}

\usepackage{mathrsfs,amsthm,graphicx,color,verbatim,bbm,amsmath,amsfonts,amssymb,newclude,nicefrac,amsfonts,graphicx,color,enumerate,hyperref,bm,geometry,dsfont,mathabx}

\usepackage[latin1]{inputenc}

\newcommand{\R}{\mathbb{R}}
\newcommand{\N}{\mathbb{N}}

\newtheorem{theorem}{Theorem}[section]
\newtheorem{lemma}[theorem]{Lemma}

\newtheorem{cor}[theorem]{Corollary}
\newtheorem{definition}[theorem]{Definition}

\newtheorem{prop}[theorem]{Proposition}

\newcommand{\smallsum}{\textstyle\sum}

\newcommand{\ANNs}{\mathbf{N}}
\newcommand{\activation}{a}
\newcommand{\activationDim}[1]{\mathfrak{M}_{\activation,#1}}
\newcommand{\functionANN}{\mathcal{R}_{\activation}}
\newcommand{\paramANN}{\mathcal{P}}

\newcommand{\lengthANN}{\mathcal{L}}
\newcommand{\inDimANN}{\mathcal{I}}
\newcommand{\compANN}[2]{{#1 \bullet #2}}

\newcommand{\concPsiANN}[2]{{#1 \odot_{\Psi} #2}}

\newcommand{\outDimANN}{\mathcal{O}}
\newcommand{\longerANN}[1]{\mathcal{E}_{#1}}

\newcommand{\idANNshort}[1]{\mathbb{I}_{#1}}
\newcommand{\dims}{\mathcal{D}}
\newcommand{\hiddenLength}{\mathcal{H}}
\newcommand{\hiddenDimId}{\mathfrak{i}}
\newcommand{\parallelization}{\operatorname{P}}
\newcommand{\parallelizationSpecial}{\mathbf{P}}

\newcommand{\affineMap}{\mathbb{A}}
\newcommand{\idMatrix}{\operatorname{I}}
\newcommand{\pa}[1]{\left({#1}\right)}

\newcommand{\indicator}[1]{\mathds{1}_{#1}}

\newcommand{\affineProcess}{Y}

\newcommand{\norm}[1]{ \left\| #1 \right\| }
\newcommand{\NOrm}[1]{\thickvert\!\thickvert #1\thickvert\!\thickvert}
\newcommand{\NORm}[1]{\left\thickvert\!\left\thickvert #1\right\thickvert\!\right\thickvert}
\newcommand{\NOrmm}[1]{\big\thickvert\!\big\thickvert #1\big\thickvert\!\big\thickvert}
\newcommand{\NOrmmm}[1]{\Big\thickvert\!\Big\thickvert #1\Big\thickvert\!\Big\thickvert}
\newcommand{\NOrmmmm}[1]{\bigg\thickvert\!\bigg\thickvert #1\bigg\thickvert\!\bigg\thickvert}

\newcommand{\qandqShort}{\quad\text{and}\quad}

\newcommand{\qandq}{\qquad\text{and}\qquad}
\newcommand{\andq}{\text{and}\qquad}
\newcommand{\andShort}{\text{ and }}
\newcommand{\forallDist}{\forall\,}

\newcommand{\brac}[1]{\left[{#1}\right]}
\newcommand{\bracbig}[1]{\big[{#1}\big]}
\newcommand{\bracLclosedRopen}[1]{\left[{#1}\right)}
\newcommand{\zeroto}[1]{\{0,1,\dots,{#1}\}}
\newcommand{\eps}{\varepsilon}

\newcommand{\abs}[1]{\left|{#1}\right|}

\newcommand{\relu}{a}

\renewcommand{\L}{\lengthANN}

\newcommand{\LogBin}{\log_2}

\newcommand{\vertiii}[1]{{\left\vert\kern-0.25ex\left\vert\kern-0.25ex\left\vert {#1} 
		\right\vert\kern-0.25ex\right\vert\kern-0.25ex\right\vert}}
\newcommand{\vertiiibig}[1]{{\big\vert\kern-0.25ex\big\vert\kern-0.25ex\big\vert {#1} 
		\big\vert\kern-0.25ex\big\vert\kern-0.25ex\big\vert}}
\newcommand{\vertiiiStandard}[1]{{\vert\kern-0.25ex\vert\kern-0.25ex\vert {#1} 
		\vert\kern-0.25ex\vert\kern-0.25ex\vert}}
\newcommand{\tripleNorm}[1]{\vertiii#1}

\begin{document}
	\title{Space-time error estimates for deep neural network approximations for differential equations}

\author{
	Philipp Grohs$^1$,
	Fabian Hornung$^{2,3}$,\\
	Arnulf Jentzen$^4$,
	and Philipp Zimmermann$^{5,6}$ 
	\bigskip
	\\
	\small{$^1$Faculty of Mathematics and Research Platform Data Science, University of Vienna}\\
	\small{Austria, e-mail:  philipp.grohs@univie.ac.at}
	\\
	\small{$^2$ Department of Mathematics, ETH Zurich, Z\"urich,}\\
	\small{Switzerland, e-mail: fabianhornung89@gmail.com}
	\\
	\small{$^3$ Faculty of Mathematics, Karlsruhe Institute of Technology, Karlsruhe,}\\
	\small{Germany, e-mail: fabianhornung89@gmail.com} 
	\\
	\small{$^4$ Department of Mathematics, ETH Zurich, Z\"urich,}\\
	\small{Switzerland, e-mail: arnulf.jentzen@sam.math.ethz.ch} 
	\\
	\small{$^5$ Department of Mathematics, ETH Zurich, Z\"urich,}\\
	\small{Switzerland, e-mail: philipp.zimmermann@math.ethz.ch} 
	\\
	\small{$^6$ Faculty of Mathematics, University of Regensburg, Regensburg,}\\
	\small{Germany, e-mail: philipp.zimmermann@ur.de}
}

\maketitle

\begin{abstract}
	Over the last few years deep artificial neural networks (DNNs) have very successfully been used in numerical simulations for a wide variety of computational problems including computer vision, image classification, speech recognition, natural language processing, as well as computational advertisement. In addition, it has recently been proposed to approximate solutions of partial differential equations (PDEs) by means of stochastic learning problems involving DNNs. There are now also a few rigorous mathematical results in the scientific literature which provide error estimates for such deep learning based approximation methods for PDEs. All of these articles provide spatial error estimates for neural network approximations for PDEs but do not provide error estimates for the entire space-time error for the considered neural network approximations. It is the subject of the main result of this article to provide space-time error estimates for DNN approximations of Euler approximations of certain perturbed differential equations. Our proof of this result is based (i) on a certain artificial neural network (ANN) calculus and (ii) on ANN approximation results for products of the form $[0,T]\times \R^d\ni (t,x)\mapsto tx\in \R^d$ where $T\in (0,\infty)$, $d\in \N$, which we both develop within this article.
\end{abstract}

\tableofcontents

\section{Introduction}
\label{sec:intro}
\input{Intro}

\section{Artificial neural network (ANN) calculus}
\label{sec:Artificial neural network (ANN) calculus}
\input{ANNIntroSec2}
\input{ANNcompositions} % subsection 2.1 (Artificial neural networks (ANNs) and their realizations), 2.2 (Compositions of ANNs) (2.2.1-2.2.2)
\input{ANNcompositionsSpecial} % subsection 2.2 (2.2.3-2.2.5)
\input{ANNparallelizations} % subsection 2.3 (Parallelizations of ANNs)

\subsection{Sums of ANNs} % New subsection 2.4
\subsubsection{Sums of ANNs with the same length} % New subsection 2.4.1 (2.25)
\input{SumsSameLengths}
\subsubsection{Sums of ANNs with different lengths} % New subsection 2.4.2 (2.26)
\input{SumsDifferentLengths}

%\subsection{ANN representation results} % subsection 2.4
%\subsubsection{ANN representation of sums} % subsubsection 2.4.1
%\input{ANNforSums} 
%\subsubsection{ANN representation of abstract Euler type schemes} % subsubsection 2.4.2
%\input{ANNsForEulerSpace}

\subsection{ANN representations for Euler approximations} % New subsection 2.5
\label{subsec:ANNrepforEulerApprox}
\subsubsection{ANN representations for one Euler step} % New subsection 2.5.1
\input{ANNrepOneEuler}

\subsubsection{ANN representations for multiple nested Euler steps} % New subsection 2.5.2
\label{subsec:ANN representations for multiple nested Euler steps}
\input{ANNrepMultNestedEuler}
\subsubsection{ANN representations for multiple perturbed nested Euler steps} % New subsection 2.5.3
\input{ANNrepMultiplePertEuler}

% \section{ANN approximation results} % section 3
% \input{ProduktLemma} % subsection 3.1 (ANN approximation of the square function)
% \input{ANNforProducts} % subsection 3.2 (ANN approximation of products)
% \subsection{ANN approximation of the Euler scheme} 
% \input{DNNeuler} % % subsection 3.3 (ANN approximation of the Euler scheme)
%\section[DNN approximations for partial differential equations (PDEs)]{DNN approximations for partial differential equations (PDEs)}
%\label{sec:DNN_PDEs}
%%\subsection{DNN Approximation of PDEs}
%\input{DNNapproximation}

\section{ANN approximation results} % New section 3
\input{ANNIntroSec3}
\subsection{ANN approximations for the square function} % New subsection 3.1
\label{subsec:ANNApproxforSquare}
\subsubsection{Explicit approximations for the square function on $[0,1]$} % New subsection 3.1.1
\label{subsubsec:Explicit approximations for the square function on $[0,1]$}
\input{ExplicitApproxSquare}
\subsubsection{ANN approximations for the square function on $[0,1]$} % New subsection 3.1.2
\label{subsubsec:ANN approximations for the square function on $[0,1]$}
\input{ANNApproxSquare}
\subsubsection{ANN approximations for the square function on $\R$} % New subsection 3.1.3
\input{ANNApproxSquareOnR}
\subsection{ANN approximations for products} % New subsection 3.2
\label{subsec:ANNApproxProducts}
\subsubsection{ANN approximations for one-dimensional products} % New subsection 3.2.1
\label{subsubsec:ANN approximations for one-dimensional products}
\input{ANNApproxOneDimProducts}

\subsubsection{ANN approximations for multi-dimensional products} % New subsection 3.2.2
\input{ANNApproxMultiDimProducts}

\subsection{Space-time ANN approximations for Euler approximations} % New subsection 3.3
\label{subsec:ANNapproxEuler}
\subsubsection{Space-time representations for Euler approximations} % New subsection 3.3.1
\input{SpaceTimeRepEuler}
\subsubsection{ANN representations for hat functions} % New subsection 3.3.2

\input{ANNrepHatFct}
\subsubsection[A posteriori error estimates for ANN approximations]{A posteriori error estimates for space-time ANN approximations} % New subsection 3.3.3
\input{PosterioriEstimates}
\subsubsection{A priori estimates for Euler approximations} % New subsection 3.3.4
\input{PrioriEstimatesEuler}
\subsubsection[A priori error estimates for ANN approximations]{A priori error estimates for space-time ANN approximations} % New subsection 3.3.5
\label{subsubsec:A priori error estimates for space-time ANN approximations}
\input{PrioriEstimatesANN}

\bibliographystyle{acm}
\bibliography{bibfileDNN}

\end{document}

%% file: Intro.tex
% Introduction :Space-time error estimates for deep neural network approximations for differential equations

Over the last few years deep artificial neural networks (DNNs) have very successfully been used in numerical simulations for a wide variety of computational problems including computer vision, image classification, speech recognition, natural language processing, as well as computational advertisement (cf., e.g., the references mentioned in \cite{goodfellow2016deep,DNNApproxTheory,Salimova2018}). In addition, the articles \cite{weinan2017deep,Han8505} suggest to approximate solutions of partial differential equations (PDEs) by means of stochastic learning problems involving DNNs. We also refer to \cite{DeepSplitting,Kolmogorov,beck2017machine,becker2018deep,TWelti,berg2018unified,MLsemiPDE,DeepRitzMethod,FarahmandNabiNikovski17,FujiiTakahashiTakahashi17,goudenege2019machine,han2018convergence,HenryLabordere17,hure2019some,jacquier2019deep,LongLuMaDong17,lye2019deep,magill2018neural,FullyNonlinPDEs,Raissi18,sirignano2017dgm} for extensions and improvements of such deep learning based approximation methods for PDEs. There are now also a few rigorous mathematical results in the scientific literature which provide error estimates for such deep learning based approximation methods for PDEs; see, e.g., \cite{berner2018analysis,ElbraechterSchwab2018,GrohsWurstemberger2018,han2018convergence,hutzenthaler2019proof,Salimova2018,kutyniok2019theoretical,RectifiedDeepNetworks,sirignano2017dgm}. The articles in this reference list all provide spatial error estimates for neural network approximations for PDEs but do not provide error estimates for the entire space-time error for the considered neural network approximations. It is the subject of Theorem~\ref{Thm:ApproxOfEulerWithGronwall} in this article, which is the main result of this article, to provide space-time error estimates for DNN approximations of Euler approximations of certain perturbed differential equations. To illustrate the findings of the main result of this article in more details, we now formulate in Theorem~\ref{Thm:ApproxOfEulerWithGronwall:SpecialCase} below a special case of Theorem~\ref{Thm:ApproxOfEulerWithGronwall}.
\begin{theorem}\label{Thm:ApproxOfEulerWithGronwall:SpecialCase}
	Let  $\mathfrak{C}, T,\mathfrak{d} \in (0,\infty)$, let $A_d\in C(\R^d,\R^d)$, $d\in\N$, satisfy for all 
	$d\in \N$, $x=(x_1,x_2,\dots,x_d)\in\R^d$ that $A_d(x)=(\max\{x_1,0\},\max\{x_2,0\},\dots,\max\{x_d,0\})$, let
	$\ANNs
	=
	\cup_{L \in \N}
	\cup_{ (l_0,l_1,\ldots, l_L) \in \N^{L+1} }
	\left(
	\times_{k = 1}^L (\R^{l_k \times l_{k-1}} \times \R^{l_k})
	\right)$,
	let $R\colon\ANNs\to \cup_{k,l\in\N}\,C(\R^k,\R^l)$ and $P\colon\ANNs\to \N$ satisfy for all $ L\in\N$, $l_0,l_1,\ldots, l_L \in \N$, 
	$
	\Phi 
	\in  \allowbreak
	( \times_{k = 1}^L\allowbreak(\R^{l_k \times l_{k-1}} \times \R^{l_k}))$, $\Psi =	\allowbreak((W_1, B_1),\allowbreak(W_2, B_2),\allowbreak \ldots, (W_L,\allowbreak B_L))	\in  \allowbreak
	( \times_{k = 1}^L\allowbreak(\R^{l_k \times l_{k-1}} \times \R^{l_k}))$, $x_0 \in \R^{l_0}, x_1 \in \R^{l_1}, \ldots, x_{L-1} \in \R^{l_{L-1}}$ 
	with $\forall \, k \in \N \cap (0,L) \colon x_k =A_{l_k}(W_k x_{k-1} + B_k)$
	that
	$P(\Phi)
	=
	\sum_{k = 1}^L l_k(l_{k-1} + 1)$,
	$R(\Psi) \in C(\R^{l_0},\R^{l_L}),\text{ and }
	( R(\Psi) ) (x_0) = W_L x_{L-1} + B_L$, let $\Phi_d\in \ANNs$, $d\in\N$, satisfy for all $d\in\N$, $x\in \R^d$ that $R(\Phi_d)\in C(\R^d,\R^d)$, $\Vert (R(\Phi_d))(x)\Vert\le \mathfrak{C}(1+\Vert x\Vert)$, and 
	$P(\Phi_d)\le \mathfrak{C}d^{\mathfrak{d}}$,
	let $Y^{d,N}= (Y^{d,N}_{t,x})_{(t,x)\in [0,T]\times \R^d} \colon\allowbreak [0,T]\times \R^d \to \R^d $, $N, d\in \N$, 
	be the functions which satisfy for all $d, N\in\N$, $n\in\{0,1,\dots,N-1\}$, $ t \in \big[\frac{nT}{N},\frac{(n+1)T}{N}\big]$, $ x \in \R^d $ that $\affineProcess^{d,N }_{0,x}=x$ and
	\begin{equation}
	\label{ApproxOfEulerWithGronwall:Y_processes}
	\begin{split}
	&\affineProcess^{d,N}_{t,x} 
	=
	\affineProcess^{d,N }_{\frac{nT}{N},x}+ \left(t-\tfrac{nT}{N}\right)(R(\Phi_d)) ( 
	\affineProcess^{d,N }_{\frac{nT}{N},x} )\,.
	\end{split}
	\end{equation}
	Then there exist $C\in\R$ and $\Psi_{\varepsilon,d,N}\in \ANNs$, $N, d\in \N$, $\varepsilon\in (0,1]$, such that 
	\begin{enumerate}[(i)]
		\item \label{cor1} it holds for all  $\varepsilon\in (0,1]$, $d,N\in\N$ that $R (\Psi_{\varepsilon,d,N})\in C(\R^{d+1},\R^d)$,
		\item  it holds for all
		$\varepsilon\in (0,1]$, $d,N\in\N$,
		$ t \in [0,T]$,
		$x\in\R^d$
		that 
		\begin{equation}
			\Vert \affineProcess^{d,N}_{t,x} -(R (\Psi_{\varepsilon,d,N}))(t,x)\Vert\le Cd^{\nicefrac{1}{2}} N^{\nicefrac{3}{2}}\varepsilon(1+\Vert x\Vert^3),
		\end{equation}
		\item it holds for all 
		$\varepsilon\in (0,1]$, $d,N\in\N$,
		$ t \in [0,T]$,
		$x\in\R^d$
		that 
		\begin{equation}
			\Vert (R (\Psi_{\varepsilon,d,N}))(t,x)\Vert \le Cd^{\nicefrac{1}{2}}N(1+\Vert x\Vert^2),
		\end{equation}
		and
		\item it holds for all  $\varepsilon\in (0,1]$, $d,N\in\N$ that 
			\begin{equation}
				P(\Psi_{\varepsilon,d,N})\le Cd^{16+8\mathfrak{d}}N^6\big[1+|\!\ln(\varepsilon)|^2\big].
			\end{equation}
	\end{enumerate}
\end{theorem}
Theorem~\ref{Thm:ApproxOfEulerWithGronwall:SpecialCase} is an immediate consequence of Corollary~\ref{Cor:ApproxOfEulerWithGronwall} in Subsection~\ref{subsubsec:A priori error estimates for space-time ANN approximations} below. Corollary~\ref{Cor:ApproxOfEulerWithGronwall}, in turn, follows
from Theorem~\ref{Thm:ApproxOfEulerWithGronwall} in Subsection~\ref{subsubsec:A priori error estimates for space-time ANN approximations}, which is the main result of this article. Our proof of Theorem~\ref{Thm:ApproxOfEulerWithGronwall:SpecialCase} and Theorem~\ref{Thm:ApproxOfEulerWithGronwall}, respectively, is based on a certain artificial neural network (ANN) calculus, which we develop in Section~\ref{sec:Artificial neural network (ANN) calculus}. Section~\ref{sec:Artificial neural network (ANN) calculus} is in parts based on several well-known concepts and results in the scientific literature (cf., e.g., \cite{ElbraechterSchwab2018,Salimova2018,petersen2017optimal,yarotsky2017error}). We refer to the beginning of Section~\ref{sec:Artificial neural network (ANN) calculus} for a more detailed comparison of the content of Section~\ref{sec:Artificial neural network (ANN) calculus} with the material in related articles in the scientific literature. Our proof of Theorem~\ref{Thm:ApproxOfEulerWithGronwall:SpecialCase} and Theorem~\ref{Thm:ApproxOfEulerWithGronwall}, respectively, is mainly inspired by \cite{Salimova2018}, \cite[Section~6]{ElbraechterSchwab2018}, and \cite[Section~3.2]{yarotsky2017error}.
\label{Activation function} Theorem~\ref{Thm:ApproxOfEulerWithGronwall:SpecialCase} and Theorem~\ref{Thm:ApproxOfEulerWithGronwall}, respectively, provide error estimates for rectified DNN approximations of Euler approximations of certain perturbed differential equations. Many of the DNN approximation and representation results of this work, however, apply to DNNs with more general activation functions than only the rectifier function (cf., e.g., Li et al.~\cite[Section~1]{RectifiedPowerUnits} and Petersen et al.~\cite[Section~2]{ParametricReLU} for further activation functions).
\label{Conclusion} The error estimates for rectified DNN approximations of Euler approximations of perturbed differential equations, which we establish in Theorem~\ref{Thm:ApproxOfEulerWithGronwall:SpecialCase} and Theorem~\ref{Thm:ApproxOfEulerWithGronwall}, respectively, can then be used to establish space-time error estimates for DNN approximations for PDEs. This will be the subject of a future research article, which will be based on this article. 

The remainder of this article is organized as follows. In Section~\ref{sec:Artificial neural network (ANN) calculus} we develop the above mentioned ANN calculus and, in particular, we establish in Subsection~\ref{subsec:ANNrepforEulerApprox} ANN representation results for Euler approximations. In Subsection~\ref{subsec:ANNApproxforSquare} we develop ANN approximation results for the square function $\R\ni x\mapsto x^2\in\R$. These ANN approximation results for the square function are then used in Subsection~\ref{subsec:ANNApproxProducts} to develop ANN approximation results for products of the form $[0,T]\times \R^d\ni (t,x)\mapsto tx\in \R^d$ where $T\in (0,\infty)$, $d\in \N$. In Subsection~\ref{subsec:ANNapproxEuler} we then combine the ANN representation results in Subsection~\ref{subsec:ANNrepforEulerApprox} with the ANN approximation results for products in Subsection~\ref{subsec:ANNApproxProducts} to establish in Theorem~\ref{Thm:ApproxOfEulerWithGronwall} the main result of this article.

%% file: ANNIntroSec2.tex
This section develops a certain calculus for ANNs. Some of the notions and results which we present here are rather elementary, but for convenience of the reader we present here all details and we include the proof of every result. The material in this section is also in parts based on several well-known concepts and results in the scientific literature. In particular, Definition~\ref{Def:ANN}, Definition~\ref{Def:multidim_version}, and Definition~\ref{Definition:ANNrealization} are slight reformulations of Petersen \& Voigtlaender \cite[Definition~2.1]{petersen2017optimal}. Moreover, Lemma~\ref{Lemma:elementaryPropertiesANN} is elementary and well-known in the scientific literature. Furthermore, Definition~\ref{Definition:ANNcomposition} is also a slight reformulation of Petersen \& Voigtlaender \cite[Definition~2.2]{petersen2017optimal}. In addition, Proposition~\ref{Lemma:PropertiesOfCompositions}, Corollary~\ref{Corollary:Composition}, and Lemma~\ref{Lemma:CompositionAssociative} are elementary and essentially well-known in the scientific literature (cf., e.g., Petersen \& Voigtlaender \cite{petersen2017optimal}). Moreover, Definition~\ref{Definition:iteratedANNcomposition} is an extension of Elbr{\"a}chter et al.~\cite[Setting~5.2]{ElbraechterSchwab2018} and Proposition~\ref{Lemma:PropertiesOfConcatenations} is in parts an extension of Elbr{\"a}chter et al.~\cite[Lemma~5.3]{ElbraechterSchwab2018}. Furthermore, Definition~\ref{Definition:simpleParallelization} and Definition~\ref{Definition:generalParallelization} extend Elbr{\"a}chter et al.~\cite[Setting~5.2]{ElbraechterSchwab2018} (cf., e.g., Petersen \& Voigtlaender \cite[Definition~2.7]{petersen2017optimal}). In addition, Proposition~\ref{Lemma:SumsOfANNSequalArchitecture} is a reformulation of \cite[Lemma~5.1]{Salimova2018}. Moreover, Lemma~\ref{Lemma:CompositionSum} and Proposition~\ref{Cor:CompositionSum} are significantly inspired by \cite[Proposition~5.3]{Salimova2018}. Furthermore, item~\eqref{CompositionSum:Params} in Lemma~\ref{Lemma:CompositionSum} and item~\eqref{CorCompositionSum:Params} in Proposition~\ref{Cor:CompositionSum}, respectively, improve the parameter estimates in \cite[Proposition~5.3]{Salimova2018}. In addition, Corollary~\ref{Cor:CompositionSumInductionTilde} in Subsection~\ref{subsec:ANN representations for multiple nested Euler steps} below is also in parts inspired by \cite[Proposition~6.1]{Salimova2018}.

%% file: ANNcompositions.tex
%\subsection[Representation properties for ANNs]{Representation properties for artificial neural networks}

\subsection{Artificial neural networks (ANNs) and their realizations}

\begin{definition}[Artificial neural networks (ANNs)]
	\label{Def:ANN}
	We denote by $\ANNs$ the set given by 
		\begin{equation}
	\begin{split}
	\ANNs
	&=
	\cup_{L \in \N}
	\cup_{ (l_0,l_1,\ldots, l_L) \in \N^{L+1} }
	\left(
	\times_{k = 1}^L (\R^{l_k \times l_{k-1}} \times \R^{l_k})
	\right)
	\end{split}
	\end{equation}
%	and we call $\ANNs$ the set of all fully connected artificial neural networks,
	and we denote by 	$
	\paramANN, 
	%	\paramNotZeroANN, 
	\lengthANN,  \inDimANN, \outDimANN \colon \ANNs \to \N
	$, $\hiddenLength \colon \ANNs \to \N_0$, and
	$\dims\colon\ANNs\to  \cup_{L=2}^\infty\, \N^{L}$
	the functions which satisfy
	for all $ L\in\N$, $l_0,l_1,\ldots, l_L \in \N$, 
	$
	\Phi 
	\in  \allowbreak
	( \times_{k = 1}^L\allowbreak(\R^{l_k \times l_{k-1}} \times \R^{l_k}))$
	that
	$\paramANN(\Phi)
	=
	\sum_{k = 1}^L l_k(l_{k-1} + 1) 
	$, $\lengthANN(\Phi)=L$,  $\inDimANN(\Phi)=l_0$,  $\outDimANN(\Phi)=l_L$, $\hiddenLength(\Phi)=L-1$, and $\dims(\Phi)= (l_0,l_1,\ldots, l_L)$.
%	we denote by $\activationDim{n}\colon \R^n\to\R^n$, $\activation\in C(\R,\R)$, $n\in\N$, the functions
%	which satisfy 
%	for all 
%	$\activation\in C(\R,\R)$,
%	$
%	n \in \N
%	$,
%	$ x = ( x_1, \dots, x_n ) \in \R^n $
%	that
%	$ 
%	\activationDim{n}(x)
%	=
%	( \activation(x_1), \ldots, \activation(x_n) )
%	$,	 
\end{definition}

\begin{definition}[Multidimensional versions]
	\label{Def:multidim_version}
	Let $d \in \N$ and let $\psi \colon \R \to \R$ be a function.
	Then we denote by $\mathfrak{M}_{\psi, d} \colon \R^d \to \R^d$ the function which satisfies for all $ x = ( x_1, \dots, x_{d} ) \in \R^{d} $ that
	\begin{equation}\label{multidim_version:Equation}
		\mathfrak{M}_{\psi, d}( x ) 
		=
		\left(
		\psi(x_1)
		,
		\ldots
		,
		\psi(x_d)
		\right).
	\end{equation}
	%	and we call $\mathfrak{M}_{\psi, d} $ the $d$-dimensional version of $\psi$.
\end{definition}

\begin{definition}[Realizations associated to ANNs]
	\label{Definition:ANNrealization}
	Let $a\in C(\R,\R)$.
	Then we denote by 
	$
	\functionANN \colon \ANNs \to \cup_{k,l\in\N}\,C(\R^k,\R^l)
	$
	the function which satisfies
	for all  $ L\in\N$, $l_0,l_1,\ldots, l_L \in \N$, 
	$
	\Phi 
	=
	((W_1, B_1),(W_2, B_2),\allowbreak \ldots, (W_L,\allowbreak B_L))
	\in  \allowbreak
	( \times_{k = 1}^L\allowbreak(\R^{l_k \times l_{k-1}} \times \R^{l_k}))
	$,
	$x_0 \in \R^{l_0}, x_1 \in \R^{l_1}, \ldots, x_{L-1} \in \R^{l_{L-1}}$ 
	with $\forall \, k \in \N \cap (0,L) \colon x_k =\activationDim{l_k}(W_k x_{k-1} + B_k)$  
	that
	\begin{equation}
	\label{setting_NN:ass2}
	\functionANN(\Phi) \in C(\R^{l_0},\R^{l_L})\qandq
	( \functionANN(\Phi) ) (x_0) = W_L x_{L-1} + B_L
	\end{equation}
	(cf.\ Definition~\ref{Def:multidim_version} and Definition~\ref{Def:ANN}).
\end{definition}

%\begin{definition}\label{Definition:ANNbullet}
%	Let $\Phi_1,\Phi_2\in\ANNs$ satisfy that $\outDimANN(\Phi_2)=\inDimANN(\Phi_1)$. Then we denote by 
%	$\compANN{\Phi_1}{\Phi_2}$ the element of $\ANNs$ which satisfies for all 
%	$i\in\{1,2\}$, $\big((W_1^i, B_1^i),\allowbreak(W_2^i, B_2^i),\allowbreak \ldots, (W_{\lengthANN(\Phi_i)}^i,\allowbreak B_{\lengthANN(\Phi_i)}^i)\big)\in\ANNs$ with
%	$\forall i\in\{1,2\}\colon \Phi_i
%	=
%	\big((W_1^i, B_1^i),\allowbreak(W_2^i, B_2^i),\allowbreak \ldots, (W_{\lengthANN(\Phi_i)}^i,\allowbreak B_{\lengthANN(\Phi_i)}^i)\big)$
% that
%	\begin{equation}\label{ANNoperations:Composition}
%	\begin{split}
%	\compANN{\Phi_1}{\Phi_2}=&\big((W_1^2, B_1^2),(W_2^2, B_2^2),\allowbreak \ldots, (W_{\lengthANN(\Phi_2)-1}^2,\allowbreak B_{\lengthANN(\Phi_2)-1}^2), (W_1^1 W_{\lengthANN(\Phi_2)}^2, W_1^1 B_{\lengthANN(\Phi_2)}^2+B_{1}^1), \\&\quad(W_2^1, B_2^1), (W_3^1, B_3^1),\ldots,(W_{\lengthANN(\Phi_1)}^1,\allowbreak B_{\lengthANN(\Phi_1)}^1) \big).
%	\end{split}
%	\end{equation}
%\end{definition}

\begin{lemma}\label{Lemma:elementaryPropertiesANN}
	Let $\Phi\in\ANNs$ (cf.\ Definition~\ref{Def:ANN}). Then 
	\begin{enumerate}[(i)]
		\item\label{elementaryPropertiesANN:ItemOne} it holds that $\dims(\Phi)\in \N^{\lengthANN(\Phi)+1}$ and 
		\item\label{elementaryPropertiesANN:ItemTwo} it holds for all $a\in C(\R,\R)$ that $\functionANN(\Phi)\in C(\R^{\inDimANN(\Phi)},\R^{\outDimANN(\Phi)})$
	\end{enumerate}
	(cf.\ Definition~\ref{Definition:ANNrealization}).
\end{lemma}

\begin{proof}[Proof of Lemma~\ref{Lemma:elementaryPropertiesANN}]	
	Note that the assumption that $\Phi\in \ANNs=\cup_{L \in \N}
	\cup_{ (l_0,l_1,\ldots, l_L) \in \N^{L+1} }
	(
	\times_{k = 1}^L (\R^{l_k \times l_{k-1}} \times \R^{l_k})
	)$ ensures that there exist $L \in \N$, $l_0,l_1,\ldots, l_L\in \N$ such that 
	\begin{equation}\label{elementaryPropertiesANN:SetForPhi}
		\Phi\in \left(
		\times_{k = 1}^L (\R^{l_k \times l_{k-1}} \times \R^{l_k})
		\right).
	\end{equation}
	Observe that \eqref{elementaryPropertiesANN:SetForPhi} assures that 
	\begin{equation}\label{elementaryPropertiesANN:Length}
		\lengthANN(\Phi)=L,\qquad \inDimANN(\Phi)=l_0,\qquad \outDimANN(\Phi)=l_L,
	\end{equation}
	\begin{equation}\label{elementaryPropertiesANN:Dims}
		\andq \dims(\Phi)=(l_0,l_1,\dots, l_L)\in \N^{L+1}=\N^{\lengthANN(\Phi)+1}.
	\end{equation}
	This establishes item~\eqref{elementaryPropertiesANN:ItemOne}. Moreover, note that \eqref{elementaryPropertiesANN:Length} and \eqref{setting_NN:ass2} show that 
	$\functionANN(\Phi) \in C(\R^{\inDimANN(\Phi)},\R^{\outDimANN(\Phi)})$. This establishes item~\eqref{elementaryPropertiesANN:ItemTwo}.
	The proof of Lemma~\ref{Lemma:elementaryPropertiesANN} is thus completed.
\end{proof}

\subsection{Compositions of ANNs}
\subsubsection{Standard compositions of ANNs}

%\begin{definition}[Standard compositions of GANNs]
%	\label{Definition:ANNcomposition}
%	We denote by $\compANN{(\cdot)}{(\cdot)}\colon\allowbreak \{(\Phi_1,\Phi_2)\allowbreak\in\ANNs\times \ANNs\colon \inDimANN(\Phi_1)=\outDimANN(\Phi_2)\}\allowbreak\to\ANNs$ the function which satisfies for all 
%	$ L,\mathfrak{L}\in\N$, $l_0,l_1,\ldots, l_L, \mathfrak{l}_0,\mathfrak{l}_1,\ldots, \mathfrak{l}_\mathfrak{L} \in \N$, 
%	$
%	\Phi_1
%	=
%	((W_1, B_1),(W_2, B_2),\allowbreak \ldots, (W_L,\allowbreak B_L))
%	\in  \allowbreak
%	( \times_{k = 1}^L\allowbreak(\R^{l_k \times l_{k-1}} \times \R^{l_k}))
%	$,
%	$
%\Phi_2
%=
%((\mathcal{W}_1, \mathfrak{B}_1),\allowbreak(\mathcal{W}_2, \mathfrak{B}_2),\allowbreak \ldots, (\mathcal{W}_\mathfrak{L},\allowbreak \mathfrak{B}_\mathfrak{L}))
%\in  \allowbreak
%( \times_{k = 1}^\mathfrak{L}\allowbreak(\R^{\mathfrak{l}_k \times \mathfrak{l}_{k-1}} \times \R^{\mathfrak{l}_k}))
%$ 
%with $l_0=\inDimANN(\Phi_1)=\outDimANN(\Phi_2)=\mathfrak{l}_{\mathfrak{L}}$
%%with $\inDimANN(N)=\outDimANN(\Phi_2)$
%	that
%	\begin{multline}\label{ANNoperations:Composition}
%			\compANN{\Phi_1}{\Phi_2}=\big((\mathcal{W}_1, \mathfrak{B}_1),(\mathcal{W}_2, \mathfrak{B}_2),\allowbreak \ldots, (\mathcal{W}_{\mathfrak{L}-1},\allowbreak \mathfrak{B}_{\mathfrak{L}-1}), (W_1 \mathcal{W}_{\mathfrak{L}}, W_1 \mathfrak{B}_{\mathfrak{L}}+B_{1}), \\(W_2, B_2), (W_3, B_3),\ldots,(W_{L},\allowbreak B_{L}) \big)
%	\end{multline}
%	(cf.\ Definition~\ref{Def:ANN}).
%\end{definition}

\begin{definition}[Standard compositions of ANNs]
	\label{Definition:ANNcomposition}
	We denote by $\compANN{(\cdot)}{(\cdot)}\colon\allowbreak \{(\Phi_1,\Phi_2)\allowbreak\in\ANNs\times \ANNs\colon \inDimANN(\Phi_1)=\outDimANN(\Phi_2)\}\allowbreak\to\ANNs$ the function which satisfies for all 
	$ L,\mathfrak{L}\in\N$, $l_0,l_1,\ldots, l_L, \mathfrak{l}_0,\mathfrak{l}_1,\ldots, \mathfrak{l}_\mathfrak{L} \in \N$, 
	$
	\Phi_1
	=
	((W_1, B_1),(W_2, B_2),\allowbreak \ldots, (W_L,\allowbreak B_L))
	\in  \allowbreak
	( \times_{k = 1}^L\allowbreak(\R^{l_k \times l_{k-1}} \times \R^{l_k}))
	$,
	$
	\Phi_2
	=
	((\mathcal{W}_1, \mathfrak{B}_1),\allowbreak(\mathcal{W}_2, \mathfrak{B}_2),\allowbreak \ldots, (\mathcal{W}_\mathfrak{L},\allowbreak \mathfrak{B}_\mathfrak{L}))
	\in  \allowbreak
	( \times_{k = 1}^\mathfrak{L}\allowbreak(\R^{\mathfrak{l}_k \times \mathfrak{l}_{k-1}} \times \R^{\mathfrak{l}_k}))
	$ 
	with $l_0=\inDimANN(\Phi_1)=\outDimANN(\Phi_2)=\mathfrak{l}_{\mathfrak{L}}$
	%with $\inDimANN(N)=\outDimANN(\Phi_2)$
	that
	\begin{equation}\label{ANNoperations:Composition}
	\begin{split}
	&\compANN{\Phi_1}{\Phi_2}=\\&
	\begin{cases} 
			\begin{array}{r}
			\big((\mathcal{W}_1, \mathfrak{B}_1),(\mathcal{W}_2, \mathfrak{B}_2),\ldots, (\mathcal{W}_{\mathfrak{L}-1},\allowbreak \mathfrak{B}_{\mathfrak{L}-1}),
			(W_1 \mathcal{W}_{\mathfrak{L}}, W_1 \mathfrak{B}_{\mathfrak{L}}+B_{1}),\\ (W_2, B_2), (W_3, B_3),\ldots,(W_{L},\allowbreak B_{L})\big)
			\end{array}
	&: L>1<\mathfrak{L} \\[3ex]
	\big( (W_1 \mathcal{W}_{1}, W_1 \mathfrak{B}_1+B_{1}), (W_2, B_2), (W_3, B_3),\ldots,(W_{L},\allowbreak B_{L}) \big)
	&: L>1=\mathfrak{L}\\[1ex]
	\big((\mathcal{W}_1, \mathfrak{B}_1),(\mathcal{W}_2, \mathfrak{B}_2),\allowbreak \ldots, (\mathcal{W}_{\mathfrak{L}-1},\allowbreak \mathfrak{B}_{\mathfrak{L}-1}),(W_1 \mathcal{W}_{\mathfrak{L}}, W_1 \mathfrak{B}_{\mathfrak{L}}+B_{1}) \big)
	&: L=1<\mathfrak{L}  \\[1ex]
	(W_1 \mathcal{W}_{1}, W_1 \mathfrak{B}_1+B_{1}) 
	&: L=1=\mathfrak{L} 
	\end{cases}
	\end{split}
	\end{equation}
	(cf.\ Definition~\ref{Def:ANN}).
\end{definition}

\begin{prop}\label{Lemma:PropertiesOfCompositions}
	Let 
	$\Phi_1,\Phi_2\in\ANNs$, 
	$l_{1,0},l_{1,1},\dots,\allowbreak l_{1,\lengthANN(\Phi_1)},l_{2,0},\allowbreak l_{2,1},\dots,\allowbreak l_{2,\lengthANN(\Phi_2)}
	\in\N$  satisfy for all $k\in\{1,2\}$ that
	$\inDimANN(\Phi_1)=\outDimANN(\Phi_2)$
    and $\dims(\Phi_k)=(l_{k,0},l_{k,1},\dots, l_{k,\lengthANN(\Phi_k)})$
	(cf.\ Definition~\ref{Def:ANN}).
	Then
	\begin{enumerate}[(i)]
\item \label{PropertiesOfCompositions:Dims} it holds that
\begin{equation}
\dims(\compANN{\Phi_1}{\Phi_2})=(l_{2,0},l_{2,1},\dots, l_{2,\lengthANN(\Phi_2)-1},l_{1,1},l_{1,2},\dots,l_{1,\lengthANN(\Phi_1)} ),
\end{equation}
\item \label{PropertiesOfCompositions:Length} it holds that
\begin{equation}\label{PropertiesOfCompositions:LengthDisplay}
[\lengthANN(\compANN{\Phi_1}{\Phi_2})-1]=[\lengthANN(\Phi_1)-1]+[\lengthANN(\Phi_2)-1],
\end{equation}
\item \label{PropertiesOfCompositions:HiddenLength} it holds that
\begin{equation}\label{PropertiesOfCompositions:HiddenLengthDisplay}
\hiddenLength(\compANN{\Phi_1}{\Phi_2})=\hiddenLength(\Phi_1)+\hiddenLength(\Phi_2),
\end{equation}
% $\lengthANN(\compANN{\Phi_1}{\Phi_2})-1=\big[\lengthANN(\Phi_1)-1\big]+\big[\lengthANN(\Phi_2)-1\big]$
\item \label{PropertiesOfCompositions:Params} it holds that
\begin{equation}
\begin{split}
\paramANN(\compANN{\Phi_1}{\Phi_2})&
= \paramANN(\Phi_1)+\paramANN(\Phi_2)
+l_{1,1}( l_{2,\lengthANN(\Phi_2)-1}+1)
\\&\quad
-l_{1,1}(l_{1,0} + 1)
-l_{2,\lengthANN(\Phi_2)}(l_{2,\lengthANN(\Phi_2)-1} + 1)
\\&\le
\paramANN(\Phi_1)+\paramANN(\Phi_2)+l_{1,1} l_{2,\lengthANN(\Phi_2)-1},
\end{split}
\end{equation}  
%\begin{equation}
%\begin{split}
%\paramANN(\compANN{\Phi_1}{\Phi_2})&
%=\left[\smallsum\limits_{m = 1}^{\lengthANN(\Phi_2)-1} l_{2,m}(l_{2,m-1} + 1) \right]+
%\left[\smallsum\limits_{m = 2}^{\lengthANN(\Phi_1)} l_{1,m}(l_{1,m-1} + 1) \right]
%\\&\quad+l_{1,1}(l_{2,\lengthANN(\Phi_2)-1} + 1)
%\\&\le
%\paramANN(\Phi_1)+\paramANN(\Phi_2)+l_{1,1}  (l_{2,\lengthANN(\Phi_2)-1}+1),
%\end{split}
%\end{equation}  
and
		\item \label{PropertiesOfCompositions:Realization} it holds
for all  $\activation\in C(\R,\R)$
that $\functionANN(\compANN{\Phi_1}{\Phi_2})\in C(\R^{\inDimANN(\Phi_2)},\R^{\outDimANN(\Phi_1)})$ and
\begin{equation}\label{PropertiesOfCompositions:RealizationEquation}
\functionANN(\compANN{\Phi_1}{\Phi_2})=[\functionANN(\Phi_1)]\circ [\functionANN(\Phi_2)]
\end{equation}
	\end{enumerate}
(cf.\  Definition~\ref{Definition:ANNrealization} and Definition~\ref{Definition:ANNcomposition}).
\end{prop}

\begin{proof}[Proof of Proposition~\ref{Lemma:PropertiesOfCompositions}]	
	Throughout this proof   let $a\in C(\R,\R)$,
	 let  $L_k\in\N$, $k\in\{1,2\}$, satisfy for all $k\in\{1,2\}$ that  
	 $L_k=\lengthANN(\Phi_k)$,
	let $\big((W_{k,1}, B_{k,1}),(W_{k,2}, B_{k,2}),\allowbreak \ldots,\allowbreak (W_{k,L_k}, B_{k,L_k})\big) 
	\in ( \times_{j = 1}^{L_k}\allowbreak(\R^{l_{k,j} \times l_{k,j-1}} \times \R^{l_{k,j}}))$, $k\in\{1,2\}$,
	%	$\big((W_{2,1}, B_{2,1}), \ldots, (W_{2,L_2}, B_{2,L_2})\big) 
	%	\in ( \times_{k = 1}^{L_2}\allowbreak(\R^{l_{2,k} \times l_{2,k-1}} \times \R^{l_{2,k}})) $,
	%	$\big((W_{3,1}, B_{3,1}),\allowbreak \ldots,\allowbreak (W_{3,L_1}, B_{3,L_1})\big) 
	%	\in ( \times_{k = 1}^{L_3}\allowbreak(\R^{l_{3,k} \times l_{3,k-1}} \times \R^{l_{3,k}}))$ 
	%	satisfy for all $k\in\{1,2,3\}$ that 
	%	\begin{equation}
	%	\Phi_1=\big((W_{1,1}, B_{1,1}),\allowbreak \ldots,\allowbreak (W_{1,L_1}, B_{1,L_1})\big), \qandq 
	%	\Phi_2=\big((W_{2,1}, B_{2,1}), \ldots, (W_{2,L_2}, B_{2,L_2})\big),
	%	\end{equation}	
	%		 $\dims(\Phi_3)=(l_{3,0},l_{3,1},\dots,\allowbreak l_{3,\lengthANN(\Phi_3)})$,
	satisfy for all $k\in\{1,2\}$ that
	\begin{equation}\label{PropertiesOfCompositions:DefOfArchitectures}
	%L_k=\lengthANN(\Phi_k)\qandq
	\Phi_k=\big((W_{k,1}, B_{k,1}),(W_{k,2}, B_{k,2}),\allowbreak \ldots,\allowbreak (W_{k,L_k}, B_{k,L_k})\big),
	\end{equation}
	let $L_3\in\N$, $l_{3,0}, l_{3,1},\dots, l_{3,L_3}\in\N$,
	$\Phi_3=\big((W_{3,1}, B_{3,1}), \ldots, (W_{3,L_3}, B_{3,L_3})\big) \in ( \times_{j = 1}^{L_3}\allowbreak(\R^{l_{3,j} \times l_{3,j-1}} \times \R^{l_{3,j}}))$
	satisfy  that  $\Phi_3=\compANN{\Phi_1}{\Phi_2}$,
	let $x_0 \in \R^{l_{2,0}}, x_1 \in \R^{l_{2,1}}, \ldots, x_{L_2-1} \in \R^{l_{2,L_2-1}}$ satisfy that
	\begin{equation}\label{PropertiesOfCompositions:NeuralSequenceOne}
	\forall\,j \in \N \cap (0,L_2)\colon\, x_j =\activationDim{l_{2,j}}(W_{2,j} x_{j-1} + B_{2,j})
	\end{equation}
	(cf.\ Definition~\ref{Def:multidim_version}),
%Falsch!
%	let $y_0 \in \R^{l_{1,0}}, y_1 \in \R^{l_{1,1}}, \ldots, y_{L_1-1} \in \R^{l_{1,L_1-1}}$ satisfy for all  $k \in \N \cap (0,L_1)$ that
%	\begin{equation}\label{PropertiesOfCompositions:NeuralSequenceTwo}
%	y_0= W_{2,L_2} x_{L_2-1} + B_{2,L_2} \qandq y_k =\activationDim{l_{1,k}}(W_{1,k} y_{k-1} + B_{1,k}),
%	\end{equation}
	let $y_0 \in \R^{l_{1,0}}, y_1 \in \R^{l_{1,1}}, \ldots, y_{L_1-1} \in \R^{l_{1,L_1-1}}$ satisfy that $y_0= W_{2,L_2} x_{L_2-1} + B_{2,L_2}$ and
	\begin{equation}\label{PropertiesOfCompositions:NeuralSequenceTwo}
	 \forall\,j \in \N \cap (0,L_1)\colon\, y_j =\activationDim{l_{1,j}}(W_{1,j} y_{j-1} + B_{1,j}),
	\end{equation}
	% $y_0= W_{L_1} x_{L_1-1} + B_{L_1}$ and $y_k =\activationDim{l_k^2}(W_k^2 y_{k-1} + B_k^2)$,
	and let $z_0 \in \R^{l_{3,0}}, z_1 \in \R^{l_{3,1}}, \ldots, z_{L_3-1} \in \R^{l_{3,L_3-1}}$ satisfy that $z_0=x_0$ and
	%  $z_k =\activationDim{l_k^4}(W_k^4 x_{k-1} + B_k^4)$, 
	\begin{equation}\label{PropertiesOfCompositions:NeuralSequenceThree}
	\forall\, j \in \N \cap (0,L_3)\colon\, z_j =\activationDim{l_{3,j}}(W_{3,j} z_{j-1} + B_{3,j}).
	\end{equation} 
	Note that \eqref{ANNoperations:Composition} ensures that
%	\begin{equation}\label{PropertiesOfCompositions:CompositionArchitecture}
%	\begin{split}
%	&\Phi_3
%	=\compANN{\Phi_1}{\Phi_2}
%	\\&=\big((W_{2,1}, B_{2,1}),(W_{2,2}, B_{2,2}),\allowbreak \ldots, (W_{2,L_2-1},\allowbreak B_{2,L_2-1}),(W_{1,1} W_{2,L_2}, W_{1,1} B_{2,L_2}+B_{1,1}), \\&\qquad (W_{1,2}, B_{1,2}), (W_{1,3}, B_{1,3}),\ldots,(W_{1,L_1},\allowbreak B_{1,L_1}) \big).
%	\end{split}
%	\end{equation}
%		\begin{equation}\label{PropertiesOfCompositions:CompositionArchitecture}
%		\begin{split}
%		&\Phi_3=\compANN{\Phi_1}{\Phi_2}=\\&
%		\begin{cases} 
%		\begin{array}{c}
%		\big((\mathcal{W}_1, \mathfrak{B}_1),(\mathcal{W}_2, \mathfrak{B}_2),\ldots, (\mathcal{W}_{L_2-1},\allowbreak \mathfrak{B}_{L_2-1}),\\
%		\,\,\,(W_1 \mathcal{W}_{L_2}, W_1 \mathfrak{B}_{L_2}+B_{1}), (W_2, B_2), (W_3, B_3),\ldots,(W_{L_1},\allowbreak B_{L_1})\big)
%		\end{array}
%		&: L_1>1<L_2 \\[3ex]
%		\big( (W_1 \mathcal{W}_{1}, W_1 \mathfrak{B}_1+B_{1}), (W_2, B_2), (W_3, B_3),\ldots,(W_{L_1},\allowbreak B_{L_1}) \big)
%		&: L_1>1=L_2\\[1ex]
%		\big((\mathcal{W}_1, \mathfrak{B}_1),(\mathcal{W}_2, \mathfrak{B}_2),\allowbreak \ldots, (\mathcal{W}_{L_2-1},\allowbreak \mathfrak{B}_{L_2-1}),(W_1 \mathcal{W}_{L_2}, W_1 \mathfrak{B}_{L_2}+B_{1}) \big)
%		&: L_1=1<L_2  \\[1ex]
%		(W_1 \mathcal{W}_{1}, W_1 \mathfrak{B}_1+B_{1}) 
%		&: L_1=1=L_2 
%		\end{cases}
%		\end{split}
%		\end{equation}
				\begin{equation}\label{PropertiesOfCompositions:CompositionArchitecture}
				\begin{split}
				&\Phi_3=\compANN{\Phi_1}{\Phi_2}=\\&
				\begin{cases} 
				\begin{array}{r}
				\big((W_{2,1}, B_{2,1}),(W_{2,2}, B_{2,2}),\ldots, (W_{2,L_2-1},\allowbreak B_{2,L_2-1}),\\
				\,\,\,(W_{1,1} W_{2,L_2}, W_{1,1} B_{2,L_2}+B_{1,1}),
				 (W_{1,2}, B_{1,2}),\\ (W_{1,3}, B_{1,3}),\ldots,(W_{1,L_1},\allowbreak B_{1,L_1})\big)
				\end{array}
				&: L_1>1<L_2 \\[5ex]
				\begin{array}{r}
				\big( (W_{1,1} W_{2,1}, W_{1,1} B_{2,1}+B_{1,1}), (W_{1,2}, B_{1,2}),\\ (W_{1,3}, B_{1,3}),\ldots,(W_{1,L_1},\allowbreak B_{1,L_1}) \big)
				\end{array}
				&: L_1>1=L_2\\[4ex]
				\begin{array}{r}
				\big((W_{2,1}, B_{2,1}),(W_{2,2}, B_{2,2}),\allowbreak \ldots, (W_{2,L_2-1},\allowbreak B_{2,L_2-1}),\\(W_{1,1} W_{2,L_2}, W_{1,1} B_{2,L_2}+B_{1,1}) \big)
				\end{array}
				&: L_1=1<L_2  \\[4ex]
				(W_{1,1} W_{2,1}, W_{1,1} B_{2,1}+B_{1,1}) 
				&: L_1=1=L_2 
				\end{cases}.
				\end{split}
				\end{equation}
%	(cf.\ Definition~\ref{Definition:ANNcomposition}).
	Hence, we obtain that 
	\begin{equation}\label{PropertiesOfCompositions:CompositionLength}
	\begin{split}
		[\lengthANN(\compANN{\Phi_1}{\Phi_2})-1]&=[(L_2-1)+1+(L_1-1)]-1\\&=[L_1-1]+[L_2-1]=[\lengthANN(\Phi_1)-1]+[\lengthANN(\Phi_2)-1]
	\end{split}
	\end{equation}
	%$\lengthANN(\compANN{\Phi_1}{\Phi_2})=(L_2-1)+1+(L_1-1)=L_1+L_2-1=\lengthANN(\Phi_1)+\lengthANN(\Phi_2)-1,$
	\begin{equation}\label{PropertiesOfCompositions:CompositionDims}
	\andq\dims(\compANN{\Phi_1}{\Phi_2})=(l_{2,0},l_{2,1},\dots, l_{2,L_2-1},l_{1,1},l_{1,2},\dots,l_{1,L_1} ).
	\end{equation}
%and 
%	\begin{equation}\label{PropertiesOfCompositions:CompositionDims}
%	\dims(\compANN{\Phi_1}{\Phi_2})
%	=				\begin{cases} 
%	(l_{2,0},l_{2,1},\dots, l_{2,L_2-1},l_{1,1},l_{1,2},\dots,l_{1,L_1} )
%	&: L_1>1<L_2 \\
%	(l_{1,1},l_{1,2},\dots,l_{1,L_1} )
%	&: L_1>1=L_2\\
%	(l_{2,0},l_{2,1},\dots, l_{2,L_2-1},l_{1,1},l_{1,2},\dots,l_{1,L_1} )
%	&: L_1=1<L_2  \\
%	(l_{2,0},l_{2,1},\dots, l_{2,L_2-1},l_{1,1},l_{1,2},\dots,l_{1,L_1} )
%	&: L_1=1=L_2 
%	\end{cases}.
%%	=(l_{2,0},l_{2,1},\dots, l_{2,L_2-1},l_{1,1},l_{1,2},\dots,l_{1,L_1} ).
%	\end{equation}
	This establishes items~\eqref{PropertiesOfCompositions:Dims}--\eqref{PropertiesOfCompositions:HiddenLength}.
	%	$\dims(\compANN{\Phi_1}{\Phi_2})=(\ell_{2,0},\ell_{2,1},\dots, \ell_{2,\lengthANN(\Phi_2)-1},\ell_{1,1},\ell_{1,2},\dots, \ell_{1,\lengthANN(\Phi_1)})$.
	In addition, observe that \eqref{PropertiesOfCompositions:CompositionDims} demonstrates that
	\begin{equation}\label{PropertiesOfCompositions:CompositionParameters}
	\begin{split}
	&\paramANN(\compANN{\Phi_1}{\Phi_2})
	=
	\smallsum\limits_{j = 1}^{L_3} l_{3,j}(l_{3,j-1} + 1)
	\\&=\left[\smallsum\limits_{j = 1}^{L_2-1} l_{3,j}(l_{3,j-1} + 1)\right]+l_{3,L_2}(l_{3,L_2-1} + 1)+\left[\smallsum\limits_{j = L_2+1}^{L_3} l_{3,j}(l_{3,j-1} + 1)\right]
	\\&=\left[\smallsum\limits_{j = 1}^{L_2-1} l_{2,j}(l_{2,j-1} + 1)\right]+l_{1,1}(l_{2,L_2-1} + 1)+\left[\smallsum\limits_{j = L_2+1}^{L_3} l_{1,j-L_2+1}(l_{1,j-L_2} + 1)\right]
	\\&=
	\left[\smallsum\limits_{j = 1}^{L_2-1} l_{2,j}(l_{2,j-1} + 1) \right]+
	\left[\smallsum\limits_{j = 2}^{L_1} l_{1,j}(l_{1,j-1} + 1) \right]
	+l_{1,1}\big(l_{2,L_2-1} + 1\big)
	\\&= \left[\smallsum\limits_{j = 1}^{L_2} l_{2,j}(l_{2,j-1} + 1) \right]+
	\left[\smallsum\limits_{j = 1}^{L_1} l_{1,j}(l_{1,j-1} + 1) \right]
	+l_{1,1}( l_{2,L_2-1}+1)
	\\&\quad-l_{2,L_2}(l_{2,L_2-1} + 1)-l_{1,1}(l_{1,0} + 1)
		\\&= \paramANN(\Phi_1)+\paramANN(\Phi_2)
	+l_{1,1}( l_{2,L_2-1}+1)-l_{2,L_2}(l_{2,L_2-1} + 1)
	\\&\quad-l_{1,1}(l_{1,0} + 1)
	\\&\le \paramANN(\Phi_1)+\paramANN(\Phi_2)
	+l_{1,1}l_{2,L_2-1}.
	\end{split}
	\end{equation}
	This establishes item~\eqref{PropertiesOfCompositions:Params}.
	%		\begin{equation}\label{PropertiesOfCompositions:CompositionDimsTwo}
	%	\inDimANN(\compANN{\Phi_1}{\Phi_2})=\inDimANN(\Phi_2), \qandq \outDimANN(\compANN{\Phi_1}{\Phi_2})=\outDimANN(\Phi_1).
	%	\end{equation}
	Moreover, observe that \eqref{PropertiesOfCompositions:CompositionArchitecture} and the fact that $a\in C(\R,\R)$ ensure that 	
	\begin{equation}\label{PropertiesOfCompositions:Continuous}
	\functionANN(\compANN{\Phi_1}{\Phi_2})\in
	C(\R^{l_{2,0}},\R^{l_{1,L_1}})=
	 C(\R^{\inDimANN(\Phi_2)},\R^{\outDimANN(\Phi_1)}).
	\end{equation}
	Next note that \eqref{PropertiesOfCompositions:CompositionLength} implies that $L_3=L_1+L_2-1$.
	This,  \eqref{PropertiesOfCompositions:CompositionArchitecture},  and \eqref{PropertiesOfCompositions:CompositionDims} ensure that
%	for all $j\in\N\cap(0, L_2)$, $m\in \N\cap(L_2, L_1+L_2)$ it holds that
	\begin{equation}\label{PropertiesOfCompositions:CompositionDimsTwo}
	(l_{3,0},l_{3,1},\dots, l_{3,L_1+L_2-1})=(l_{2,0},l_{2,1},\dots, l_{2,L_2-1},l_{1,1},l_{1,2},\dots,l_{1,L_1}),
	\end{equation}
	\begin{equation}\label{PropertiesOfCompositions:CompositionSchnittstelleOne}
	\begin{split}
	\big[\forall\,j\in\N\cap(0, L_2)\colon\, (W_{3,j}, B_{3,j})=(W_{2,j}, B_{2,j})\big],
	\end{split}
	\end{equation}
		\begin{equation}\label{PropertiesOfCompositions:CompositionSchnittstelleTwo}
		(W_{3,L_2}, B_{3,L_2})=(W_{1,1} W_{2,L_2}, W_{1,1} B_{2,L_2}+B_{1,1}),
		\end{equation}
	\begin{equation}\label{PropertiesOfCompositions:CompositionSchnittstelleThree}
		\andq \big[\forall\,j\in \N\cap(L_2, L_1+L_2)\colon\,	(W_{3,j}, B_{3,j})=(W_{1,j+1-L_2}, B_{1,j+1-L_2})\big].
	\end{equation}
	This, \eqref{PropertiesOfCompositions:NeuralSequenceOne}, \eqref{PropertiesOfCompositions:NeuralSequenceThree},  and induction imply that for all $j\in\N_0\cap [0,L_2)$ it holds that $z_j=x_j$. Combining this with \eqref{PropertiesOfCompositions:CompositionSchnittstelleTwo} and the fact that $y_0=W_{2,L_2} x_{{L_2}-1} + B_{2,L_2}$ ensures that
		\begin{equation}\label{PropertiesOfCompositions:CompositionCalc}
		\begin{split}
		W_{3,L_2} z_{{L_2}-1} + B_{3,L_2}
		&=W_{3,L_2} x_{{L_2}-1} + B_{3,L_2}
		\\&=W_{1,1} W_{2,L_2} x_{{L_2}-1} +W_{1,1} B_{2,L_2}+B_{1,1}
		\\&=W_{1,1} (W_{2,L_2} x_{{L_2}-1} + B_{2,L_2})+B_{1,1}
		=W_{1,1} y_0+B_{1,1}.
		\end{split}
		\end{equation}
			Next we claim that for all $j\in\N\cap [L_2,L_1+L_2)$ it holds that
			\begin{equation}\label{PropertiesOfCompositions:AuxiliaryClaim}
			W_{3,j} z_{j-1} + B_{3,j}
			=W_{1,j+1-L_2} y_{j-L_2} + B_{1,j+1-L_2}.
			\end{equation}
%			To prove \eqref{PropertiesOfCompositions:AuxiliaryClaim} we distinguish between the cases $L_1=1$ and $L_1>1$. 
%			%	We first prove \eqref{PropertiesOfCompositions:AuxiliaryClaim} in the case $L_1=1$. 
%			Observe that \eqref{PropertiesOfCompositions:CompositionCalc} establishes  \eqref{PropertiesOfCompositions:AuxiliaryClaim} in the case $L_1=1$.
%%			We now prove \eqref{PropertiesOfCompositions:AuxiliaryClaim} in the case $L_1>1$. 
%			In the case $L_1>1$ 
			We prove \eqref{PropertiesOfCompositions:AuxiliaryClaim} by induction on $j\in\N\cap [L_2,L_1+L_2)$.
			Note that \eqref{PropertiesOfCompositions:CompositionCalc} establishes  \eqref{PropertiesOfCompositions:AuxiliaryClaim} in the base case $j=L_2$.
			For the induction step 
%			$[L_2,L_1+L_2-1)\ni m\to m+1\in [L_2+1,L_1+L_2)$ assume that there exists $m\in\N\cap [L_2,L_1+L_2-1)$ which satisfies that
%					\begin{equation}\label{PropertiesOfCompositions:AuxiliaryClaimInduction}
%					W_{3,j} z_{j-1} + B_{3,j}
%					=W_{1,mj+1-L_2} y_{j-L_2} + B_{1,j+1-L_2}.
%					\end{equation}		
			note that the fact that $L_3=L_1+L_2-1$, \eqref{PropertiesOfCompositions:NeuralSequenceTwo}, \eqref{PropertiesOfCompositions:NeuralSequenceThree},  \eqref{PropertiesOfCompositions:CompositionDimsTwo}, and \eqref{PropertiesOfCompositions:CompositionSchnittstelleThree} imply that
			for all $j\in\N\cap [L_2,\infty)\cap (0,L_1+L_2-1)$ with 
								\begin{equation}\label{PropertiesOfCompositions:AuxiliaryClaimInduction}
								W_{3,j} z_{j-1} + B_{3,j}
								=W_{1,j+1-L_2} y_{j-L_2} + B_{1,j+1-L_2}
								\end{equation}
				it holds that 
			\begin{equation}
				\begin{split}
				&W_{3,j+1} z_{j} + B_{3,j+1}
				=W_{3,j+1}\activationDim{l_{3,j}}(W_{3,j} z_{j-1} + B_{3,j}) + B_{3,j+1}
				\\&=W_{1,j+2-L_2}\activationDim{l_{1,j+1-L_2}}(W_{1,j+1-L_2} y_{j-L_2} + B_{1,j+1-L_2}) + B_{1,j+2-L_2}
				\\&=W_{1,j+2-L_2} y_{j+1-L_2} + B_{1,j+2-L_2}.
				\end{split}
			\end{equation}
			Induction hence proves  \eqref{PropertiesOfCompositions:AuxiliaryClaim}.
			Next observe that \eqref{PropertiesOfCompositions:AuxiliaryClaim} and the fact that $L_3=L_1+L_2-1$ assure that
			\begin{equation}
				W_{3,L_3} z_{L_3-1} + B_{3,L_3}
				=W_{3,L_1+L_2-1} z_{L_1+L_2-2} + B_{3,L_1+L_2-1}
				=W_{1,L_1} y_{L_1-1} + B_{1,L_1}.
			\end{equation}
The fact that $\Phi_3=\compANN{\Phi_1}{\Phi_2}$, \eqref{PropertiesOfCompositions:NeuralSequenceOne},
	\eqref{PropertiesOfCompositions:NeuralSequenceTwo}, and
	\eqref{PropertiesOfCompositions:NeuralSequenceThree}
	therefore prove that
	\begin{equation}
	\begin{split}
	[\functionANN(\compANN{\Phi_1}{\Phi_2})](x_0)
	&=[\functionANN(\Phi_3)](x_0)
	=[\functionANN(\Phi_3)](z_0)
	=W_{3,L_3} z_{{L_3}-1} + B_{3,L_3}
	%	\\&=W_{3,L_3} z_{{L_3}-1} + B_{3,L_3}
	\\&=W_{1,L_1} y_{L_1-1} + B_{1,L_1}
	=[\functionANN(\Phi_1)](y_0)
	\\&=[\functionANN(\Phi_1)]\big(W_{2,L_2} x_{L_2-1} + B_{2,L_2}\big)
	\\&=[\functionANN(\Phi_1)]\big([\functionANN(\Phi_2)](x_0)\big)
	=[(\functionANN(\Phi_1))\circ (\functionANN(\Phi_2))](x_0).
	\end{split}
	\end{equation}
%	Combining this with the fact that $\functionANN(\Phi_1)\in C(\R^{\inDimANN(\Phi_1)},\R^{\outDimANN(\Phi_1)})$ and 
%	$\functionANN(\Phi_2)\in C(\R^{\inDimANN(\Phi_2)},\R^{\outDimANN(\Phi_2)})$
%	%	 $\functionANN(\compANN{\Phi_1}{\Phi_2})=(\functionANN(\Phi_1))\circ (\functionANN(\Phi_2))$.
%	demonstrates that
%	\begin{equation}\label{PropertiesOfCompositions:BulletComp}
%	\functionANN(\compANN{\Phi_1}{\Phi_2})=(\functionANN(\Phi_1))\circ (\functionANN(\Phi_2))\in C(\R^{\inDimANN(\Phi_2)},\R^{\outDimANN(\Phi_1)}).
%	\end{equation}	
%Hence, we obtain that $\functionANN(\compANN{\Phi_1}{\Phi_2})=(\functionANN(\Phi_1))\circ (\functionANN(\Phi_2))$.
	Combining this with \eqref{PropertiesOfCompositions:Continuous} establishes item~\eqref{PropertiesOfCompositions:Realization}.
%	This, \eqref{PropertiesOfCompositions:CompositionLength}, \eqref{PropertiesOfCompositions:CompositionDims}, and \eqref{PropertiesOfCompositions:CompositionParameters} establish item~\eqref{PropertiesOfCompositions:ItemComp}.
	The proof of Proposition~\ref{Lemma:PropertiesOfCompositions} is thus completed.
\end{proof}

\begin{cor}\label{Corollary:Composition}
	Let $L_1,L_2,L_3\in\N$, $l_{1,0},l_{1,1},\dots, l_{1,L_1},l_{2,0},l_{2,1},\dots,l_{2,L_2},l_{3,0},l_{3,1},\dots,\allowbreak l_{3,L_3}\in\N$  satisfy that $l_{1,0}=l_{2,L_2}$ and let 
	$\Phi_k=\big((W_{k,1}, B_{k,1}),(W_{k,2}, B_{k,2}),\allowbreak \ldots,\allowbreak (W_{k,L_k}, B_{k,L_k})\big) 
	\in ( \times_{j = 1}^{L_k}\allowbreak(\R^{l_{k,j} \times l_{k,j-1}} \times \R^{l_{k,j}}))$, $k\in\{1,2, 3\}$, satisfy that 
	$\Phi_3=\compANN{\Phi_1}{\Phi_2}$ (cf.\ Definition~\ref{Def:ANN} and Definition~\ref{Definition:ANNcomposition}). Then 
	\begin{enumerate}[(i)]
		\item \label{Corollary:Composition:ItemOne} 
		it holds that 
		\begin{equation}
		L_3=\lengthANN(\Phi_3)=\lengthANN(\Phi_1)+\lengthANN(\Phi_2)-1=L_1+L_2-1\ge \max\{L_1,L_2\},
		\end{equation}
		\item \label{Corollary:Composition:ItemTwo}  it holds for all $j\in\N\cap (0,L_2)$ that
		\begin{equation}
		(W_{3,j}, B_{3,j})=(W_{2,j}, B_{2,j}),
		\end{equation}	
		\item \label{Corollary:Composition:ItemThree}  it holds that 
		\begin{equation}
		(W_{3,L_2}, B_{3,L_2})=(W_{1,1}W_{2,L_2},W_{1,1}B_{2,L_2}+B_{1,1}),
		\end{equation}
		and
		\item \label{Corollary:Composition:ItemFour}  it holds for all $j\in\N\cap (L_2,L_1+L_2)=\N\cap (L_2,\infty)\cap [1,L_3]$ that
		\begin{equation}
		(W_{3,j}, B_{3,j})=(W_{1,j-L_2+1}, B_{1,j-L_2+1}).
		\end{equation}
	\end{enumerate} 
\end{cor}

\begin{proof}[Proof of Corollary~\ref{Corollary:Composition}]	
	Observe that item~\eqref{PropertiesOfCompositions:Length} in Proposition~\ref{Lemma:PropertiesOfCompositions} proves item~\eqref{Corollary:Composition:ItemOne}.
	Moreover, note that \eqref{ANNoperations:Composition} establishes items~\eqref{Corollary:Composition:ItemTwo}--\eqref{Corollary:Composition:ItemFour}.	
	The proof of Corollary~\ref{Corollary:Composition} is thus completed.
\end{proof}

\subsubsection{Associativity of standard compositions of ANNs}

\begin{lemma}\label{Lemma:CompositionAssociative}
	Let 
	$\Phi_1,\Phi_2,\Phi_3\in\ANNs$
	satisfy that
	$\inDimANN(\Phi_1)=\outDimANN(\Phi_2)$ and
	$\inDimANN(\Phi_2)=\outDimANN(\Phi_3)$ 
	(cf.\ Definition~\ref{Def:ANN}).
	Then
	it holds that 
	\begin{equation}
	\compANN{(\compANN{\Phi_1}{\Phi_2})}{\Phi_3}=\compANN{\Phi_1}{(\compANN{\Phi_2}{\Phi_3})}
	\end{equation}
	(cf.\  Definition~\ref{Definition:ANNcomposition}).
\end{lemma}

\begin{proof}[Proof of Lemma~\ref{Lemma:CompositionAssociative}]	
	Throughout this proof   
	let $\Phi_4,\Phi_5,\Phi_6,\Phi_7\in \ANNs$ satisfy that $\Phi_4=\compANN{\Phi_1}{\Phi_2}$,
	$\Phi_5=\compANN{\Phi_2}{\Phi_3}$, $\Phi_6=\compANN{\Phi_4}{\Phi_3}$, and $\Phi_7=\compANN{\Phi_1}{\Phi_5}$,
	let
	$L_k\in\N$, $k\in\{1,2,\dots, 7\}$, satisfy 
	for all $k\in\{1,2,\dots, 7\}$ that $L_k=\lengthANN(\Phi_k)$,
	let $l_{k,0},l_{k,1},\dots,\allowbreak l_{k,L_k}\in\N$, $k\in\{1,2,\dots, 7\}$,
	and let $\big((W_{k,1}, B_{k,1}),(W_{k,2}, B_{k,2}),\allowbreak \ldots,\allowbreak (W_{k,L_k}, B_{k,L_k})\big) 
	\in ( \times_{j = 1}^{L_k}\allowbreak(\R^{l_{k,j} \times l_{k,j-1}} \times \R^{l_{k,j}}))$, $k\in\{1,2,\dots, 7\}$, 
	%	$\big((W_{2,1}, B_{2,1}), \ldots, (W_{2,L_2}, B_{2,L_2})\big) 
	%	\in ( \times_{k = 1}^{L_2}\allowbreak(\R^{l_{2,k} \times l_{2,k-1}} \times \R^{l_{2,k}})) $,
	%	$\big((W_{3,1}, B_{3,1}),\allowbreak \ldots,\allowbreak (W_{3,L_1}, B_{3,L_1})\big) 
	%	\in ( \times_{k = 1}^{L_3}\allowbreak(\R^{l_{3,k} \times l_{3,k-1}} \times \R^{l_{3,k}}))$ 
	%	satisfy for all $k\in\{1,2,3\}$ that 
	%	\begin{equation}
	%	\Phi_1=\big((W_{1,1}, B_{1,1}),\allowbreak \ldots,\allowbreak (W_{1,L_1}, B_{1,L_1})\big), \qandq 
	%	\Phi_2=\big((W_{2,1}, B_{2,1}), \ldots, (W_{2,L_2}, B_{2,L_2})\big),
	%	\end{equation}
	satisfy for all $k\in\{1,2,\dots, 7\}$ that  
	\begin{equation}\label{CompositionAssociative:DefOfArchitectures}
	%L_k=\lengthANN(\Phi_k)\qandq
	\Phi_k=\big((W_{k,1}, B_{k,1}),(W_{k,2}, B_{k,2}),\allowbreak \ldots,\allowbreak (W_{k,L_k}, B_{k,L_k})\big).
	\end{equation}
	%	and let $L_3\in\N$, $l_{3,0}, l_{3,1},\dots, l_{3,L_3}\in\N$,
	%	$\Phi_3=\big((W_{3,1}, B_{3,1}), \ldots, (W_{3,L_3}, B_{3,L_3})\big) \in ( \times_{j = 1}^{L_3}\allowbreak(\R^{l_{3,j} \times l_{3,j-1}} \times \R^{l_{3,j}}))$
	%	satisfy  that  $\Phi_3=\compANN{\Phi_1}{\Phi_2}$.
	Observe that item~\eqref{PropertiesOfCompositions:Length} in Proposition~\ref{Lemma:PropertiesOfCompositions} and the fact that for all $k\in \{1,2,3\}$ it holds that $\lengthANN(\Phi_k)=L_k$ proves that
	\begin{equation}\label{associative:Length}
	\begin{split}
	\lengthANN(\Phi_6)&=\lengthANN(\compANN{(\compANN{\Phi_1}{\Phi_2})}{\Phi_3})
	=\lengthANN(\compANN{\Phi_1}{\Phi_2})+\lengthANN(\Phi_3)-1
	\\&=\lengthANN(\Phi_1)+\lengthANN(\Phi_2)+\lengthANN(\Phi_3)-2
	=L_1+L_2+L_3-2
	\\&=\lengthANN(\Phi_1)+\lengthANN(\compANN{\Phi_2}{\Phi_3})-1
	=\lengthANN(\compANN{\Phi_1}{(\compANN{\Phi_2}{\Phi_3})})
	=\lengthANN(\Phi_7).
	\end{split}
	\end{equation}
	Next note that Corollary~\ref{Corollary:Composition}, \eqref{CompositionAssociative:DefOfArchitectures}, and the fact that $\Phi_4=\compANN{\Phi_1}{\Phi_2}$ imply that
	\begin{equation}\label{FourONE}
	\big[\forall\, j\in\N\cap (0,L_2)\colon\, (W_{4,j}, B_{4,j})=(W_{2,j}, B_{2,j})\big],
	\end{equation}	
	\begin{equation}\label{FourTWO}
	(W_{4,L_2}, B_{4,L_2})=(W_{1,1}W_{2,L_2},W_{1,1}B_{2,L_2}+B_{1,1}),
	\end{equation}
	%$(W_{4,L_2}, B_{4,L_2})=(W_{1,1}W_{2,L_2},W_{1,1}B_{2,L_2}+B_{1,1})$, and 
	\begin{equation}\label{FourTHREE}
	\andq\big[\forall\, j\in\N\cap (L_2,L_1+L_2)\colon\, (W_{4,j}, B_{4,j})=(W_{1,j+1-L_2}, B_{1,j+1-L_2})\big].
	\end{equation}
	Hence, we obtain that 
	\begin{multline}\label{FourOne}
	\big[\forall\, j\in\N\cap (L_3-1,L_2+L_3-1)\colon\\ (W_{4,j+1-L_3}, B_{4,j+1-L_3})=(W_{2,j+1-L_3}, B_{2,j+1-L_3})\big],
	\end{multline}	
	\begin{equation}
	(W_{4,L_2}, B_{4,L_2})=(W_{1,1}W_{2,L_2},W_{1,1}B_{2,L_2}+B_{1,1}),
	\end{equation}
	and
	\begin{multline}\label{FourTwo}
	\big[\forall\, j\in\N\cap (L_2+L_3-1,L_1+L_2+L_3-1)\colon\\ (W_{4,j+1-L_3}, B_{4,j+1-L_3})=(W_{1,j+2-L_2-L_3}, B_{1,j+2-L_2-L_3})\big].
	\end{multline}
	%		\begin{equation}\label{FourTwo}
	%\forall\, j\in\N\cap (L_2+L_3-1,L_1+L_2+L_3-1)\colon\, (W_{4,j+1-L_3}, B_{4,j+1-L_3})=(W_{1,j+2-L_2-L_3}, B_{1,j+2-L_2-L_3}).
	%		\end{equation}
	In addition, observe that Corollary~\ref{Corollary:Composition}, \eqref{CompositionAssociative:DefOfArchitectures}, and the fact that $\Phi_5=\compANN{\Phi_2}{\Phi_3}$ demonstrate that
	\begin{equation}\label{FiveOne}
	\big[\forall\, j\in\N\cap (0,L_3)\colon\, (W_{5,j}, B_{5,j})=(W_{3,j}, B_{3,j})\big],
	\end{equation}	
	\begin{equation}\label{FiveTwo}
	(W_{5,L_3}, B_{5,L_3})=(W_{2,1}W_{3,L_3},W_{2,1}B_{3,L_3}+B_{2,1}),
	\end{equation}
	\begin{equation}\label{FiveThree}
	\andq\big[\forall\, j\in\N\cap (L_3,L_2+L_3)\colon\, (W_{5,j}, B_{5,j})=(W_{2,j+1-L_3}, B_{2,j+1-L_3})\big].
	\end{equation}
	Moreover, note that Corollary~\ref{Corollary:Composition}, \eqref{CompositionAssociative:DefOfArchitectures}, and the fact that $\Phi_6=\compANN{\Phi_4}{\Phi_3}$ ensure that 
	\begin{equation}\label{SixOne}
	\big[\forall\, j\in\N\cap (0,L_3)\colon\, (W_{6,j}, B_{6,j})=(W_{3,j}, B_{3,j})\big],
	\end{equation}	
	\begin{equation}\label{SixTwo}
	(W_{6,L_3}, B_{6,L_3})=(W_{4,1}W_{3,L_3},W_{4,1}B_{3,L_3}+B_{4,1}),
	\end{equation}
	\begin{equation}\label{SixThree}
	\andq\big[\forall\, j\in\N\cap (L_3,L_4+L_3)\colon\, (W_{6,j}, B_{6,j})=(W_{4,j+1-L_3}, B_{4,j+1-L_3})\big].
	\end{equation}
	Furthermore, observe that Corollary~\ref{Corollary:Composition}, \eqref{CompositionAssociative:DefOfArchitectures}, and the fact that $\Phi_7=\compANN{\Phi_1}{\Phi_5}$ show that 
	\begin{equation}\label{SevenOne}
	\big[\forall\, j\in\N\cap (0,L_5)\colon\, (W_{7,j}, B_{7,j})=(W_{5,j}, B_{5,j})\big],
	\end{equation}	
	\begin{equation}\label{SevenTwo}
	(W_{7,L_5}, B_{7,L_5})=(W_{1,1}W_{5,L_5},W_{1,1}B_{5,L_5}+B_{1,1}),
	\end{equation}
	\begin{equation}\label{SevenThree}
	\andq\big[\forall\, j\in\N\cap (L_5,L_1+L_5)\colon\, (W_{7,j}, B_{7,j})=(W_{1,j+1-L_5}, B_{1,j+1-L_5})\big].
	\end{equation}
	This, the fact that $L_3\le L_2+L_3-1=L_5$, \eqref{FiveOne}, and \eqref{SixOne} imply that for all $j\in \N\cap (0,L_3)$ it holds that 
	\begin{equation}\label{associative:One}
	(W_{6,j}, B_{6,j})=(W_{3,j}, B_{3,j})=(W_{5,j}, B_{5,j})=(W_{7,j}, B_{7,j}).
	\end{equation}
	In addition, observe that \eqref{FourONE}, \eqref{FourTWO}, \eqref{FiveOne}, \eqref{FiveTwo}, \eqref{SixTwo}, \eqref{SevenOne},  \eqref{SevenTwo}, and the fact that $L_5=L_2+L_3-1$ demonstrate that
	\begin{equation}\label{associative:Two}
	\begin{split}
	&(W_{6,L_3}, B_{6,L_3})=(W_{4,1}W_{3,L_3},W_{4,1}B_{3,L_3}+B_{4,1})
	\\&=\begin{cases} 
	(W_{2,1}W_{3,L_3},W_{2,1}B_{3,L_3}+B_{2,1})
	&: L_2>1\\
	(W_{1,1}W_{2,1}W_{3,L_3},W_{1,1}W_{2,1}B_{3,L_3}+W_{1,1}B_{2,1}+B_{1,1})
	&: L_2=1
	\end{cases}
		\\&=\begin{cases} 
		(W_{2,1}W_{3,L_3},W_{2,1}B_{3,L_3}+B_{2,1})
		&: L_2>1\\
		(W_{1,1}(W_{2,1}W_{3,L_3}),W_{1,1}(W_{2,1}B_{3,L_3}+B_{2,1})+B_{1,1})
		&: L_2=1
		\end{cases}
	\\&=\begin{cases} 
	(W_{5,L_3},B_{5,L_3})
	&: L_2>1\\
	(W_{1,1}W_{5,L_3},W_{1,1}B_{5,L_3}+B_{1,1})
	&: L_2=1
	\end{cases}
	\\&=(W_{7,L_3}, B_{7,L_3}).
	\end{split}
	\end{equation}
	Next note that the fact that $L_5=L_2+L_3-1<L_1+L_2+L_3-1=L_3+L_4$, \eqref{SixThree}, \eqref{FourOne}, \eqref{FiveThree}, and \eqref{SevenOne} ensure that for all $j\in \N$ with $L_3<j<L_5$  it holds that  
	\begin{equation}\label{associative:Three}
	\begin{split}
	(W_{6,j}, B_{6,j})&=(W_{4,j+1-L_3}, B_{4,j+1-L_3})=(W_{2,j+1-L_3}, B_{2,j+1-L_3})
	\\&=	(W_{5,j}, B_{5,j})=	(W_{7,j}, B_{7,j}).
	\end{split}
	\end{equation}
	Moreover, observe that the fact that $L_5=L_2+L_3-1<L_1+L_2+L_3-1=L_3+L_4$, \eqref{SixThree}, \eqref{associative:Two}, \eqref{FourTWO}, \eqref{FiveThree}, and \eqref{SevenTwo} prove that
	\begin{equation}\label{associative:Four}
	\begin{split}
	&(W_{6,L_5}, B_{6,L_5})=\begin{cases} 
	(W_{4,L_5+1-L_3}, B_{4,L_5+1-L_3})
	&: L_2>1\\
	(W_{6,L_3}, B_{6,L_3})
	&: L_2=1
	\end{cases}
	\\&=\begin{cases} 
	(W_{4,L_2}, B_{4,L_2})
	&: L_2>1\\
	(W_{7,L_3}, B_{7,L_3})
	&: L_2=1
	\end{cases}
	\\&=\begin{cases} 
	(W_{1,1}W_{2,L_2},W_{1,1}B_{2,L_2}+B_{1,1})
	&: L_2>1\\
	(W_{7,L_5}, B_{7,L_5})
	&: L_2=1
	\end{cases}
	\\&=\begin{cases} 
	(W_{1,1}W_{5,L_5},W_{1,1}B_{5,L_5}+B_{1,1})
	&: L_2>1\\
	(W_{7,L_5}, B_{7,L_5})
	&: L_2=1
	\end{cases}
	\\&=(W_{7,L_5}, B_{7,L_5}).
	\end{split}
	\end{equation}
	Furthermore, note that \eqref{SixThree}, \eqref{FourTwo}, \eqref{SevenThree}, and the fact that $L_5=L_2+L_3-1\ge L_3$ assure that for all $j\in \N$ with $L_5<j\le L_6$ it holds that 
	\begin{equation}\label{associative:Five}
	\begin{split}
	(W_{6,j}, B_{6,j})&=(W_{4,j+1-L_3}, B_{4,j+1-L_3})=(W_{1,j+2-L_2-L_3}, B_{1,j+2-L_2-L_3})
	\\&=(W_{1,j+1-L_5}, B_{1,j+1-L_5})=(W_{7,j}, B_{7,j}).
	\end{split}
	\end{equation}
	Combining this with \eqref{associative:Length}, \eqref{associative:One}, \eqref{associative:Two}, \eqref{associative:Three}, and \eqref{associative:Four} establishes that 
	\begin{equation}
	\compANN{(\compANN{\Phi_1}{\Phi_2})}{\Phi_3}
	=\compANN{\Phi_4}{\Phi_3}=\Phi_6=\Phi_7=\compANN{\Phi_1}{\Phi_5}
	=\compANN{\Phi_1}{(\compANN{\Phi_2}{\Phi_3})}.
	\end{equation}	
	The proof of Lemma~\ref{Lemma:CompositionAssociative} is thus completed.
\end{proof}

%% file: ANNcompositionsSpecial.tex
\subsubsection{Compositions of ANNs and affine linear transformations}

\begin{cor}\label{Lemma:PropertiesOfCompositionsWithAffineMaps}
	Let $\Phi\in\ANNs$
	(cf.\ Definition~\ref{Def:ANN}).
	Then
	\begin{enumerate}[(i)]
		\item \label{PropertiesOfCompositionsWithAffineMaps:ItemFront}
		it holds for all $\affineMap\in\ANNs$ with $\lengthANN(\affineMap)=1$ and  $\inDimANN(\affineMap)=\outDimANN(\Phi)$ that
		\begin{equation}
		\paramANN(\compANN{\affineMap}{\Phi})\le   \left[\max\!\left\{1,\tfrac{\outDimANN(\affineMap)}{\outDimANN(\Phi)}\right\}\right] \paramANN(\Phi)
		\end{equation}
		and
		\item \label{PropertiesOfCompositionsWithAffineMaps:ItemBehind}
		it holds for all $\affineMap\in\ANNs$ with $\lengthANN(\affineMap)=1$ and  $\inDimANN(\Phi)=\outDimANN(\affineMap)$ that
		\begin{equation}
		\paramANN(\compANN{\Phi}{\affineMap})\le   \left[\max\!\left\{1,\tfrac{\inDimANN(\affineMap)+1}{\inDimANN(\Phi)+1}\right\}\right] \paramANN(\Phi)
		\end{equation}
	\end{enumerate}
		(cf.\   Definition~\ref{Definition:ANNcomposition}).
\end{cor}

\begin{proof}[Proof of Corollary~\ref{Lemma:PropertiesOfCompositionsWithAffineMaps}]	
	Throughout this proof let $L\in\N$, $l_{0},l_{1},\dots, l_L\in\N$, $ \affineMap_1,\affineMap_2\in\ANNs$ satisfy that $\lengthANN(\affineMap_1)=\lengthANN(\affineMap_2)=1$, $\inDimANN(\affineMap_1)=\outDimANN(\Phi)$, $\inDimANN(\Phi)=\outDimANN(\affineMap_2)$, and $\dims(\Phi)=(l_{0},l_{1},\dots, l_{L})$.
%Observe that the fact that 
%	$\dims(\Phi)=(l_{0},l_{1},\dots, l_{L})$,
%	 the fact that for all $k\in\{1,2\}$ it holds that $\dims(\affineMap_k)=(\inDimANN(\affineMap_k),\outDimANN(\affineMap_k))$, and item~\eqref{PropertiesOfCompositions:Dims} in Proposition~\ref{Lemma:PropertiesOfCompositions} ensure that 
%	\begin{equation}
%	\dims(\compANN{\affineMap_1}{\Phi})=(l_{0},l_{1},\dots,l_{L-1},\outDimANN(\affineMap_1))
%	\qandq 
%	\dims(\compANN{\Phi}{\affineMap_2})=(\inDimANN(\affineMap_2), l_{1},l_{2},\dots,l_{L_2}).
%	\end{equation}		
%	Hence, we obtain that
Observe that item~\eqref{PropertiesOfCompositions:Params} in Proposition~\ref{Lemma:PropertiesOfCompositions}, the fact that $\outDimANN(\Phi)=l_L$, the fact that $\inDimANN(\Phi)=l_0$,
and the fact that for all $k\in\{1,2\}$ it holds that $\dims(\affineMap_k)=(\inDimANN(\affineMap_k),\outDimANN(\affineMap_k))$ ensure that
	\begin{equation}
	\begin{split}
	&\paramANN(\compANN{\affineMap_1}{\Phi})
	=
	\left[\smallsum\limits_{m = 1}^{L-1} l_{m}(l_{m-1} + 1) \right]+\big[\outDimANN(\affineMap_1)\big] (l_{L-1} + 1)
	\\&=\left[\smallsum\limits_{m = 1}^{L-1} l_{m}(l_{m-1} + 1) \right]+\left[\tfrac{\outDimANN(\affineMap_1)}{l_{L}}\right]l_L (l_{L-1} + 1)
	\\&\le \left[\max\!\left\{1,\tfrac{\outDimANN(\affineMap_1)}{l_{L}}
	%	\big[\outDimANN(\affineMap_1)\big]/l_{L}^1
	\right\}\right] \left[\smallsum\limits_{m = 1}^{L-1} l_{m}(l_{m-1} + 1) \right]
	+\left[\max\!\left\{1,\tfrac{\outDimANN(\affineMap_1)}{l_{L}}
	%	\big[\outDimANN(\affineMap_1)\big]/l_{L}^1
	\right\}\right]l_L (l_{L-1} + 1)
	\\&=
	 \left[\max\!\left\{1,\tfrac{\outDimANN(\affineMap_1)}{l_{L}}
	\right\}\right]\left[\smallsum\limits_{m = 1}^{L} l_{m}(l_{m-1} + 1) \right]
	= \left[\max\!\left\{1,\tfrac{\outDimANN(\affineMap_1)}{\outDimANN(\Phi)}\right\}\right] \paramANN(\Phi)
	\end{split}
	\end{equation}
	and
	\begin{equation}
	\begin{split}
	&\paramANN(\compANN{\Phi}{\affineMap_2})
	=
	\left[\smallsum\limits_{m = 2}^{L} l_{m}(l_{m-1} + 1) \right]+ l_{1} \big[\inDimANN(\affineMap_2)+1\big] 
	\\&=\left[\smallsum\limits_{m = 2}^{L} l_{m}(l_{m-1} + 1) \right]+\left[\tfrac{\inDimANN(\affineMap_2)+1}{l_{0}+1}\right] l_{1} (l_0+1) 
	\\&\le \left[\max\!\left\{1,\tfrac{\inDimANN(\affineMap_2)+1}{l_{0}+1}\right\}\right] \left[\smallsum\limits_{m = 2}^{L} l_{m}(l_{m-1} + 1)\right]
	+ \left[\max\!\left\{1,\tfrac{\inDimANN(\affineMap_2)+1}{l_{0}+1}\right\}\right]  l_{1} (l_0+1) 
	\\&= \left[\max\!\left\{1,\tfrac{\inDimANN(\affineMap_2)+1}{l_{0}+1}
	%	\big[\outDimANN(\affineMap_2)\big]/l_{L}^2
	\right\}\right]\left[\smallsum\limits_{m = 1}^{L} l_{m}(l_{m-1} + 1) \right]
	= \left[\max\!\left\{1,\tfrac{\inDimANN(\affineMap_2)+1}{\inDimANN(\Phi)+1}\right\}\right] \paramANN(\Phi).
	\end{split}
	\end{equation}
	This establishes items~\eqref{PropertiesOfCompositionsWithAffineMaps:ItemFront}--\eqref{PropertiesOfCompositionsWithAffineMaps:ItemBehind}.
	The proof of Corollary~\ref{Lemma:PropertiesOfCompositionsWithAffineMaps} is thus completed.
	%Moreover, observe that 	
	%	\begin{equation}
	%\compANN{\affineMap_1}{\Phi_1}=\big((W_1^1, B_1^1),(W_2^1, B_2^1),\allowbreak \ldots, (W_{L_{k-1}}^1,\allowbreak B_{L_{k-1}}^1),(A_1 W_{L_k}^1,A_1B_{L_k}^1+b_1)\big)
	%\end{equation}
	%and 
	%		\begin{equation}
	%\compANN{\affineMap_2}{\Phi_2}=\big((W_{1}^2A_2 , W_{1}^2 b_2+B_{1}^2), (W_2^2, B_2^2),(W_3^2, B_3^2),\allowbreak \ldots, (W_{L_{k}}^2,\allowbreak B_{L_{k}}^2)\big).
	%\end{equation}
\end{proof}

\subsubsection{Powers and extensions of ANNs}
\begin{definition}\label{Definition:identityMatrix}
	Let $d\in\N$. Then we denote by $\idMatrix_{d}\in \R^{d\times d}$ the identity matrix in $\R^{d\times d}$.
\end{definition}

\begin{definition}\label{Definition:iteratedANNcomposition}
	We denote by $(\cdot)^{\bullet n}\colon \{\Phi\in \ANNs\colon \inDimANN(\Phi)=\outDimANN(\Phi)\}\allowbreak\to\ANNs$, $n\in\N_0$, the functions
%	 which satisfy for all $n\in\N$, $\Phi\in\ANNs$ with $\inDimANN(\Phi)=\outDimANN(\Phi)$ that 
%		\begin{equation}
%	\begin{split}
%	\Phi^{\bullet 0}=\big(\idMatrix_{\outDimANN(\Phi)},(0,0,\dots, 0)\big)\in\R^{\outDimANN(\Phi)\times \outDimANN(\Phi)}\times \R^{\outDimANN(\Phi)}
%	\qandq
%	\Phi^{\bullet n}=\compANN{\Phi}{\Phi^{\bullet (n-1)}}
%	\end{split}
%	\end{equation}
		 which satisfy for all $n\in\N_0$, $\Phi\in\ANNs$ with $\inDimANN(\Phi)=\outDimANN(\Phi)$ that 
	\begin{equation}\label{iteratedANNcomposition:equation}
		\begin{split}
		\Phi^{\bullet n}=
		\begin{cases} \big(\idMatrix_{\outDimANN(\Phi)},(0,0,\dots, 0)\big)\in\R^{\outDimANN(\Phi)\times \outDimANN(\Phi)}\times \R^{\outDimANN(\Phi)}
		&: n=0 \\
		\,\compANN{\Phi}{(\Phi^{\bullet (n-1)})} &: n\in\N
		\end{cases}
		\end{split}
	\end{equation}	
	(cf.\ Definition~\ref{Def:ANN},  Definition~\ref{Definition:ANNcomposition}, and Definition~\ref{Definition:identityMatrix}).
\end{definition}

\begin{definition}[Extension of ANNs]\label{Definition:ANNenlargement}
	Let $L\in\N$, $\Psi\in \ANNs$ satisfy that $\inDimANN(\Psi)=\outDimANN(\Psi)$.
	Then
	we denote by $\longerANN{L,\Psi}\colon \{\Phi\in\ANNs\colon (\lengthANN(\Phi)\le L \andShort \outDimANN(\Phi)=\inDimANN(\Psi)) \}\to \ANNs$ the function which satisfies for all $\Phi\in\ANNs$ with $\lengthANN(\Phi)\le L$ and $\outDimANN(\Phi)=\inDimANN(\Psi)$ that
%	\begin{equation}
%	\longerANN{L,\Psi}(\Phi)=	 \compANN{\Psi^{\bullet (L-\lengthANN(\Phi))}}{\Phi}
%	\end{equation}
		\begin{equation}\label{ANNenlargement:Equation}
	\longerANN{L,\Psi}(\Phi)=	 \compANN{(\Psi^{\bullet (L-\lengthANN(\Phi))})}{\Phi}
	\end{equation}
	%	\begin{equation}
	%	\longerANN{L}(\Phi)=	\begin{cases} \compANN{\idANN{\outDimANN(\Phi)}{L-\lengthANN(\Phi)+1}}{\Phi}
	%	&: \lengthANN(\Phi)< L \\
	%	\Phi &: \lengthANN(\Phi)= L
	%	\end{cases}
	%	\end{equation}
	(cf.\ Definition~\ref{Def:ANN},   Definition~\ref{Definition:ANNcomposition}, and Definition~\ref{Definition:iteratedANNcomposition}).
\end{definition}

\begin{lemma}\label{Lemma:PropertiesOfANNenlargementGeometry}
	Let  $d,\hiddenDimId\in\N$, $\Psi\in \ANNs$ satisfy  that  $\dims(\Psi)=(d,\hiddenDimId,d)$
	(cf.\ Definition~\ref{Def:ANN}).
	Then
	\begin{enumerate}[(i)]
	\item \label{PropertiesOfANNenlargementGeometry:BulletPower} it holds for all $n\in\N_0$ that
	$\lengthANN(\Psi^{\bullet n})=n+1$,  $\dims(\Psi^{\bullet n})\in \N^{n+2}$, and 
%	\begin{equation}\label{BulletPower:Dimensions}
%	\lengthANN(\Psi^{\bullet n})=n+1
%	\qandq \dims(\Psi^{\bullet n}) = (d,\hiddenDimId,\hiddenDimId,\dots,\hiddenDimId,d)\in\N^{n+2}
%	\end{equation}
	\begin{equation}\label{BulletPower:Dimensions}
%	\lengthANN(\Psi^{\bullet n})=n+1
%	\qandq 
	\dims(\Psi^{\bullet n}) = 
	\begin{cases}
	(d,d) &: n=0\\
	(d,\hiddenDimId,\hiddenDimId,\dots,\hiddenDimId,d) &: n\in\N
	\end{cases}
	\end{equation}
	 and
	\item \label{PropertiesOfANNenlargementGeometry:ItemLonger}
	it holds 
	for all $\Phi\in\ANNs$,
% $L\in\{\lengthANN(\Phi),\lengthANN(\Phi)+1,\dots\}$
	$L\in\N\cap [\lengthANN(\Phi),\infty)$ with $\outDimANN(\Phi)=d$
	that $\lengthANN\big(\longerANN{L,\Psi}(\Phi)\big)=L$ and
	\begin{equation}\label{PropertiesOfANNenlargementGeometry:ParamsLonger}
	\begin{split}
	&\paramANN(\longerANN{L,\Psi}(\Phi))
%	\le 
%\paramANN(\Phi)\, \indicator{\{\lengthANN(\Phi)\}}(L)
%\\&+\left[\big(\!\max\!\big\{1,\tfrac{\hiddenDimId}{d}\big\}\big)\paramANN(\Phi)+
%(L-\lengthANN(\Phi)-1) \,\hiddenDimId\,(\hiddenDimId+1)
%+d\,(\hiddenDimId+1)\right]
%\indicator{(\lengthANN(\Phi),\infty)}(L)
\\&\le 
\begin{cases} \paramANN(\Phi)
		&: \lengthANN(\Phi)=L \\
		\left[\big(\!\max\!\big\{1,\tfrac{\hiddenDimId}{d}\big\}\big)\paramANN(\Phi)+
		\big((L-\lengthANN(\Phi)-1) \,\hiddenDimId+d\big)(\hiddenDimId+1)
		\right] &: \lengthANN(\Phi)<L
		\end{cases}
	\end{split}
	\end{equation}
\end{enumerate}
	(cf.\  Definition~\ref{Definition:iteratedANNcomposition} and Definition~\ref{Definition:ANNenlargement}).
\end{lemma}

\begin{proof}[Proof of Lemma~\ref{Lemma:PropertiesOfANNenlargementGeometry}]
	Throughout this proof let  $\Phi\in\ANNs$, $l_0, l_1,\dots, l_{\lengthANN(\Phi)}\in\N$ satisfy that  $\outDimANN(\Phi)=d$ and
	$\dims(\Phi)=(l_0,l_1,\dots, l_{\lengthANN(\Phi)})\in \N^{\lengthANN(\Phi)+1}$ 
	and let $a_{L,k}\in \N$,  $k\in \N_0\cap [0,L]$, $L\in \N\cap [\lengthANN(\Phi),\infty)$, satisfy 
%	for all $L\in \N\cap (\lengthANN(\Phi),\infty)$, $k\in \{0,1,\dots,L\}$ that
%	\begin{equation}\label{PropertiesOfANNenlargementGeometry:extendedDimensions}
%	a_{L,k}=	\begin{cases} l_k
%	&: k\in\{0,1,\dots,\lengthANN(\Phi)-1\} \\
%	\hiddenDimId &: k\in \{\lengthANN(\Phi),\lengthANN(\Phi)+1,\dots, L-1\}\\
%	l_{\lengthANN(\Phi)} &: k= L
%	\end{cases}.
%	\end{equation}
	for all $L\in \N\cap [\lengthANN(\Phi),\infty)$, $k\in \N_0\cap [0,L]$ that
	\begin{equation}\label{PropertiesOfANNenlargementGeometry:extendedDimensions}
	 a_{L,k}=	
	 \begin{cases} 
	 l_k &: k<\lengthANN(\Phi) \\
	\hiddenDimId &: \lengthANN(\Phi)\le k<L\\
	d &: k= L
	\end{cases}.
	\end{equation}
	We claim that for all $n\in\N_0$ it holds that 
		\begin{equation}\label{BulletPower:DimensionsProof}
		\lengthANN(\Psi^{\bullet n})=n+1
		\qandq \N^{n+2}\ni \dims(\Psi^{\bullet n}) = 
		\begin{cases}
		(d,d) &: n=0\\
		(d,\hiddenDimId,\hiddenDimId,\dots,\hiddenDimId,d) &: n\in\N
		\end{cases}.
		\end{equation}
	We now prove \eqref{BulletPower:DimensionsProof}  by induction on $n\in\N_0$.
Note that the fact  that $\Psi^{\bullet 0}=(\idMatrix_{d},0)\in \R^{d\times d}\times \R^d$ (cf.\ Definition~\ref{Definition:identityMatrix})  
establishes \eqref{BulletPower:Dimensions} in the base case $n=0$.
For the induction step $\N_0\ni n\to n+1\in \N$ assume that there exists $n\in\N_0$ such that 
\begin{equation}\label{PropertiesOfANNenlargementGeometry:DimensionsInduction}
\lengthANN(\Psi^{\bullet n})=n+1
\qandq 
\N^{n+2}\ni\dims(\Psi^{\bullet n}) = 
\begin{cases}
(d,d) &: n=0\\
(d,\hiddenDimId,\hiddenDimId,\dots,\hiddenDimId,d) &: n\in\N
\end{cases}.
%\dims(\Psi^{\bullet n}) = (d,\hiddenDimId,\hiddenDimId,\dots,\hiddenDimId,d)\in\N^{n+2}.
%	\qandqShort (\functionANN(\mathbb{I}^{\bullet n}))(x)=x
\end{equation}
Observe that Lemma~\ref{Lemma:elementaryPropertiesANN}, 
\eqref{iteratedANNcomposition:equation}, items~\eqref{PropertiesOfCompositions:Dims}--\eqref{PropertiesOfCompositions:Length} in Proposition~\ref{Lemma:PropertiesOfCompositions}, \eqref{PropertiesOfANNenlargementGeometry:DimensionsInduction}, and the hypothesis that $\dims(\Psi) = (d,\hiddenDimId,d)$ imply that 
\begin{equation}\label{BulletPower:DimensionsInduction2}
\begin{split}
\lengthANN(\Psi^{\bullet (n+1)})&=\lengthANN(\compANN{\Psi}{(\Psi^{\bullet n})})=\lengthANN(\Psi)+\lengthANN(\Psi^{\bullet n})-1=2+(n+1)-1=(n+1)+1
\\&\andq \dims(\Psi^{\bullet (n+1)})=\dims(\compANN{\Psi}{(\Psi^{\bullet n})}) = (d,\hiddenDimId,\hiddenDimId,\dots,\hiddenDimId,d)\in \N^{n+3}.
\end{split}
\end{equation}
Induction thus proves \eqref{BulletPower:DimensionsProof}.
Next note that \eqref{BulletPower:DimensionsProof} establishes item~\eqref{PropertiesOfANNenlargementGeometry:BulletPower}.
In addition, observe that 
	items~\eqref{PropertiesOfCompositions:Dims}--\eqref{PropertiesOfCompositions:Length} in
	Proposition~\ref{Lemma:PropertiesOfCompositions}, item~\eqref{PropertiesOfANNenlargementGeometry:BulletPower},  \eqref{ANNenlargement:Equation}, and \eqref{PropertiesOfANNenlargementGeometry:extendedDimensions} ensure that
	for all $L\in \N\cap [\lengthANN(\Phi),\infty)$ it holds 
	that
	\begin{equation}\label{PropertiesOfANNenlargementGeometryLengthLonger}
	\begin{split}
	\lengthANN\big(\longerANN{L,\Psi}(\Phi)\big)
	&=	 	\lengthANN\big(\compANN{(\Psi^{\bullet (L-\lengthANN(\Phi))})}{\Phi}\big)
	=\lengthANN\big( \Psi^{\bullet (L-\lengthANN(\Phi))} \big)+\lengthANN(\Phi)-1
	\\&=(L-\lengthANN(\Phi)+1)+\lengthANN(\Phi)-1
	=L
	\end{split}
	\end{equation}
	and
	\begin{equation}\label{PropertiesOfANNenlargementGeometryDimsLonger}
	\begin{split}
	 \dims\big(\longerANN{L,\Psi}(\Phi)\big)
	&=\dims\big(\compANN{(\Psi^{\bullet (L-\lengthANN(\Phi))})}{\Phi}\big)
	=(a_{L,0},a_{L,1},\dots, a_{L,L}).
	%	 \\&=\big(\inDimANN(\Phi),l_1,\dots, l_{\lengthANN(\Phi)-1},\hiddenDimId,\hiddenDimId,\dots, \hiddenDimId,\outDimANN(\Phi)\big).
	\end{split}
	\end{equation}
Combining this with \eqref{PropertiesOfANNenlargementGeometry:extendedDimensions} demonstrates that 
\begin{equation}
	\lengthANN\big(\longerANN{\lengthANN(\Phi),\Psi}(\Phi)\big)
	=\lengthANN(\Phi)
\end{equation}
%$\lengthANN\big(\longerANN{\lengthANN(\Phi),\Psi}(\Phi)\big)
%=\lengthANN(\Phi)$
	and
	\begin{equation}
	\begin{split}
	\dims\big(\longerANN{\lengthANN(\Phi),\Psi}(\Phi)\big)
	&=(a_{\lengthANN(\Phi),0},a_{\lengthANN(\Phi),1},\dots, a_{\lengthANN(\Phi),\lengthANN(\Phi)})
	\\&=(l_0,l_1,\dots, l_{\lengthANN(\Phi)})=\dims(\Phi).
	\end{split}
	\end{equation}
	Hence, we obtain that 
	\begin{equation}\label{PropertiesOfANNenlargementGeometry:NotLongerParams}
		\paramANN\big(\longerANN{\lengthANN(\Phi),\Psi}(\Phi)\big)=\paramANN(\Phi).
	\end{equation}
%	
%	$\paramANN\big(\longerANN{\lengthANN(\Phi),\Psi}(\Phi)\big)=\paramANN(\Phi)$.
Next note that \eqref{PropertiesOfANNenlargementGeometry:extendedDimensions}, \eqref{PropertiesOfANNenlargementGeometryDimsLonger}, and the fact that $l_{\lengthANN(\Phi)}=\outDimANN(\Phi)=d$ imply that for all $L\in\N\cap (\lengthANN(\Phi),\infty)$ it holds that 
\begin{equation}
	\begin{split}
&\paramANN\big(\longerANN{L,\Psi}(\Phi)\big)
=\smallsum\limits_{k = 1}^{L} a_{L,k}(a_{L,k-1} + 1)
\\&=\left[\smallsum\limits_{k = 1}^{\lengthANN(\Phi)-1} a_{L,k}(a_{L,k-1} + 1)\right]+\left[\smallsum\limits_{k = \lengthANN(\Phi)}^{L} a_{L,k}(a_{L,k-1} + 1)\right]
\\&=\left[\smallsum\limits_{k = 1}^{\lengthANN(\Phi)-1} l_{k}(l_{k-1} + 1)\right]+\left[\smallsum\limits_{k = \lengthANN(\Phi)}^{\lengthANN(\Phi)} a_{L,k}(a_{L,k-1} + 1)\right]
\\&\quad+\left[\smallsum\limits_{k = \lengthANN(\Phi)+1}^{L} a_{L,k}(a_{L,k-1} + 1)\right]
\\&=\left[\smallsum\limits_{k = 1}^{\lengthANN(\Phi)-1} l_{k}(l_{k-1} + 1)\right]+ a_{L,\lengthANN(\Phi)}(a_{L,\lengthANN(\Phi)-1} + 1)
\\&\quad+\left[\smallsum\limits_{k = \lengthANN(\Phi)+1}^{L-1} a_{L,k}(a_{L,k-1} + 1)\right]+\left[\smallsum\limits_{k = L}^{L} a_{L,k}(a_{L,k-1} + 1)\right]
\\&=\left[\smallsum\limits_{k = 1}^{\lengthANN(\Phi)-1} l_{k}(l_{k-1} + 1)\right]+ \hiddenDimId(l_{\lengthANN(\Phi)-1} + 1)
\\&\quad+ \big(L-1-(\lengthANN(\Phi)+1)+1\big)\hiddenDimId(\hiddenDimId+1)
+ a_{L,L}(a_{L,L-1} + 1)
\\&=\left[\smallsum\limits_{k = 1}^{\lengthANN(\Phi)-1} l_{k}(l_{k-1} + 1)\right]+ \tfrac{\hiddenDimId}{d} \big[l_{\lengthANN(\Phi)}(l_{\lengthANN(\Phi)-1} + 1)\big]
\\&\quad+ \big(L-\lengthANN(\Phi)-1\big)\hiddenDimId(\hiddenDimId+1)
+ d(\hiddenDimId+1)
\\&\le \left[\max\{1,\tfrac{\hiddenDimId}{d}\}\right]\left[\smallsum\limits_{k = 1}^{\lengthANN(\Phi)} l_{k}(l_{k-1} + 1)\right]
+ \big(L-\lengthANN(\Phi)-1\big)\hiddenDimId(\hiddenDimId+1)
+ d(\hiddenDimId+1)
\\&=\left[\max\{1,\tfrac{\hiddenDimId}{d}\}\right]\paramANN(\Phi)
+ \big(L-\lengthANN(\Phi)-1\big)\hiddenDimId(\hiddenDimId+1)
+ d(\hiddenDimId+1).
	\end{split}
\end{equation}
%
%
%
%
%	 This proves \eqref{PropertiesOfANNenlargementGeometry:ParamsLonger} in the case $\lengthANN(\Phi)=L$.
%	We now prove \eqref{PropertiesOfANNenlargementGeometry:ParamsLonger} in the case $\lengthANN(\Phi)<L$.
%Observe that the fact that $l_{\lengthANN(\Phi)}=\outDimANN(\Phi)=d$, \eqref{PropertiesOfANNenlargementGeometry:extendedDimensions}, and 
%\eqref{PropertiesOfANNenlargementGeometryDimsLonger}
% imply that
%	\begin{equation}\label{PropertiesOfANNenlargementGeometryParametersLonger}
%	\begin{split}
%	\paramANN(\longerANN{L,\Psi}(\Phi))
%	&= \bigg[\smallsum\limits_{m = 1}^{\lengthANN(\Phi)-1} l_m(l_{m-1} + 1) \bigg]+\hiddenDimId\, (l_{\lengthANN(\Phi)-1}+1)
%	\\&\quad+ 
%	(L-\lengthANN(\Phi)-1) \,\hiddenDimId\,(\hiddenDimId+1)
%	+d\,(\hiddenDimId+1)
%	\\&\le  \big(\!\max\!\big\{1,\tfrac{\hiddenDimId}{d}\big\}\big)\,\paramANN(\Phi)+
%	(L-\lengthANN(\Phi)-1) \,\hiddenDimId\,(\hiddenDimId+1)
%	+d\,(\hiddenDimId+1) .
%	\end{split}
%	\end{equation}
%This establishes \eqref{PropertiesOfANNenlargementGeometry:ParamsLonger} in the case $\lengthANN(\Phi)<L$.
Combining this with \eqref{PropertiesOfANNenlargementGeometry:NotLongerParams} establishes \eqref{PropertiesOfANNenlargementGeometry:ParamsLonger}.
	The proof of Lemma~\ref{Lemma:PropertiesOfANNenlargementGeometry} is thus completed.
\end{proof}

\begin{lemma}\label{Lemma:PropertiesOfANNenlargementRealization}
	Let $\activation\in C(\R,\R)$,  $ \mathbb{I}\in \ANNs$ satisfy for all  $x\in\R^{\inDimANN(\mathbb{I})}$ that $\inDimANN(\mathbb{I})=\outDimANN(\mathbb{I})$ 
%	$\functionANN( \mathbb{I})\in C(\R^d,\R^d)$
	and $(\functionANN( \mathbb{I}))(x)=x$
	(cf.\ Definition~\ref{Def:ANN} and Definition~\ref{Definition:ANNrealization}).
	Then
	\begin{enumerate}[(i)]
		\item \label{PropertiesOfANNenlargementRealization: itemOne} it holds for all $n\in\N_0$, $x\in\R^{\inDimANN(\mathbb{I})}$ that
		\begin{equation}\label{BulletPowerRealization:Dimensions}
			\functionANN(\mathbb{I}^{\bullet n})\in C(\R^{\inDimANN(\mathbb{I})},\R^{\inDimANN(\mathbb{I})}) \qandq (\functionANN(\mathbb{I}^{\bullet n}))(x)=x
		\end{equation}
%		$\functionANN(\mathbb{I}^{\bullet n})\in C(\R^d,\R^d)$
%		and $(\functionANN(\mathbb{I}^{\bullet n}))(x)=x$ 	
		 and	
		\item 	\label{PropertiesOfANNenlargementRealization:ItemIdentityLonger}
			it holds 
		for all $\Phi\in\ANNs$,
		% $L\in\{\lengthANN(\Phi),\lengthANN(\Phi)+1,\dots\}$
		$L\in\N\cap [\lengthANN(\Phi),\infty)$, $x\in\R^{\inDimANN(\Phi)}$ with $\outDimANN(\Phi)={\inDimANN(\mathbb{I})}$
		that	
%		$\functionANN(\longerANN{L,\mathbb{I}}(\Phi))\in C(\R^{\inDimANN(\Phi)},\R^{\outDimANN(\Phi)})$ and
		\begin{equation}
		\functionANN(\longerANN{L,\mathbb{I}}(\Phi))\in C(\R^{\inDimANN(\Phi)},\R^{\outDimANN(\Phi)})
		\qandqShort
		\big(\functionANN(\longerANN{L,\mathbb{I}}(\Phi))\big)(x)=\big(\functionANN(\Phi)\big)(x)
		\end{equation}
	\end{enumerate}
(cf.\  Definition~\ref{Definition:iteratedANNcomposition} and Definition~\ref{Definition:ANNenlargement}).
\end{lemma}

\begin{proof}[Proof of Lemma~\ref{Lemma:PropertiesOfANNenlargementRealization}]
	Throughout this proof let  $\Phi\in\ANNs$, $L,d\in\N$ satisfy that $\lengthANN(\Phi)\le L$ and $\inDimANN(\mathbb{I})=\outDimANN(\Phi)=d$.
	We claim that for all $n\in\N_0$ it holds that 
			\begin{equation}\label{BulletPowerRealization:DimensionsProof}
			\functionANN(\mathbb{I}^{\bullet n})\in C(\R^{d},\R^{d}) \qandq \forall\, x\in\R^d\colon\,(\functionANN(\mathbb{I}^{\bullet n}))(x)=x.
			\end{equation}
		We now prove \eqref{BulletPowerRealization:DimensionsProof}  by induction on $n\in\N_0$.
		Note that \eqref{iteratedANNcomposition:equation} and the fact that $\outDimANN(\mathbb{I})=d$  demonstrate that
		$\functionANN(\mathbb{I}^{\bullet 0})\in C(\R^d,\R^d)$
		and $\forall\,x\in\R^d\colon(\functionANN(\mathbb{I}^{\bullet 0}))(x)=x$.
		This 
		establishes \eqref{BulletPowerRealization:DimensionsProof} in the base case $n=0$.
		For the induction step observe that for all $n\in\N_0$ with $\functionANN(\mathbb{I}^{\bullet n})\in C(\R^d,\R^d)$
		and $\forall\,x\in\R^d\colon(\functionANN(\mathbb{I}^{\bullet n}))(x)=x$ it holds that 
		\begin{equation}
			\functionANN(\mathbb{I}^{\bullet (n+1)})=\functionANN(\compANN{\mathbb{I}}{(\mathbb{I}^{\bullet n})})
			=(\functionANN(\mathbb{I}))\circ (\functionANN(\mathbb{I}^{\bullet n}))\in C(\R^d,\R^d)
		\end{equation}
		and
		\begin{equation}
		\begin{split}
				\forall\, x\in\R^d\colon\, \big(\functionANN(\mathbb{I}^{\bullet (n+1)})\big)(x)&=\big([\functionANN(\mathbb{I})]\circ [\functionANN(\mathbb{I}^{\bullet n})]\big)(x)
				\\&=(\functionANN(\mathbb{I}))\big(\big(\functionANN(\mathbb{I}^{\bullet n})\big)(x)\big)
				=(\functionANN(\mathbb{I}))(x)=x.
		\end{split}
		\end{equation}
		Induction thus proves \eqref{BulletPowerRealization:DimensionsProof}. Next observe that \eqref{BulletPowerRealization:DimensionsProof} establishes 
		item~\eqref{PropertiesOfANNenlargementRealization: itemOne}.
	Moreover, note that \eqref{ANNenlargement:Equation},
	 item~\eqref{PropertiesOfCompositions:Realization} in Proposition~\ref{Lemma:PropertiesOfCompositions}, item~\eqref{PropertiesOfANNenlargementRealization: itemOne},
	 and the fact that $\inDimANN(\mathbb{I})=\outDimANN(\Phi)$
	  ensure that 
%	  for all $x\in\R^{\inDimANN(\Phi)}$ it holds that
	  \begin{equation}
	  \begin{split}
	  	  	&\functionANN(\longerANN{L,\mathbb{I}}(\Phi))=	\functionANN( \compANN{(\mathbb{I}^{\bullet (L-\lengthANN(\Phi))})}{\Phi})
	  	  	\\&\in C(\R^{\inDimANN(\Phi)},\R^{\outDimANN(\mathbb{I})})=C(\R^{\inDimANN(\Phi)},\R^{\inDimANN(\mathbb{I})})=C(\R^{\inDimANN(\Phi)},\R^{\outDimANN(\Phi)})
	  \end{split}
	  \end{equation}
and
	\begin{equation}
	\begin{split}
		\forall\,x\in\R^{\inDimANN(\Phi)}\colon\,	\big(\functionANN(\longerANN{L,\mathbb{I}}(\Phi))\big)(x)&=\big(\functionANN(\mathbb{I}^{\bullet (L-\lengthANN(\Phi))})\big)\big((\functionANN(\Phi))(x)\big)\\&=(\functionANN(\Phi))(x).
	\end{split}
%	\forall\,x\in\R^{\inDimANN(\Phi)}\colon\,	\big(\functionANN(\longerANN{L,\mathbb{I}}(\Phi))\big)(x)=\big(\functionANN(\mathbb{I}^{\bullet (L-\lengthANN(\Phi))})\big)\big((\functionANN(\Phi))(x)\big)=(\functionANN(\Phi))(x).
	\end{equation}
	This establishes item~\eqref{PropertiesOfANNenlargementRealization:ItemIdentityLonger}.
	The proof of Lemma~\ref{Lemma:PropertiesOfANNenlargementRealization} is thus completed.
\end{proof}

\subsubsection{Compositions of ANNs involving artificial identities}

\begin{definition}[Composition of ANNs involving artificial identities]
	\label{Definition:ANNconcatenation}
	Let $\Psi\in \ANNs$.
	Then	
	we denote by 
	\begin{equation}
		\concPsiANN{(\cdot)}{(\cdot)}\colon \{(\Phi_1,\Phi_2)\in\ANNs\times \ANNs\colon \inDimANN(\Phi_1)=\outDimANN(\Psi) \text{ and } \outDimANN(\Phi_2)=\inDimANN(\Psi)\}\allowbreak\to\ANNs
	\end{equation}
	the function which satisfies for all $\Phi_1,\Phi_2\in\ANNs$ with $\inDimANN(\Phi_1)=\outDimANN(\Psi)$ and $\outDimANN(\Phi_2)=\inDimANN(\Psi)$ that 
	\begin{equation}
	\begin{split}
	\concPsiANN{\Phi_1}{\Phi_2}
	= \compANN{\Phi_1}
	{
		(\compANN
		{\Psi}
		{\Phi_2})}
	=\compANN{(\compANN{\Phi_1}{\Psi})}{\Phi_2}
	\end{split}
	\end{equation}
	(cf.\ Definition~\ref{Def:ANN}, Definition~\ref{Definition:ANNcomposition}, and  Lemma~\ref{Lemma:CompositionAssociative}).
\end{definition}

%Definition with concatenation with more special \Psi

%\begin{definition}\label{Definition:ANNconcatenation}
%		Let $\Psi\in \ANNs$ satisfy that $\inDimANN(\Psi)=\outDimANN(\Psi)$.
%	Then	
%	we denote by $\concPsiANN{(\cdot)}{(\cdot)}\colon \{(\Phi_1,\Phi_2)\in\ANNs\times \ANNs\colon \inDimANN(\Phi_1)=\outDimANN(\Phi_2)=\inDimANN(\Psi)\}\allowbreak\to\ANNs$ the function which satisfies for all $\Phi_1,\Phi_2\in\ANNs$ with $\inDimANN(\Phi_1)=\outDimANN(\Phi_2)=\inDimANN(\Psi)$ that 
%	\begin{equation}
%	\begin{split}
%	\concPsiANN{\Phi_1}{\Phi_2}
%	= \compANN{\Phi_1}
%	{
%		(\compANN
%			{\Psi}
%				{\Phi_2})}
%	\end{split}
%	\end{equation}
%	(cf.\ Definition~\ref{Def:ANN}, Definition~\ref{Definition:ANNrealization}, and Definition~\ref{Definition:ANNcomposition}).
%\end{definition}

%Definition of concWith family

%\begin{definition}\label{Definition:ANNconcatenation}
%	Let $L\in\N$, let $\varPsi=\Psi\subseteq \ANNs$ satisfy for all $d\in\N$ that $\inDimANN(\Psi_d)=\outDimANN(\Psi_d)=d$.
%	Then	
%	we denote by $\concPsiANN{(\cdot)}{(\cdot)}\colon \{(\Phi_1,\Phi_2)\in\ANNs\times \ANNs\colon \inDimANN(\Phi_1)=\outDimANN(\Phi_2)\}\allowbreak\to\ANNs$ the function which satisfies for all $\Phi_1,\Phi_2\in\ANNs$ with $\inDimANN(\Phi_1)=\outDimANN(\Phi_2)$ that 
%	\begin{equation}
%	\begin{split}
%	\concPsiANN{\Phi_1}{\Phi_2}
%	= \compANN{\Phi_1}
%	{
%		(\compANN
%		{\Psi_{\outDimANN(\Phi_2)}}
%		{\Phi_2})}
%	\end{split}
%	\end{equation}
%	(cf.\ Definition~\ref{Def:ANN}, Definition~\ref{Definition:ANNrealization}, and Definition~\ref{Definition:ANNcomposition}).
%\end{definition}

\begin{prop}\label{Lemma:PropertiesOfConcatenations}
	Let   $\Psi, \Phi_1,\Phi_2\in \ANNs$,  $\hiddenDimId,l_{1,0},l_{1,1},\dots,\allowbreak l_{1,\lengthANN(\Phi_1)},l_{2,0},\allowbreak l_{2,1},\dots, l_{2,\lengthANN(\Phi_2)}\in\N$ satisfy for all $k\in\{1,2\}$ that $\dims(\Psi)=(\inDimANN(\Psi),\hiddenDimId,\outDimANN(\Psi))$, $\inDimANN(\Phi_1)=\outDimANN(\Psi)$, $\outDimANN(\Phi_2)=\inDimANN(\Psi)$,
	and
	$\dims(\Phi_k)=(l_{k,0},l_{k,1},\dots, l_{k,\lengthANN(\Phi_k)})$
	(cf.\ Definition~\ref{Def:ANN}).
	Then
	\begin{enumerate}[(i)]
		\item \label{PropertiesOfConcatenations:Dims}
		it holds that	\begin{equation}
		\dims(\concPsiANN{\Phi_1}{\Phi_2})=(l_{2,0},l_{2,1},\dots, l_{2,\lengthANN(\Phi_2)-1},\hiddenDimId,l_{1,1},l_{1,2},\allowbreak\dots,l_{1,\lengthANN(\Phi_1)} ),
		\end{equation}
		\item \label{PropertiesOfConcatenations:Length} it holds that 
		\begin{equation}
		\lengthANN(\concPsiANN{\Phi_1}{\Phi_2})=\lengthANN(\Phi_1)+\lengthANN(\Phi_2),
		\end{equation}
		\item \label{PropertiesOfConcatenations:Params}
		it holds that 
		\begin{equation}
		\paramANN(\concPsiANN{\Phi_1}{\Phi_2})\le \left[\max\!\left\{1,\tfrac{\hiddenDimId}{\inDimANN(\Psi)}, \tfrac{\hiddenDimId}{{\outDimANN(\Psi)}}\right\}\right] \big(\paramANN(\Phi_1)+\paramANN(\Phi_2)\big),
		\end{equation} 
		and
		\item \label{PropertiesOfConcatenations:Realization}
		it holds for all $a\in C(\R,\R)$ that 
		$\functionANN(\concPsiANN{\Phi_1}{\Phi_2})\in C(\R^{\inDimANN(\Phi_2)},\R^{\outDimANN(\Phi_1)})$ and
		\begin{equation}
		\functionANN(\concPsiANN{\Phi_1}{\Phi_2})=[\functionANN(\Phi_1)]\circ [\functionANN(\Psi)]\circ [\functionANN(\Phi_2)]
		\end{equation}
	\end{enumerate}
		(cf.\   Definition~\ref{Definition:ANNrealization} and Definition~\ref{Definition:ANNconcatenation}).
\end{prop}

\begin{proof}[Proof of Propositions~\ref{Lemma:PropertiesOfConcatenations}]
	Throughout this proof let $a\in C(\R,\R)$, $L_1,L_2\in\N$ satisfy that $L_1=\lengthANN(\Phi_1)$ and $L_2=\lengthANN(\Phi_2)$.	
	Note that item~\eqref{PropertiesOfCompositions:Dims} in Proposition~\ref{Lemma:PropertiesOfCompositions}, the hypothesis that $\dims(\Phi_2)=(l_{2,0},l_{2,1},\dots,\allowbreak l_{2,L_2})$, the hypothesis that
	$\dims(\Psi) = (\inDimANN(\Psi),\hiddenDimId,\outDimANN(\Psi))$,
	and the hypothesis that $\inDimANN(\Psi)=\outDimANN(\Phi_2)$
	show that 
	\begin{equation}
	\dims(\compANN{\Psi}{\Phi_2})=(l_{2,0},l_{2,1},\dots, l_{2,L_2-1},\hiddenDimId,\outDimANN(\Psi))
	\end{equation}
	(cf.\ Definition~\ref{Definition:ANNcomposition}).
	Combining this with item~\eqref{PropertiesOfCompositions:Dims} in Proposition~\ref{Lemma:PropertiesOfCompositions}, the hypothesis that $\dims(\Phi_1)=(l_{1,0},l_{1,1},\dots, l_{1,L_1})$, 
	and the hypothesis that $\inDimANN(\Phi_1)=\outDimANN(\Psi)$ proves that 
	\begin{equation}\label{PropertiesOfConcatenations:ConcDims}
	\begin{split}
	\dims(\concPsiANN{\Phi_1}{\Phi_2})=
	\dims\big(\compANN{\Phi_1}{(\compANN{\Psi}{\Phi_2})}\big)=
	(l_{2,0},l_{2,1},\dots,l_{2,L_2-1}, \hiddenDimId,l_{1,1},l_{1,2},\dots, l_{1,L_1}).
	\end{split}
	\end{equation}
	This establishes item~\eqref{PropertiesOfConcatenations:Dims}.
	%	Next note that 
	%\begin{equation}\label{PropertiesOfConcatenations:ConcMatrices}
	%\begin{split}
	%\concPsiANN{\Phi_1}{\Phi_2}&= \big((W_1^2, B_1^2),(W_2^2, B_2^2),\allowbreak \ldots, (W_{L_2-1}^2,\allowbreak B_{L_2-1}^2),
	%\\&\qquad(W_{3,1} W_{L_2}^2,\allowbreak W_{3,1} B_{L_2}^2+B_{3,1}),
	%(W_{1}^1 W_{3,2},\allowbreak B_{1}^1+ W_{1}^1 B_{3,2}),
	%\\&\qquad (W_2^1, B_2^1),(W_3^1, B_3^1),\allowbreak \ldots, (W_{L_1}^1,B_{L_1}^1)
	%\big).
	%\end{split}
	%\end{equation}
	Moreover, observe that item~\eqref{PropertiesOfCompositions:Length} in Proposition~\ref{Lemma:PropertiesOfCompositions} and the fact that  $\lengthANN(\Psi)=2$  ensure that
	\begin{equation}\label{PropertiesOfConcatenations:ConcLength}
	\begin{split}
	\lengthANN(\concPsiANN{\Phi_1}{\Phi_2})&=
	\lengthANN\big(\compANN{\Phi_1}{(\compANN{\Psi}{\Phi_2})}\big)
	= \lengthANN({\Phi_1})+\lengthANN(\compANN{\Psi}{\Phi_2})-1
	\\&=\lengthANN({\Phi_1})+\lengthANN(\Psi)+\lengthANN({\Phi_2})-2=\lengthANN({\Phi_1})+\lengthANN({\Phi_2}).
	\end{split}
	\end{equation}	
	This establishes item~\eqref{PropertiesOfConcatenations:Length}.
	In addition, observe that \eqref{PropertiesOfConcatenations:ConcDims}, the fact that $\inDimANN(\Psi)=\outDimANN(\Phi_2)=l_{2,L_2}$, and the fact that $\outDimANN(\Psi)=\inDimANN(\Phi_1)=l_{1,0}$ demonstrate that
	\begin{equation}
	\begin{split}
	\paramANN(\concPsiANN{\Phi_1}{\Phi_2})
	&=
	\left[\smallsum\limits_{m = 1}^{L_2-1} l_{2,m}(l_{2,m-1} + 1) \right]+
	\left[\smallsum\limits_{m = 2}^{L_1} l_{1,m}(l_{1,m-1} + 1) \right]
	\\&\quad+\hiddenDimId\big(l_{2,L_2-1} + 1\big)
	+l_{1,1} (\hiddenDimId+1)
	\\&=
	\left[\smallsum\limits_{m = 1}^{L_2-1} l_{2,m}(l_{2,m-1} + 1) \right]+
	\left[\smallsum\limits_{m = 2}^{L_1} l_{1,m}(l_{1,m-1} + 1) \right]
	\\&\quad+\tfrac{\hiddenDimId}{\inDimANN(\Psi)}\, l_{2,L_2} \big(l_{2,L_2-1} + 1\big)
	+l_{1,1} \big(\tfrac{\hiddenDimId}{\outDimANN(\Psi)}\,l_{1,0}+1\big)
	\\&\le 
	\left[\max\!\left\{1,\tfrac{\hiddenDimId}{\inDimANN(\Psi)}\right\}\right]
	\left[\smallsum\limits_{m = 1}^{L_2} l_{2,m}(l_{2,m-1} + 1) \right]
	\\&\quad+\left[\max\!\left\{1,\tfrac{\hiddenDimId}{\outDimANN(\Psi)}\right\}\right]
	\left[\smallsum\limits_{m = 1}^{L_1} l_{1,m}(l_{1,m-1} + 1) \right]
	\\&
	\le\left[\max\!\left\{1,\tfrac{\hiddenDimId}{\inDimANN(\Psi)},\tfrac{\hiddenDimId}{\outDimANN(\Psi)}\right\}\right] \big(\paramANN(\Phi_1)+\paramANN(\Phi_2)\big).
	\end{split}
	\end{equation}
	%	This, \eqref{PropertiesOfConcatenations:ConcFunction},  \eqref{PropertiesOfConcatenations:ConcLength}, and \eqref{PropertiesOfConcatenations:ConcDims} establish item~\eqref{PropertiesOfConcatenations:ItemConc}.
	This establishes item~\eqref{PropertiesOfConcatenations:Params}.
	Next	note that item~\eqref{PropertiesOfCompositions:Realization} in Proposition~\ref{Lemma:PropertiesOfCompositions} implies that
%		$\functionANN(\concPsiANN{\Phi_1}{\Phi_2})\in C(\R^{\inDimANN(\Phi_2)},\R^{\outDimANN(\Phi_1)})$ and
	\begin{equation}\label{PropertiesOfConcatenations:ConcFunction}
	\begin{split}
	\functionANN(\concPsiANN{\Phi_1}{\Phi_2})
	&=	\functionANN\big(\compANN{\Phi_1}{(\compANN{\Psi}{\Phi_2})}\big)
	\\&=\big[\functionANN(\Phi_1)\big]\circ \big[\functionANN(\compANN{\Psi}{\Phi_2})\big]
	\\&=\big(\big[\functionANN(\Phi_1)\big]\circ \big[\functionANN(\Psi)\big]\circ
	\big[\functionANN({\Phi_2})\big]\big)
	\in C(\R^{\inDimANN(\Phi_2)},\R^{\outDimANN(\Phi_1)}).
	%	\\&=\big[\functionANN(\Phi_1)\big]\circ
	%	\big[\functionANN({\Phi_2})\big]
	%	=\functionANN(\compANN{\Phi_1}{\Phi_2}).
	\end{split}
	\end{equation}	
	This establishes item~\eqref{PropertiesOfConcatenations:Realization}.
	The proof of Proposition~\ref{Lemma:PropertiesOfConcatenations} is thus completed.
\end{proof}

%% file: ANNparallelizations.tex
\subsection{Parallelizations of ANNs}
\subsubsection{Parallelizations of ANNs with the same length}

\begin{definition}[Parallelization of ANNs with the same length]
	\label{Definition:simpleParallelization}
	Let $n\in\N$. Then we denote by 
%	$\parallelizationSpecial_{n}\colon \{(\Phi_1,\Phi_2,\dots, \Phi_n)\in\ANNs^n\colon \lengthANN(\Phi_1)= \lengthANN(\Phi_2)=\dots =\lengthANN(\Phi_n) \}\to \ANNs$ 
	\begin{equation}
		\parallelizationSpecial_{n}\colon \big\{(\Phi_1,\Phi_2,\dots, \Phi_n)\in\ANNs^n\colon \lengthANN(\Phi_1)= \lengthANN(\Phi_2)=\ldots =\lengthANN(\Phi_n) \big\}\to \ANNs
	\end{equation}
	the function which satisfies for all  $L\in\N$,
	$(l_{1,0},l_{1,1},\dots, l_{1,L}), (l_{2,0},l_{2,1},\dots, l_{2,L}),\dots,\allowbreak (l_{n,0},\allowbreak l_{n,1},\allowbreak\dots, l_{n,L})\in\N^{L+1}$, 
	$\Phi_1=((W_{1,1}, B_{1,1}),(W_{1,2}, B_{1,2}),\allowbreak \ldots, (W_{1,L},\allowbreak B_{1,L}))\in ( \times_{k = 1}^L\allowbreak(\R^{l_{1,k} \times l_{1,k-1}} \times \R^{l_{1,k}}))$, 
		$\Phi_2=((W_{2,1}, B_{2,1}),(W_{2,2}, B_{2,2}),\allowbreak \ldots, (W_{2,L},\allowbreak B_{2,L}))\in ( \times_{k = 1}^L\allowbreak(\R^{l_{2,k} \times l_{2,k-1}} \times \R^{l_{2,k}}))$,
%	$\Phi_2=((W_1^2, B_1^2),\allowbreak(W_2^2, B_2^2),\allowbreak \ldots, (W_L^2,\allowbreak B_L^2))\in ( \times_{k = 1}^L\allowbreak(\R^{l_k^2 \times l_{k-1}^2} \times \R^{l_k^2}))$, 
\dots,  
		$\Phi_n=((W_{n,1}, B_{n,1}),(W_{n,2}, B_{n,2}),\allowbreak \ldots, (W_{n,L},\allowbreak B_{n,L}))\in ( \times_{k = 1}^L\allowbreak(\R^{l_{n,k} \times l_{n,k-1}} \times \R^{l_{n,k}}))$
%	
%	
%	$(\Phi_1,\Phi_2,\dots,\Phi_n)=\big(((W_1^k, B_1^k),(W_2^k, B_2^k),\allowbreak \ldots, (W_L^k,\allowbreak B_L^k))\big)_{k\in \{0,1,\dots, n\}}\allowbreak\in \{\Phi\in\ANNs\colon \lengthANN(\Phi)= L \}^n$
	that
	%\begin{equation}\label{parallelisationSameLengthDef}
	%\begin{split}
	%&\parallelizationSpecial_{n,L}(\Phi_1,\Phi_2,\dots,\Phi_n)=
	%\left[
	%%\pa{\begin{pmatrix}W^1_1\\W^2_1\\\vdots\\ W^n_1\end{pmatrix},\begin{pmatrix}B^1_1\\B^2_1\\\vdots\\ B^n_1\end{pmatrix}},
	%\pa{\begin{pmatrix}W^1_1&&&\\& W^2_1 &&\\&&
	%	\ddots&\\&&&W^n_1\end{pmatrix},\begin{pmatrix}B^1_1\\B^2_1\\\vdots\\ B^n_1\end{pmatrix}},\right.
	%\\&
	%\left. \pa{\begin{pmatrix}W^1_2&&&\\& W^2_2 &&\\&&
	%	\ddots&\\&&&W^n_2\end{pmatrix},\begin{pmatrix}B^1_2\\B^2_2\\\vdots\\ B^n_2\end{pmatrix}},\dots,
	%\pa{\begin{pmatrix}W^1_L&&&\\& W^2_L &&\\&&\ddots&\\&&&W^n_L\end{pmatrix},\begin{pmatrix}B^1_L\\B^2_L\\\vdots\\ B^n_L\end{pmatrix}}\right],
	%\end{split}
	%\end{equation}	
	\begin{equation}\label{parallelisationSameLengthDef}
	\begin{split}
	\parallelizationSpecial_{n}(\Phi_1,\Phi_2,\dots,\Phi_n)&=
	\left(
	%\pa{\begin{pmatrix}W^1_1\\W^2_1\\\vdots\\ W^n_1\end{pmatrix},\begin{pmatrix}B^1_1\\B^2_1\\\vdots\\ B^n_1\end{pmatrix}},
	\pa{\begin{pmatrix}
		W_{1,1}& 0& 0& \cdots& 0\\
		0& W_{2,1}& 0&\cdots& 0\\
		0& 0& W_{3,1}&\cdots& 0\\
		\vdots& \vdots&\vdots& \ddots& \vdots\\
		0& 0& 0&\cdots& W_{n,1}
		\end{pmatrix} ,\begin{pmatrix}B_{1,1}\\B_{2,1}\\B_{3,1}\\\vdots\\ B_{n,1}\end{pmatrix}},\right.
	\\&\quad
		\pa{\begin{pmatrix}
		W_{1,2}& 0& 0& \cdots& 0\\
		0& W_{2,2}& 0&\cdots& 0\\
		0& 0& W_{3,2}&\cdots& 0\\
		\vdots& \vdots&\vdots& \ddots& \vdots\\
		0& 0& 0&\cdots& W_{n,2}
		\end{pmatrix} ,\begin{pmatrix}B_{1,2}\\B_{2,2}\\B_{3,2}\\\vdots\\ B_{n,2}\end{pmatrix}}
	,\dots,
	\\&\quad\left.
		\pa{\begin{pmatrix}
	W_{1,L}& 0& 0& \cdots& 0\\
	0& W_{2,L}& 0&\cdots& 0\\
	0& 0& W_{3,L}&\cdots& 0\\
	\vdots& \vdots&\vdots& \ddots& \vdots\\
	0& 0& 0&\cdots& W_{n,L}
	\end{pmatrix} ,\begin{pmatrix}B_{1,L}\\B_{2,L}\\B_{3,L}\\\vdots\\ B_{n,L}\end{pmatrix}}\right)
	\end{split}
	\end{equation}
	(cf.\ Definition~\ref{Def:ANN}).
\end{definition}

%\begin{lemma}\label{Lemma:ParallelizationElementary}
%	Let $n,L\in\N$, $\alpha_0,\alpha_1,\dots,\alpha_L\in\N$,
%	$(l_{1,0},l_{1,1},\dots, l_{1,L}), (l_{2,0},l_{2,1},\dots,\allowbreak l_{2,L}),\allowbreak\dots,\allowbreak (l_{n,0},\allowbreak l_{n,1},\allowbreak\dots, l_{n,L})\in\N^{L+1}$ satisfy for all $k\in\{0,1,\dots, L\}$ that  
%	$\alpha_k=\sum_{j=1}^nl_{j,k}$ and let
%	$\Phi_1=((W_{1,1}, B_{1,1}),(W_{1,2}, B_{1,2}),\allowbreak \ldots, (W_{1,L},\allowbreak B_{1,L}))\in ( \times_{k = 1}^L\allowbreak(\R^{l_{1,k} \times l_{1,k-1}}\allowbreak \times \R^{l_{1,k}}))$, 
%	$\Phi_2=((W_{2,1}, B_{2,1}),(W_{2,2}, B_{2,2}),\allowbreak \ldots, (W_{2,L},\allowbreak B_{2,L}))\in ( \times_{k = 1}^L\allowbreak(\R^{l_{2,k} \times l_{2,k-1}} \times \R^{l_{2,k}}))$,
%	\dots,  
%	$\Phi_n=((W_{n,1}, B_{n,1}),(W_{n,2}, B_{n,2}),\allowbreak \ldots, (W_{n,L},\allowbreak B_{n,L}))\in ( \times_{k = 1}^L\allowbreak(\R^{l_{n,k} \times l_{n,k-1}} \allowbreak\times \R^{l_{n,k}}))$. 
%	Then it holds that 
%	\begin{equation}\label{ParallelizationElementary:Display}
%		\parallelizationSpecial_{n}(\Phi_1,\Phi_2,\dots,\Phi_n)\in \big(\!\!\times_{k = 1}^L\allowbreak(\R^{\alpha_{k} \times \alpha_{k-1}} \times \R^{\alpha_{k}})\big)
%	\end{equation}
%		\begin{equation}\label{ParallelizationElementary:Display}
%		\parallelizationSpecial_{n}(\Phi_1,\Phi_2,\dots,\Phi_n)\in \Big(\!\!\times_{k = 1}^L\allowbreak\big(\R^{(\sum_{j=1}^nl_{j,k}) \times (\sum_{j=1}^nl_{j,k-1})} \times \R^{(\sum_{j=1}^nl_{j,k})}\big)\Big)
%		\end{equation}
%	(cf.\ Definition~\ref{Definition:simpleParallelization}).
%\end{lemma}

\begin{lemma}\label{Lemma:ParallelizationElementary}
	Let $n,L\in\N$,
	$(l_{1,0},l_{1,1},\dots, l_{1,L}), (l_{2,0},l_{2,1},\dots,\allowbreak l_{2,L}),\allowbreak\dots,\allowbreak (l_{n,0},\allowbreak l_{n,1},\allowbreak\dots, l_{n,L})\in\N^{L+1}$,
	$\Phi_1=((W_{1,1}, B_{1,1}),(W_{1,2}, B_{1,2}),\allowbreak \ldots, (W_{1,L},\allowbreak B_{1,L}))\in ( \times_{k = 1}^L\allowbreak(\R^{l_{1,k} \times l_{1,k-1}}\allowbreak \times \R^{l_{1,k}}))$, 
	$\Phi_2=((W_{2,1}, B_{2,1}),(W_{2,2}, B_{2,2}),\allowbreak \ldots, (W_{2,L},\allowbreak B_{2,L}))\in ( \times_{k = 1}^L\allowbreak(\R^{l_{2,k} \times l_{2,k-1}} \times \R^{l_{2,k}}))$,
	\dots,  
	$\Phi_n=((W_{n,1}, B_{n,1}),(W_{n,2}, B_{n,2}),\allowbreak \ldots, (W_{n,L},\allowbreak B_{n,L}))\in ( \times_{k = 1}^L\allowbreak(\R^{l_{n,k} \times l_{n,k-1}} \allowbreak\times \R^{l_{n,k}}))$. 
	Then it holds that 
	\begin{equation}\label{ParallelizationElementary:Display}
	\parallelizationSpecial_{n}(\Phi_1,\Phi_2,\dots,\Phi_n)\in \Big(\!\!\times_{k = 1}^L\allowbreak\big(\R^{(\sum_{j=1}^nl_{j,k}) \times (\sum_{j=1}^nl_{j,k-1})} \times \R^{(\sum_{j=1}^nl_{j,k})}\big)\Big)
	\end{equation}
	(cf.\ Definition~\ref{Definition:simpleParallelization}).
\end{lemma}

\begin{proof}[Proof of Lemma~\ref{Lemma:ParallelizationElementary}]	
	Note that \eqref{parallelisationSameLengthDef} establishes  \eqref{ParallelizationElementary:Display}.
	The proof of Lemma~\ref{Lemma:ParallelizationElementary} is thus completed.
\end{proof}

%\begin{lemma}\label{Lemma:PropertiesOfParallelizationEqualLength}
%	Let $a\in C(\R,\R)$, let
%	$n\in\N$, 
%	$\Phi_1,\Phi_2,\allowbreak\dots,\allowbreak \Phi_n\in\ANNs$,
%	$D,L,\mathfrak{D},d_1,d_2,\allowbreak\dots,\allowbreak d_n,\mathfrak{d}_1,\mathfrak{d}_2,\allowbreak\dots,\allowbreak \mathfrak{d}_n
%	\in\N$ satisfy for all $j\in\{1,2,\dots,n\}$ that $\lengthANN(\Phi_j)=L$, 
%	%	assume that for all $j\in\{1,2,\dots,n\}$ it holds that
%	$\inDimANN(\Phi_j)=d_j$, $\outDimANN(\Phi_j)=\mathfrak{d}_j$,
%	$D=\smallsum\nolimits_{k=1}^n d_k$, and $\mathfrak{D}=\smallsum\nolimits_{k=1}^n \mathfrak{d}_k$
%	(cf.\ Definition~\ref{Def:ANN}).	
%	Then
%		it holds for all   $x\in\R^D, x_1\in\R^{d_1},x_2\in\R^{d_2},\dots, x_n\in\R^{d_n}$ with $x=(x_1,x_2,\dots, x_n)$ that $\functionANN(\parallelizationSpecial_{n}(\Phi_1,\Phi_2,\dots,\Phi_n))\in C(\R^{D},\R^{\mathfrak{D}})$
%		%		\begin{equation}
%		%		\functionANN(\parallelizationSpecial_{n}(\Phi_1,\Phi_2,\dots,\Phi_n))\in C(\R^{D},\R^{\mathfrak{D}})
%		%		\end{equation}
%		and
%		\begin{equation}\label{PropertiesOfParallelizationEqualLengthFunction}
%		%\functionANN(\Phi) \in C(\R^{\alpha^j_0},\R^{\alpha^j_L})\qandq
%		%\qquad
%		\begin{split}
%		&\big[ \functionANN\big(\parallelizationSpecial_{n}(\Phi_1,\Phi_2,\dots,\Phi_n)\big) \big] (x) 
%		=\big([\functionANN(\Phi_1)](x_1), [\functionANN(\Phi_2)](x_2),\dots,
%		[\functionANN(\Phi_n)](x_n) \big)
%		\end{split}
%		\end{equation}
%		(cf.\  Definition~\ref{Definition:ANNrealization} and Definition~\ref{Definition:simpleParallelization}).
%\end{lemma}

\begin{prop}\label{Lemma:PropertiesOfParallelizationEqualLength}
%	Let $a\in C(\R,\R)$, 
%	$n\in\N$, 
%	$\Phi=(\Phi_1,\Phi_2,\allowbreak\dots,\allowbreak \Phi_n)\in\ANNs^n$,
%	$L,d,\mathfrak{d}\in\N$ satisfy for all $j\in\{1,2,\dots,n\}$ that $\lengthANN(\Phi_j)=L$, 
%	$d=\smallsum\nolimits_{k=1}^n \inDimANN(\Phi_k)$, and $\mathfrak{d}=\smallsum\nolimits_{k=1}^n \outDimANN(\Phi_k)$
%	(cf.\ Definition~\ref{Def:ANN}).	
	Let $a\in C(\R,\R)$, 
	$n\in\N$, 
	$\Phi=(\Phi_1,\Phi_2,\allowbreak\dots,\allowbreak \Phi_n)\in\ANNs^n$
	satisfy that $\lengthANN(\Phi_1)=\lengthANN(\Phi_2)=\ldots=\lengthANN(\Phi_n)$
	(cf.\ Definition~\ref{Def:ANN}).
	Then
	\begin{enumerate}[(i)]
		\item\label{PropertiesOfParallelizationEqualLength:ItemOne} it holds that 
		\begin{equation}
%			\functionANN(\parallelizationSpecial_{n}(\Phi_1,\Phi_2,\dots,\Phi_n))\in C(\R^{\smallsum\nolimits_{k=1}^n \inDimANN(\Phi_k)},\R^{\smallsum\nolimits_{k=1}^n \outDimANN(\Phi_k)})
\functionANN(\parallelizationSpecial_{n}(\Phi))\in C\big(\R^{[\sum_{j=1}^n \inDimANN(\Phi_j)]},\R^{[\sum_{j=1}^n \outDimANN(\Phi_j)]}\big)
		\end{equation}
		and
		\item\label{PropertiesOfParallelizationEqualLength:ItemTwo} it holds for all    $x_1\in\R^{\inDimANN(\Phi_1)},x_2\in\R^{\inDimANN(\Phi_2)},\dots, x_n\in\R^{\inDimANN(\Phi_n)}$ that 
			\begin{equation}\label{PropertiesOfParallelizationEqualLengthFunction}
			%\functionANN(\Phi) \in C(\R^{\alpha^j_0},\R^{\alpha^j_L})\qandq
			%\qquad
			\begin{split}
			&\big( \functionANN\big(\parallelizationSpecial_{n}(\Phi)\big) \big) (x_1,x_2,\dots, x_n) 
			\\&=\big((\functionANN(\Phi_1))(x_1), (\functionANN(\Phi_2))(x_2),\dots,
			(\functionANN(\Phi_n))(x_n) \big)\in \R^{[\sum_{j=1}^n \outDimANN(\Phi_j)]}
			\end{split}
			\end{equation}
	\end{enumerate}
	(cf.\  Definition~\ref{Definition:ANNrealization} and Definition~\ref{Definition:simpleParallelization}).
\end{prop}

\begin{proof}[Proof of Proposition~\ref{Lemma:PropertiesOfParallelizationEqualLength}]	
	Throughout this proof  
	let $L\in\N$ satisfy that $L=\lengthANN(\Phi_1)$,
	%	let $d,\mathfrak{d}\in\N$ satisfy that $d=\smallsum\nolimits_{k=1}^n \inDimANN(\Phi_k)$ and $\mathfrak{d}=\smallsum\nolimits_{k=1}^n \outDimANN(\Phi_k)$,
	let $l_{j,0},l_{j,1},\dots, l_{j,L}\in\N$, $j\in\{1,2,\dots,\allowbreak n\}$, satisfy for all $j\in\{1,2,\dots, n\}$ that $\dims(\Phi_j)=(l_{j,0},l_{j,1},\dots, l_{j,L})$,
	%	 let $(l^j)_{j\in\{1,2,\dots, n\}}=\big((l_0^j,\allowbreak l_1^j,\allowbreak\dots,\allowbreak l_L^j)\big)_{j\in\{1,2,\dots, n\}}\allowbreak\subseteq\N^{L+1}$ be the multi-indices which satisfy for all $j\in\{1,2,\dots, n\}$ that 
	%	$\dims(\Phi_j)=l^j$, 
	let $\big((W_{j,1},B_{j,1}),\allowbreak (W_{j,2},B_{j,2}),\dots,\allowbreak (W_{j,L},B_{j,L})\big)\in ( \times_{k = 1}^L\allowbreak(\R^{l_{j,k} \times l_{j,k-1}} \times \R^{l_{j,k}}))$, $j\in\{1,2,\dots, n\}$, satisfy for all $j\in\{1,2,\dots, n\}$ that 
	\begin{equation}\label{matricesForPsiJ}
	\Phi_j=\big((W_{j,1},B_{j,1}), (W_{j,2},B_{j,2}),\dots, (W_{j,L},B_{j,L})\big),
	\end{equation}
	let $\alpha_k\in\N$, $k\in\{0,1,\dots, L\}$, satisfy for all $k\in\{0,1,\dots, L\}$ that 
	$\alpha_k=\sum_{j=1}^nl_{j,k}$,
	let $\big((A_1,b_1), (A_2,b_2),\dots, (A_L,b_L)\big)\in ( \times_{k = 1}^L\allowbreak(\R^{\alpha_{k} \times \alpha_{k-1}} \times \R^{\alpha_{k}}))$ satisfy  that 
	\begin{equation}\label{matricesForParallelization}
	\parallelizationSpecial_{n}(\Phi)=\big((A_1,b_1), (A_2,b_2),\dots, (A_L,b_L)\big)
	\end{equation}
	(cf.\ Lemma~\ref{Lemma:ParallelizationElementary}),
	let	$(x_{j,0},x_{j,1},\ldots,x_{j,L-1}) \in (\R^{l_{j,0}}\times \R^{l_{j,1}}\times\ldots\times \R^{l_{j,{L-1}}})$, $j\in\{1,2,\dots, n\}$, 
	satisfy for all $j\in\{1,2,\dots, n\}$, $k \in \N \cap (0,L)$ that 
	\begin{equation}\label{realizationRecursionPhij}
			 x_{j,k} =\activationDim{l_{j,k}}(W_{j,k}\allowbreak x_{j,k-1} + B_{j,k})
	\end{equation}	
%	 $x_{j,k} =\activationDim{l_{j,k}}(W_{j,k}\allowbreak x_{j,k-1} + B_{j,k})$,
	 	(cf.\ Definition~\ref{Def:multidim_version}),
	and let   $\mathfrak{x}_0 \in \R^{\alpha_{0}}, \mathfrak{x}_1 \in \R^{\alpha_{1}}, \ldots, \mathfrak{x}_{L-1} \in \R^{\alpha_{{L-1}}}$ 
	satisfy for all $k\in\{0,1,\dots, L-1\}$ that $\mathfrak{x}_k=(x_{1,k},x_{2,k},\dots, x_{n,k})$.
	Observe that \eqref{matricesForParallelization} demonstrates that $\inDimANN(\parallelizationSpecial_{n}(\Phi))=\alpha_0$ and $\outDimANN(\parallelizationSpecial_{n}(\Phi))=\alpha_L$.
	Combining this with item~\eqref{elementaryPropertiesANN:ItemTwo} in Lemma~\ref{Lemma:elementaryPropertiesANN}, the fact that for all $k\in\{0,1,\dots, L\}$ it holds that 
	$\alpha_k=\sum_{j=1}^nl_{j,k}$, the fact that for all $j\in \{1,2,\dots,n\}$ it holds that $\inDimANN(\Phi_j)=l_{j,0}$, and the fact that for all $j\in \{1,2,\dots,n\}$ it holds that $\outDimANN(\Phi_j)=l_{j,L}$ ensures that 
	\begin{equation}
	\begin{split}
		\functionANN(\parallelizationSpecial_{n}(\Phi))
		&\in C(\R^{\alpha_0},\R^{\alpha_L})
		=C\big(\R^{[\sum_{j=1}^n l_{j,0}]},\R^{[\sum_{j=1}^n l_{j,L}]}\big)
		\\&=C\big(\R^{[\sum_{j=1}^n \inDimANN(\Phi_j)]},\R^{[\sum_{j=1}^n \outDimANN(\Phi_j)]}\big).
	\end{split}
	\end{equation}
	This proves item~\eqref{PropertiesOfParallelizationEqualLength:ItemOne}.
	Moreover, observe that \eqref{parallelisationSameLengthDef} and \eqref{matricesForParallelization} demonstrate that for all $k\in\{1,2,\dots,L\}$ it holds that 
	\begin{equation}\label{matricesForParallelizationTwo}
		\begin{split}
		A_k=\begin{pmatrix}
			W_{1,k}& 0& 0& \cdots& 0\\
			0& W_{2,k}& 0&\cdots& 0\\
			0& 0& W_{3,k}&\cdots& 0\\
			\vdots& \vdots&\vdots& \ddots& \vdots\\
			0& 0& 0&\cdots& W_{n,k}
			\end{pmatrix} 
			\qandq
			b_k=\begin{pmatrix}B_{1,k}\\B_{2,k}\\B_{3,k}\\\vdots\\ B_{n,k}\end{pmatrix}.
		\end{split}
	\end{equation}
	Combining this with \eqref{multidim_version:Equation}, \eqref{realizationRecursionPhij}, and the fact that for all $k\in\N\cap [0,L)$ it holds that $\mathfrak{x}_k=(x_{1,k},x_{2,k},\dots, x_{n,k})$  implies that for all $k\in\N\cap (0,L)$ it holds that
	\begin{equation}
	\begin{split}
	\activationDim{\alpha_k}(A_k \mathfrak{x}_{k-1} + b_k)
	=
	\begin{pmatrix}
	\activationDim{l_{1,k}}(W_{1,k} x_{1,k-1} + B_{1,k})\\\activationDim{l_{2,k}}(W_{2,k} x_{2,k-1} + B_{2,k})
	\\ \vdots \\
	\activationDim{l_{n,k}}(W_{n,k} x_{n,k-1} + B_{n,k})
	\end{pmatrix}
	=\begin{pmatrix}
	x_{1,k}\\x_{2,k}
	\\ \vdots \\
	x_{n,k}
	\end{pmatrix}
	=\mathfrak{x}_k.
	\end{split}
	\end{equation}
	This, \eqref{setting_NN:ass2}, \eqref{matricesForPsiJ}, \eqref{matricesForParallelization}, \eqref{realizationRecursionPhij}, \eqref{matricesForParallelizationTwo},
	the fact that $\mathfrak{x}_{0}=(x_{1,{0}},x_{2,{0}},\dots,\allowbreak x_{n,{0}})$,
	  and the fact that $\mathfrak{x}_{L-1}=(x_{1,{L-1}},x_{2,{L-1}},\dots,\allowbreak x_{n,{L-1}})$ ensure that
	\begin{equation}\label{PropertiesOfParallelizationEqualLengthFunctionInProof}
	%\functionANN(\Phi) \in C(\R^{l_{j,0},\R^{l_{j,L})\qandq
	%\qquad
	\begin{split}
	&\big( \functionANN\big(\parallelizationSpecial_{n}(\Phi)\big) \big) (x_{1,0},x_{2,0},\dots, x_{n,0})=
	\big( \functionANN\big(\parallelizationSpecial_{n}(\Phi)\big) \big) (\mathfrak{x}_0) 
	\\&= A_L \mathfrak{x}_{L-1} + b_L
	=\begin{pmatrix}
	W_{1,L} x_{1,L-1}+B_{1,L}
	\\W_{2,L} x_{2,L-1}+B_{2,L}
	\\ \vdots \\
	W_{n,L} x_{n,L-1}+B_{n,L}
	\end{pmatrix}
	=\begin{pmatrix}
	(\functionANN(\Phi_1))(x_{1,0})
	\\
	(\functionANN(\Phi_2))(x_{2,0})
	\\ \vdots \\
	(\functionANN(\Phi_n))(x_{n,0})
	\end{pmatrix}.
	\end{split}
	\end{equation}
	This establishes item~\eqref{PropertiesOfParallelizationEqualLength:ItemTwo}.
	The proof of Proposition~\ref{Lemma:PropertiesOfParallelizationEqualLength} is thus completed.
\end{proof}

\begin{prop}\label{Lemma:PropertiesOfParallelizationEqualLengthDims}
	Let 
	$n,L\in\N$, 
	$\Phi=(\Phi_1,\Phi_2,\allowbreak\dots,\allowbreak \Phi_n)\in\ANNs^n$,
	$(l_{1,0},l_{1,1},\dots, l_{1,L}),\allowbreak (l_{2,0},l_{2,1},\allowbreak\dots,\allowbreak l_{2,L}),\allowbreak \dots, (l_{n,0},l_{n,1},\dots, l_{n,L})
	\in\N^{L+1}$ satisfy for all $j\in\{1,2,\dots,n\}$ that 
	$\dims(\Phi_j)=(l_{j,0},l_{j,1},\dots,\allowbreak l_{j,L})$
	(cf.\ Definition~\ref{Def:ANN}).	
	Then
	\begin{enumerate}[(i)]
		\item \label{PropertiesOfParallelizationEqualLengthDims:Dims} it holds that
		\begin{equation}
		\dims\big(\parallelizationSpecial_{n}(\Phi)\big)=\big(\smallsum_{j=1}^n l_{j,0}, \smallsum_{j=1}^n l_{j,1},\dots, \smallsum_{j=1}^n l_{j,L}\big)
		\end{equation}
		and
		\item \label{PropertiesOfParallelizationEqualLengthDims:Params}
		it holds that
		\begin{equation}
		\paramANN(\parallelizationSpecial_{n}(\Phi))\le \tfrac{1}{2} \big[\smallsum\nolimits_{j=1}^n \paramANN(\Phi_j)\big]^2
		%		\qandq \paramNotZeroANN(\Phi)= \smallsum\nolimits_{j=1}^n \paramNotZeroANN(\Phi_j),
		\end{equation}
	\end{enumerate}
	(cf.\   Definition~\ref{Definition:simpleParallelization}).
\end{prop}

\begin{proof}[Proof of Proposition~\ref{Lemma:PropertiesOfParallelizationEqualLengthDims}]	
%Note that  the hypothesis that for all $j\in\{1,2,\dots,\allowbreak n\}$ it holds that 
%$\dims(\Phi_j)=(l_{j,0},l_{j,1},\dots, l_{j,L})$ and Lemma~\ref{Lemma:ParallelizationElementary}  assure that 
Note that  the hypothesis that $\forallDist j\in\{1,2,\dots,\allowbreak n\}\colon \dims(\Phi_j)=(l_{j,0},l_{j,1},\dots, l_{j,L})$  and Lemma~\ref{Lemma:ParallelizationElementary}  assure that
		\begin{equation}\label{PropertiesOfParallelizationEqualLengthDims:DimsProof}
\dims\big(\parallelizationSpecial_{n}(\Phi)\big)=\big(\smallsum_{j=1}^n l_{j,0}, \smallsum_{j=1}^n l_{j,1},\dots, \smallsum_{j=1}^n l_{j,L}\big).
\end{equation}
This establishes item~\eqref{PropertiesOfParallelizationEqualLengthDims:Dims}. Moreover, observe that \eqref{PropertiesOfParallelizationEqualLengthDims:DimsProof}
demonstrates that
	\begin{equation}
	\begin{split}
	&\paramANN(\parallelizationSpecial_{n}(\Phi))
		=\sum_{k=1}^L \Big[\smallsum\nolimits_{i=1}^n l_{i,k}\Big] \Big[\big(\smallsum\nolimits_{i=1}^n l_{i,k-1}\big)+1\Big]
\\&=\sum_{k=1}^L \Big[\smallsum\nolimits_{i=1}^n l_{i,k}\Big] \Big[\big(\smallsum\nolimits_{j=1}^n l_{j,k-1}\big)+1\Big]
	\\&\le \sum_{i=1}^n \sum_{j=1}^n\sum_{k=1}^L l_{i,k} (l_{j,k-1}+1)
	\le \sum_{i=1}^n \sum_{j=1}^n\sum_{k,\ell=1}^L l_{i,k} (l_{j,\ell-1}+1)
	%		\\&\le \sum_{i=1}^n \sum_{j=1}^n\left[\smallsum\nolimits_{k=1}^L \vert l_k^j\vert^2\right]^{\nicefrac{1}{2}}
	%		\left[\smallsum\nolimits_{k=1}^L \vert l_{k-1}^m+1\vert^2\right]^{\nicefrac{1}{2}}
	\\&= \sum_{i=1}^n \sum_{j=1}^n\Big[\smallsum\nolimits_{k=1}^L  l_{i,k}\Big]
	\Big[\smallsum\nolimits_{\ell=1}^L  (l_{j,\ell-1}+1)\Big]
	\\&\le \sum_{i=1}^n \sum_{j=1}^n\Big[\smallsum\nolimits_{k=1}^L  \tfrac{1}{2}l_{i,k} (l_{i,k-1}+1)\Big]
	\Big[\smallsum\nolimits_{\ell=1}^L  l_{j,\ell}(l_{j,\ell-1}+1)\Big]
	\\&=  \sum_{i=1}^n \sum_{j=1}^n \tfrac{1}{2}\paramANN(\Phi_i) \paramANN(\Phi_j)
	=\tfrac{1}{2}\Big[\smallsum\nolimits_{i=1}^n \paramANN(\Phi_i)\Big]^2.
	\end{split}
	\end{equation}
	%	Combining this with \eqref{PropertiesOfParallelizationEqualLengthDimsFunctionInProof} and the fact that 	$\dims(\parallelizationSpecial_{n}(\Phi_1,\Phi_2,\dots,\Phi_n))=(\alpha_0,\alpha_1,\dots, \alpha_L)$ establishes item~\eqref{PropertiesOfParallelizationEqualLengthDims:General}.
	The proof of Proposition~\ref{Lemma:PropertiesOfParallelizationEqualLengthDims} is thus completed.
\end{proof}

\begin{cor}\label{Lemma:ParallelizationImprovedBoundsOne}
	Let   $n\in\N$, $\Phi=( \Phi_1,\Phi_2,\dots,\allowbreak \Phi_n)\in\ANNs^n$ satisfy  that 
	$\dims(\Phi_1)=\dims(\Phi_2)=\ldots=\dims(\Phi_n)$
	(cf.\ Definition~\ref{Def:ANN}).	
	Then it holds that $\paramANN(\parallelizationSpecial_{n}(\Phi))\le n^2 \paramANN(\Phi_1)$
	(cf.\   Definition~\ref{Definition:simpleParallelization}).
\end{cor}

\begin{proof}[Proof of Corollary~\ref{Lemma:ParallelizationImprovedBoundsOne}]	
	Throughout this proof let $L\in\N$,
	$l_0,l_1,\dots, l_L
	\in\N$ satisfy  that 
	$\dims(\Phi_1)=(l_0,l_1,\dots, l_{L})$.
Note that item~\eqref{PropertiesOfParallelizationEqualLengthDims:Dims} in Proposition~\ref{Lemma:PropertiesOfParallelizationEqualLengthDims} and the fact that $\forallDist j\in \{1,2,\dots,n\}\colon\dims(\Phi_j)=(l_0,l_1,\dots,\allowbreak l_{L})$ demonstrate that
\begin{equation}
\begin{split}
\paramANN(\parallelizationSpecial_{n}(\Phi_1,\Phi_2,\dots,\Phi_n))&=\smallsum\limits_{j=1}^L  (nl_j)\big((nl_{j-1})+1\big)\le \smallsum\limits_{j=1}^L  (nl_j)\big((nl_{j-1})+n\big)
\\&= n^2 \bigg[\smallsum\limits_{j=1}^L  l_j(l_{j-1}+1)\bigg]
= n^2 \paramANN(\Phi_1).
\end{split}
\end{equation}
	The proof of Corollary~\ref{Lemma:ParallelizationImprovedBoundsOne} is thus completed.
\end{proof}

\subsubsection{Parallelizations of ANNs with different lengths}

%\begin{definition}[Parallelization of GANNs with different length]\label{Definition:generalParallelization}
%	Let $n\in\N$,  $\Psi=(\Psi_1,\Psi_2,\dots, \Psi_n)\in \ANNs^n$ satisfy for all $j\in\{1,2,\dots, n\}$ that $\lengthANN(\Psi_j)=2$ and $\inDimANN(\Psi_j)=\outDimANN(\Psi_j)$.
%	Then
%	we denote by 
%	\begin{equation}
%		\parallelization_{n,\Psi}\colon \{(\Phi_1,\Phi_2,\dots, \Phi_n)\in\ANNs^n\colon \,(\forallDist j\in\{1,2,\dots, n\}\colon \outDimANN(\Phi_j)= \inDimANN(\Psi_j))\}\to \ANNs
%	\end{equation}
%%	$\parallelization_{n,\Psi}\colon \{(\Phi_1,\Phi_2,\dots, \Phi_n)\in\ANNs^n\colon \,(\forallDist j\in\{1,2,\dots, n\}\colon \outDimANN(\Phi_j)= \inDimANN(\Psi_j))\}\to \ANNs$
%	the function which satisfies for all $L\in\N$,  $(\Phi_1,\Phi_2,\dots, \Phi_n)\in\ANNs^n$ with $\max_{k\in\{1,2,\dots,n\}} \allowbreak\lengthANN(\Phi_k)=L$
%	and $\forallDist j\in\{1,2,\dots, n\}\colon \outDimANN(\Phi_j)= \inDimANN(\Psi_j)$
%	that
%	\begin{equation}\label{generalParallelization:DefinitionFormula}
%	\parallelization_{n,\Psi}(\Phi)=\parallelizationSpecial_{n}\big(\longerANN{L,\Psi_1}({\Phi_1}),\longerANN{L,\Psi_2}({\Phi_2}),\dots,\longerANN{L,\Psi_n}({\Phi_n})\big)
%	\end{equation}
%	(cf.\ Definition~\ref{Def:ANN}, Definition~\ref{Definition:ANNenlargement}, Lemma~\ref{Lemma:PropertiesOfANNenlargementGeometry},  and Definition~\ref{Definition:simpleParallelization}).
%\end{definition}

\begin{definition}[Parallelization of ANNs with different length]\label{Definition:generalParallelization}
	Let $n\in\N$,  $\Psi=(\Psi_1,\Psi_2,\dots, \Psi_n)\in \ANNs^n$ satisfy for all $j\in\{1,2,\dots, n\}$ that $\hiddenLength(\Psi_j)=1$ and $\inDimANN(\Psi_j)=\outDimANN(\Psi_j)$.
	Then
	we denote by 
	\begin{equation}
		\parallelization_{n,\Psi}\colon \{(\Phi_1,\Phi_2,\dots, \Phi_n)\in\ANNs^n\colon \,(\forallDist j\in\{1,2,\dots, n\}\colon \outDimANN(\Phi_j)= \inDimANN(\Psi_j))\}\to \ANNs
	\end{equation}
%	$\parallelization_{n,\Psi}\colon \{(\Phi_1,\Phi_2,\dots, \Phi_n)\in\ANNs^n\colon \,(\forallDist j\in\{1,2,\dots, n\}\colon \outDimANN(\Phi_j)= \inDimANN(\Psi_j))\}\to \ANNs$
	the function which satisfies for all $\Phi=(\Phi_1,\Phi_2,\dots, \Phi_n)\in\ANNs^n$ with $\forallDist j\in\{1,2,\dots, n\}\colon\allowbreak \outDimANN(\Phi_j)= \inDimANN(\Psi_j)$
	that
%	\begin{multline}\label{generalParallelization:DefinitionFormula}
%			\parallelization_{n,\Psi}(\Phi)
%			\\=\parallelizationSpecial_{n}\big(\longerANN{\max_{1\le k\le n}\lengthANN(\Phi_k),\Psi_1}({\Phi_1}),
%			\longerANN{\max_{1\le k\le n}\lengthANN(\Phi_k),\Psi_2}({\Phi_2}),\dots,\longerANN{\max_{1\le k\le n}\lengthANN(\Phi_k),\Psi_n}({\Phi_n})\big)
%	\end{multline}
		\begin{multline}\label{generalParallelization:DefinitionFormula}
		\parallelization_{n,\Psi}(\Phi)=\parallelizationSpecial_{n}\big(\longerANN{\max_{k\in\{1,2,\dots,n\}}\lengthANN(\Phi_k),\Psi_1}({\Phi_1}),\dots,\longerANN{\max_{k\in\{1,2,\dots,n\}}\lengthANN(\Phi_k),\Psi_n}({\Phi_n})\big)
		\end{multline}
%	\begin{multline}\label{generalParallelization:DefinitionFormula}
%	\parallelization_{n,\Psi}(\Phi)=\parallelizationSpecial_{n}\big(\longerANN{\max_{k\in\{1,2,\dots,n\}}\lengthANN(\Phi_k),\Psi_1}({\Phi_1}),\longerANN{\max_{k\in\{1,2,\dots,n\}}\lengthANN(\Phi_k),\Psi_2}({\Phi_2}),\\\dots,\longerANN{\max_{k\in\{1,2,\dots,n\}}\lengthANN(\Phi_k),\Psi_n}({\Phi_n})\big)
%	\end{multline}
	(cf.\ Definition~\ref{Def:ANN}, Definition~\ref{Definition:ANNenlargement}, Lemma~\ref{Lemma:PropertiesOfANNenlargementGeometry},  and Definition~\ref{Definition:simpleParallelization}).
\end{definition}

\begin{cor}\label{Lemma:PropertiesOfParallelizationRealization}
	Let $\activation\in C(\R,\R)$, $n\in\N$, $\mathbb{I}=(\idANNshort{1},\idANNshort{2},\dots, \idANNshort{n})$, $\Phi=(\Phi_1,\Phi_2,\dots,\allowbreak \Phi_n)\in\ANNs^n$
	satisfy for all $j\in\{1,2,\dots, n\}$, $x\in\R^{\outDimANN(\Phi_j)}$ that $\hiddenLength(\idANNshort{j}) =1$, $ \inDimANN(\mathbb{I}_j)=\outDimANN(\mathbb{I}_j)=\outDimANN(\Phi_j)$, and
	$(\functionANN(\idANNshort{j}))(x)=x$
	(cf.\ Definition~\ref{Def:ANN} and Definition~\ref{Definition:ANNrealization}).
	Then
	\begin{enumerate}[(i)]
		\item \label{PropertiesOfParallelizationRealization:ItemOne} it holds that
		\begin{equation}
			\functionANN\big(\!\parallelization_{n,\mathbb{I}}(\Phi)\big)\in C\big(\R^{[\sum_{j=1}^n \inDimANN(\Phi_j)]},\R^{[\sum_{j=1}^n \outDimANN(\Phi_j)]}\big)
		\end{equation}
		and
		\item \label{PropertiesOfParallelizationRealization:ItemTwo}
		it holds for all $x_1\in\R^{\inDimANN(\Phi_1)},x_2\in\R^{\inDimANN(\Phi_2)},\dots, x_n\in\R^{\inDimANN(\Phi_n)}$ that 
		\begin{equation}\label{PropertiesOfParallelizationRealizationEqualLengthFunction}
		%\functionANN(\Phi) \in C(\R^{\alpha_{j,0},\R^{\alpha_{j,L})\qandq
		%\qquad
		\begin{split}
		&\big( \functionANN(\parallelization_{n,\mathbb{I}}(\Phi)) \big) (x_1,x_2,\dots, x_n) 
		\\&=\big((\functionANN(\Phi_1))(x_1), (\functionANN(\Phi_2))(x_2),\dots,
		(\functionANN(\Phi_n))(x_n) \big)\in \R^{[\sum_{j=1}^n \outDimANN(\Phi_j)]}
		\end{split}
		\end{equation}
	\end{enumerate}
(cf.\ Definition~\ref{Definition:generalParallelization}).
\end{cor}

%\begin{lemma}\label{Lemma:PropertiesOfParallelization}
%		Let $n\in\N$, $\Psi=(\Psi_j)_{j\in\{1,2,\dots, n\}}\subseteq \ANNs$,
%		$\Phi_1,\Phi_2,\dots, \Phi_n\in\ANNs$,  $L,\allowbreak\mathfrak{d}_1,\mathfrak{d}_2,\dots, \mathfrak{d}_n,\hiddenDimId_1,\hiddenDimId_2,\dots, \hiddenDimId_n\in\N$
%		 satisfy for all $j\in\{1,2,\dots, n\}$, $x\in\R^{\mathfrak{d}_j}$ that $\dims(\Psi_j) = (\mathfrak{d}_j,\hiddenDimId_j,\mathfrak{d}_j)$, 
%		 $(\functionANN(\idANNshort{j}))(x)=x$,
%		  $\outDimANN(\Phi_j)=\mathfrak{d}_j$, and
%%		$\functionANN(\idANNshort{j})\in C(\R^{\mathfrak{d}_j},\R^{\mathfrak{d}_j})$
%  $L=\max_{j\in\{1,2,\dots,n\}} \lengthANN(\Phi_j)$
%			(cf.\ Definition~\ref{Def:ANN} and Definition~\ref{Definition:ANNrealization}).
%	Then
%		it holds that
%		\begin{equation}
%		\begin{split}
%		&\paramANN\big(\!\parallelization_{n,\Psi}(\Phi_1,\Phi_2,\dots,\Phi_n)\big)
%		\\&\le \tfrac{1}{2} \bigg(\left[\smallsum\nolimits_{j=1}^n\left(
%		\max\big\{1,\tfrac{\hiddenDimId_j}{\mathfrak{d}_j}\}\,\paramANN(\Phi_j)
%		+(L-\lengthANN(\Phi_j)-1) \,\hiddenDimId_j\,(\hiddenDimId_j+1)
%		+\mathfrak{d}_j\,(\hiddenDimId_j+1)\right)
%		\indicator{(\lengthANN(\Phi_j),\infty)}(L)\right]
%		\\&\qquad+\left[\smallsum\nolimits_{j=1}^n  \paramANN(\Phi_j)\indicator{\{\lengthANN(\Phi_j)\}}(L)\right]\bigg)^{\!2}
%		\end{split}
%		\end{equation}
%		(cf.\  Definition~\ref{Definition:generalParallelization}).
%\end{lemma}

\begin{proof}[Proof of Corollary~\ref{Lemma:PropertiesOfParallelizationRealization}]
	Throughout this proof let $L\in\N$ satisfy that $L=	\allowbreak\max_{j\in\{1,2,\dots,n\}} \allowbreak\lengthANN(\Phi_j)$.
Note that item~\eqref{PropertiesOfANNenlargementGeometry:ItemLonger} in Lemma~\ref{Lemma:PropertiesOfANNenlargementGeometry}, the hypothesis that
for all $j\in\{1,2,\dots, n\}$ it holds that $\hiddenLength(\idANNshort{j}) =1$, \eqref{ANNenlargement:Equation}, \eqref{PropertiesOfCompositions:LengthDisplay}, and item~\eqref{PropertiesOfANNenlargementRealization:ItemIdentityLonger} in Lemma~\ref{Lemma:PropertiesOfANNenlargementRealization} demonstrate 
\begin{enumerate}[(I)]
	\item that for all $j\in\{1,2,\dots, n\}$ it holds that $\lengthANN(\longerANN{L,\mathbb{I}_j}(\Phi_j))=L$ and
	 $\functionANN(\longerANN{L,\mathbb{I}_j}(\Phi_j))\in C(\R^{\inDimANN(\Phi_j)},\R^{\outDimANN(\Phi_j)})$
	 and
	 \item that for all
	  $j\in\{1,2,\dots, n\}$, $x\in\R^{\inDimANN(\Phi_j)}$	it holds that
	 \begin{equation}
	 \big(\functionANN(\longerANN{L,\mathbb{I}_j}(\Phi_j))\big)(x)=(\functionANN(\Phi_j))(x)
	 \end{equation}
\end{enumerate}
	(cf.\  Definition~\ref{Definition:ANNenlargement}).
	Items~\eqref{PropertiesOfParallelizationEqualLength:ItemOne}--\eqref{PropertiesOfParallelizationEqualLength:ItemTwo} in	
	 Proposition~\ref{Lemma:PropertiesOfParallelizationEqualLength} therefore imply
	 \begin{enumerate}[(A)]
	 	\item that 	 
	 	\begin{equation}
	 	\functionANN\big(\parallelizationSpecial_n\big(\longerANN{L,\mathbb{I}_1}({\Phi_1}),\longerANN{L,\mathbb{I}_2}({\Phi_2}),\dots,\longerANN{L,\mathbb{I}_n}({\Phi_n})\big)
	 	\in C\big(\R^{[\sum_{j=1}^n \inDimANN(\Phi_j)]},\R^{[\sum_{j=1}^n \outDimANN(\Phi_j)]}\big)
	 	\end{equation}
	 	and 
	 	\item that for all $x_1\in\R^{\inDimANN(\Phi_1)},x_2\in\R^{\inDimANN(\Phi_2)},\dots, x_n\in\R^{\inDimANN(\Phi_n)}$ it holds that
	 		 	\begin{equation}
	 		 	%\functionANN(\Phi) \in C(\R^{\alpha_{j,0},\R^{\alpha_{j,L})\qandq
	 		 	%\qquad
	 		 	\begin{split}
	 		 	&\big(\functionANN\big(\parallelizationSpecial_n\big(\longerANN{L,\mathbb{I}_1}({\Phi_1}),\longerANN{L,\mathbb{I}_2}({\Phi_2}),\dots,\longerANN{L,\mathbb{I}_n}({\Phi_n})\big)\big)\big)(x_1,x_2,\dots, x_n)
	 		 	\\&=\Big(\big(\functionANN\big(\longerANN{L,\mathbb{I}_1}({\Phi_1})\big)\big)(x_1), \big(\functionANN\big(\longerANN{L,\mathbb{I}_2}({\Phi_2})\big)\big)(x_2),\dots,
	 		 	\big(\functionANN\big(\longerANN{L,\mathbb{I}_n}({\Phi_n})\big)\big)(x_n) \Big)
	 		 	\\&=\Big((\functionANN(\Phi_1))(x_1),(\functionANN(\Phi_2))(x_2),\dots, (\functionANN(\Phi_n))(x_n)\Big)
	 		 	\end{split}
	 		 	\end{equation}
	 \end{enumerate}
	 (cf.\ Definition~\ref{Definition:simpleParallelization}).
	 Combining this with  \eqref{generalParallelization:DefinitionFormula} and the fact that $L=\allowbreak\max_{j\in\{1,2,\dots,n\}} \allowbreak\lengthANN(\Phi_j)$ ensures
	 \begin{enumerate}
	 	\item[(C)] that 
	 		 \begin{equation}
	 		 \functionANN\big(\!\parallelization_{n,\mathbb{I}}(\Phi)\big)\in C\big(\R^{[\sum_{j=1}^n \inDimANN(\Phi_j)]},\R^{[\sum_{j=1}^n \outDimANN(\Phi_j)]}\big)
	 		 \end{equation}
	 		 and
	 	\item[(D)] 	 that for all $x_1\in\R^{\inDimANN(\Phi_1)},x_2\in\R^{\inDimANN(\Phi_2)},\dots,\allowbreak x_n\in\R^{\inDimANN(\Phi_n)}$ it holds that
	 		\begin{equation}
	 		%\functionANN(\Phi) \in C(\R^{\alpha_{j,0},\R^{\alpha_{j,L})\qandq
	 		%\qquad
	 		\begin{split}
	 		&\big( \functionANN\big(\!\parallelization_{n,\mathbb{I}}(\Phi)\big) \big) (x_1,x_2,\dots,x_n) 
	 		\\&=\big(\functionANN\big(\parallelizationSpecial_n\big(\longerANN{L,\mathbb{I}_1}({\Phi_1}),\longerANN{L,\mathbb{I}_2}({\Phi_2}),\dots,\longerANN{L,\mathbb{I}_n}({\Phi_n})\big)\big)\big)(x_1,x_2,\dots,x_n)
	 		%	\\&=\Big(\big(\functionANN\big(\longerANN{L,\mathbb{I}_1}({\Phi_1})\big)\big)(x_1), \big(\functionANN\big(\longerANN{L,\mathbb{I}_2}({\Phi_2})\big)\big)(x_2),\dots,
	 		%	\big(\functionANN\big(\longerANN{L,\mathbb{I}_n}({\Phi_n})\big)\big)(x_n) \Big)
	 		\\&=\Big((\functionANN(\Phi_1))(x_1),(\functionANN(\Phi_2))(x_2),\dots, (\functionANN(\Phi_n))(x_n)\Big).
	 		\end{split}
	 		\end{equation}
	 \end{enumerate}
This establishes items \eqref{PropertiesOfParallelizationRealization:ItemOne}--\eqref{PropertiesOfParallelizationRealization:ItemTwo}.
	The proof of Corollary~\ref{Lemma:PropertiesOfParallelizationRealization} is thus completed.
\end{proof}

\begin{cor}\label{Lemma:PropertiesOfParallelization}
	Let $n,L\in\N$, $\hiddenDimId_1,\hiddenDimId_2,\dots, \hiddenDimId_n\in\N$,  $\Psi=(\Psi_1,\Psi_2,\dots,\allowbreak \Psi_n),\allowbreak\Phi=(\Phi_1,\Phi_2,\allowbreak\dots, \Phi_n)\in \ANNs^n$
	satisfy for all $j\in\{1,2,\dots, n\}$ that $\dims(\Psi_j) = (\outDimANN(\Phi_j),\hiddenDimId_j,\allowbreak\outDimANN(\Phi_j))$ 
 and
	%		$\functionANN(\idANNshort{j})\in C(\R^{\outDimANN(\Phi_j)},\R^{\outDimANN(\Phi_j)})$
	$L=\max_{k\in\{1,2,\dots,n\}} \lengthANN(\Phi_k)$
	(cf.\ Definition~\ref{Def:ANN}).
	Then
	it holds that
	\begin{equation}
	\begin{split}
	&\paramANN\big(\!\parallelization_{n,\Psi}(\Phi)\big)
	\\&\le \tfrac{1}{2} \bigg(\left[\smallsum\nolimits_{j=1}^n
	\big[\max\big\{1,\tfrac{\hiddenDimId_j}{\outDimANN(\Phi_j)}\big\}\big]\,\paramANN(\Phi_j)
	\,\indicator{(\lengthANN(\Phi_j),\infty)}(L)\right]
	\\&\qquad+\left[\smallsum\nolimits_{j=1}^n\big(
(L-\lengthANN(\Phi_j)-1) \,\hiddenDimId_j\,(\hiddenDimId_j+1)
	+\outDimANN(\Phi_j)\,(\hiddenDimId_j+1)\big)
	\,\indicator{(\lengthANN(\Phi_j),\infty)}(L)\right]
	\\&\qquad+\left[\smallsum\nolimits_{j=1}^n  \paramANN(\Phi_j)\,\indicator{\{\lengthANN(\Phi_j)\}}(L)\right]\bigg)^{\!2}
	\end{split}
	\end{equation}
	(cf.\  Definition~\ref{Definition:generalParallelization}).
\end{cor}

\begin{proof}[Proof of Corollary~\ref{Lemma:PropertiesOfParallelization}]	
Observe that \eqref{generalParallelization:DefinitionFormula}, item~\eqref{PropertiesOfParallelizationEqualLengthDims:Params} in  Proposition~\ref{Lemma:PropertiesOfParallelizationEqualLengthDims}, and
 item~\eqref{PropertiesOfANNenlargementGeometry:ItemLonger} in Lemma~\ref{Lemma:PropertiesOfANNenlargementGeometry} assure that
\begin{equation}
\begin{split}
&\paramANN\big(\!\parallelization_{n,\Psi}(\Phi)\big)
\\&=\paramANN\big(\parallelizationSpecial_n\big(\longerANN{L,\Psi_1}({\Phi_1}),\longerANN{L,\Psi_2}({\Phi_2}),\dots,\longerANN{L,\Psi_n}({\Phi_n})\big)\big)
\\&\le \tfrac{1}{2} \left[\smallsum\nolimits_{j=1}^n \paramANN(\longerANN{L,\Psi_j}({\Phi_j}))\right]^2
%&\le \tfrac{1}{2}  \bigg(\!\left[\smallsum\nolimits_{j=1}^n\left(
%\max\big\{1,\tfrac{\hiddenDimId_j}{\outDimANN(\Phi_j)}\}\,\paramANN(\Phi_j)
%+(L-\lengthANN(\Phi_j)-1) \,\hiddenDimId_j\,(\hiddenDimId_j+1)
%+\outDimANN(\Phi_j)\,(\hiddenDimId_j+1)\right)
%\indicator{(\lengthANN(\Phi_j),\infty)}(L)\right]
%\\&\qquad+\left[\smallsum\nolimits_{j=1}^n  \paramANN(\Phi_j)\indicator{\{\lengthANN(\Phi_j)\}}(L)\right]\!\bigg)^{\!2}.
	\\&\le \tfrac{1}{2} \bigg(\left[\smallsum\nolimits_{j=1}^n
	\big[\max\big\{1,\tfrac{\hiddenDimId_j}{\outDimANN(\Phi_j)}\big\}\big]\,\paramANN(\Phi_j)
	\,\indicator{(\lengthANN(\Phi_j),\infty)}(L)\right]
	\\&\qquad+\left[\smallsum\nolimits_{j=1}^n\big(
	(L-\lengthANN(\Phi_j)-1) \,\hiddenDimId_j\,(\hiddenDimId_j+1)
	+\outDimANN(\Phi_j)\,(\hiddenDimId_j+1)\big)
	\,\indicator{(\lengthANN(\Phi_j),\infty)}(L)\right]
	\\&\qquad+\left[\smallsum\nolimits_{j=1}^n  \paramANN(\Phi_j)\,\indicator{\{\lengthANN(\Phi_j)\}}(L)\right]\bigg)^{\!2}
\end{split}
\end{equation}	
 (cf.\ Definition~\ref{Definition:ANNenlargement} and Definition~\ref{Definition:simpleParallelization}).
%
%	\begin{equation}
%\begin{split}
%&\paramANN(\longerANN{L,\Psi}(\Phi))
%%	\le 
%%\paramANN(\Phi)\, \indicator{\{\lengthANN(\Phi)\}}(L)
%%\\&+\left[\big(\!\max\!\big\{1,\tfrac{\hiddenDimId}{d}\big\}\big)\paramANN(\Phi)+
%%(L-\lengthANN(\Phi)-1) \,\hiddenDimId\,(\hiddenDimId+1)
%%+d\,(\hiddenDimId+1)\right]
%%\indicator{(\lengthANN(\Phi),\infty)}(L)
%\\&\le 
%\begin{cases} \tfrac{1}{2} \big[\smallsum\nolimits_{j=1}^n \paramANN({\Phi_j})\big]^2
%&: \lengthANN(\Phi)=L \\
%\left[\big(\!\max\!\big\{1,\tfrac{\hiddenDimId}{d}\big\}\big)\paramANN(\Phi)+
%(L-\lengthANN(\Phi)-1) \,\hiddenDimId\,(\hiddenDimId+1)
%+d\,(\hiddenDimId+1)\right] &: \lengthANN(\Phi)<L
%\end{cases}
%\end{split}
%\end{equation}
The proof of Corollary~\ref{Lemma:PropertiesOfParallelization} is thus completed.
\end{proof}

%% file: SumsSameLengths.tex
\begin{prop}\label{Lemma:SumsOfANNSequalArchitecture}
	Let $a\in C(\R,\R)$, $M\in\N$,  $h_1,h_2,\dots,h_M\in\R$,
	$\Phi_1,\Phi_2,\dots,\allowbreak \Phi_M\in\ANNs$ satisfy that
	$\dims(\Phi_{1})=\dims(\Phi_{2})=\ldots=\dims(\Phi_{M})$
	(cf.\ Definition~\ref{Def:ANN}).
	Then 
	there exists $\Psi\in \ANNs$ such that
	\begin{enumerate}[(i)]
		\item\label{SumsOfANNSequalArchitecture:ItemOne} it holds that $\paramANN(\Psi)
		\le M^2
		\paramANN(\Phi_1)$,
		\item\label{SumsOfANNSequalArchitecture:ItemTwo} it holds that $\functionANN (\Psi)\in C(\R^{\inDimANN(\Phi_1)},\R^{\outDimANN(\Phi_1)})$, and
		\item\label{SumsOfANNSequalArchitecture:ItemThree} it holds for all 
		$x\in\R^{\inDimANN(\Phi_1)}$ 
		that 
		\begin{equation}
		\begin{split}
		(\functionANN (\Psi))(x)=\smallsum\limits_{m=1}^M h_m \,(\functionANN (\Phi_{m}))(x)
		\end{split}
		\end{equation}
	\end{enumerate}
	(cf.\ Definition~\ref{Definition:ANNrealization}).
\end{prop}

\begin{proof}[Proof of Proposition~\ref{Lemma:SumsOfANNSequalArchitecture}]	
	Throughout this proof 
	%	for all $n\in\N$ let $\operatorname{I}_{n}\in \R^{n\times n}$ be the $n$-dimensional identity matrix,
%	let $\mathfrak{I}=(\mathbb{I}_m)_{m\in\{1,2,\dots, M\}}\subseteq \ANNs$ satisfy for all $m\in\{1,2,\dots, M\}$ that $\mathbb{I}_m=\mathbb{I}$, 
let $d,\mathfrak{d}\in\N$ satisfy that $\inDimANN(\Phi_1)=d$ and $\outDimANN(\Phi_1)=\mathfrak{d}$,
	let $(A_1,b_1)\in \R^{\mathfrak{d}\times (M\mathfrak{d})}\times \R^{\mathfrak{d}}$, $(A_2,b_2)\in \R^{(Md)\times d}\times \R^{Md}$ satisfy that 
	\begin{equation}\label{SumsOfANNSequalArchitectureSigmas}
	A_1=\begin{pmatrix}
	h_1 \operatorname{I}_{\mathfrak{d}} & h_2 \operatorname{I}_{\mathfrak{d}} & \dots & h_M \operatorname{I}_{\mathfrak{d}} 
	\end{pmatrix},\quad 
	A_2=\begin{pmatrix}
	\operatorname{I}_{d}\\\operatorname{I}_{d}\\\vdots\\\operatorname{I}_{d}
	\end{pmatrix}
	, \quad b_1=0, \qandqShort b_2=0
	\end{equation}
	(cf.\ Definition~\ref{Definition:identityMatrix}),
	let $\affineMap_1,\affineMap_2\in\ANNs$ satisfy that $\affineMap_1=(A_1,b_1)$ and  $\affineMap_2=(A_2,b_2)$,
	%	\begin{equation}\label{SumsOfANNSequalArchitectureDefinitionSigmas}
	%	\affineMap_1=\left(\begin{pmatrix}
	%	h_1 & h_2& \dots & h_M
	%	\end{pmatrix}, \begin{pmatrix}
	%	0
	%	\end{pmatrix}\right)
	%	\qandq 
	%	\affineMap_2=(A,b),
	%	\end{equation}
	and let $\Psi\in\ANNs$ satisfy that 
	\begin{equation}\label{SumsOfANNSequalArchitectureNetworkDefinitionPsi}
	\Psi=\compANN{\affineMap_1}{{\compANN{\big[\mathbf{P}_M(\Phi_1,\Phi_2,\dots, \Phi_M)\big]}{\affineMap_2}}}
	\end{equation}
	(cf.\ Definition~\ref{Definition:ANNcomposition}, Definition~\ref{Definition:simpleParallelization}, Lemma~\ref{Lemma:CompositionAssociative}, and Proposition~\ref{Lemma:PropertiesOfParallelizationEqualLength}).
	Note that  \eqref{SumsOfANNSequalArchitectureNetworkDefinitionPsi} and items~\eqref{PropertiesOfCompositionsWithAffineMaps:ItemFront}--\eqref{PropertiesOfCompositionsWithAffineMaps:ItemBehind} in Corollary~\ref{Lemma:PropertiesOfCompositionsWithAffineMaps} demonstrate  that 
	\begin{equation}
	\begin{split}
	\paramANN(\Psi)&
	\le \left[\max\!\left\{1,\tfrac{\outDimANN(\affineMap_1)}{\outDimANN(\mathbf{P}_M(\Phi_1,\Phi_2,\dots, \Phi_M))}\right\}\right] \paramANN\big(\compANN{\big[\mathbf{P}_M(\Phi_1,\Phi_2,\dots, \Phi_M)\big]}{\affineMap_2}\big)
	\\&\le \left[\max\!\left\{1,\tfrac{\outDimANN(\affineMap_1)}{\outDimANN(\mathbf{P}_M(\Phi_1,\Phi_2,\dots, \Phi_M))}\right\}\right]
	\left[\max\!\left\{1,\tfrac{\inDimANN(\affineMap_2)+1}{\inDimANN(\mathbf{P}_M(\Phi_1,\Phi_2,\dots, \Phi_M))+1}\right\}\right]
	\\&\quad\cdot\paramANN\big(\mathbf{P}_M(\Phi_1,\Phi_2,\dots, \Phi_M)\big)
	\\&= \left[\max\!\big\{1,\tfrac{\mathfrak{d}}{M\mathfrak{d}}\big\}\right]
	\left[\max\!\big\{1,\tfrac{d+1}{Md+1}\big\}\right]
	\paramANN\big(\mathbf{P}_M(\Phi_1,\Phi_2,\dots, \Phi_M)\big)
	\\&= 
	\paramANN\big(\mathbf{P}_M(\Phi_1,\Phi_2,\dots, \Phi_M)\big).
	\end{split}	
	\end{equation}
	Corollary~\ref{Lemma:ParallelizationImprovedBoundsOne}  and the hypothesis that 
	for all $m\in\{1,2,\dots, M\}$ it holds that $\dims(\Phi_{m})=\dims(\Phi_{1})$ hence prove that 
	\begin{equation}\label{SumsOfANNSequalArchitectureParamsTwo}
	\begin{split}
	\paramANN(\Psi)
	&\le 
	\paramANN\big(\mathbf{P}_M(\Phi_1,\Phi_2,\dots, \Phi_M)\big)
	\le M^2
	\paramANN(\Phi_1).
	\end{split}	
	\end{equation}	
	Next note that  \eqref{SumsOfANNSequalArchitectureSigmas} and the fact that $\affineMap_2=(A_2,b_2)$ prove that for all $x\in\R^d$ it holds that $\functionANN(\affineMap_2)\in C(\R^{d},\R^{Md})$ and $(\functionANN(\affineMap_2))(x)=(x,x,\dots,x)\in \R^{Md}$.
	Proposition~\ref{Lemma:PropertiesOfParallelizationEqualLength} and 
	item~\eqref{PropertiesOfCompositions:Realization} in Proposition~\ref{Lemma:PropertiesOfCompositions}
	therefore ensure that 
	\begin{equation}\label{SumsOfANNSequalArchitectureCalculationPhiTwoCont}
	\begin{split}
			\functionANN\big(\compANN{(\mathbf{P}_M(\Phi_1,\Phi_2,\dots,\allowbreak \Phi_M))}{\affineMap_2}\big)
			&=\big(\functionANN(\mathbf{P}_M(\Phi_1,\Phi_2,\dots,\allowbreak \Phi_M))\big)\circ (\functionANN({\affineMap_2}))
			\\&\in C\big(\R^{\inDimANN(\affineMap_2)},\R^{\outDimANN(\mathbf{P}_M(\Phi_1,\Phi_2,\dots,\allowbreak \Phi_M))}\big)
			\\&=C\big(\R^{d},\R^{\outDimANN(\Phi_1)+\outDimANN(\Phi_2)+\ldots+\outDimANN(\Phi_M)}\big)
			\\&=C(\R^d,\R^{M\mathfrak{d}})
	\end{split}
	\end{equation}
	and
	\begin{equation}\label{SumsOfANNSequalArchitectureCalculationPhiTwo}
	\begin{split}
	\forallDist x\in\R^d\colon\,
	&\big(\functionANN\big(\compANN{(\mathbf{P}_M(\Phi_1,\Phi_2,\dots, \Phi_M))}{\affineMap_2}\big)\big)(x)
	\\&=\big(\big[\functionANN(\mathbf{P}_M(\Phi_1,\Phi_2,\dots,\allowbreak \Phi_M))\big]\circ [\functionANN({\affineMap_2})]\big)(x)
	\\&=\big(\functionANN(\mathbf{P}_M(\Phi_1,\Phi_2,\dots,\allowbreak \Phi_M))\big)(x,x,\dots,x)
	\\&=\big( (\functionANN(\Phi_1))(x),(\functionANN(\Phi_2))(x),\dots, (\functionANN(\Phi_M))(x) \big).
	\end{split}
	\end{equation}
	Furthermore, observe that \eqref{SumsOfANNSequalArchitectureSigmas} and the fact that $\affineMap_1=(A_1,b_1)$ assure that for all $y_1,y_2,\dots, y_M\in\R^{\mathfrak{d}}$ 
	%$y=(y_1,y_2,\dots, y_M)\in\R^{(M\mathfrak{d})}$ 
	it holds that $\functionANN(\affineMap_1)\in C(\R^{M\mathfrak{d}},\R^{\mathfrak{d}})$ and
	\begin{equation}
	(\functionANN(\affineMap_1))(y_1,y_2,\dots, y_M) 
	%		=\begin{pmatrix}
	%		h_1 & h_2& \dots & h_M
	%		\end{pmatrix} y
	=\smallsum\limits_{m=1}^M h_m y_m.
	\end{equation}
	Combining this and item~\eqref{PropertiesOfCompositions:Realization} in Proposition~\ref{Lemma:PropertiesOfCompositions}
	with  \eqref{SumsOfANNSequalArchitectureNetworkDefinitionPsi},  \eqref{SumsOfANNSequalArchitectureCalculationPhiTwoCont}, and \eqref{SumsOfANNSequalArchitectureCalculationPhiTwo} 
	 demonstrates that for all $x\in\R^d$ it holds that $\functionANN(\Psi)\in C(\R^d,\R^{\mathfrak{d}})$ and
	\begin{equation}
	\begin{split}
	(\functionANN(\Psi))(x)=\smallsum\limits_{m=1}^M h_m(\functionANN(\Phi_m))(x).
	\end{split}
	\end{equation}
	This and \eqref{SumsOfANNSequalArchitectureParamsTwo} establish items~\eqref{SumsOfANNSequalArchitecture:ItemOne}--\eqref{SumsOfANNSequalArchitecture:ItemThree}.
	The proof of Proposition~\ref{Lemma:SumsOfANNSequalArchitecture} is thus completed.
\end{proof}

%% file: SumsDifferentLengths.tex
\begin{prop}\label{Lemma:SumsOfANNS}
	Let $a\in C(\R,\R)$, $M,d,\mathfrak{d},\hiddenDimId,L\in\N$,  $h_1,h_2,\dots,h_M\in\R$,
	$\mathbb{I}, \Phi_1,\Phi_2,\dots,\allowbreak \Phi_M\in\ANNs$ satisfy  for all $m\in\{1,2,\dots, M\}$, $x\in\R^{\mathfrak{d}}$ that 
	$\dims(\mathbb{I}) = (\mathfrak{d},\hiddenDimId,\mathfrak{d})$, $(\functionANN(\mathbb{I}))(x)=x$, 
	$\inDimANN(\Phi_m)=d$, $\outDimANN(\Phi_m)=\mathfrak{d}$, and $L=\max_{m\in\{1,2,\dots,M\}} \lengthANN(\Phi_m)$
	(cf.\ Definition~\ref{Def:ANN} and Definition~\ref{Definition:ANNrealization}).
	Then 
	there exists $\Psi\in \ANNs$ 
 such that
 \begin{enumerate}[(i)]
 	\item \label{SumsOfANNS:Continuity} it holds that $\functionANN (\Psi)\in C(\R^{d},\R^{\mathfrak{d}})$,
 	\item \label{SumsOfANNS:Realization} it holds  for all
 	$x\in\R^d$ that  
 		\begin{equation}
 	\begin{split}
 	(\functionANN (\Psi))(x)=\smallsum\limits_{m=1}^M h_m \,(\functionANN (\Phi_{m}))(x),
 	\end{split}
 	\end{equation}
 	and
 	\item \label{SumsOfANNS:Params} it holds that \begin{equation}
 	\begin{split}
 	\paramANN(\Psi)
 	&\le \tfrac{1}{2} \bigg(\left[\smallsum\nolimits_{m=1}^M
 	\big[\max\big\{1,\tfrac{\hiddenDimId}{\mathfrak{d}}\big\}\big]\,\paramANN(\Phi_m)
 	\,\indicator{(\lengthANN(\Phi_m),\infty)}(L)\right]
 	\\&\qquad+\left[\smallsum\nolimits_{m=1}^M\big(
 	(L-\lengthANN(\Phi_m)-1) \,\hiddenDimId\,(\hiddenDimId+1)
 	+\mathfrak{d}\,(\hiddenDimId+1)\big)
 	\,\indicator{(\lengthANN(\Phi_m),\infty)}(L)\right]
 	\\&\qquad+\left[\smallsum\nolimits_{m=1}^M  \paramANN(\Phi_m)\,\indicator{\{\lengthANN(\Phi_m)\}}(L)\right]\bigg)^{\!2}.
 	\end{split}
 	\end{equation}
 \end{enumerate}		
\end{prop}

%\begin{lemma}\label{Lemma:SumsOfANNS}
%	Let $a\in C(\R,\R)$, $M,\hiddenDimId\in\N$,  $h_1,h_2,\dots,h_M\in\R$,
%	$\mathbb{I}, \Phi_1,\Phi_2,\dots,\allowbreak \Phi_M\in\ANNs$ satisfy
%	for all $x\in\R^{\outDimANN(\Phi_1)}$
%	that 
%	$\dims(\Phi_{1})=\dims(\Phi_{2})=\ldots=\dims(\Phi_{M})$,
%	$\dims(\mathbb{I}) = (\outDimANN(\Phi_1),\hiddenDimId,\outDimANN(\Phi_1))$, and $(\functionANN(\mathbb{I}))(x)=x$
%	(cf.\ Definition~\ref{Def:ANN} and Definition~\ref{Definition:ANNrealization}).
%	Then 
%	there exists $\Psi\in \ANNs$ 
%	which satisfies
%	\begin{enumerate}[(i)]
%		\item that $\functionANN (\Psi)\in C(\R^{d},\R^{\mathfrak{d}})$,
%		\item that for all
%		$x\in\R^d$ it holds that 
%		\begin{equation}
%		\begin{split}
%		(\functionANN (\Psi))(x)=\smallsum\limits_{m=1}^M h_m \,(\functionANN (\Phi_{m}))(x),
%		\end{split}
%		\end{equation}
%		\item and \begin{equation}
%		\begin{split}
%		\paramANN(\Psi)
%		&\le \tfrac{1}{2} \bigg(\left[\smallsum\nolimits_{m=1}^M
%		\max\big\{1,\tfrac{\hiddenDimId}{\mathfrak{d}}\}\,\paramANN(\Phi_m)
%		\indicator{(\lengthANN(\Phi_m),\infty)}(L)\right]
%		\\&\qquad+\left[\smallsum\nolimits_{m=1}^M\big(
%		(L-\lengthANN(\Phi_m)-1) \,\hiddenDimId\,(\hiddenDimId+1)
%		+\mathfrak{d}\,(\hiddenDimId+1)\big)
%		\indicator{(\lengthANN(\Phi_m),\infty)}(L)\right]
%		\\&\qquad+\left[\smallsum\nolimits_{j=1}^n  \paramANN(\Phi_m)\indicator{\{\lengthANN(\Phi_m)\}}(L)\right]\bigg)^{\!2}.
%		\end{split}
%		\end{equation}
%	\end{enumerate}		
%\end{lemma}

\begin{proof}[Proof of Proposition~\ref{Lemma:SumsOfANNS}]	
	Throughout this proof 
	%	for all $n\in\N$ let $\operatorname{I}_{n}\in \R^{n\times n}$ be the $n$-dimensional identity matrix,
	let $\mathfrak{I}=(\mathfrak{I}_1,\mathfrak{I}_2,\dots, \mathfrak{I}_M)\in \ANNs^M$ satisfy for all $m\in\{1,2,\dots, M\}$ that $\mathfrak{I}_m=\mathbb{I}$, 
	let $(A_1,b_1)\in \R^{\mathfrak{d}\times (M\mathfrak{d})}\times \R^{\mathfrak{d}}$, $(A_2,b_2)\in \R^{(Md)\times d}\times \R^{Md}$ satisfy that 
	\begin{equation}\label{SumsOfANNsSigmas}
	A_1=\begin{pmatrix}
	h_1 \operatorname{I}_{\mathfrak{d}} & h_2 \operatorname{I}_{\mathfrak{d}} & \dots & h_M \operatorname{I}_{\mathfrak{d}} 
	\end{pmatrix},\quad 
	A_2=\begin{pmatrix}
	\operatorname{I}_{d}\\\operatorname{I}_{d}\\\vdots\\\operatorname{I}_{d}
	\end{pmatrix}
	, \quad b_1=0, \qandqShort b_2=0
	\end{equation}
	(cf.\ Definition~\ref{Definition:identityMatrix}),
	let $\affineMap_1,\affineMap_2\in\ANNs$ satisfy that $\affineMap_1=(A_1,b_1)$ and  $\affineMap_2=(A_2,b_2)$,
	%	\begin{equation}\label{SumsOfANNSDefinitionSigmas}
	%	\affineMap_1=\left(\begin{pmatrix}
	%	h_1 & h_2& \dots & h_M
	%	\end{pmatrix}, \begin{pmatrix}
	%	0
	%	\end{pmatrix}\right)
	%	\qandq 
	%	\affineMap_2=(A,b),
	%	\end{equation}
	and let $\Psi\in\ANNs$ satisfy that 
	\begin{equation}\label{SumsOfANNSNetworkDefinitionPhi}
	\Psi=\compANN{\affineMap_1}{{\compANN{\big(\!\operatorname{P}_{M,\mathfrak{I}}(\Phi_1,\Phi_2,\dots, \Phi_M)\big)}{\affineMap_2}}}
	\end{equation}
	(cf.\ Definition~\ref{Definition:ANNcomposition}, Definition~\ref{Definition:generalParallelization},  Lemma~\ref{Lemma:CompositionAssociative}, and Corollary~\ref{Lemma:PropertiesOfParallelizationRealization}).
	Note that \eqref{SumsOfANNSNetworkDefinitionPhi} and items~\eqref{PropertiesOfCompositionsWithAffineMaps:ItemFront}--\eqref{PropertiesOfCompositionsWithAffineMaps:ItemBehind} in Corollary~\ref{Lemma:PropertiesOfCompositionsWithAffineMaps} demonstrate  that 
	\begin{equation}
	\begin{split}
	\paramANN(\Psi)&\le \left[\max\!\left\{1,\tfrac{\outDimANN(\affineMap_1)}{\outDimANN(\operatorname{P}_{M,\mathfrak{I}}(\Phi_1,\Phi_2,\dots, \Phi_M))}\right\}\right]
	\left[\max\!\left\{1,\tfrac{\inDimANN(\affineMap_2)+1}{\inDimANN(\operatorname{P}_{M,\mathfrak{I}}(\Phi_1,\Phi_2,\dots, \Phi_M))+1}\right\}\right]
	\\&\quad\cdot\paramANN\big(\!\operatorname{P}_{M,\mathfrak{I}}(\Phi_1,\Phi_2,\dots, \Phi_M)\big)
	\\&= \left[\max\!\big\{1,\tfrac{\mathfrak{d}}{M\mathfrak{d}}\big\}\right]
	\left[\max\!\big\{1,\tfrac{d+1}{Md+1}\big\}\right]
	\paramANN\big(\!\operatorname{P}_{M,\mathfrak{I}}(\Phi_1,\Phi_2,\dots, \Phi_M)\big)
	\\&=
	\paramANN\big(\!\operatorname{P}_{M,\mathfrak{I}}(\Phi_1,\Phi_2,\dots, \Phi_M)\big).
	\end{split}	
	\end{equation}
	Corollary~\ref{Lemma:PropertiesOfParallelization} hence proves that 
	\begin{equation}\label{SumsOfANNSParamsOne}
	\begin{split}
	&\paramANN(\Psi)
	\le 
	\paramANN\big(\!\operatorname{P}_{M,\mathfrak{I}}(\Phi_1,\Phi_2,\dots, \Phi_M)\big)
%	\\&\le \tfrac{1}{2}
%	\bigg(\!\left[\smallsum\nolimits_{m=1}^M\left(
%	\max\big\{1,\tfrac{\hiddenDimId}{\mathfrak{d}}\}\,\paramANN(\Phi_m)
%	+(L-\lengthANN(\Phi_m)-1) \,\hiddenDimId\,(\hiddenDimId+1)
%	+\mathfrak{d}\,(\hiddenDimId+1)\right)
%	\,\indicator{(\lengthANN(\Phi_m),\infty)}(L)\right]
%	\\&\qquad+\left[\smallsum\nolimits_{j=1}^n  \paramANN(\Phi_m)\,\indicator{\{\lengthANN(\Phi_m)\}}(L)\right]\!\bigg)^{\!2}.
	\\&\le \tfrac{1}{2} \bigg(\left[\smallsum\nolimits_{m=1}^M
	\big[\max\big\{1,\tfrac{\hiddenDimId}{\mathfrak{d}}\big\}\big]\,\paramANN(\Phi_m)
	\,\indicator{(\lengthANN(\Phi_m),\infty)}(L)\right]
	\\&\qquad+\left[\smallsum\nolimits_{m=1}^M\big(
	(L-\lengthANN(\Phi_m)-1) \,\hiddenDimId\,(\hiddenDimId+1)
	+\mathfrak{d}\,(\hiddenDimId+1)\big)
	\,\indicator{(\lengthANN(\Phi_m),\infty)}(L)\right]
	\\&\qquad+\left[\smallsum\nolimits_{m=1}^M  \paramANN(\Phi_m)\,\indicator{\{\lengthANN(\Phi_m)\}}(L)\right]\bigg)^{\!2}.
	\end{split}	
	\end{equation}	
%	This establishes item~\eqref{SumsOfANNS:Params}.
	Next note that \eqref{SumsOfANNsSigmas} and the fact that $\affineMap_2=(A_2,b_2)$ prove that  $\functionANN(\affineMap_2)\in C(\R^{d},\R^{Md})$ and
	\begin{equation}
	\forallDist x\in\R^d\colon\,	(\functionANN(\affineMap_2))(x)=(x,x,\dots,x)\in \R^{Md}.
	\end{equation}
%	 $(\functionANN(\affineMap_2))(x)=(x,x,\dots,x)\in \R^{Md}$.
	Corollary~\ref{Lemma:PropertiesOfParallelizationRealization} and item~\eqref{PropertiesOfCompositions:Realization} in Proposition~\ref{Lemma:PropertiesOfCompositions} 
	therefore ensure that 
	\begin{equation}\label{SumsOfANNSCalculationPhiOneCont}
	\begin{split}
			\functionANN\big(\compANN{\operatorname{P}_{M,\mathfrak{I}}(\Phi_1,\Phi_2,\dots, \Phi_M)}{\affineMap_2}\big)
			&\in C\big(\R^{\inDimANN(\affineMap_2)},\R^{\outDimANN(\operatorname{P}_{M,\mathfrak{I}}(\Phi_1,\Phi_2,\dots, \Phi_M))}\big)
			\\&=C(\R^d,\R^{M\mathfrak{d}})
	\end{split}
	\end{equation}
	and
	\begin{equation}\label{SumsOfANNSCalculationPhiOne}
	\begin{split}
	\forallDist x\in\R^d\colon\,&\big(\functionANN\big(\compANN{\operatorname{P}_{M,\mathfrak{I}}(\Phi_1,\Phi_2,\dots, \Phi_M)}{\affineMap_2}\big)\big)(x)
	\\&=\big( (\functionANN(\Phi_1))(x),(\functionANN(\Phi_2))(x),\dots, (\functionANN(\Phi_M))(x) \big).
	\end{split}
	\end{equation}
	In addition, observe that \eqref{SumsOfANNsSigmas}  and the fact that $\affineMap_1=(A_1,b_1)$ assure that 
	%$y=(y_1,y_2,\dots, y_M)\in\R^{(M\mathfrak{d})}$ 
	 $\functionANN(\affineMap_1)\in C(\R^{M\mathfrak{d}},\R^{\mathfrak{d}})$ and
	\begin{equation}
	\forallDist y_1,y_2,\dots, y_M\in\R^{\mathfrak{d}}\colon\,(\functionANN(\affineMap_1))(y_1,y_2,\dots, y_M) 
	%		=\begin{pmatrix}
	%		h_1 & h_2& \dots & h_M
	%		\end{pmatrix} y
	=\smallsum\limits_{m=1}^M h_m y_m.
	\end{equation}
	Combining this, \eqref{SumsOfANNSCalculationPhiOneCont}, \eqref{SumsOfANNSCalculationPhiOne}, and \eqref{SumsOfANNSNetworkDefinitionPhi} with item~\eqref{PropertiesOfCompositions:Realization} in Proposition~\ref{Lemma:PropertiesOfCompositions} demonstrates that
	\begin{equation}
		\functionANN(\Psi)\in C\big(\R^{\inDimANN((\compANN{\operatorname{P}_{M,\mathfrak{I}}(\Phi_1,\Phi_2,\dots, \Phi_M))}{\affineMap_2})},\R^{\outDimANN(\affineMap_1)}\big)
		= C(\R^d,\R^{\mathfrak{d}})
	\end{equation}
	and
	\begin{equation}
	\forallDist x\in\R^d\colon\,(\functionANN(\Psi))(x)=\smallsum\limits_{m=1}^M h_m(\functionANN(\Phi_m))(x).
	\end{equation}
	This and \eqref{SumsOfANNSParamsOne} establish items~\eqref{SumsOfANNS:Continuity}--\eqref{SumsOfANNS:Params}.
	The proof of Proposition~\ref{Lemma:SumsOfANNS} is thus completed.
\end{proof}

%% file: ANNrepOneEuler.tex
\begin{lemma}\label{Lemma:CompositionSum}
	Let $a\in C(\R,\R)$, $L_1\in \N\cap [2,\infty)$, $L_2\in\N$, $d,\hiddenDimId, l_{1,0},l_{1,1},\dots,\allowbreak l_{1,L_1},l_{2,0},l_{2,1},\dots,\allowbreak l_{2,L_2}\in\N$,
	$\mathbb{I}, \Phi_1,\Phi_2\in\ANNs$
	satisfy  for all $k\in\{1,2\}$, $x\in\R^{d}$ that  $\dims(\mathbb{I}) = (d,\hiddenDimId,d)$, $(\functionANN(\mathbb{I}))(x)=x$, $\inDimANN(\Phi_k)=\outDimANN(\Phi_k)=d$, and
	$\dims(\Phi_k)=(l_{k,0},l_{k,1},\dots,\allowbreak l_{k,L_k})$
	%	\begin{equation}\label{CorCompositionSum:HypothesisDims}
	%	\dims(\Phi_k)=(l_0^k,l_1^k,\dots, l_{L_k}^k),\qandq \dims(\mathbb{I}) = (d,\hiddenDimId,d)
	%	\end{equation}
	(cf.\ Definition~\ref{Def:ANN} and Definition~\ref{Definition:ANNrealization}).
	Then there exists $\Psi\in \ANNs$ such that
	 \begin{enumerate}[(i)]
	 	\item \label{CompositionSum:Continuity} it holds  that
	 	$\functionANN(\Psi)\in C(\R^d,\R^d)$,
	 	\item \label{CompositionSum:Realization} it holds for all $x\in\R^d$ that
	 	\begin{equation} (\functionANN(\Psi))(x)=(\functionANN(\Phi_2))(x)+\big((\functionANN(\Phi_1))\circ (\functionANN(\Phi_2))\big)(x),
	 	\end{equation}
	 	\item \label{CompositionSum:Dims} it holds that
	 			\begin{equation}
	 			\dims(\Psi)=
	 			(l_{2,0},l_{2,1},\dots, l_{2,L_2-1},l_{1,1}+\hiddenDimId,l_{1,2}+\hiddenDimId,\dots,l_{1,L_1-1}+\hiddenDimId,l_{1,L_1} ),
	 			\end{equation}
	 			and
	 	\item \label{CompositionSum:Params}
%	 		\begin{equation}
%	 		\begin{split}
%	 		\paramANN(\Psi)
%	 		&\le   
%	 		\paramANN(\Phi_2)+
%	 		\paramANN(\Phi_1)+\paramANN(\mathbb{I})+(L_1-2)\hiddenDimId(\hiddenDimId+1)
%	 		+l_{1,1} d 
%	 		\\&\quad+\hiddenDimId \bigg[\smallsum\limits_{m=2}^{L_1} l_{1,m} \bigg]
%	 		+ \hiddenDimId \bigg[\smallsum\limits_{m=0}^{L_1-2} l_{1,m} \bigg]
%	 		+d \,l_{1,L_1-1}
%	 		+ (l_{1,1}+\hiddenDimId) ( l_{2,L_2-1}+1).
%	 		\end{split}
%	 		\end{equation}
it holds that
%\begin{multline} 
%		\paramANN(\Psi)
%			 		\le   
%			 		\paramANN(\Phi_2)+
%			 		\paramANN(\Phi_1)+\paramANN(\mathbb{I})+(L_1-2)\hiddenDimId(\hiddenDimId+1)
%			 		+l_{1,1} d 
%			 		\\+\hiddenDimId \bigg[\smallsum\limits_{m=2}^{L_1} l_{1,m} \bigg]
%			 		+ \hiddenDimId \bigg[\smallsum\limits_{m=0}^{L_1-2} l_{1,m} \bigg]
%			 		+d \,l_{1,L_1-1}
%			 		+ (l_{1,1}+\hiddenDimId) ( l_{2,L_2-1}+1).
%\end{multline}
\begin{equation}
\begin{split}
\paramANN(\Psi)
&=
\paramANN(\Phi_1)+\paramANN(\Phi_2)
+ (\hiddenDimId-d)(l_{2,L_2-1}+1)+l_{1,1} (l_{2,L_2-1}-d)
\\&\quad+(L_1-2)\hiddenDimId(\hiddenDimId+1)
+\hiddenDimId\bigg[\smallsum\limits_{m=2}^{L_1} l_{1,m} \bigg]
+\hiddenDimId\bigg[\smallsum\limits_{m=1}^{L_1-2} l_{1,m} \bigg]
.
\end{split}
\end{equation}
	 	\end{enumerate}
\end{lemma}

\begin{proof}[Proof of Lemma~\ref{Lemma:CompositionSum}]
	Throughout this proof 
%	let $\idMatrix_{d}\in \R^{d\times d}$ be the $d$-dimensional identity matrix,
	let $A_1\in\R^{d\times 2d}$, $A_2\in\R^{2d\times d}$, $b_1\in\R^d$, $b_2\in\R^{2d}$ satisfy that 
	\begin{equation}\label{CompositionSum:Matrices}
	A_1=\begin{pmatrix}
	\idMatrix_{d}& \idMatrix_{d}
	\end{pmatrix}, \qquad
	A_2=\begin{pmatrix}
	\idMatrix_{d}\\ \idMatrix_{d}
	\end{pmatrix},\qquad
	b_1=0,\qandq 
	b_2=0
	\end{equation}
	(cf.\ Definition~\ref{Definition:identityMatrix})
	and let
	$\affineMap_1\in (\R^{d\times 2d}\times \R^d)\subseteq \ANNs$, $\affineMap_2\in (\R^{2d\times d}\times \R^{2d})\subseteq \ANNs$, $\Psi\in\ANNs$ satisfy that $\affineMap_1=(A_1,b_1)$, $\affineMap_2=(A_2,b_2)$, and
	%	\begin{equation}
	%		\Psi=\compANN{\affineMap_1}{
	%			\Big(\compANN{\big[\parallelizationSpecial_{2}(\Phi_1,\idANN{d}{L_1})\big]}
	%			{\big[\compANN{\affineMap_2}
	%			{\Phi_2}\big]}\Big)
	%		}
	%	\end{equation}
	\begin{equation}\label{CompositionSum:ANNdefinition}
	\Psi=
	\compANN
	{\affineMap_1}
	{\compANN
		{\compANN
			{\big[\parallelizationSpecial_{2}\big(\Phi_1,\mathbb{I}^{\bullet (L_1-1)}\big)\big]}
			{\affineMap_2}}
		{\Phi_2}
	}
	\end{equation}
	(cf.\ Definition~\ref{Definition:ANNcomposition}, Definition~\ref{Definition:iteratedANNcomposition}, Definition~\ref{Definition:simpleParallelization},   Lemma~\ref{Lemma:CompositionAssociative}, and item~\eqref{PropertiesOfANNenlargementGeometry:BulletPower} in Lemma~\ref{Lemma:PropertiesOfANNenlargementGeometry}).	
	Observe that \eqref{CompositionSum:Matrices} and the fact that $\affineMap_2=(A_2,b_2)$ ensure that for all $x\in\R^d$ it holds that 
	\begin{equation}
		\functionANN(\affineMap_2)\in C(\R^d,\R^{2d}) \qandq (\functionANN(\affineMap_2))(x)=(x,x).
	\end{equation}
%	$\functionANN(\affineMap_1)\in C(\R^d,\R^{2d})$ and $(\functionANN(\affineMap_1))(x)=(x,x)$. 
	Item~\eqref{PropertiesOfCompositions:Realization} in Proposition~\ref{Lemma:PropertiesOfCompositions}, item~\eqref{PropertiesOfANNenlargementRealization: itemOne} in Lemma~\ref{Lemma:PropertiesOfANNenlargementRealization}, and Proposition~\ref{Lemma:PropertiesOfParallelizationEqualLength} hence imply that for all $x\in\R^d$ it holds that
	$\functionANN([\compANN{\parallelizationSpecial_{2}(\Phi_1,\mathbb{I}^{\bullet (L_1-1)})]}{\affineMap_2})\allowbreak\in C(\R^d,\R^{2d})$ and
	\begin{equation}
		\begin{split}
		&\big(\functionANN\big(\big[\compANN{\parallelizationSpecial_{2}\big(\Phi_1,\mathbb{I}^{\bullet (L_1-1)}\big)\big]}{\affineMap_2}\big)\big)(x)
		=\big(\functionANN\big(\parallelizationSpecial_{2}\big(\Phi_1,\mathbb{I}^{\bullet (L_1-1)}\big)\big)\big)(x,x)
		\\&=\big((\functionANN(\Phi_1))(x),(\functionANN(\mathbb{I}^{\bullet (L_1-1)}))(x)\big)
		=\big((\functionANN(\Phi_1))(x),x\big).
		\end{split}
	\end{equation}
	Item~\eqref{PropertiesOfCompositions:Realization} in Proposition~\ref{Lemma:PropertiesOfCompositions} therefore demonstrates that for all $x\in\R^d$ it holds that
		$\functionANN\big(\compANN
		{\compANN
			{\big[\parallelizationSpecial_{2}\big(\Phi_1,\mathbb{I}^{\bullet (L_1-1)}\big)\big]}
			{\affineMap_2}}
		{\Phi_2}\big)\in C(\R^d,\R^{2d})$ and
		\begin{equation}\label{CompositionSum:RealizationOfParall}
		\begin{split}
		&\big(\functionANN\big(\compANN
		{\compANN
			{\big[\parallelizationSpecial_{2}\big(\Phi_1,\mathbb{I}^{\bullet (L_1-1)}\big)\big]}
			{\affineMap_2}}
		{\Phi_2}\big)\big)(x)
		=\Big((\functionANN(\Phi_1))\big((\functionANN(\Phi_2))(x)\big),(\functionANN(\Phi_2))(x)\Big).
		\end{split}
		\end{equation}
		In addition, note that \eqref{CompositionSum:Matrices} and the fact that $\affineMap_1=(A_1,b_1)$ ensure that for all $y=(y_1,y_2)\in\R^{d}\times \R^d$ it holds that
		\begin{equation}
			\functionANN(\affineMap_1)\in C(\R^{2d},\R^d)\qandq (\functionANN(\affineMap_1))(y)=y_1+y_2.
		\end{equation}
%		 $\functionANN(\affineMap_2)\in C(\R^{2d},\R^d)$ and $(\functionANN(\affineMap_2))(y)=y_1+y_2$.
	Item~\eqref{PropertiesOfCompositions:Realization} in Proposition~\ref{Lemma:PropertiesOfCompositions}, \eqref{CompositionSum:ANNdefinition}, and \eqref{CompositionSum:RealizationOfParall} hence prove that 
	 for all $x\in\R^d$ it holds that
	 $\functionANN(\Psi)\in C(\R^d,\R^{d})$ and
	 \begin{equation}\label{CompositionSum:RealizationProof}
	 \begin{split}
	 &(\functionANN(\Psi))(x)=
	 (\functionANN(\Phi_1))\big((\functionANN(\Phi_2))(x)\big)+(\functionANN(\Phi_2))(x).
	 \end{split}
	 \end{equation}
%	This establishes item~\eqref{CompositionSum:Realization}.
	Next note that item~\eqref{PropertiesOfANNenlargementGeometry:BulletPower} in Lemma~\ref{Lemma:PropertiesOfANNenlargementGeometry} and item~\eqref{PropertiesOfParallelizationEqualLengthDims:Dims}	in Proposition~\ref{Lemma:PropertiesOfParallelizationEqualLengthDims} demonstrate that 
	\begin{equation}
	\dims\big(\parallelizationSpecial_{2}\big(\Phi_1,\mathbb{I}^{\bullet (L_1-1)}\big)\big)
	=
	(2d,l_{1,1}+\hiddenDimId,l_{1,2}+\hiddenDimId,\dots, l_{1,L_1-1}+\hiddenDimId,2d).	
	\end{equation} 
%	and
%	\begin{equation}\label{CompositionSum:ParallelizationEstimate}
%	\begin{split}
%	\paramANN\big(\parallelizationSpecial_{2}\big(\Phi_1,\mathbb{I}^{\bullet (L_1-1)}\big)\big)&=\paramANN(\Phi_1)+\paramANN\big(\mathbb{I}^{\bullet (L_1-1)}\big)
%	+l_{1,1} d 
%	+
%	\hiddenDimId \bigg[\smallsum\limits_{m=2}^{L_1} l_{1,m} \bigg]
%	\\&\quad+ \hiddenDimId \bigg[\smallsum\limits_{m=0}^{L_1-2} l_{1,m} \bigg]
%	+d \,l_{1,L_1-1}.
%	%	\\&\le\paramANN(\Phi_1)+\paramANN\big(\mathbb{I}^{\bullet (L_1-1)}\big)
%	%	+
%	%	\max\{d,\hiddenDimId\}\, \bigg[\smallsum\limits_{m=1}^{L_1} l_{1,m} \bigg]
%	%	\\&\quad+ \max\{d,\hiddenDimId\}\, \bigg[\smallsum\limits_{m=0}^{L_1-1} l_{1,m} \bigg]
%	%		\\&\le \big(1+2\max\{d,\hiddenDimId\}\big)\,\paramANN(\Phi_1)+\paramANN\big(\mathbb{I}^{\bullet (L_1-1)}\big).
%	\end{split}
%	\end{equation}
	Item~\eqref{PropertiesOfCompositions:Dims} in Proposition~\ref{Lemma:PropertiesOfCompositions} therefore ensures that 
	%\begin{equation}
	%	\begin{split}
	%	\paramANN\big(\compANN{\affineMap_2}{\Phi_2}\big)
	%	\le   \max\!\left\{1,\tfrac{\outDimANN(\affineMap_2)}{\outDimANN(\Phi_2)}\right\} \paramANN(\Phi_2)
	%	=2 \paramANN(\Phi_2)
	%	\end{split}
	%\end{equation}
	\begin{multline}
		\dims\big(\compANN{\affineMap_1}{\big[\compANN{\parallelizationSpecial_{2}\big(\Phi_1,\mathbb{I}^{\bullet (L_1-1)}\big)\big]}{\affineMap_2}}\big)
		= 
		(d,l_{1,1}+\hiddenDimId,l_{1,2}+\hiddenDimId,\dots, l_{1,L_1-1}+\hiddenDimId,d).	
	\end{multline} 
	Combining this with item~\eqref{PropertiesOfCompositions:Dims} in Proposition~\ref{Lemma:PropertiesOfCompositions}, \eqref{CompositionSum:ANNdefinition}, and  the fact that $\outDimANN(\Phi_2)=l_{2,L_2}=d$ shows that
	 \begin{equation}\label{CompositionSum:DimsProof}
	 	 \dims(\Psi)=
	 	 (l_{2,0},l_{2,1},\dots, l_{2,L_2-1},l_{1,1}+\hiddenDimId,l_{1,2}+\hiddenDimId,\dots,l_{1,L_1-1}+\hiddenDimId,d ).
	 \end{equation} 
%	Next we claim that 
%		 			\begin{equation}\label{CompositionSum:itemIIIproof}
%		 			\paramANN(\Psi)= 
%		 			\begin{cases}
%		 			\begin{array}{r}
%		 			\paramANN(\Phi_1)+\paramANN(\Phi_2)
%		 			+ (\hiddenDimId-d)(l_{2,L_2-1}+1)+l_{1,1} (l_{2,L_2-1}-d)
%		 			\\+(L_1-2)\hiddenDimId(\hiddenDimId+1)
%		 			+\hiddenDimId\bigg[\smallsum\limits_{m=2}^{L_1} l_{1,m} \bigg]
%		 			+\hiddenDimId\bigg[\smallsum\limits_{m=1}^{L_1-2} l_{1,m} \bigg]
%		 			\end{array}
%		 			&: L_1>1\\[5ex]
%		 			\paramANN(\Phi_2) &: L_1=1
%		 			\end{cases}.
%		 			\end{equation}
%	To prove \eqref{CompositionSum:itemIIIproof} we distinguish between the cases $L_1=1$ and $L_1>1$.
%	Note that \eqref{CompositionSum:DimsProof} implies \eqref{CompositionSum:itemIIIproof} in the case $L_1=1$.
%	We now prove \eqref{CompositionSum:itemIIIproof} in the case $L_1>1$. 
	The fact that $l_{1,L_1}=\outDimANN(\Phi_1)=d$ hence ensures that 
	\begin{equation}
	\begin{split}
	\paramANN(\Psi)
	&=  \bigg[\smallsum\limits_{m=1}^{L_2-1} l_{2,m}(l_{2,m-1}+1) \bigg]
	+ (l_{1,1}+\hiddenDimId)(l_{2,L_2-1}+1)
	\\&\quad+\bigg[\smallsum\limits_{m=2}^{L_1-1} (l_{1,m}+\hiddenDimId)(l_{1,m-1}+\hiddenDimId+1) \bigg]
	+d(l_{1,L_1-1}+\hiddenDimId+1)
	\\&=\paramANN(\Phi_2)-l_{2,L_2}(l_{2,L_2-1}+1)
	+ (l_{1,1}+\hiddenDimId)(l_{2,L_2-1}+1)
	\\&\quad+\hiddenDimId\bigg[\smallsum\limits_{m=2}^{L_1-1} l_{1,m} \bigg]
	+\hiddenDimId\bigg[\smallsum\limits_{m=2}^{L_1-1} l_{1,m-1} \bigg]
	+\bigg[\smallsum\limits_{m=2}^{L_1-1} l_{1,m}(l_{1,m-1}+1) \bigg]
	\\&\quad+(L_1-2)\hiddenDimId(\hiddenDimId+1)+l_{1,L_1}(l_{1,L_1-1}+1)+l_{1,L_1} \hiddenDimId.
%	\\&=\paramANN(\Phi_2)
%	+ (l_{1,1}-d+\hiddenDimId)(l_{2,L_2-1}+1)
%	+(L_1-2)\hiddenDimId(\hiddenDimId+1)
%	\\&\quad+\hiddenDimId\bigg[\smallsum\limits_{m=2}^{L_1} l_{1,m} \bigg]
%	+\hiddenDimId\bigg[\smallsum\limits_{m=1}^{L_1-2} l_{1,m} \bigg]
%	+\paramANN(\Phi_1)
%	-l_{1,1}(l_{1,0}+1)
%	\\&=\paramANN(\Phi_1)+\paramANN(\Phi_2)
%	+ (\hiddenDimId-d)(l_{2,L_2-1}+1)+l_{1,1} (l_{2,L_2-1}-d)
%	\\&\quad+(L_1-2)\hiddenDimId(\hiddenDimId+1)
%	+\hiddenDimId\bigg[\smallsum\limits_{m=2}^{L_1} l_{1,m} \bigg]
%	+\hiddenDimId\bigg[\smallsum\limits_{m=1}^{L_1-2} l_{1,m} \bigg].
	\end{split}
	\end{equation}
	This, the fact that $l_{2,L_2}=\outDimANN(\Phi_2)=d$, and the fact that $l_{1,0}=\inDimANN(\Phi_1)=d$ demonstrate that 
		\begin{equation}
		\begin{split}
		\paramANN(\Psi)
%		\\&=  \bigg[\smallsum\limits_{m=1}^{L_2-1} l_{2,m}(l_{2,m-1}+1) \bigg]
%		+ (l_{1,1}+\hiddenDimId)(l_{2,L_2-1}+1)
%		\\&\quad+\bigg[\smallsum\limits_{m=2}^{L_1-1} (l_{1,m}+\hiddenDimId)(l_{1,m-1}+\hiddenDimId+1) \bigg]
%		+d(l_{1,L_1-1}+\hiddenDimId+1)
%		\\&=\paramANN(\Phi_2)-l_{2,L_2}(l_{2,L_2-1}+1)
%		+ (l_{1,1}+\hiddenDimId)(l_{2,L_2-1}+1)
%		+(L_1-2)\hiddenDimId(\hiddenDimId+1)
%		\\&\quad+\hiddenDimId\bigg[\smallsum\limits_{m=2}^{L_1-1} l_{1,m} \bigg]
%		+\hiddenDimId\bigg[\smallsum\limits_{m=2}^{L_1-1} l_{1,m-1} \bigg]
%		+\bigg[\smallsum\limits_{m=2}^{L_1-1} l_{1,m}(l_{1,m-1}+1) \bigg]
%		\\&\quad+l_{1,L_1}(l_{1,L_1-1}+1)+l_{1,L_1} \hiddenDimId
		&=\paramANN(\Phi_2)
		+ (l_{1,1}-d+\hiddenDimId)(l_{2,L_2-1}+1)
		+(L_1-2)\hiddenDimId(\hiddenDimId+1)
		\\&\quad+\hiddenDimId\bigg[\smallsum\limits_{m=2}^{L_1} l_{1,m} \bigg]
		+\hiddenDimId\bigg[\smallsum\limits_{m=1}^{L_1-2} l_{1,m} \bigg]
		+\paramANN(\Phi_1)
		-l_{1,1}(l_{1,0}+1)
		\\&=\paramANN(\Phi_1)+\paramANN(\Phi_2)
		+ (\hiddenDimId-d)(l_{2,L_2-1}+1)+l_{1,1} (l_{2,L_2-1}-d)
		\\&\quad+(L_1-2)\hiddenDimId(\hiddenDimId+1)
		+\hiddenDimId\bigg[\smallsum\limits_{m=2}^{L_1} l_{1,m} \bigg]
		+\hiddenDimId\bigg[\smallsum\limits_{m=1}^{L_1-2} l_{1,m} \bigg].
		\end{split}
		\end{equation}
%	The fact that $l_{1,L_1}=\outDimANN(\Phi_1)=d=\outDimANN(\Phi_2)=l_{2,L_2}$ and the fact that $l_{1,0}=\inDimANN(\Phi_1)=d$ hence ensure that 
%	\begin{equation}
%		\begin{split}
%			&\paramANN(\Psi)
%			\\&=  \bigg[\smallsum\limits_{m=1}^{L_2-1} l_{2,m}(l_{2,m-1}+1) \bigg]
%			+ (l_{1,1}+\hiddenDimId)(l_{2,L_2-1}+1)
%			\\&\quad+\bigg[\smallsum\limits_{m=2}^{L_1-1} (l_{1,m}+\hiddenDimId)(l_{1,m-1}+\hiddenDimId+1) \bigg]
%			+d(l_{1,L_1-1}+\hiddenDimId+1)
%			\\&=\paramANN(\Phi_2)-l_{2,L_2}(l_{2,L_2-1}+1)
%			+ (l_{1,1}+\hiddenDimId)(l_{2,L_2-1}+1)
%			+(L_1-2)\hiddenDimId(\hiddenDimId+1)
%			\\&\quad+\hiddenDimId\bigg[\smallsum\limits_{m=2}^{L_1-1} l_{1,m} \bigg]
%			+\hiddenDimId\bigg[\smallsum\limits_{m=2}^{L_1-1} l_{1,m-1} \bigg]
%			+\bigg[\smallsum\limits_{m=2}^{L_1-1} l_{1,m}(l_{1,m-1}+1) \bigg]
%			\\&\quad+l_{1,L_1}(l_{1,L_1-1}+1)+l_{1,L_1} \hiddenDimId
%						\\&=\paramANN(\Phi_2)
%						+ (l_{1,1}-d+\hiddenDimId)(l_{2,L_2-1}+1)
%						+(L_1-2)\hiddenDimId(\hiddenDimId+1)
%						\\&\quad+\hiddenDimId\bigg[\smallsum\limits_{m=2}^{L_1} l_{1,m} \bigg]
%						+\hiddenDimId\bigg[\smallsum\limits_{m=1}^{L_1-2} l_{1,m} \bigg]
%						+\paramANN(\Phi_1)
%						-l_{1,1}(l_{1,0}+1)
%								\\&=\paramANN(\Phi_1)+\paramANN(\Phi_2)
%								+ (\hiddenDimId-d)(l_{2,L_2-1}+1)+l_{1,1} (l_{2,L_2-1}-d)
%								\\&\quad+(L_1-2)\hiddenDimId(\hiddenDimId+1)
%								+\hiddenDimId\bigg[\smallsum\limits_{m=2}^{L_1} l_{1,m} \bigg]
%								+\hiddenDimId\bigg[\smallsum\limits_{m=1}^{L_1-2} l_{1,m} \bigg].
%		\end{split}
%	\end{equation}
Combining this with \eqref{CompositionSum:RealizationProof} and \eqref{CompositionSum:DimsProof} establishes items~\eqref{CompositionSum:Continuity}--\eqref{CompositionSum:Params}.
	The proof of Lemma~\ref{Lemma:CompositionSum} is thus completed.
\end{proof}

%\begin{cor}\label{Cor:CompositionSum}
%			Let $a\in C(\R,\R)$, $L_1,L_2\in\N$, $d,\hiddenDimId, l_{1,0},l_{1,1},\dots,\allowbreak l_{1,L_1},l_{2,0},l_{2,1},\dots,\allowbreak l_{2,L_2}\in\N$,
%	$\mathbb{I}, \Phi_1,\Phi_2\in\ANNs$
%	satisfy  for all $k\in\{1,2\}$, $x\in\R^{d}$ that $2\le\hiddenDimId\le 2d$,
%	$l_{2,L_2-1}\le l_{1,L_1-1}+\hiddenDimId$,
%	 $\dims(\mathbb{I}) = (d,\hiddenDimId,d)$, $(\functionANN(\mathbb{I}))(x)=x$, 
%%	 $\functionANN(\Phi_k)\in C(\R^d,\R^d)$, 
%$\inDimANN(\Phi_k)=\outDimANN(\Phi_k)=d$,
%	 and
%	$\dims(\Phi_k)=(l_{k,0},l_{k,1},\dots, l_{k,L_k})$
%	%	\begin{equation}\label{CorCompositionSum:HypothesisDims}
%	%	\dims(\Phi_k)=(l_0^k,l_1^k,\dots, l_{L_k}^k),\qandq \dims(\mathbb{I}) = (d,\hiddenDimId,d)
%	%	\end{equation}
%	(cf.\ Definition~\ref{Def:ANN} and Definition~\ref{Definition:ANNrealization}).
%	Then there exists $\Psi\in \ANNs$ which satisfies for all $x\in\R^d$ that
%$\functionANN(\Psi)\in C(\R^d,\R^d)$ and 
%	\begin{equation}
% [\functionANN(\Psi)](x)=[\functionANN(\Phi_2)](x)+[(\functionANN(\Phi_1))\circ (\functionANN(\Phi_2))](x),
%\end{equation}
%which satisfies that
%	\begin{equation}
%	 \dims(\Psi)=(l_{2,0},l_{2,1},\dots, l_{2,L_2-1},l_{1,1}+\hiddenDimId,l_{1,2}+\hiddenDimId,\dots,l_{1,L_1-1}+\hiddenDimId,l_{1,L_1} ),
%	\end{equation}
%and which satisfies that
%	\begin{equation}
%\begin{split}
%\paramANN(\Psi)
%&\le   
%\paramANN(\Phi_2)+[\paramANN(\mathbb{I})+\paramANN(\Phi_1)]^2.
%\end{split}
%\end{equation}
%\end{cor}

\begin{prop}\label{Cor:CompositionSum}
	Let $a\in C(\R,\R)$, $L_1\in \N\cap [2,\infty)$, $L_2\in\N$, $\mathbb{I}, \Phi_1,\Phi_2\in\ANNs$, $d,\hiddenDimId, l_{1,0},l_{1,1},\dots,\allowbreak l_{1,L_1},l_{2,0},l_{2,1},\allowbreak\dots,\allowbreak l_{2,L_2}\in\N$	
	satisfy  for all $k\in\{1,2\}$, $x\in\R^{d}$ that $2\le\hiddenDimId\le 2d$,
	$l_{2,L_2-1}\le l_{1,L_1-1}+\hiddenDimId$,
	$\dims(\mathbb{I}) = (d,\hiddenDimId,d)$, $(\functionANN(\mathbb{I}))(x)=x$, 
	%	 $\functionANN(\Phi_k)\in C(\R^d,\R^d)$, 
	$\inDimANN(\Phi_k)=\outDimANN(\Phi_k)=d$,
	and
	$\dims(\Phi_k)=(l_{k,0},l_{k,1},\dots, l_{k,L_k})$
	%	\begin{equation}\label{CorCompositionSum:HypothesisDims}
	%	\dims(\Phi_k)=(l_0^k,l_1^k,\dots, l_{L_k}^k),\qandq \dims(\mathbb{I}) = (d,\hiddenDimId,d)
	%	\end{equation}
	(cf.\ Definition~\ref{Def:ANN} and Definition~\ref{Definition:ANNrealization}).
	Then there exists $\Psi\in \ANNs$ such that
	\begin{enumerate}[(i)]
				\item \label{CorCompositionSum:Continuity}
				it holds that
				$\functionANN(\Psi)\in C(\R^d,\R^d)$,
		\item \label{CorCompositionSum:Realization}
		it holds for all $x\in\R^d$  that
		\begin{equation}
		(\functionANN(\Psi))(x)=(\functionANN(\Phi_2))(x)+\big((\functionANN(\Phi_1))\circ (\functionANN(\Phi_2))\big)(x),
		\end{equation}
		\item \label{CorCompositionSum:Dims}
		it holds that
%		\begin{equation}
%		\dims(\Psi)=
%		\begin{cases}
%		(l_{2,0},l_{2,1},\dots, l_{2,L_2-1},l_{1,1}+\hiddenDimId,l_{1,2}+\hiddenDimId,\dots,l_{1,L_1-1}+\hiddenDimId,l_{1,L_1} ) &: L_1=1\\
%		\dims(\Phi_2) &: L_1>1
%		\end{cases},
%		\end{equation}
			 	\begin{equation}
			 	\dims(\Psi)=(l_{2,0},l_{2,1},\dots, l_{2,L_2-1},l_{1,1}+\hiddenDimId,l_{1,2}+\hiddenDimId,\dots,l_{1,L_1-1}+\hiddenDimId,l_{1,L_1} ),
			 	\end{equation} 
			 	and
		\item \label{CorCompositionSum:Params}
		it holds that 
		\begin{equation}
		\begin{split}
					\paramANN(\Psi)
					&\le  
					\paramANN(\Phi_2)+\paramANN(\Phi_1)\big[\tfrac{1}{4}\paramANN(\Phi_1)+\paramANN(\mathbb{I})-1\big] 
					\\&\le \paramANN(\Phi_2)+\big[\tfrac{1}{2}\paramANN(\mathbb{I})+\paramANN(\Phi_1)\big]^{\!2}.
		\end{split}
		\end{equation}
%		$\paramANN(\Psi)
%		\le   
%		\paramANN(\Phi_2)+\big[\tfrac{1}{2}\paramANN(\mathbb{I})+\paramANN(\Phi_1)\big]^{\!2}$.
%			\begin{equation}
%			\begin{split}
%			\paramANN(\Psi)
%			&\le   
%			\paramANN(\Phi_2)+[\paramANN(\mathbb{I})+\paramANN(\Phi_1)]^2.
%			\end{split}
%			\end{equation}
	\end{enumerate}
\end{prop}

\begin{proof}[Proof of Proposition~\ref{Cor:CompositionSum}]
	Throughout this proof 
let $\Psi\in \ANNs$ satisfy that 
		\begin{enumerate}[(I)]
			\item \label{CorCompositionSum:ContinuityLoad}  it holds that
			$\functionANN(\Psi)\in C(\R^d,\R^d)$,
			\item \label{CorCompositionSum:RealizationLoad} it holds  for all $x\in\R^d$  that
			\begin{equation} (\functionANN(\Psi))(x)=(\functionANN(\Phi_2))(x)+\big((\functionANN(\Phi_1))\circ (\functionANN(\Phi_2))\big)(x),
			\end{equation}
			\item \label{CorCompositionSum:DimsLoad}  it holds that
			\begin{equation}
			\dims(\Psi)=(l_{2,0},l_{2,1},\dots, l_{2,L_2-1},l_{1,1}+\hiddenDimId,l_{1,2}+\hiddenDimId,\dots,l_{1,L_1-1}+\hiddenDimId,l_{1,L_1} ),
			\end{equation} 
			and
			\item \label{CorCompositionSum:ParamsLoad} it holds that
%			\begin{multline}
%			\paramANN(\Psi)
%			\le   
%			\paramANN(\Phi_2)+
%			\paramANN(\Phi_1)+\paramANN(\mathbb{I})+(L_1-2)\hiddenDimId(\hiddenDimId+1)
%			+l_{1,1} d 
%			\\+\hiddenDimId \bigg[\smallsum\limits_{m=2}^{L_1} l_{1,m} \bigg]
%			+ \hiddenDimId \bigg[\smallsum\limits_{m=0}^{L_1-2} l_{1,m} \bigg]
%			+d \,l_{1,L_1-1}
%			+ (l_{1,1}+\hiddenDimId) ( l_{2,L_2-1}+1)
%			\end{multline}
			\begin{equation}
			\begin{split}
						\paramANN(\Psi)
						&=
						\paramANN(\Phi_1)+\paramANN(\Phi_2)
						+ (\hiddenDimId-d)(l_{2,L_2-1}+1)+l_{1,1} (l_{2,L_2-1}-d)
						\\&\quad+(L_1-2)\hiddenDimId(\hiddenDimId+1)
						+\hiddenDimId\bigg[\smallsum\limits_{m=2}^{L_1} l_{1,m} \bigg]
						+\hiddenDimId\bigg[\smallsum\limits_{m=1}^{L_1-2} l_{1,m} \bigg]
			\end{split}
			\end{equation}
		\end{enumerate}
		(cf.\ Lemma~\ref{Lemma:CompositionSum}). 
	Note that the fact that $l_{1,0}=\inDimANN(\Phi_1)=d=\outDimANN(\Phi_1)=l_{1,L_1}$ implies that
	\begin{equation}\label{CorCompositionSum:EstimateOfFirstSum}
	\begin{split}
	\hiddenDimId \bigg[\smallsum\limits_{m=2}^{L_1} l_{1,m} \bigg]
	&\le \tfrac{1}{2}	\hiddenDimId \bigg[\smallsum\limits_{m=2}^{L_1} l_{1,m} (l_{1,m-1}+1) \bigg]
		= \tfrac{1}{2}\hiddenDimId \big[\paramANN(\Phi_1)-l_{1,1}(d+1)\big]
	\end{split}
	\end{equation}
	and
	\begin{equation}\label{CorCompositionSum:EstimateOfSecondSum}
	\begin{split}
	\hiddenDimId \bigg[\smallsum\limits_{m=1}^{L_1-2} l_{1,m} \bigg]
	&\le \tfrac{1}{2}	\hiddenDimId \bigg[\smallsum\limits_{m=1}^{L_1-2} l_{1,m} (l_{1,m-1}+1) \bigg]
	\\&=\tfrac{1}{2}\hiddenDimId\big[\paramANN(\Phi_1)-d(l_{1,L_1-1}+1)-l_{1,L_1-1}(l_{1,L_1-2}+1)\big].
	\end{split}
	\end{equation}
	Combining this with \eqref{CorCompositionSum:ParamsLoad} and the hypothesis that 
	$l_{2,L_2-1}\le l_{1,L_1-1}+\hiddenDimId$
	ensures that
\begin{equation}\label{CorCompositionSum:FirstParameterEstimateFirst}
\begin{split}
\paramANN(\Psi)&\le    
					[1+\hiddenDimId]\,	\paramANN(\Phi_1)+\paramANN(\Phi_2)
						+ (\hiddenDimId-d)(l_{2,L_2-1}+1)+l_{1,1} (l_{2,L_2-1}-d)
						\\&\quad+(L_1-2)\hiddenDimId(\hiddenDimId+1)
					-\tfrac{1}{2}\hiddenDimId\, l_{1,1}(d+1)
					\\&\quad-\tfrac{1}{2}\hiddenDimId\,d(l_{1,L_1-1}+1)-\tfrac{1}{2}\hiddenDimId\,l_{1,L_1-1}(l_{1,L_1-2}+1)
					\\&\le    
					[1+\hiddenDimId]\,	\paramANN(\Phi_1)+\paramANN(\Phi_2)
					+ [\max\{\hiddenDimId-d,0\}](l_{1,L_1-1}+\hiddenDimId+1)
					\\&\quad+l_{1,1} \big[l_{1,L_1-1}+\hiddenDimId-d-\tfrac{1}{2}\hiddenDimId\, (d+1)\big]
					+(L_1-2)\hiddenDimId(\hiddenDimId+1)
					\\&\quad-\tfrac{1}{2}\hiddenDimId\,d(l_{1,L_1-1}+1)-\tfrac{1}{2}\hiddenDimId\,l_{1,L_1-1}(l_{1,L_1-2}+1).
%					\\&\le    
%					[1+\hiddenDimId]\,	\paramANN(\Phi_1)+\paramANN(\Phi_2)
%					+ (\hiddenDimId-d)l_{1,L_1-1}
%					+l_{1,1} l_{1,L_1-1}
%					\\&\quad+(L_1-1)\hiddenDimId(\hiddenDimId+1)
%					-\tfrac{3}{2}\hiddenDimId\,l_{1,L_1-1}
%					\\&\le    
%					[1+\hiddenDimId]\,	\paramANN(\Phi_1)+\paramANN(\Phi_2)
%					+l_{1,1} l_{1,L_1-1}+(L_1-1)\hiddenDimId(\hiddenDimId+1).
\end{split}
\end{equation}	
Moreover, observe that the hypothesis that $2\le \hiddenDimId\le 2d$ shows that
\begin{equation}
	\begin{split}
	&l_{1,1} \big[\hiddenDimId-d-\tfrac{1}{2}\hiddenDimId\, (d+1)\big]-\tfrac{1}{2}\hiddenDimId\,d(l_{1,L_1-1}+1)-\tfrac{1}{2}\hiddenDimId\,l_{1,L_1-1}(l_{1,L_1-2}+1)
	\\&\le l_{1,1} \big[2d-d- (d+1)\big]-\tfrac{1}{2}\hiddenDimId \,l_{1,L_1-1}-\hiddenDimId\,l_{1,L_1-1}\le -\tfrac{3}{2}\hiddenDimId \,l_{1,L_1-1}.
	\end{split}
\end{equation}
This and \eqref{CorCompositionSum:FirstParameterEstimateFirst} prove that 
\begin{equation}\label{CorCompositionSum:FirstParameterEstimate}
\begin{split}
\paramANN(\Psi)&\le    
[1+\hiddenDimId]\,	\paramANN(\Phi_1)+\paramANN(\Phi_2)
+ [\max\{\hiddenDimId-d,0\}]\,l_{1,L_1-1}
\\&\quad+ [\max\{\hiddenDimId-d,0\}]\,(\hiddenDimId+1)
+l_{1,1} l_{1,L_1-1}-\tfrac{3}{2}\hiddenDimId\, l_{1,L_1-1}
+(L_1-2)\hiddenDimId(\hiddenDimId+1)
\\&\le    
[1+\hiddenDimId]\,	\paramANN(\Phi_1)+\paramANN(\Phi_2)
+ [\max\{\hiddenDimId-d,0\}]\,l_{1,L_1-1}
\\&\quad+\hiddenDimId (\hiddenDimId+1)
+l_{1,1} l_{1,L_1-1}
+(L_1-2)\hiddenDimId(\hiddenDimId+1)-\tfrac{3}{2}\hiddenDimId\, l_{1,L_1-1}
%\\&\le    
%[1+\hiddenDimId]\,	\paramANN(\Phi_1)+\paramANN(\Phi_2)
%+ [\max\{\hiddenDimId-d,0\}]\,(l_{1,L_1-1}+\hiddenDimId+1)
%\\&\quad+l_{1,1} \big[l_{1,L_1-1}+\hiddenDimId-d-\tfrac{1}{2}\hiddenDimId\, (d+1)\big]
%+(L_1-2)\hiddenDimId(\hiddenDimId+1)
%\\&\quad-\tfrac{1}{2}\hiddenDimId\,d(l_{1,L_1-1}+1)-\tfrac{1}{2}\hiddenDimId\,l_{1,L_1-1}(l_{1,L_1-2}+1)
\\&\le    
[1+\hiddenDimId]\,	\paramANN(\Phi_1)+\paramANN(\Phi_2)
+ [\max\{\hiddenDimId-d,0\}]\,l_{1,L_1-1}
+l_{1,1} l_{1,L_1-1}
\\&\quad+(L_1-1)\hiddenDimId(\hiddenDimId+1)
-\tfrac{3}{2}\hiddenDimId\,l_{1,L_1-1}
\\&\le    
[1+\hiddenDimId]\,	\paramANN(\Phi_1)+\paramANN(\Phi_2)
+l_{1,1} l_{1,L_1-1}+(L_1-1)\hiddenDimId(\hiddenDimId+1).
\end{split}
\end{equation}	
Moreover, observe that
	\begin{equation}
	\begin{split}
		L_1-1&\le \bigg[\smallsum\limits_{m=1}^{L_1} l_{1,m}\bigg]-1
		\le \tfrac{1}{2} \bigg[\smallsum\limits_{m=1}^{L_1} l_{1,m}(l_{1,m-1}+1)\bigg]-1
		\\&\le\tfrac{1}{2} \paramANN(\Phi_1)-1
		\le\tfrac{1}{2} \paramANN(\Phi_1).
	\end{split}
	\end{equation}
Combining this and \eqref{CorCompositionSum:FirstParameterEstimate} with the fact that $\forall\, k\in \N\cap [1,L_1]\colon l_{1,k}\le \tfrac{1}{2} l_{1,k} (l_{1,k-1}+1)\le \tfrac{1}{2}\paramANN(\Phi_1)$ demonstrates that
	\begin{equation}\label{CorCompositionSum:SecondParameterEstimate}
	\begin{split}
	\paramANN(\Psi)
	&\le [1+\hiddenDimId]\paramANN(\Phi_1)+\paramANN(\Phi_2)+l_{1,1}l_{1,L_1-1}+\tfrac{1}{2}\paramANN(\Phi_1) \hiddenDimId(\hiddenDimId+1)
	\\&=\paramANN(\Phi_2)+\big[1+\hiddenDimId+\tfrac{1}{2}\hiddenDimId(\hiddenDimId+1)\big] \paramANN(\Phi_1)+l_{1,1}l_{1,L_1-1}
	\\&\le   
	\paramANN(\Phi_2)+[1+\hiddenDimId+\tfrac{1}{2}\hiddenDimId(\hiddenDimId+1)]\,\paramANN(\Phi_1)
	+\tfrac{1}{4}[\paramANN(\Phi_1)]^2.
	\end{split}
	\end{equation}
	Furthermore, note that  the hypothesis that $2\le \hiddenDimId\le 2d$ and the hypothesis that $\dims(\mathbb{I}) = (d,\hiddenDimId,d)$ prove that 
%	\begin{equation}
%		\paramANN(\mathbb{I})-\tfrac{1}{2}\hiddenDimId^2=\hiddenDimId (d+1)+d(\hiddenDimId+1)-\tfrac{1}{2}\hiddenDimId^2=2 d \hiddenDimId+d+\hiddenDimId-\tfrac{1}{2}\hiddenDimId^2
%		\ge \tfrac{1}{2}\hiddenDimId^2+\tfrac{3}{2}\hiddenDimId=\hiddenDimId+\tfrac{1}{2}\hiddenDimId(\hiddenDimId+1).
%	\end{equation}
	\begin{equation}
	\begin{split}
		\hiddenDimId+\tfrac{1}{2}\hiddenDimId(\hiddenDimId+1)
	&=\hiddenDimId^2+\hiddenDimId+\tfrac{1}{2}\hiddenDimId-\tfrac{1}{2}\hiddenDimId^2
	\le 2 d \hiddenDimId+\hiddenDimId+d-\tfrac{1}{2}\hiddenDimId^2
	\\&=\hiddenDimId (d+1)+d(\hiddenDimId+1)-\tfrac{1}{2}\hiddenDimId^2
	=\paramANN(\mathbb{I})-\tfrac{1}{2}\hiddenDimId^2
	\le \paramANN(\mathbb{I})-2.
	\end{split}
\end{equation}
	Combining this and \eqref{CorCompositionSum:SecondParameterEstimate} implies that 
	\begin{equation}
	\begin{split}
	\paramANN(\Psi)
	&\le   
	\paramANN(\Phi_2)+[1+\paramANN(\mathbb{I})-2]\,\paramANN(\Phi_1)
	+\tfrac{1}{4}[\paramANN(\Phi_1)]^2
	\\&=\paramANN(\Phi_2)+\big[\tfrac{1}{4}\paramANN(\Phi_1)+\paramANN(\mathbb{I})-1\big] \paramANN(\Phi_1)
	\\&\le   
	\paramANN(\Phi_2)+\paramANN(\mathbb{I})\,\paramANN(\Phi_1)+\tfrac{1}{4}[\paramANN(\mathbb{I})]^2
	+[\paramANN(\Phi_1)]^2
	\\&=\paramANN(\Phi_2)+\big[\tfrac{1}{2}\paramANN(\mathbb{I})+\paramANN(\Phi_1)\big]^{\!2}.
	\end{split}
	\end{equation}
%		Observe that \eqref{CorCompositionSum:FunctionFromLemma} and \eqref{CorCompositionSum:DimsFromLemma} establish items~\eqref{CompositionSum:Realization} and \eqref{CompositionSum:Dims}.
	This, \eqref{CorCompositionSum:ContinuityLoad}, \eqref{CorCompositionSum:RealizationLoad}, and \eqref{CorCompositionSum:DimsLoad} establish items~\eqref{CorCompositionSum:Continuity}--\eqref{CorCompositionSum:Params}.
	The proof of
 Proposition~\ref{Cor:CompositionSum} is thus completed.
\end{proof}

%% file: ANNrepMultNestedEuler.tex
\begin{cor}\label{Cor:CompositionSumInduction}
	Let $a\in C(\R,\R)$, $d,\hiddenDimId,\mathfrak{L}\in\N$, $\ell_0,\ell_1,\dots, \ell_\mathfrak{L}\in\N$, $(L_n)_{n\in\N_0}\subseteq \N\cap [2,\infty)$,
%	$(l_{n,k})_{n\in\N_0, k\in\{1,2,\dots, L_n\}}\subseteq\N$,
	$\mathbb{I}, \psi\in\ANNs$, $(\phi_n)_{n\in\N_0}\subseteq \ANNs$,
	let $l_{n,k}\in\N$, $k\in\{0,1,\dots, L_n\}$, $n\in\N_0$,  assume
	 for all  $n\in\N_0$, $x\in\R^{d}$  that $2\le\hiddenDimId\le 2d$,
	$\ell_{\mathfrak{L}-1}\le l_{0,L_0-1}+\hiddenDimId$, $l_{n,L_n-1}\le l_{n+1,L_{n+1}-1}$,
	$\dims(\mathbb{I}) = (d,\hiddenDimId,d)$, $(\functionANN(\mathbb{I}))(x)=x$, 
	%	 $\functionANN(\phi_k)\in C(\R^d,\R^d)$, 
	$\inDimANN(\phi_n)=\outDimANN(\phi_n)=\inDimANN(\psi)=\outDimANN(\psi)=d$,
	$\dims(\phi_n)=(l_{n,0},l_{n,1},\dots, l_{n,L_n})$, and
	$\dims(\psi)=(\ell_{0},\ell_{1},\dots, \ell_{\mathfrak{L}})$, and
	%	\begin{equation}\label{CorCompositionSum:HypothesisDims}
	%	\dims(\phi_n_k)=(l_0^k,l_1^k,\dots, l_{L_k}^k),\qandq \dims(\mathbb{I}) = (d,\hiddenDimId,d)
	%	\end{equation}
	let $f_n\colon \R^d\to \R^d$, $n\in\N_0$, be the functions which satisfy for all  $n\in\N_0$, $x\in\R^d$ that 
	\begin{equation}\label{CompositionSumInduction:Recursion}
	f_0(x)=(\functionANN(\psi))(x)
	\qandq f_{n+1}(x)=f_n(x)+\big([\functionANN(\phi_n)]\circ f_n\big)(x)
	\end{equation}
	(cf.\ Definition~\ref{Def:ANN} and Definition~\ref{Definition:ANNrealization}).
	Then for every $n\in\N$ there exists $\Psi\in \ANNs$
	such that
	\begin{enumerate}[(i)]
%		\item \label{CompositionSumInduction:InitialValue}
%		that $\Psi_0=\Psi$,
	\item \label{CompositionSumInduction:Continuity}
		it holds that 
		$\functionANN(\Psi)\in C(\R^d,\R^d)$,
		\item \label{CompositionSumInduction:Realization}
		it holds for all $x\in\R^d$ that 
		$(\functionANN(\Psi))(x)=f_n(x)$,
		\item \label{CompositionSumInduction:Dims}
		it holds that $\hiddenLength(\Psi)={\hiddenLength(\psi)+\sum_{k=0}^{n-1} \hiddenLength(\phi_k)}$, 
		\item it holds that 
				\begin{multline}\label{CompositionSumInduction:ClaimTwo}
		\dims(\Psi)=
		\big(\ell_{0},\ell_{1},\dots, \ell_{\mathfrak{L}-1},
		l_{0,1}+\hiddenDimId, l_{0,2}+\hiddenDimId,\dots, l_{0,L_0-1}+\hiddenDimId,
		\\l_{1,1}+\hiddenDimId, l_{1,2}+\hiddenDimId,\dots, l_{1,L_1-1}+\hiddenDimId,
		\dots, l_{n-1,1}+\hiddenDimId, l_{n-1,2}+\hiddenDimId,\dots, l_{n-1,L_{n-1}-1}+\hiddenDimId,d
		\big),
%		\in \N^{\mathfrak{L}+1+\sum_{k=0}^{n-1} (L_{k}-1)}
		\end{multline}
		and
		\item \label{CompositionSumInduction:Params} 
		it holds that $\paramANN(\Psi)\le \paramANN(\psi)+ \smallsum_{k=0}^{n-1} \big[\tfrac{1}{2}\paramANN(\mathbb{I})+\paramANN(\phi_k)\big]^2$. 	
%		\begin{equation}\label{CompositionSumInduction:ClaimThree}
%		\paramANN(\Psi)\le \paramANN(\Psi)+ \smallsum_{k=0}^{n-1} [\paramANN(\mathbb{I})+\paramANN(\phi_k)]^2.
%		\end{equation}
	\end{enumerate}
\end{cor}

\begin{proof}[Proof of Corollary~\ref{Cor:CompositionSumInduction}]
	%	Throughout this proof let $\Psi_0\in\ANNs$ satisfy that $\Psi_0=\Psi$.
	We prove items~\eqref{CompositionSumInduction:Continuity}--\eqref{CompositionSumInduction:Params} by induction on $n\in\N$.
%	Note that \eqref{CompositionSumInduction:Recursion}, the fact that $\functionANN(\psi)\in C(\R^d,\R^d)$, the fact that $\lengthANN(\psi)=\mathfrak{L}={\mathfrak{L}+\allowbreak\sum_{k=0}^{0-1} (L_{k}-1)}$,
%	the fact that $\dims(\psi)=(\ell_{0},\ell_{1},\ell_{2},\dots, \ell_{\mathfrak{L}-1},d)$,
%	and the fact that $\paramANN(\psi)= \paramANN(\psi)+ \smallsum_{k=0}^{0-1} \big[\tfrac{1}{2}\paramANN(\mathbb{I})+\paramANN(\phi_k)\big]^2$ prove items~\eqref{CompositionSumInduction:Continuity}--\eqref{CompositionSumInduction:Params} in the base case $n=0$.
		Note that 
			the  hypothesis that $\dims(\psi)=(\ell_{0},\ell_{1},\dots, \ell_{\mathfrak{L}})$, the fact that $\ell_{0}=\inDimANN(\psi)=\ell_{\mathfrak{L}}=\outDimANN(\psi)=d$, the hypothesis that $\dims(\phi_0)=(l_{0,0},l_{0,1},\dots, l_{0,L_0})$,
			the hypothesis that 	$\ell_{\mathfrak{L}-1}\le l_{0,L_0-1}+\hiddenDimId$, the hypothesis that $L_0\in\N\cap [2,\infty)$, Proposition~\ref{Cor:CompositionSum} (with $a=a$, $L_1=L_0$, $L_2=\mathfrak{L}$, $\mathbb{I}=\mathbb{I}$,
			$\Phi_1=\phi_0$, $\Phi_2=\psi$, $d=d$, $\hiddenDimId=\hiddenDimId$, $l_{1,v}=l_{0,v}$, $l_{2,w}=\ell_w$ for $v\in\{0,1,\dots,L_0\}$, $w\in \{0,1,\dots,\mathfrak{L}\}$ in the notation of Proposition~\ref{Cor:CompositionSum}), and \eqref{CompositionSumInduction:Recursion} imply that there exists 
			$\Upsilon\in\ANNs$	which satisfies that
			\begin{enumerate}[(I)]
				\item \label{CompositionSumInduction:ContinuityBase}
				it holds that 
				$\functionANN(\Upsilon)\in C(\R^d,\R^d)$,
				\item \label{CompositionSumInduction:RealizationBase}
				it holds for all $x\in\R^d$  that
			\begin{equation}
			\begin{split}
			(\functionANN(\Upsilon))(x)&=(\functionANN(\psi))(x)+\big([\functionANN(\phi_0)]\circ [\functionANN(\psi)]\big)(x)
			\\&=f_0(x)+([\functionANN(\phi_0)]\circ f_0)(x)
			=f_{1}(x),
			\end{split}
			\end{equation}
				\item \label{CompositionSumInduction:DimsBase}
			it holds	that 
			\begin{equation}\label{CompositionSumInduction:BaseCase}
			\dims(\Upsilon)=
			\big(\ell_{0},\ell_{1},\dots, \ell_{\mathfrak{L}-1},
			l_{0,1}+\hiddenDimId, l_{0,2}+\hiddenDimId,\dots, l_{0,L_0-1}+\hiddenDimId,
			l_{0,L_0}\big),
			%	\in \N^{\mathfrak{L}+1+ (L_{0}-1)},
			\end{equation}
			and
				\item \label{CompositionSumInduction:ParamsBase}	
				it holds that  $\paramANN(\Upsilon)\le \paramANN(\psi)+   \big[\tfrac{1}{2}\paramANN(\mathbb{I})+\paramANN(\phi_0)\big]^2$.
				%			\begin{equation}
				%			\begin{split}
				%			\paramANN(\Pi)
				%			&\le   
				%			\paramANN(\phi_2)+[\paramANN(\mathbb{I})+\paramANN(\phi_1)]^2.
				%			\end{split}
				%			\end{equation}
			\end{enumerate}
		Observe that \eqref{CompositionSumInduction:DimsBase} shows that $\lengthANN(\Upsilon)=\mathfrak{L}+ L_{0}-1$. Hence, we obtain that 
		\begin{equation}
			\hiddenLength(\Upsilon)=\lengthANN(\Upsilon)-1=(\mathfrak{L}-1)+ (L_{0}-1)=\hiddenLength(\psi)+\hiddenLength(\phi_0).
		\end{equation}
		Combining this with \eqref{CompositionSumInduction:ContinuityBase}--\eqref{CompositionSumInduction:ParamsBase} establishes items~\eqref{CompositionSumInduction:Continuity}--\eqref{CompositionSumInduction:Params} in the base case $n=1$.
	%	
	%	Note that \eqref{CompositionSumInduction:Recursion} and the hypothesis that 
	%	$\dims(\Psi)=(\ell_{0},\ell_{1},\dots, \ell_{\mathfrak{L}})$ prove 
	%	\eqref{CompositionSumInduction:ClaimTwo} and \eqref{CompositionSumInduction:ClaimThree}
	%	in the base case $n=0$.
	For the induction step $\N\ni n\to n+1 \in \N\cap [2,\infty)$ let $n\in \N$, $\Psi\in \ANNs$, $\mathfrak{l}_0,\mathfrak{l}_1,\dots, \mathfrak{l}_{\mathfrak{L}+\sum_{k=0}^{n-1} (L_{k}-1)}\in\N$  satisfy that
	\begin{enumerate}[(a)] 
		\item \label{CompositionSumInduction:ContinuityInduction}
		it holds that
		$\functionANN(\Psi)\in C(\R^d,\R^d)$,
		\item \label{CompositionSumInduction:RealizationInduction}
		it holds for all $x\in\R^d$ that 
		$(\functionANN(\Psi))(x)=f_n(x)$,
		\item \label{CompositionSumInduction:HiddenInduction}
		it holds that $\hiddenLength(\Psi)={\hiddenLength(\psi)+\sum_{k=0}^{n-1} \hiddenLength(\phi_k)}$, 
		\item \label{CompositionSumInduction:DimsInduction} it holds that 
		\begin{equation}\label{CompositionSumInduction:ClaimTwoInduction}
		\begin{split}
		\dims(\Psi)&=
		\big(\ell_{0},\ell_{1},\dots, \ell_{\mathfrak{L}-1},
		l_{0,1}+\hiddenDimId, l_{0,2}+\hiddenDimId,\dots, l_{0,L_0-1}+\hiddenDimId,
		l_{1,1}+\hiddenDimId, l_{1,2}+\hiddenDimId,
		\\& \qquad\dots, l_{1,L_1-1}+\hiddenDimId, \dots,
		l_{n-1,1}+\hiddenDimId, l_{n-1,2}+\hiddenDimId,\dots, l_{n-1,L_{n-1}-1}+\hiddenDimId,d
		\big)
		\\&=\big(\mathfrak{l}_0,\mathfrak{l}_1,\dots, \mathfrak{l}_{\mathfrak{L}+\sum_{k=0}^{n-1} (L_{k}-1)}\big),
		%		\in \N^{\mathfrak{L}+1+\sum_{k=0}^{n-1} (L_{k}-1)},
		\end{split}
		\end{equation}
		and
		\item \label{CompositionSumInduction:ParamsInduction} 
		it holds that	$\paramANN(\Psi)\le \paramANN(\psi)+ \smallsum_{k=0}^{n-1}  \big[\tfrac{1}{2}\paramANN(\mathbb{I})+\paramANN(\phi_k)\big]^2$.
		%			\begin{equation}\label{CompositionSumInduction:ClaimThree}
		%			\paramANN(\Psi)\le \paramANN(\Psi)+ \smallsum_{k=0}^{n-1} [\paramANN(\mathbb{I})+\paramANN(\phi_k)]^2.
		%			\end{equation}
	\end{enumerate}		
	Observe that \eqref{CompositionSumInduction:DimsInduction} and the hypothesis that $\forallDist k\in\N_0\colon l_{k,L_{k}-1}\le l_{k+1,L_{k+1}-1}$ demonstrate that
	\begin{equation}
		\mathfrak{l}_{\lengthANN(\Psi)-1}=\mathfrak{l}_{\mathfrak{L}-1+\sum_{k=0}^{n-1} (L_{k}-1)}=l_{n-1,L_{n-1}-1}+\hiddenDimId\le l_{n,L_n-1}+\hiddenDimId.
	\end{equation}
%		 $\mathfrak{l}_{\lengthANN(\Psi)-1}=\mathfrak{l}_{\mathfrak{L}-1+\sum_{k=0}^{n-1} (L_{k}-1)}\in \{\ell_{\mathfrak{L}-1},l_{L-1}+\hiddenDimId\}$. 
	%	The hypothesis that for all $m\in\N_0$ it holds that $\ell_{\mathfrak{L}-1}\le l_{0,L-1}+\hiddenDimId$ and 
	%	$l_{m,L_m-1}\le l_{m+1,L_{m+1}-1}$
	%	therefore ensures that $\mathfrak{l}_{\lengthANN(\Psi)-1}\le l_{n,L_n-1}+\hiddenDimId$.
	The hypothesis that $\dims(\phi_n)=(l_{n,0},l_{n,1},\dots, l_{n,L_n})$, \eqref{CompositionSumInduction:DimsInduction}, the hypothesis that $L_n\in\N\cap [2,\infty)$, and Proposition~\ref{Cor:CompositionSum} 
	(with $a=a$, $L_1=L_n$, $L_{2}=\mathfrak{L}+\sum_{k=0}^{n-1}(L_{k}-1)$, $\mathbb{I}=\mathbb{I}$, $\Phi_1=\phi_n$, $\Phi_2=\Psi$, $d=d$, $\mathfrak{i}=\mathfrak{i}$, $l_{1,v}=l_{n,v}$, $l_{2,w}=\ell_w$ for $v\in\{0,1,\dots,L_n\}$, 
	$w\in \{0,1,\dots,\mathfrak{L}+\sum_{k=0}^{n-1}(L_k-1)\}$ in the notation of Proposition~\ref{Cor:CompositionSum}) hence prove that there exists $\Phi\in\ANNs$ which satisfies that
	\begin{enumerate}[(A)] 
		\item \label{CompositionSumInduction:ContinuityStep}
		it holds that 
		$\functionANN(\Phi)\in C(\R^d,\R^d)$,
		\item \label{CompositionSumInduction:RealizationStep}
		it holds for all $x\in\R^d$  that 
		\begin{equation}
		(\functionANN(\Phi))(x)=(\functionANN(\Psi))(x)+\big([\functionANN(\phi_n)]\circ [\functionANN(\Psi)]\big)(x),
		\end{equation}
		\item \label{CompositionSumInduction:DimsStep}
		%	that  $\lengthANN(\Phi)={\mathfrak{L}+\sum_{k=0}^{n} (L_{k}-1)}$ and
		it holds that
		\begin{multline}\label{CompositionSumInduction:InductionStepDims}
		\dims(\Phi)=
		\big(\ell_0,\ell_{1},\dots, \ell_{\mathfrak{L}-1},
		l_{0,1}+\hiddenDimId, l_{0,2}+\hiddenDimId,\dots, l_{0,L_0-1}+\hiddenDimId,
		l_{1,1}+\hiddenDimId, l_{1,2}+\hiddenDimId,
		\\\dots, 
		l_{1,L_1-1}+\hiddenDimId, \dots,
		l_{n,1}+\hiddenDimId, l_{n,2}+\hiddenDimId,\dots, l_{n,L_{n}-1}+\hiddenDimId,l_{n,L_n}
		\big),
		%	\in \N^{\mathfrak{L}+1+\sum_{k=0}^{n} (L_{k}-1)},
		\end{multline}
		and
		\item \label{CompositionSumInduction:ParamsStep} 
		it holds that 	$\paramANN(\Phi)\le   
		\paramANN(\Psi)+ \big[\tfrac{1}{2}\paramANN(\mathbb{I})+\paramANN(\phi_n)\big]^2$.
		%			\begin{equation}\label{CompositionSumInduction:ClaimThree}
		%			\paramANN(\Psi)\le \paramANN(\Psi)+ \smallsum_{k=0}^{n-1} [\paramANN(\mathbb{I})+\paramANN(\phi_k)]^2.
		%			\end{equation}
	\end{enumerate}
	%	Proposition~\ref{Cor:CompositionSum} hence proves that there exists $\Psi_{n+1}\in\ANNs$ which satisfies that for all $x\in\R^d$ it holds that 
	%	\begin{multline}\label{CompositionSumInduction:InductionStepDims}
	%	\dims(\Psi_{n+1})=
	%	\big(d,\ell_{1},\ell_{2},\dots, \ell_{\mathfrak{L}-1},
	%	l_{0,1}+\hiddenDimId, l_{0,2}+\hiddenDimId,\dots, l_{0,L_0-1}+\hiddenDimId,
	%	l_{1,1}+\hiddenDimId, l_{1,2}+\hiddenDimId,\dots, 
	%	\\ l_{1,L_1-1}+\hiddenDimId, \dots,
	%	l_{n,1}+\hiddenDimId, l_{n,2}+\hiddenDimId,\dots, l_{n,L_{n}-1}+\hiddenDimId,d
	%	\big)\in \N^{\mathfrak{L}+1+\sum_{k=0}^{n} (L_{k}-1)},
	%	\end{multline}
	%	\begin{equation}
	%	\begin{split}
	%	\paramANN(\Psi_{n+1})
	%	&\le   
	%	\paramANN(\Psi_n)+[\paramANN(\mathbb{I})+\paramANN(\phi_n)]^2,
	%	\end{split}
	%	\end{equation}
	%	$\lengthANN(\Psi_{n+1})={\mathfrak{L}+\sum_{k=0}^{n} (L_{k}-1)}$,  $\functionANN(\Psi_{n+1})\in C(\R^d,\R^d)$, and
	%	\begin{equation}
	%	(\functionANN(\Psi_{n+1}))(x)=(\functionANN(\Psi_n))(x)+(\functionANN(\phi_n)\circ \functionANN(\Psi_n))(x).
	%	\end{equation}
	Next note that \eqref{CompositionSumInduction:DimsStep} implies that $\lengthANN(\Phi)={\mathfrak{L}+\sum_{k=0}^{n} (L_{k}-1)}$. Hence, we obtain that 
	\begin{equation}\label{CompositionSumInduction:hiddenStep}
	\hiddenLength(\Phi)=\lengthANN(\Phi)-1=(\mathfrak{L}-1)+\smallsum_{k=0}^{n} (L_{k}-1)=\hiddenLength(\psi)+\smallsum_{k=0}^{n} \hiddenLength(\phi_k).
	\end{equation}
	Moreover, observe that \eqref{CompositionSumInduction:RealizationStep}, \eqref{CompositionSumInduction:Recursion}, and \eqref{CompositionSumInduction:RealizationInduction} demonstrate that for all $x\in\R^d$ it holds that
	\begin{equation}\label{CompositionSumInduction:StepCalculation}
	\begin{split}
	(\functionANN(\Phi))(x)&=(\functionANN(\Psi))(x)+\big([\functionANN(\phi_n)]\circ [\functionANN(\Psi)]\big)(x)
	\\&=f_n(x)+([\functionANN(\phi_n)]\circ f_n)(x)
	=f_{n+1}(x).
	\end{split}
	\end{equation}
	In addition, note that \eqref{CompositionSumInduction:ParamsStep}  and \eqref{CompositionSumInduction:ParamsInduction} ensure that
	\begin{equation}
	\begin{split}
	\paramANN(\Phi)
	&\le   
	\paramANN(\psi)+\left[\smallsum\limits_{k=0}^{n-1}  \big[\tfrac{1}{2}\paramANN(\mathbb{I})+\paramANN(\phi_k)\big]^2\right]+ \big[\tfrac{1}{2}\paramANN(\mathbb{I})+\paramANN(\phi_n)\big]^2
	\\&=\paramANN(\psi)+\smallsum\limits_{k=0}^{n}  \big[\tfrac{1}{2}\paramANN(\mathbb{I})+\paramANN(\phi_k)\big]^2.
	\end{split}
	\end{equation}
	This, \eqref{CompositionSumInduction:ContinuityStep}, \eqref{CompositionSumInduction:DimsStep}, \eqref{CompositionSumInduction:hiddenStep}, and \eqref{CompositionSumInduction:StepCalculation}
	prove items~\eqref{CompositionSumInduction:Continuity}--\eqref{CompositionSumInduction:Params} in the case $n+1$. Induction thus establishes items~\eqref{CompositionSumInduction:Continuity}--\eqref{CompositionSumInduction:Params}.
	The proof of
	Corollary~\ref{Cor:CompositionSumInduction} is thus completed.
\end{proof}

\begin{prop}\label{Cor:CompositionSumInductionAffine}
	Let $a\in C(\R,\R)$, $d,\mathfrak{L}\in\N$, $\ell_0,\ell_1,\dots, \ell_\mathfrak{L}\in\N$, 
	$ \psi\in\ANNs$, $(\phi_n)_{n\in\N_0}\subseteq \ANNs$
	satisfy  for all  $n\in\N_0$ that 
	$\inDimANN(\phi_n)=\outDimANN(\phi_n)=\inDimANN(\psi)=\outDimANN(\psi)=d$,
	$\lengthANN(\phi_n)=1$, and
	$\dims(\psi)=(\ell_{0},\ell_{1},\dots, \ell_{\mathfrak{L}})$ and
	%	\begin{equation}\label{CorCompositionSum:HypothesisDims}
	%	\dims(\phi_n_k)=(l_0^k,l_1^k,\dots, l_{L_k}^k),\qandq \dims(\mathbb{I}) = (d,\hiddenDimId,d)
	%	\end{equation}
	let $f_n\colon \R^d\to \R^d$, $n\in\N_0$, be the functions which satisfy for all  $n\in\N_0$, $x\in\R^d$ that 
	\begin{equation}\label{CompositionSumInductionAffine:Recursion}
	f_0(x)=(\functionANN(\psi))(x)
	\qandq f_{n+1}(x)=f_n(x)+\big([\functionANN(\phi_n)]\circ f_n\big)(x)
	\end{equation}
	(cf.\ Definition~\ref{Def:ANN} and Definition~\ref{Definition:ANNrealization}).
	Then for every $n\in\N_0$ there exists $\Psi\in\ANNs$
	such that
	\begin{enumerate}[(i)]
		%		\item \label{CompositionSumInductionAffine:InitialValue}
		%		that $\Psi_0=\psi$,
		\item \label{CompositionSumInductionAffine:Continuity}
		it holds that 
		$\functionANN(\Psi)\in C(\R^d,\R^d)$,
		\item \label{CompositionSumInductionAffine:Realization}
		it holds for all  $x\in\R^d$ that $(\functionANN(\Psi))(x)=f_n(x)$,
		and
		\item \label{CompositionSumInductionAffine:Dims}
		it holds  that $\dims(\Psi)=\dims(\psi)$.
	\end{enumerate}
\end{prop}

\begin{proof}[Proof of Proposition~\ref{Cor:CompositionSumInductionAffine}]
%	Throughout this proof let $\Psi_0\in\ANNs$ satisfy that $\Psi_0=\psi$.
	We prove items~\eqref{CompositionSumInductionAffine:Continuity}--\eqref{CompositionSumInductionAffine:Dims} by induction on $n\in\N_0$.
	Note that \eqref{CompositionSumInductionAffine:Recursion} and the fact that $\functionANN(\psi)\in C(\R^d,\R^d)$
	 establish items~\eqref{CompositionSumInductionAffine:Continuity}--\eqref{CompositionSumInductionAffine:Dims} in the base case $n=0$.
	For the induction step $\N_0\ni n\to n+1 \in \N$ let $n\in \N_0$, $\Psi\in \ANNs$ satisfy that
		\begin{enumerate}[(I)]
			%		\item \label{CompositionSumInductionAffine:InitialValue}
			%		that $\Psi_0=\psi$,
			\item \label{CompositionSumInductionAffineInduction:Continuity}
			it holds that 
			$\functionANN(\Psi)\in C(\R^d,\R^d)$,
			\item \label{CompositionSumInductionAffineInduction:Realization}
			it holds for all  $x\in\R^d$ that $(\functionANN(\Psi))(x)=f_n(x)$,
			and
			\item \label{CompositionSumInductionAffineInduction:Dims}
			it holds  that $\dims(\Psi)=\dims(\psi)$,
		\end{enumerate}
	 and let  $(A,b)\in (\R^{d\times d}\times \R^d)\subseteq \ANNs$, $\affineMap\in (\R^{d\times d}\times \R^d)\subseteq \ANNs$, $\Phi\in\ANNs$ satisfy that $\phi_n=(A,b)$, $\affineMap=(A+\idMatrix_d,b)$, and $\Phi=\compANN{\affineMap}{\Psi}$ (cf.\ Definition~\ref{Definition:ANNcomposition} and Definition~\ref{Definition:identityMatrix}).
	Observe that item~\eqref{PropertiesOfCompositions:Realization} in Proposition~\ref{Lemma:PropertiesOfCompositions} demonstrates that for all $x\in\R^d$ it holds that 
	$\functionANN(\Phi)\in C(\R^d,\R^d)$ and
	\begin{equation}
	\begin{split}
			(\functionANN(\Phi))(x)&=(\functionANN(\affineMap))\big((\functionANN(\Psi))(x)\big)
			\\&=(A+\idMatrix_d)\big((\functionANN(\Psi))(x)\big)+b
						\\&=A\big((\functionANN(\Psi))(x)\big)+b +(\functionANN(\Psi))(x)
			\\&=(\functionANN(\phi_n))\big((\functionANN(\Psi))(x)\big) +(\functionANN(\Psi))(x).
	\end{split}
	\end{equation}
	Combining this with \eqref{CompositionSumInductionAffine:Recursion} and \eqref{CompositionSumInductionAffineInduction:Realization} proves that 
	for all $x\in\R^d$ it holds that 
%	$\functionANN(\Phi)\in C(\R^d,\R^d)$ and
	\begin{equation}\label{CompositionSumInductionAffine:InductiveCalc}
			\begin{split}
			(\functionANN(\Phi))(x)&=(\functionANN(\phi_n))\big(f_n(x)\big) +f_n(x)=f_{n+1}(x).
			\end{split}
	\end{equation}
	In addition, note that \eqref{CompositionSumInductionAffineInduction:Dims}, the fact that $\Phi=\compANN{\affineMap}{\Psi}$, the fact that $\lengthANN(\affineMap)=1$, the fact that $\inDimANN(\affineMap)=\outDimANN(\affineMap)=\outDimANN(\Psi)=d$,  and item~\eqref{PropertiesOfCompositions:Dims} in Proposition~\ref{Lemma:PropertiesOfCompositions} imply that $\dims(\Phi)=\dims(\Psi)=\dims(\psi)$.
	Combining this and the fact that $\functionANN(\Phi)\in C(\R^d,\R^d)$ with \eqref{CompositionSumInductionAffine:InductiveCalc} 
	proves items~\eqref{CompositionSumInductionAffine:Continuity}--\eqref{CompositionSumInductionAffine:Dims} in the case $n+1$.
	Induction thus establishes items~\eqref{CompositionSumInductionAffine:Continuity}--\eqref{CompositionSumInductionAffine:Dims}.
		The proof of
		Proposition~\ref{Cor:CompositionSumInductionAffine} is thus completed.
%	
%	Observe that the fact that $\lengthANN(\phi_n)=1$ and Lemma~\ref{Lemma:CompositionSum} prove that there exists $\Psi_{n+1}\in\ANNs$ which satisfies 
%		that for all $x\in\R^d$ it holds that $\dims(\Psi_{n+1})=\dims(\Psi_n)$,
%		$\functionANN(\Psi_{n+1})\in C(\R^d,\R^d)$, and
%		\begin{equation}
%		(\functionANN(\Psi_{n+1}))(x)=(\functionANN(\Psi_n))(x)+(\functionANN(\phi_n)\circ \functionANN(\Psi_n))(x).
%		\end{equation}
%	Combining this, the fact that  for all  $x\in\R^d$ it holds that 
%	$[\functionANN(\Psi_n)](x)=f_n(x)$, \eqref{CompositionSumInductionAffine:Recursion}, and the fact that $\dims(\Psi_n)=\dims(\Psi)$ with induction
%	establishes items~\eqref{CompositionSumInductionAffine:Continuity}--\eqref{CompositionSumInductionAffine:Dims}.
%	The proof of
%	Proposition~\ref{Cor:CompositionSumInductionAffine} is thus completed.
\end{proof}

\begin{cor}\label{Cor:CompositionSumInductionTilde}
	Let $a\in C(\R,\R)$, $d,\hiddenDimId,L,\mathfrak{L}\in\N$, $ \ell_0, \ell_1,\dots, \ell_\mathfrak{L}\in\N$,
	%	$(l_{n,k})_{n\in\N_0, k\in\{1,2,\dots, L_n\}}\subseteq\N$,
	$\mathbb{I}, \psi\in\ANNs$, $(\phi_n)_{n\in\N_0}\subseteq \ANNs$,
		let $l_{n,k}\in\N$, $k\in\{0,1,\dots, L\}$, $n\in\N_0$,
	assume
	for all  $n\in\N_0$, $x\in\R^{d}$  that $2\le\hiddenDimId\le 2d$,
	$\ell_{\mathfrak{L}-1}\le l_{0,L-1}+\hiddenDimId$, $l_{n,L-1}\le l_{n+1,L-1}$, 
	$\dims(\mathbb{I}) = (d,\hiddenDimId,d)$, $(\functionANN(\mathbb{I}))(x)=x$, 
	%	 $\functionANN(\phi_k)\in C(\R^d,\R^d)$, 
	$\inDimANN(\phi_n)=\outDimANN(\phi_n)=\inDimANN(\psi)=\outDimANN(\psi)=d$,
	$\dims(\phi_n)=(l_{n,0},l_{n,1},\dots, l_{n,L})$, and
	$\dims(\psi)=(\ell_{0},\ell_{1},\dots, \ell_{\mathfrak{L}})$, and
	%	\begin{equation}\label{CorCompositionSum:HypothesisDims}
	%	\dims(\phi_n_k)=(l_0^k,l_1^k,\dots, l_{L_k}^k),\qandq \dims(\mathbb{I}) = (d,\hiddenDimId,d)
	%	\end{equation}
	let $f_n\colon \R^d\to \R^d$, $n\in\N_0$, be the functions which satisfy for all  $n\in\N_0$, $x\in\R^d$ that 
	\begin{equation}\label{CompositionSumInductionTilde:Recursion}
	f_0(x)=(\functionANN(\psi))(x)
	\qandq f_{n+1}(x)=f_n(x)+\big([\functionANN(\phi_n)]\circ f_n\big)(x)
	\end{equation}
	(cf.\ Definition~\ref{Def:ANN} and Definition~\ref{Definition:ANNrealization}).
	Then for every $n\in\N_0$ there exists $\Psi\in \ANNs$
	such that
	\begin{enumerate}[(i)]
		%		\item \label{CompositionSumInductionTilde:InitialValue}
		%		that $\Psi_0=\Psi$,
		\item \label{CompositionSumInductionTilde:Continuity}
		it holds that 
		$\functionANN(\Psi)\in C(\R^d,\R^d)$,
		\item \label{CompositionSumInductionTilde:Realization}
		it holds for all $x\in\R^d$ that 
		$(\functionANN(\Psi))(x)=f_n(x)$,
		\item \label{CompositionSumInductionTilde:Dims}
		it holds that $\hiddenLength(\Psi)={\hiddenLength(\psi)+\sum_{k=0}^{n-1} \hiddenLength(\phi_k)}={\hiddenLength(\psi)+n \,\hiddenLength(\phi_0)}$, 
		%		\item it holds that 
		%		\begin{multline}\label{CompositionSumInductionTilde:ClaimTwo}
		%		\dims(\Psi)=
		%		\big(\ell_{0},\ell_{1},\ell_{2},\dots, \ell_{\mathfrak{L}-1},
		%		l_{0,1}+\hiddenDimId, l_{0,2}+\hiddenDimId,\dots, l_{0,L_0-1}+\hiddenDimId,
		%		\\l_{1,1}+\hiddenDimId, l_{1,2}+\hiddenDimId,\dots, l_{1,L_1-1}+\hiddenDimId,
		%		\dots, l_{n-1,1}+\hiddenDimId, l_{n-1,2}+\hiddenDimId,\dots, l_{n-1,L_{n-1}-1}+\hiddenDimId,d
		%		\big),
		%		%		\in \N^{\mathfrak{L}+1+\sum_{k=0}^{n-1} (L_{k}-1)}
		%		\end{multline}
		and
		\item \label{CompositionSumInductionTilde:Params} 
		it holds that 
				$\paramANN(\Psi)\le 
				\paramANN(\psi)+\sum_{k=0}^{n-1} \big[\tfrac{1}{2}\paramANN(\mathbb{I})+\paramANN(\phi_k)\big]^2$.
		%		=\paramANN(\psi)+n \big[\tfrac{1}{2}\paramANN(\mathbb{I})+\paramANN(\phi_0)\big]^2$. 	
%		\begin{equation}
%		\begin{split}
%		\paramANN(\Psi)\le 
%		\paramANN(\psi)+\sum_{k=0}^{n-1} \big[\tfrac{1}{2}\paramANN(\mathbb{I})+\paramANN(\phi_k)\big]^2.
%		\end{split}
%		\end{equation}
	\end{enumerate}
\end{cor}

\begin{proof}[Proof of Corollary~\ref{Cor:CompositionSumInductionTilde}]
	To prove items~\eqref{CompositionSumInductionTilde:Continuity}--\eqref{CompositionSumInductionTilde:Params} we distinguish between the case $L=1$ and the case $L\in \N\cap [2,\infty)$.
	We first prove  items~\eqref{CompositionSumInductionTilde:Continuity}--\eqref{CompositionSumInductionTilde:Params} in the case $L=1$.
	Observe that Proposition~\ref{Cor:CompositionSumInductionAffine} ensures that there exist $\Psi_n\in \ANNs$, $n\in \N_0$, which satisfy that
	\begin{enumerate}[(I)]
		%		\item \label{CompositionSumInductionAffine:InitialValue}
		%		that $\Psi_0=\psi$,
		\item \label{CompositionSumInductionAffineTilde:Continuity}
		it holds for all $n\in \N_0$ that 
		$\functionANN(\Psi_n)\in C(\R^d,\R^d)$,
		\item \label{CompositionSumInductionAffineTilde:Realization}
		it holds for all $n\in\N_0$, $x\in\R^d$ that $(\functionANN(\Psi_n))(x)=f_n(x)$,
		and
		\item \label{CompositionSumInductionAffineTilde:Dims}
		it holds for all $n\in\N_0$ that $\dims(\Psi_n)=\dims(\psi)$.
	\end{enumerate}
	Next note that the hypothesis that $L=1$ demonstrates that for all $n\in\N_0$ it holds that $\hiddenLength(\phi_n)=0$. Combining this with
	\eqref{CompositionSumInductionAffineTilde:Dims} implies that for all $n\in\N_0$ it holds that
	\begin{equation}\label{CompositionSumInductionTilde:HiddenAffineCase}
		\hiddenLength(\Psi_n)=\hiddenLength(\psi)={\hiddenLength(\psi)+\smallsum_{k=0}^{n-1} \hiddenLength(\phi_k)}=\hiddenLength(\psi)+n \hiddenLength(\phi_0).
	\end{equation}
	In addition, observe that \eqref{CompositionSumInductionAffineTilde:Dims} shows that for all $n\in\N_0$ it holds that
	\begin{equation}
		\paramANN(\Psi_n)=\paramANN(\psi)\le 
		\paramANN(\psi)+\smallsum_{k=0}^{n-1} \big[\tfrac{1}{2}\paramANN(\mathbb{I})+\paramANN(\phi_k)\big]^2.
	\end{equation}
	Combining this and \eqref{CompositionSumInductionTilde:HiddenAffineCase} with \eqref{CompositionSumInductionAffineTilde:Continuity}--\eqref{CompositionSumInductionAffineTilde:Realization} establishes items~\eqref{CompositionSumInductionTilde:Continuity}--\eqref{CompositionSumInductionTilde:Params} in the case $L=1$.
	We now prove  items~\eqref{CompositionSumInductionTilde:Continuity}--\eqref{CompositionSumInductionTilde:Params} in the case $L\in \N\cap [2,\infty)$.
	Note that \eqref{CompositionSumInductionTilde:Recursion}, the fact that $\functionANN(\psi)\in C(\R^d,\R^d)$, the fact that  $\hiddenLength(\psi)={\hiddenLength(\psi)+\sum_{k=0}^{-1} \hiddenLength(\phi_k)}={\hiddenLength(\psi)+0\cdot  \hiddenLength(\phi_0)}$, and the fact that $\paramANN(\psi)= 
	\paramANN(\psi)+\sum_{k=0}^{-1} \big[\tfrac{1}{2}\paramANN(\mathbb{I})+\paramANN(\phi_k)\big]^2$ prove that there exists $\Psi\in \ANNs$
	such that
	\begin{enumerate}[(a)]
		\item \label{CompositionSumInductionTildeZero:Continuity}
		it holds that 
		$\functionANN(\Psi)\in C(\R^d,\R^d)$,
		\item \label{CompositionSumInductionTildeZero:Realization}
		it holds for all $x\in\R^d$ that 
		$(\functionANN(\Psi))(x)=f_0(x)$,
		\item \label{CompositionSumInductionTildeZero:Dims}
		it holds that $\hiddenLength(\Psi)={\hiddenLength(\psi)+\sum_{k=0}^{-1} \hiddenLength(\phi_k)}={\hiddenLength(\psi)+0\cdot \hiddenLength(\phi_0)}$, 
		and
		\item \label{CompositionSumInductionTildeZero:Params} 
		it holds that 
		$\paramANN(\Psi)\le 
		\paramANN(\psi)+\sum_{k=0}^{-1} \big[\tfrac{1}{2}\paramANN(\mathbb{I})+\paramANN(\phi_k)\big]^2$.
	\end{enumerate}
	Moreover, observe that Corollary~\ref{Cor:CompositionSumInduction} and the fact that for all $k\in\N_0$ it holds that $\hiddenLength(\phi_k)=L-1=\hiddenLength(\phi_0)$  ensure that 
	for every $n\in\N$ there exists $\Psi\in \ANNs$
		such that
		\begin{enumerate}[(A)]
			%		\item \label{CompositionSumInduction:InitialValue}
			%		that $\Psi_0=\Psi$,
			\item \label{CompositionSumInductionTildeNonZero:Continuity}
			it holds that 
			$\functionANN(\Psi)\in C(\R^d,\R^d)$,
			\item \label{CompositionSumInductionTildeNonZero:Realization}
			it holds for all $x\in\R^d$ that 
			$(\functionANN(\Psi))(x)=f_n(x)$,
			\item \label{CompositionSumInductionTildeNonZero:Dims}
			it holds that $\hiddenLength(\Psi)={\hiddenLength(\psi)+\sum_{k=0}^{n-1} \hiddenLength(\phi_k)}=\hiddenLength(\psi)+n\,\hiddenLength(\phi_0)$, 
%			\item it holds that 
%			\begin{multline}\label{CompositionSumInductionTildeNonZero:ClaimTwo}
%			\dims(\Psi)=
%			\big(\ell_{0},\ell_{1},\ell_{2},\dots, \ell_{\mathfrak{L}-1},
%			l_{0,1}+\hiddenDimId, l_{0,2}+\hiddenDimId,\dots, l_{0,L_0-1}+\hiddenDimId,
%			\\l_{1,1}+\hiddenDimId, l_{1,2}+\hiddenDimId,\dots, l_{1,L_1-1}+\hiddenDimId,
%			\dots, l_{n-1,1}+\hiddenDimId, l_{n-1,2}+\hiddenDimId,\dots, l_{n-1,L_{n-1}-1}+\hiddenDimId,d
%			\big),
%			%		\in \N^{\mathfrak{L}+1+\sum_{k=0}^{n-1} (L_{k}-1)}
%			\end{multline}
			and
			\item \label{CompositionSumInductionTildeNonZero:Params} 
			it holds that $\paramANN(\Psi)\le \paramANN(\psi)+ \smallsum_{k=0}^{n-1} \big[\tfrac{1}{2}\paramANN(\mathbb{I})+\paramANN(\phi_k)\big]^2$. 	
			%		\begin{equation}\label{CompositionSumInduction:ClaimThree}
			%		\paramANN(\Psi)\le \paramANN(\Psi)+ \smallsum_{k=0}^{n-1} [\paramANN(\mathbb{I})+\paramANN(\phi_k)]^2.
			%		\end{equation}
		\end{enumerate}
		Combining this with \eqref{CompositionSumInductionTildeZero:Continuity}--\eqref{CompositionSumInductionTildeZero:Params} proves items~\eqref{CompositionSumInductionTilde:Continuity}--\eqref{CompositionSumInductionTilde:Params} in the case $L\in \N\cap [2,\infty)$.
	The proof of
	Corollary~\ref{Cor:CompositionSumInductionTilde} is thus completed.
\end{proof}

%% file: ANNrepMultiplePertEuler.tex
\begin{prop}\label{Lemma:DNNrepresentationEulerSpace}
	Let $a\in C(\R,\R)$, $N, d,\hiddenDimId \in \N$, $\mathbb{I},\Phi\in\ANNs$, $A_1,A_2,\dots, A_N\in\R^{d\times d}$
	satisfy  for all  $x\in\R^{d}$ that $2\le\hiddenDimId\le 2d$, 
	$\dims(\mathbb{I}) = (d,\hiddenDimId,d)$, $(\functionANN(\mathbb{I}))(x)=x$,  and $\inDimANN(\Phi)=\outDimANN(\Phi)=d$
%	let $(A_n)_{n\in\{1,2,\dots, N\}}\subseteq\R^{d\times d}$,
	and let 
	$Y_n= (Y^{x,y }_n)_{ (x,y)\in\R^d\times(\R^d)^N} \colon\allowbreak  \R^d\times(\R^d)^N \to \R^d $, $n\in\{0,1,\dots, N\}$,
	be the
	functions 	
	which satisfy for all 
	$n\in\{0,1,\dots, N-1\}$,
	$ x \in \R^d $,
	$y=(y_1,y_2,\dots, y_N)\in (\R^d)^N$  
	that $\affineProcess^{x,y} _0=x$ and
	\begin{equation}
	\label{DNNrepresentationEulerSpace:Y_processes}
	\begin{split}
	&\affineProcess^{x,y} _{n+1} 
	=
	\affineProcess^{x,y} _{n}+ A_{n+1}\big((\functionANN(\Phi))( 
	\affineProcess^{x,y} _{n})\big)
	+
	y_{n+1}
	\end{split}
	\end{equation}
	(cf.\ Definition~\ref{Def:ANN} and Definition~\ref{Definition:ANNrealization}).
	Then there exists $(\Psi_{n,y})_{ (n,y)\in\{0,1,\dots,N\}\times (\R^d)^N} \subseteq \ANNs$ such that
	\begin{enumerate}[(i)]
		\item \label{itemRepresentation:DNNrepresentationEulerSpaceFunction} 
		it holds for all $n\in\{0,1,\dots,N\}$, $y\in(\R^d)^N$ that
		$\functionANN(\Psi_{n,y}) \in C(\R^d, \R^{d})$,
		\item it holds for all $n\in\{0,1,\dots,N\}$, $y\in(\R^d)^N$, $x \in \R^d$ that $(\functionANN(\Psi_{n,y}))(x)=Y_n^{x,y}$,
		\item 
		it holds for all $n\in\{0,1,\dots,N\}$, $y\in(\R^d)^N$ that 
		\begin{equation}
		\hiddenLength(\Psi_{n,y})=\hiddenLength(\mathbb{I})+n \,\hiddenLength(\Phi)=1+n\, \hiddenLength(\Phi),
		\end{equation}
		\item \label{itemRepresentation:DNNrepresentationEulerSpaceParams} 
		it holds for all $n\in\{0,1,\dots,N\}$, $y\in(\R^d)^N$ that
		\begin{equation}
		\paramANN(\Psi_{n,y})\le \paramANN(\mathbb{I})+n \big[\tfrac{1}{2}\paramANN(\mathbb{I})+\paramANN(\Phi)\big]^{\!2},
		\end{equation}
		\item \label{itemContinuity:DNNrepresentationEulerSpaceFunction} 
		it holds for all $n\in\{0,1,\dots,N\}$, $x\in\R^d$ that 
		\begin{equation}
		\big[(\R^d)^N\ni y\mapsto (\functionANN (\Psi_{n,y}))(x)\in\R^d\big]\in C\big((\R^d)^N,\R^d\big),
		\end{equation}
		and
		\item \label{itemAdaptedness:DNNrepresentationEulerSpaceFunction} 
		it holds 
		for all $n\in\{0,1,\dots, N\}$, $m\in \N_0\cap [0,n]$,  $x\in\R^d$, $y=(y_1,y_2,\dots, y_N), \allowbreak z=(z_1,z_2,\dots, z_N)\in (\R^d)^N$  with $\forall\, k\in \N\cap [0,n]\colon y_k=z_k$ that 
		\begin{equation}
	(\functionANN(\Psi_{m,y}))(x)=(\functionANN(\Psi_{m,z}))(x).
		\end{equation}
	\end{enumerate}
\end{prop}

\begin{proof}[Proof of Proposition~\ref{Lemma:DNNrepresentationEulerSpace}]	
	Throughout this proof 
	let $l_0,l_1,\dots, l_{\lengthANN(\Phi)}\in\N$ satisfy that $\dims(\Phi)=(l_0,l_1,\dots, l_{\lengthANN(\Phi)})$, let  
%	$(\affineMap_{n,b})_{n\in\{1,2,\dots,N\},b\in\R^d}\subseteq\ANNs$ 
	$\affineMap_{n,b}\in (\R^{d\times d}\times \R^d)\subseteq \ANNs$, $n\in\{1,2,\dots,N\}$, $b\in\R^d$,
	satisfy for all $n\in\{1,2,\dots,\allowbreak N\}$, $b\in\R^d$ that 
%			$\affineMap_{n,b}=\big(A_n,b\big)$,
	%	$\rho_{n,y}=\compANN{\affineMap_{n,y}}{\Phi}$.
	\begin{equation}
	\affineMap_{n,b}=(A_n,b)\in (\R^{d\times d}\times \R^d),
	\end{equation}
%	let $(\rho_{n,y})_{n\in\{1,2,\dots,N\},y\in(\R^d)^N}\subseteq \ANNs$ satisfy for all $n\in \{1,2,\dots, N\}$, $y=(y_1,y_2,\dots, y_N)\allowbreak\in(\R^d)^N$ that
%	$\rho_{n,y}=\compANN{\affineMap_{n,y_{n}}}{\Phi}$
%	(cf.\ Definition~\ref{Definition:ANNcomposition}).
	let $\rho_{n,y}\in \ANNs$, $n\in\N$, $y\in(\R^d)^N$, satisfy for all $n\in \N$, $y=(y_1,y_2,\dots, y_N)\allowbreak\in(\R^d)^N$ that
\begin{equation}
	\rho_{n,y}=\compANN{\affineMap_{\min\{n,N\},y_{\min\{n,N\}}}}{\Phi}
\end{equation}
	(cf.\ Definition~\ref{Definition:ANNcomposition}), 
	and let 	$\mathcal{Y}_n= (\mathcal{Y}^{x,y }_n)_{ (x,y)\in\R^d\times(\R^d)^N} \colon\allowbreak  \R^d\times(\R^d)^N \to \R^d $, $n\in\N_0$,
	be the
	functions 	
	which satisfy for all 
	$n\in\N_0$,
	$ x \in \R^d $,
	$y=(y_1,y_2,\dots, y_N)\in (\R^d)^N$  
	that $\mathcal{Y}^{x,y} _0=x$ and
	\begin{equation}
	\label{DNNrepresentationEulerSpace:Z_processes}
	\begin{split}
	&\mathcal{Y}^{x,y}_{n+1} 
	=
	\mathcal{Y}^{x,y}_{n}+ \big(\functionANN(\rho_{n+1,y})\big)( 
	\mathcal{Y}^{x,y}_{n}).
	\end{split}
	\end{equation}
	Observe that  item~\eqref{PropertiesOfCompositions:Dims} in Proposition~\ref{Lemma:PropertiesOfCompositions} and the fact that for all $n\in \{1,2,\dots, N\}$, $y=(y_1,y_2,\dots, y_N)\allowbreak\in(\R^d)^N$ it holds that $\rho_{n,y}=\compANN{\affineMap_{n,y_n}}{\Phi}$ prove that for all $n\in \{1,2,\dots, N\}$, $y\in(\R^d)^N$ it holds that
	$\dims(\rho_{n,y})=\dims(\Phi)=(l_0,l_1,\dots, l_{\lengthANN(\Phi)})$.
%	Combining this with the hypothesis that $\dims(\mathbb{I}) = (d,\hiddenDimId,d)$ ensures that for all $n\in \N\cap (0,N)$, $y\in(\R^d)^N$ it holds that
%	\begin{equation}
%	(\dims(\mathbb{I}))_{\lengthANN(\mathbb{I})-1}=\hiddenDimId\le l_{\lengthANN(\Phi)-1}+\hiddenDimId,
%	\quad(\dims(\rho_{n,y}))_{\lengthANN(\rho_{n,y})-1}\le (\dims(\rho_{n,y+1}))_{\lengthANN(\rho_{n,y+1})-1},
%	\end{equation}
%	and $\paramANN(\rho_{n,y})=\paramANN(\Phi)$.
	Corollary~\ref{Cor:CompositionSumInductionTilde} (with $a=a$, $d=d$, $\hiddenDimId=\hiddenDimId$, $L=\lengthANN(\Phi)$, $\mathfrak{L}=2$, $\ell_0=d$, $\ell_1=\hiddenDimId$, $\ell_2=d$, $\mathbb{I}=\mathbb{I}$, $\psi=\mathbb{I}$, $(\N_0\ni n\mapsto \phi_n\in\ANNs)=(\N_0\ni n\mapsto \rho_{n+1,y}\in\ANNs)$, $(\N_0\times\{0,1,\dots,\mathcal{L}(\Phi)\}\ni(n,k)\mapsto l_{n,k}\in \N)=(\N_0\times\{0,1,\dots,\mathcal{L}(\Phi)\}\ni(n,k)\mapsto l_k\in\N)$, $(\N_0\ni n\mapsto f_n\in C(\R^d,\R^d))=(\N_0\ni n\mapsto (\R^d\ni x\mapsto \mathcal{Y}^{x,y}_n\in \R^d)\in C(\R^d,\R^d))$ for $y\in (\R^d)^N$ in the notation of Corollary~\ref{Cor:CompositionSumInductionTilde})   and the fact that for all $x\in\R^d$, $y\in(\R^d)^N$ it holds that $(\functionANN(\mathbb{I}))(x)=x=\mathcal{Y}_0^{x,y}=Y_0^{x,y}$  hence prove that 
	there exist $\Psi_{n,y} \in \ANNs$, $(n,y)\in\{0,1,\dots,N\}\times(\R^d)^N$, which satisfy that 
	\begin{enumerate}[(I)]
		\item\label{DNNrepresentationEulerSpace:EnuContinuity} it holds  for all $n\in \{0,1,\dots, N\}$, $y\in(\R^d)^N$ that $\functionANN(\Psi_{n,y})\in C(\R^d,\R^d)$,
		\item\label{DNNrepresentationEulerSpace:EnuOne} it holds for all $n\in \{0,1,\dots, N\}$, $y\in(\R^d)^N$, $x\in\R^d$ that  $(\functionANN(\Psi_{n,y}))(x)=\mathcal{Y}_n^{x,y}=Y_n^{x,y}$,
		\item\label{DNNrepresentationEulerSpace:EnuTwo} it holds  for all $n\in \{0,1,\dots, N\}$, $y\in(\R^d)^N$ that
		\begin{equation}
		\hiddenLength(\Psi_{n,y})=\hiddenLength(\mathbb{I})+\smallsum\limits_{k=0}^{n-1} \hiddenLength(\rho_{k+1,y})=1+n\hiddenLength(\Phi),
		\end{equation}
		and
		\item\label{DNNrepresentationEulerSpace:EnuThree} it holds for all $n\in \{0,1,\dots, N\}$, $y\in(\R^d)^N$ that
		\begin{equation}
		\paramANN(\Psi_{n,y})\le \paramANN(\mathbb{I})+ \smallsum\limits_{k=0}^{n-1} \big[\tfrac{1}{2}\paramANN(\mathbb{I})+\paramANN(\rho_{k+1,y})\big]^{\!2}
		=\paramANN(\mathbb{I})+n \big[\tfrac{1}{2}\paramANN(\mathbb{I})+\paramANN(\Phi)\big]^{\!2}.
		\end{equation}
	\end{enumerate}	
	Next we claim that for all $n\in\{0,1,\dots,N\}$ it holds that 
	\begin{equation}\label{DNNrepresentationEulerSpace:continuousInduction}
	\forallDist x\in\R^d\colon\,\Big[\big((\R^d)^N\ni y\mapsto Y^{x,y}_n\in\R^d\big)\in C\big((\R^d)^N,\R^d\big)\Big].
	\end{equation}
	We now prove \eqref{DNNrepresentationEulerSpace:continuousInduction} by induction on $n\in\{0,1,\dots,N\}$.
	Note that the fact that for all $x\in\R^d$, $y\in (\R^d)^N$ it holds that $Y^{x,y}_0=x$ proves \eqref{DNNrepresentationEulerSpace:continuousInduction} in the base case $n=0$. For the induction step observe that \eqref{DNNrepresentationEulerSpace:Y_processes} and the fact that $\functionANN(\Phi)\in C(\R^d,\R^d)$ ensure that
	for all $n\in\{0,1,\dots,N-1\}$ with 
	\begin{equation}
	\forall\, x\in\R^d\colon\,\Big[\big((\R^d)^N\ni y\mapsto Y^{x,y}_n\in\R^d\big)\in C\big((\R^d)^N,\R^d\big)\Big]
	\end{equation}
	it holds that 
	\begin{equation}
	\forall\, x\in\R^d\colon\,\Big[\big((\R^d)^N\ni y\mapsto Y^{x,y}_{n+1}\in\R^d\big)\in C\big((\R^d)^N,\R^d\big)\Big].
	\end{equation}
	Induction thus proves \eqref{DNNrepresentationEulerSpace:continuousInduction}.
	In addition, observe that \eqref{DNNrepresentationEulerSpace:continuousInduction} and \eqref{DNNrepresentationEulerSpace:EnuOne} imply that 
	for all $n\in\{0,1,\dots,N\}$, $x\in\R^d$ it holds that 
	\begin{equation}\label{DNNrepresentationEulerSpace:continuousANN}
	\big((\R^d)^N\ni y\mapsto (\functionANN(\Psi_{n,y}))(x)\in\R^d\big)\in C\big((\R^d)^N,\R^d\big).
	\end{equation}
	Next let $n\in\{0,1,\dots, N\}$,  $x\in\R^d$, $y=(y_1,y_2,\dots, y_N),\allowbreak z=(z_1,z_2,\dots, z_N)\in (\R^d)^N$ satisfy for all $k\in\N\cap [0,n]$ that $y_k=z_k$. 
	We claim that for all $m\in\N_0\cap [0,n]$ it holds that 
	\begin{equation}\label{DNNrepresentationEulerSpace:adaptedInduction}
Y^{x,y}_m=Y^{x,z}_m.
	\end{equation}
	We now prove \eqref{DNNrepresentationEulerSpace:adaptedInduction} by induction on $m\in \N_0\cap [0,n]$.
	Note that the fact that $Y^{x,y}_0=x=Y^{x,z}_0$ implies \eqref{DNNrepresentationEulerSpace:adaptedInduction} in the base case $m=0$. 
	For the induction step 
%	$\N_0\cap (-\infty,n)\ni k\to k+1\in \N_0\cap (-\infty,n]$ 
	observe that \eqref{DNNrepresentationEulerSpace:Y_processes} and
	the fact  that for all $k\in \N\cap [0,n]$ it holds
	that $y_k=z_k$ ensure that 
	for all  $m\in \N_0\cap (-\infty,n)$ with $Y^{x,y}_m=Y^{x,z}_m$ it holds that 
	\begin{equation}
	\begin{split}
	\affineProcess^{x,y}_{m+1} 
	&=
	\affineProcess^{x,y}_{m}+ A_{m+1}\big((\functionANN(\Phi)) ( 
	\affineProcess^{x,y}_{m} )\big)
	+
	y_{m+1}
	\\&=
	\affineProcess^{x,z}_{m}+ A_{m+1}\big((\functionANN(\Phi)) ( 
	\affineProcess^{x,z}_{m} )\big)
	+
	z_{m+1}
	= \affineProcess^{x,z} _{m+1}. 
	\end{split}
	\end{equation}
	Induction thus proves \eqref{DNNrepresentationEulerSpace:adaptedInduction}. 
	Note that \eqref{DNNrepresentationEulerSpace:adaptedInduction} and \eqref{DNNrepresentationEulerSpace:EnuOne} assure that 
	for all $n\in\{0,1,\dots, N\}$, $m\in \N_0\cap [0,n]$,  $x\in\R^d$, $y=(y_1,y_2,\dots, y_N)$, $z=(z_1,z_2,\dots, z_N)\in (\R^d)^N$ with $\forall\, k\in \N\cap [0,n]\colon y_k=z_k$ it holds that 
	\begin{equation}\label{DNNrepresentationEulerSpace:adaptedANN}
	(\functionANN(\Psi_{m,y}))(x)=(\functionANN(\Psi_{m,z}))(x).
	\end{equation}
	Combining this with \eqref{DNNrepresentationEulerSpace:continuousANN} and \eqref{DNNrepresentationEulerSpace:EnuContinuity}--\eqref{DNNrepresentationEulerSpace:EnuThree} establishes items~\eqref{itemRepresentation:DNNrepresentationEulerSpaceFunction}--\eqref{itemAdaptedness:DNNrepresentationEulerSpaceFunction}. 	
	The proof of Proposition~\ref{Lemma:DNNrepresentationEulerSpace} is thus completed.
\end{proof}

%% file: ANNIntroSec3.tex
This section establishes in Theorem~\ref{Thm:ApproxOfEulerWithGronwall} in Subsection~\ref{subsec:ANNapproxEuler} below the main result of this article. Some of the material presented in Subsection~\ref{subsec:ANNApproxforSquare} and Subsection~\ref{subsec:ANNApproxProducts} are well-known concepts and results in the scientific literature.
In particular, the material in Subsection~\ref{subsubsec:Explicit approximations for the square function on $[0,1]$} and Subsection~\ref{subsubsec:ANN approximations for the square function on $[0,1]$} consists mainly of reformulations of concepts and results in Elbr\"achter et al.~\cite[Appendix~A.3 and Appendix~A.4]{ElbraechterSchwab2018}. Moreover, our proof of Proposition~\ref{Lemma:ApproxOfProduct} in Subsection~\ref{subsubsec:ANN approximations for one-dimensional products} below is inspired by Elbr\"achter et al.~\cite[Section~6]{ElbraechterSchwab2018} and Yarotsky \cite[Section~3.2]{yarotsky2017error}. Furthermore, Lemma~\ref{Lem:EulerCont} and Lemma~\ref{Lemma:DNNhatFunction} are elementary and essentially well-known in the scientific literature. In addition, our proof of Lemma~\ref{Lemma:Gronwall} is based on a well-known Gronwall argument. 

%% file: ExplicitApproxSquare.tex
\begin{lemma}\label{Lemma:productAuxiliaryOne}
	Let $g_n\colon \R\to[0,1]$, $n\in\N$, be the functions which satisfy for all $n\in\N$, $x\in\R$ that
%	\begin{align}\label{NNsgDef}
%	g_n(x)=\begin{cases}
%	2x & \colon n=1,x<\tfrac{1}{2}\\
%	2-2x & \colon n=1,x\geq\tfrac{1}{2}\\
%	g_1(g_{n-1}(x)) & \colon n\geq 1
%	\end{cases}.
%	\end{align}
	\begin{align}\label{NNsgDef}
	g_1(x)=\begin{cases}
	2x & \colon x\in [0,\tfrac{1}{2})\\
	2-2x & \colon x\in [\tfrac{1}{2},1]\\
	0 & \colon x\in \R\backslash[0,1]\\
	\end{cases}
	\end{align}
	and $g_{n+1}(x)=g_1(g_{n}(x))$.
	Then
	\begin{enumerate}[(i)]
		\item \label{productAuxiliaryOne:interior}
		it holds for all $n\in\N$, $k\in\zeroto{2^{n-1}-1}$, $x\in\brac{\tfrac{k}{2^{n-1}},{\tfrac{k+1}{2^{n-1}}}}$ that
		\begin{equation}\label{NNsClaimLemma}
		g_n(x)=\begin{cases}
		2^n(x-\tfrac{2k}{2^n}) & \colon x\in\brac{\tfrac{2k}{2^n},{\tfrac{2k+1}{2^n}}}\\
		2^n(\tfrac{2k+2}{2^n}-x) & \colon x\in\brac{\tfrac{2k+1}{2^n},{\tfrac{2k+2}{2^n}}}
		\end{cases}
		\end{equation}
		and
		\item \label{productAuxiliaryOne:exterior}
		it holds for all $n\in\N$, $x\in \R\backslash[0,1]$ that $g_n(x)=0$. 
	\end{enumerate}
\end{lemma}

\begin{proof}[Proof of Lemma~\ref{Lemma:productAuxiliaryOne}]	
First, we claim that for all $n\in\N$ it holds that
%\begin{equation}\label{NNsClaim}
%\forallDist k\in\zeroto{2^{n-1}-1}, x\in\brac{\tfrac{k}{2^{n-1}},{\tfrac{k+1}{2^{n-1}}}}\colon
% g_n(x)=\begin{cases}
%2^n(x-\tfrac{2k}{2^n}) & \colon x\in\brac{\tfrac{2k}{2^n},{\tfrac{2k+1}{2^n}}}\\
%2^n(\tfrac{2k+2}{2^n}-x) & \colon x\in\brac{\tfrac{2k+1}{2^n},{\tfrac{2k+2}{2^n}}}
%\end{cases}.
%\end{equation}
\begin{multline}\label{NNsClaim}
\Bigg(\forallDist k\in\zeroto{2^{n-1}-1}, x\in\brac{\tfrac{k}{2^{n-1}},{\tfrac{k+1}{2^{n-1}}}}\colon
\\ g_n(x)=\begin{cases}
2^n(x-\tfrac{2k}{2^n}) & \colon x\in\brac{\tfrac{2k}{2^n},{\tfrac{2k+1}{2^n}}}\\
2^n(\tfrac{2k+2}{2^n}-x) & \colon x\in\brac{\tfrac{2k+1}{2^n},{\tfrac{2k+2}{2^n}}}
\end{cases}\Bigg).
\end{multline}
%\begin{multline}\label{NNsClaim}
%	\left(\forallDist k\in\zeroto{2^{n-1}-1}, x\in\brac{\tfrac{k}{2^{n-1}},{\tfrac{k+1}{2^{n-1}}}}\colon\right.
%	\\\left. g_n(x)=\begin{cases}
%	2^n(x-\tfrac{2k}{2^n}) & \colon x\in\brac{\tfrac{2k}{2^n},{\tfrac{2k+1}{2^n}}}\\
%	2^n(\tfrac{2k+2}{2^n}-x) & \colon x\in\brac{\tfrac{2k+1}{2^n},{\tfrac{2k+2}{2^n}}}
%	\end{cases}\right).
%\end{multline}
%We prove item~\eqref{productAuxiliaryOne:interior} by induction on $n\in\N$. 
We now prove \eqref{NNsClaim} by induction on $n\in\N$. 
Note that \eqref{NNsgDef} establishes \eqref{NNsClaim} in the base case $n=1$. 
For the induction step $\N\ni n\to n+1\in\N\cap [2,\infty)$ 
assume that there exists $n\in\N$ such that for all $k\in\zeroto{2^{n-1}-1}$, $x\in\bracbig{\tfrac{k}{2^{n-1}},{\tfrac{k+1}{2^{n-1}}}}$ it holds that
\begin{equation}\label{NNsClaimInduction}
g_n(x)=\begin{cases}
2^n(x-\tfrac{2k}{2^n}) & \colon x\in\brac{\tfrac{2k}{2^n},{\tfrac{2k+1}{2^n}}}\\
2^n(\tfrac{2k+2}{2^n}-x) & \colon x\in\brac{\tfrac{2k+1}{2^n},{\tfrac{2k+2}{2^n}}}
\end{cases}.
\end{equation}
Observe that \eqref{NNsgDef} and \eqref{NNsClaimInduction}
imply that for all  $l\in\zeroto{2^{n-1}-1}$, $x\in\bracbig{\tfrac{2l}{2^n},\tfrac{2l+(\nicefrac{1}{2})}{2^n}}$ it holds that
\begin{equation}\begin{split}\label{NNsEq1}
g_{n+1}(x)&=g_1(g_n(x))=g_1(2^n(x-\tfrac{2l}{2^n}))=2\brac{2^n(x-\tfrac{2l}{2^n})}=2^{n+1}(x-\tfrac{2l}{2^n}).
%=2^{n+1}(x-\tfrac{4l}{2^{n+1}}).
\end{split}\end{equation}
In addition, note that \eqref{NNsgDef}  and \eqref{NNsClaimInduction}
ensure that for all  $l\in\zeroto{2^{n-1}-1}$, $x\in\bracbig{\tfrac{2l+(\nicefrac{1}{2})}{2^n},\tfrac{2l+1}{2^n}}$  it holds that
\begin{equation}\begin{split}\label{NNsEq2}
	g_{n+1}(x)&=g_1(g_n(x))=g_1(2^n(x-\tfrac{2l}{2^n}))=2-2\brac{2^n(x-\tfrac{2l}{2^n})}\\
	&=2-2^{n+1}x+4l=2^{n+1}(\tfrac{4l+2}{2^{n+1}}-x).
%	\\
%	&=2^{n+1}(\tfrac{2(2l+1)}{2^{n+1}}-x).
\end{split}\end{equation}  
Moreover, observe that \eqref{NNsgDef} and \eqref{NNsClaimInduction}
demonstrate that for all  $l\in\zeroto{2^{n-1}-1}$, $x\in\bracbig{\tfrac{2l+1}{2^n},\tfrac{2l+(\nicefrac{3}{2})}{2^n}}$   it holds that
\begin{equation}\begin{split}\label{NNsEq3}
g_{n+1}(x)&=g_1(g_n(x))=g_1(2^n(\tfrac{2l+2}{2^n}-x))=2-2\brac{2^n(\tfrac{2l+2}{2^n}-x)}\\
&=2-2(2l+2)+2^{n+1}x=2^{n+1}x-4l-2\\
&=2^{n+1}(x-\tfrac{4l+2}{2^{n+1}}).
\end{split}\end{equation}
Next note that \eqref{NNsgDef} and \eqref{NNsClaimInduction}
prove that for all  $l\in\zeroto{2^{n-1}-1}$, $x\in\bracbig{\tfrac{2l+(\nicefrac{3}{2})}{2^n},\tfrac{2l+2}{2^n}}$    it holds that
\begin{equation}\begin{split}\label{NNsEq4}
g_{n+1}(x)&=g_1(g_n(x))=g_1(2^n(\tfrac{2l+2}{2^n}-x))=2\brac{2^n(\tfrac{2l+2}{2^n}-x)}=2^{n+1}(\tfrac{2l+2}{2^n}-x).
%=2^{n+1}(\tfrac{4l+4}{2^{n+1}}-x).
\end{split}\end{equation} 
%Next observe that for all  $k\in\zeroto{2^n-1}$ there exists $l\in\zeroto{2^{n-1}-1}$ 
%which satisfies that
%\begin{equation}\label{NNsS1}
%\brac{\tfrac{2k}{2^{n+1}},{\tfrac{2k+1}{2^{n+1}}}}=\brac{\tfrac{2l}{2^n},\tfrac{2l+(\nicefrac{1}{2})}{2^n}}\quad\mathrm{or}\quad\brac{\tfrac{2k}{2^{n+1}},{\tfrac{2k+1}{2^{n+1}}}}=\brac{\tfrac{2l+1}{2^n},\tfrac{2l+(\nicefrac{3}{2})}{2^n}}.
%\end{equation}
%Furthermore, note that for all  $k\in\zeroto{2^n-1}$ there exists $l\in\zeroto{2^{n-1}-1}$ which satisfies that
%\begin{equation}\label{NNsS2}
%\brac{\tfrac{2k+1}{2^{n+1}},{\tfrac{2k+2}{2^{n+1}}}}=\brac{\tfrac{2l+(\nicefrac{1}{2})}{2^n},\tfrac{2l+1}{2^n}}\quad\mathrm{or}\quad\brac{\tfrac{2k+1}{2^{n+1}},{\tfrac{2k+2}{2^{n+1}}}}=\brac{\tfrac{2l+(\nicefrac{3}{2})}{2^n},\tfrac{2l+2}{2^n}}.
%\end{equation}
Moreover, observe that for all $k\in\{0,2,4,6,\dots\}\cap [0,2^n-2]$ it holds that 
\begin{equation}
	\brac{\tfrac{2k}{2^{n+1}},{\tfrac{2k+1}{2^{n+1}}}}=\brac{\tfrac{2(\nicefrac{k}{2})}{2^n},\tfrac{2(\nicefrac{k}{2})+(\nicefrac{1}{2})}{2^n}},
	\qquad
	\brac{\tfrac{2k+1}{2^{n+1}},{\tfrac{2k+2}{2^{n+1}}}}=\brac{\tfrac{2(\nicefrac{k}{2})+(\nicefrac{1}{2})}{2^n},\tfrac{2(\nicefrac{k}{2})+1}{2^n}},
\end{equation}
%$\brac{\tfrac{2k}{2^{n+1}},{\tfrac{2k+1}{2^{n+1}}}}=\brac{\tfrac{2(\nicefrac{k}{2})}{2^n},\tfrac{2(\nicefrac{k}{2})+(\nicefrac{1}{2})}{2^n}}$
and $\nicefrac{k}{2}\in\zeroto{2^{n-1}-1}$.
This, \eqref{NNsEq1}, and \eqref{NNsEq2} demonstrate that for all $k\in\{0,2,4,6,\dots\}\cap [0,2^n-2]$, $x\in\bracbig{\tfrac{k}{2^{n}},{\tfrac{k+1}{2^{n}}}}$ it holds that 
\begin{equation}\label{inductionStepEvenNumbers}
	\begin{split}
	g_{n+1}(x)
	&=\begin{cases}
	2^{n+1}(x-\tfrac{2(\nicefrac{k}{2})}{2^{n}}) & \colon x\in\brac{\tfrac{2(\nicefrac{k}{2})}{2^n},\tfrac{2(\nicefrac{k}{2})+(\nicefrac{1}{2})}{2^n}}\\
	2^{n+1}(\tfrac{4(\nicefrac{k}{2})+2}{2^{n+1}}-x) & \colon x\in\brac{\tfrac{2(\nicefrac{k}{2})+(\nicefrac{1}{2})}{2^n},\tfrac{2(\nicefrac{k}{2})+1}{2^n}}
	\end{cases}
	\\&=\begin{cases}
	2^{n+1}(x-\tfrac{2k}{2^{n+1}}) & \colon x\in\brac{\tfrac{2k}{2^{n+1}},{\tfrac{2k+1}{2^{n+1}}}}\\
	2^{n+1}(\tfrac{2k+2}{2^{n+1}}-x) & \colon x\in\brac{\tfrac{2k+1}{2^{n+1}},{\tfrac{2k+2}{2^{n+1}}}}
	\end{cases}.
	\end{split}
\end{equation}
In addition, observe  
that for all $k\in\{1,3,5,7,\dots\}\cap [1,2^n-1]$ it holds that 
\begin{equation}
\begin{split}
\brac{\tfrac{2k}{2^{n+1}},{\tfrac{2k+1}{2^{n+1}}}}=\brac{\tfrac{2(\nicefrac{(k-1)}{2})+1}{2^n},\tfrac{2(\nicefrac{(k-1)}{2})+(\nicefrac{3}{2})}{2^n}},
\\
\brac{\tfrac{2k+1}{2^{n+1}},{\tfrac{2k+2}{2^{n+1}}}}=\brac{\tfrac{2(\nicefrac{(k-1)}{2})+(\nicefrac{3}{2})}{2^n},\tfrac{2(\nicefrac{(k-1)}{2})+2}{2^n}},
\end{split}
\end{equation}
%$\brac{\tfrac{2k}{2^{n+1}},{\tfrac{2k+1}{2^{n+1}}}}=\brac{\tfrac{2(\nicefrac{k}{2})}{2^n},\tfrac{2(\nicefrac{k}{2})+(\nicefrac{1}{2})}{2^n}}$
and $\nicefrac{(k-1)}{2}\in\zeroto{2^{n-1}-1}$.
This, \eqref{NNsEq3}, and \eqref{NNsEq4} demonstrate that for all $k\in\{1,3,5,7,\dots\}\cap [1,2^n-1]$, $x\in\bracbig{\tfrac{k}{2^{n}},{\tfrac{k+1}{2^{n}}}}$ it holds that 
\begin{equation}
\begin{split}
g_{n+1}(x)
&=\begin{cases}
2^{n+1}(x-\tfrac{4(\nicefrac{(k-1)}{2})+2}{2^{n+1}}) & \colon x\in\brac{\tfrac{2(\nicefrac{(k-1)}{2})+1}{2^n},\tfrac{2(\nicefrac{(k-1)}{2})+(\nicefrac{3}{2})}{2^n}}\\
2^{n+1}(\tfrac{2(\nicefrac{(k-1)}{2})+2}{2^{n}}-x) & \colon x\in\brac{\tfrac{2(\nicefrac{(k-1)}{2})+(\nicefrac{3}{2})}{2^n},\tfrac{2(\nicefrac{(k-1)}{2})+2}{2^n}}
\end{cases}
\\&=\begin{cases}
2^{n+1}(x-\tfrac{2k}{2^{n+1}}) & \colon x\in\brac{\tfrac{2k}{2^{n+1}},{\tfrac{2k+1}{2^{n+1}}}}\\
2^{n+1}(\tfrac{2k+2}{2^{n+1}}-x) & \colon x\in\brac{\tfrac{2k+1}{2^{n+1}},{\tfrac{2k+2}{2^{n+1}}}}
\end{cases}.
\end{split}
\end{equation}
Combining this with \eqref{inductionStepEvenNumbers} ensures that 
for all $k\in\zeroto{2^n-1}$, $x\in\bracbig{\tfrac{k}{2^{n}},{\tfrac{k+1}{2^{n}}}}$ it holds that 
\begin{equation}
\begin{split}
g_{n+1}(x)
=\begin{cases}
2^{n+1}(x-\tfrac{2k}{2^{n+1}}) & \colon x\in\brac{\tfrac{2k}{2^{n+1}},{\tfrac{2k+1}{2^{n+1}}}}\\
2^{n+1}(\tfrac{2k+2}{2^{n+1}}-x) & \colon x\in\brac{\tfrac{2k+1}{2^{n+1}},{\tfrac{2k+2}{2^{n+1}}}}
\end{cases}.
\end{split}
\end{equation}
Induction thus proves \eqref{NNsClaim}.
Observe that \eqref{NNsClaim} establishes item~\eqref{productAuxiliaryOne:interior}.
 Next we claim that for all $n\in\N$ it holds that 
 \begin{equation}\label{NNsClaimTwo}
 	\forallDist x\in \R\backslash[0,1]\colon\, g_n(x)=0.
 \end{equation} 
We now  
 prove \eqref{NNsClaimTwo} by induction on $n\in\N$. 
 Note that \eqref{NNsgDef} establishes  \eqref{NNsClaimTwo} in the base case $n=1$. 
% For the induction step $\N\ni n\to n+1\in\N\cap [2,\infty)$ 
%assume that there exists $n\in\N$ such that for all $x\in\R\backslash [0,1]$ it holds that $g_n(x)=0$.
 For the induction step observe that \eqref{NNsgDef} ensures that for all $n\in\N$ with 
 $(\forall x\in \R\backslash[0,1]\colon g_n(x)=0)$ it holds that 
 \begin{equation}
 	\big(\forallDist x\in \R\backslash[0,1]\colon\,g_{n+1}(x)=g_1(g_{n}(x))=g_1(0)=0\big).
 \end{equation}
 Induction thus proves \eqref{NNsClaimTwo}.
 Note that \eqref{NNsClaimTwo} establishes item~\eqref{productAuxiliaryOne:exterior}.
	The proof of Lemma~\ref{Lemma:productAuxiliaryOne} is thus completed.
\end{proof}

\begin{lemma}\label{Lemma:productAuxiliaryTwo}
	Let $g_n\colon [0,1]\to[0,1]$, $n\in\N$, be the functions which satisfy for all $n\in\N$, $x\in[0,1]$ that
%	\begin{equation}\label{productAuxiliaryTwo:NNsgDef}
%	g_n(x)=\begin{cases}
%	2x & \colon n=1,x<\tfrac{1}{2}\\
%	2-2x & \colon n=1,x\geq\tfrac{1}{2}\\
%	g_1(g_{n-1}(x)) & \colon n\geq 1
%	\end{cases}
%	\end{equation}
		\begin{equation}\label{productAuxiliaryTwo:NNsgDef}
		g_1(x)=\begin{cases}
		2x & \colon x\in [0,\tfrac{1}{2})\\
		2-2x & \colon x\in [\tfrac{1}{2},1]\\
		\end{cases}
		\end{equation}
		and $g_{n+1}(x)=g_1(g_n(x))$,
	and let $f_n\colon [0,1]\to[0,1]$, $n\in\N_0$, 
	be the functions which satisfy for all
	$n\in\N_0$, $k\in\zeroto{2^n-1}$, $x\in\bracLclosedRopen{\tfrac{k}{2^n},\tfrac{k+1}{2^n}}$ 
	that $f_n(1)=1$ and 
	\begin{equation}\label{productAuxiliaryTwo:NNsfmDef}
	f_n(x)=\brac{\tfrac{2k+1}{2^n}}x-\tfrac{(k^2+k)}{2^{2n}}.
	\end{equation}
	Then 
%	\begin{enumerate}[(i)]
%		\item it holds for all $n\in\N$, $x\in[0,1]$ that
%		\begin{equation}\label{NNsfmSum}
%		f_n(x)=x-\sum_{s=1}^m 2^{-2s} g_s(x)
%		\end{equation}
%		and 
%		\item it holds for all $n\in\N$ that
%		\begin{equation}
%		\sup_{x\in[0,1]}\abs{x^2-f_n(x)}
%		\le 2^{-2m-2}.
%		\end{equation}
%	\end{enumerate}
it holds for all $n\in\N_0$, $x\in[0,1]$ that
\begin{equation}\label{NNsfmSum}
f_n(x)=x-\left[\smallsum\limits_{m=1}^n\big( 2^{-2m} g_m(x)\big)\right]
\qandq \abs{x^2-f_n(x)}
\le 2^{-2n-2}.
\end{equation}
%
%it holds for all $n\in\N$, $x\in[0,1]$ that
%		\begin{equation}\label{NNsfmSum}
%		f_n(x)=x-\sum_{s=1}^m 2^{-2s} g_s(x)
%		\end{equation}
%		and 
%		\begin{equation}
%				\sup_{x\in[0,1]}\abs{x^2-f_n(x)}
%				\le 2^{-2m-2}.
%				\end{equation}
\end{lemma}

\begin{proof}[Proof of Lemma~\ref{Lemma:productAuxiliaryTwo}]
	Note that \eqref{productAuxiliaryTwo:NNsfmDef} proves that for all $n\in\N_0$, $l\in\zeroto{2^n-1}$ it holds that 
	\begin{equation}
%	\label{productAuxiliaryTwo:interpolationProperty}
		\begin{split}
		f_n(\tfrac{l}{2^n})=\brac{\tfrac{2l+1}{2^n}}\tfrac{l}{2^n}-\tfrac{(l^2+l)}{2^{2n}}
		=\tfrac{(2l+1)l-(l^2+l)}{2^{2n}}
		=\tfrac{l^2}{2^{2n}}=\left[\tfrac{l}{2^n}\right]^{\!2}.
		\end{split}
	\end{equation}
	The hypothesis that for all $n\in\N_0$ it holds that $f_n(1)=1$ hence ensures that for all $n\in\N_0$, $l\in\zeroto{2^n}$ it holds that
		\begin{equation}\label{productAuxiliaryTwo:interpolationProperty}
		\begin{split}
		f_n(\tfrac{l}{2^n})=\left[\tfrac{l}{2^n}\right]^{\!2}.
		\end{split}
		\end{equation}
	This and Lemma~\ref{Lemma:productAuxiliaryOne} demonstrate that for all $n\in\N$, $k\in\zeroto{2^{n-1}}$ it holds that
	\begin{equation}\label{productAuxiliaryTwo:NNsfDIFF1}
	\begin{split}
	f_{n-1}(\tfrac{2k}{2^n})-f_n(\tfrac{2k}{2^n})&=f_{n-1}(\tfrac{k}{2^{n-1}})-f_n(\tfrac{2k}{2^n})
	=\brac{\tfrac{k}{2^{n-1}}}^2-\brac{\tfrac{2k}{2^n}}^2 \\&=0=2^{-2n}g_n(\tfrac{2k}{2^n}).
	\end{split}\end{equation}
	In addition, note that \eqref{productAuxiliaryTwo:NNsfmDef} and \eqref{productAuxiliaryTwo:interpolationProperty} imply that 
	for all $n\in\N$, $k\in\zeroto{2^{n-1}-1}$ it holds that
	\begin{equation}\begin{split}\label{productAuxiliaryTwo:NNsT3}
	f_{n-1}(\tfrac{2k+1}{2^n})
	=f_{n-1}\big(\tfrac{k+(\nicefrac{1}{2})}{2^{n-1}}\big)
	&=\brac{\tfrac{2k+1}{2^{n-1}}}\brac{\tfrac{2k+1}{2^n}}-\tfrac{(k^2+k)}{2^{2(n-1)}}=\tfrac{(4k^2+4k+1)}{2^{2n-1}}-\tfrac{(2k^2+2k)}{2^{2n-1}}
	\\&=\tfrac{2k^2+2k+1}{2^{2n-1}} =\tfrac{4k^2+4k+2}{2^{2n}}
	\end{split}\end{equation}   
	and 
	\begin{equation}\begin{split}\label{productAuxiliaryTwo:NNsfDIFF2}
	 f_n(\tfrac{2k+1}{2^n})
	&=\brac{\tfrac{2k+1}{2^n}}^2
	=\tfrac{4k^2+4k+1}{2^{2n}}.
	\end{split}\end{equation}	
	Lemma~\ref{Lemma:productAuxiliaryOne} hence assures that for all $n\in\N$, $k\in\zeroto{2^{n-1}-1}$ it holds that
	\begin{equation}
	f_{n-1}(\tfrac{2k+1}{2^n})-f_n(\tfrac{2k+1}{2^n})=\tfrac{(4k^2+4k+2)}{2^{2n}}-\tfrac{(4k^2+4k+1)}{2^{2n}}={2^{-2n}}=2^{-2n}g_n(\tfrac{2k+1}{2^n}).
	\end{equation}
%	Furthermore, note that  proves that 
%	for all $n\in\N$, $k\in\zeroto{2^{n-1}-1}$, $x\in\brac{\tfrac{k}{2^{n-1}},{\tfrac{k+1}{2^{n-1}}}}$ it holds that
%	\begin{equation}\label{NNsClaim}
%	2^{-2n} g_n(x)=\begin{cases}
%	2^{-n}(x-\tfrac{2k}{2^n}) & \colon x\in\brac{\tfrac{2k}{2^n},{\tfrac{2k+1}{2^n}}}\\
%	2^{-n}(\tfrac{2k+2}{2^n}-x) & \colon x\in\brac{\tfrac{2k+1}{2^n},{\tfrac{2k+2}{2^n}}}
%	\end{cases}.
%	\end{equation}
Combining this with \eqref{productAuxiliaryTwo:NNsfDIFF1} shows that for all $n\in\N$, $l\in\zeroto{2^{n}}$ it holds that
	\begin{equation}\label{interpolationPropertyFnDifference}
	f_{n-1}(\tfrac{l}{2^n})-f_n(\tfrac{l}{2^n})=2^{-2n}g_n(\tfrac{l}{2^n}).
	\end{equation}
		Furthermore, observe that \eqref{productAuxiliaryTwo:interpolationProperty} demonstrates that for all $n\in\N_0$, $l\in\{0,1,\dots,\allowbreak 2^n-1\}$ it holds that 
		\begin{equation}
		\begin{split}
		\brac{\tfrac{2l+1}{2^n}} 	\brac{\tfrac{l+1}{2^n}}-\tfrac{(l^2+l)}{2^{2n}}
		= \tfrac{(2l+1)(l+1)-l(l+1)}{2^{2n}}=\tfrac{(l+1)^2}{2^{2n}}=\left[\tfrac{l+1}{2^n}\right]^{\!2}=f_n(\tfrac{l+1}{2^n}).
		\end{split}
		\end{equation}
		Combining this with \eqref{productAuxiliaryTwo:NNsfmDef} implies that for all	$n\in\N_0$ it holds that $f_n\in C([0,1],\R)$ and
		\begin{equation}\label{fnClosedInterval}
		\forallDist l\in\zeroto{2^n-1}, x\in\bracbig{\tfrac{l}{2^n},\tfrac{l+1}{2^n}}\colon\, f_n(x)=\brac{\tfrac{2l+1}{2^n}}x-\tfrac{(l^2+l)}{2^{2n}}.
		\end{equation}
	The fact that for all $n\in\N$, $k\in\{0,1,\dots,\allowbreak 2^{n-1}-1\}$ it holds that $\bracbig{\tfrac{k}{2^{n-1}},\tfrac{k+1}{2^{n-1}}}=\bracbig{\tfrac{2k}{2^n},\tfrac{2k+1}{2^n}}\cup \bracbig{\tfrac{2k+1}{2^n},\tfrac{2k+2}{2^n}}$ 
%	and the fact that  $f_n\colon [0,1]\to \R$, $n\in\N_0$, are continuous functions 
	hence ensures that there exist $(a_{n,k},b_{n,k},c_{n,k})\in\R^3$, $k\in\{0,1,\dots,\allowbreak 2^{n-1}-1\}$, $n\in\N$, such that for all $n\in\N$, $k\in\{0,1,\dots,\allowbreak 2^{n-1}-1\}$, $x\in \bracbig{\tfrac{k}{2^{n-1}},{\tfrac{k+1}{2^{n-1}}}}$	it holds that
				\begin{equation}
				f_{n-1}(x)-f_n(x)
				=\begin{cases}
				a_{n,k}\big(x-\tfrac{(2k+1)}{2^n}\big)+b_{n,k} & \colon x\in\brac{\tfrac{2k}{2^n},{\tfrac{2k+1}{2^n}}}\\
				c_{n,k}\big(x-\tfrac{(2k+1)}{2^n}\big)+b_{n,k} & \colon x\in\brac{\tfrac{2k+1}{2^n},{\tfrac{2k+2}{2^n}}}
				\end{cases}.
				\end{equation}
	 Lemma~\ref{Lemma:productAuxiliaryOne} and \eqref{interpolationPropertyFnDifference} therefore
	prove 
	that for all $n\in\N$, $k\in\zeroto{2^{n-1}-1}$, $x\in \bracbig{\tfrac{k}{2^{n-1}},{\tfrac{k+1}{2^{n-1}}}}$ it holds that
	\begin{equation}
	f_{n-1}(x)-f_n(x)=2^{-2n}g_n(x).
	\end{equation}
	Hence, we obtain 
	that for all $n\in\N$, $x\in [0,1]$ it holds that
	\begin{equation}\label{productAuxiliaryTwo:NNsClaim2}
	f_{n-1}(x)-f_n(x)=2^{-2n}g_n(x).
	\end{equation}
	Next note that \eqref{productAuxiliaryTwo:NNsfmDef} ensures that 
	 for all $x\in[0,1]$ it holds that $f_0(x)=x$. Combining this with \eqref{productAuxiliaryTwo:NNsClaim2} implies that 
	for all $m\in\N_0$, $x\in[0,1]$ it holds that
	\begin{equation}\label{productAuxiliaryTwo:NNsfmSum}
	\begin{split}
		f_m(x)&=
		f_0(x)+\Big[\smallsum\limits_{n=1}^m (f_n(x)-f_{n-1}(x))\Big]
		\\&= f_0(x)-\Big[\smallsum\limits_{n=1}^m (f_{n-1}(x)-f_n(x))\Big]
		=x-\Big[\smallsum\limits_{n=1}^m 2^{-2n} g_n(x)\Big].
	\end{split}
	\end{equation}
%Moreover, observe that  \eqref{productAuxiliaryTwo:interpolationProperty}, \eqref{fnClosedInterval}, and the fact that $([0,1]\ni x\mapsto x^2 \in [0,1])$ is a convex function ensure that for all $m\in\N_0$, $l\in\zeroto{2^m-1}$, $x\in\bracbig{\tfrac{l}{2^m},\tfrac{l+1}{2^m}}$ it holds that $f_m(x)-x^2\ge 0$. Combining this with \eqref{productAuxiliaryTwo:interpolationProperty} demonstrates that for all $m\in\N_0$, $l\in\zeroto{2^m-1}$, $x\in\bracbig{\tfrac{l}{2^m},\tfrac{l+1}{2^m}}$ it holds that
%\begin{equation}
%0\le f_m(x)-x^2=\left(x-\tfrac{l}{2^m}\right)\left(\tfrac{l+1}{2^m}-x\right)\!.
%\end{equation}
%The fact that for all $a,b,r\in\R$ with $a<b$ and $r\in [a,b]$ it holds that 
Moreover, observe that \eqref{fnClosedInterval} demonstrates that for all $m\in\N_0$, $l\in\zeroto{2^m-1}$, $x\in\bracbig{\tfrac{l}{2^m},\tfrac{l+1}{2^m}}$ it holds that
\begin{equation}
\begin{split}
f_m(x)-x^2
&=\brac{\tfrac{2l+1}{2^m}}x-\tfrac{(l^2+l)}{2^{2m}}-x^2
=\brac{\tfrac{l+1}{2^m}}x+\brac{\tfrac{l}{2^m}}x-\brac{\tfrac{l+1}{2^{m}}}\brac{\tfrac{l}{2^{m}}}-x^2
\\&=\left(x-\tfrac{l}{2^m}\right)\left(\tfrac{l+1}{2^m}-x\right)\ge 0.
\end{split}
\end{equation}
The fact that for all $a\in\R$, $b\in (a,\infty)$, $r\in [a,b]$ it holds that 
$(r-a)(b-r)\le \tfrac{1}{4}(b-a)^{2}$ hence proves that for all $m\in\N_0$, $l\in\zeroto{2^m-1}$, $x\in\brac{\tfrac{l}{2^m},\tfrac{l+1}{2^m}}$ it holds that 
\begin{equation}
\begin{split}
\abs{f_m(x)-x^2} &=f_m(x)-x^2=\left(x-\tfrac{l}{2^m}\right)\left(\tfrac{l+1}{2^m}-x\right)
\\&\le \tfrac{1}{4} \left(\tfrac{l+1}{2^m}-\tfrac{l}{2^m}\right)^2
=\tfrac{1}{4} \left(\tfrac{1}{2^m}\right)^2
=\tfrac{1}{2^2} \left(\tfrac{1}{2^{2m}}\right)
=\tfrac{1}{2^{2m+2}}
=2^{-2m-2}.
\end{split}
\end{equation}
Therefore, we obtain that for all $m\in\N_0$, $x\in[0,1]$ it holds that 
\begin{equation}
	\abs{f_m(x)-x^2} \le 2^{-2m-2}.
\end{equation}
%$\abs{f_m(x)-x^2} \le 2^{-2m-2}$.
 Combining this with \eqref{productAuxiliaryTwo:NNsfmSum} establishes \eqref{NNsfmSum}.
The proof of Lemma~\ref{Lemma:productAuxiliaryTwo} is thus completed.
\end{proof}

%\begin{prop}\label{NNsquare}
%Let $a\colon\R\to\R$ satisfy for all $x\in\R$ that $a(x)=\max\{0,x\}$.
%Then there exists $(\Phi_{\eps})_{\eps\in(0,1]}\subseteq\ANNs$ 
%such that for all $\eps\in(0,1]$, $x\in \R\backslash [0,1]$ it holds that
%\begin{equation}
%	 \paramANN(\Phi\leps)
%	\le 10\LogBin(\eps^{-1})+13,\qquad \L(\Phi\leps)\leq\max\{\tfrac{1}{2}\LogBin(\eps^{-1})+1,2\},
%\end{equation}
% \begin{equation}
% 	\sup_{y\in[0,1]}\abs{y^2-\brac{\functionANN(\Phi_{\eps})}\!(y)}\leq \eps,\qandq [\functionANN(\Phi_{\eps})](x) = a(x)
% \end{equation}     
% (cf.\ Definition~\ref{Def:ANN} and Definition~\ref{Definition:ANNrealization}).
%  \end{prop}

%% file: ANNApproxSquare.tex
\begin{prop}\label{Prop:NNsquare}
	Let $\varepsilon\in (0,1]$, $a\in C(\R,\R)$ satisfy for all $x\in\R$ that $a(x)=\max\{x,0\}$.
	Then there exists $\Phi\in \ANNs$
	such that
	\begin{enumerate}[(i)] 
		\item \label{NNsquare:RealizationCont} it holds  that $\functionANN(\Phi)\in C(\R,\R)$,
		\item \label{NNsquare:Realization} it holds for all  $x\in \R\backslash [0,1]$ that  $(\functionANN(\Phi))(x)\allowbreak = a(x)$,
		\item \label{NNsquare:Estimate} it holds for all  $x\in [0,1]$ that $\abs{x^2-(\functionANN(\Phi))(x)}\leq \eps$, 
		\item \label{NNsquare:Params} it holds  that 
%		$\paramANN(\Phi)\le 10\LogBin(\eps^{-1})+13$, 
		$\paramANN(\Phi)\le \max\{10\LogBin(\eps^{-1})-7,13\}$, 
		and 
		\item \label{NNsquare:Length} it holds that $\L(\Phi)\leq\max\{\tfrac{1}{2}\LogBin(\eps^{-1})+1,2\}$
	\end{enumerate}  
	(cf.\ Definition~\ref{Def:ANN} and Definition~\ref{Definition:ANNrealization}).
\end{prop}

\begin{proof}[Proof of Proposition~\ref{Prop:NNsquare}] 
	Throughout this proof 
	let $M\in\N$ satisfy that
	\begin{equation}\label{NNsquare:M_epsDef}
		M=\min\!\Big(\N\cap [2,\infty)\cap \big[\tfrac{1}{2}\LogBin(\varepsilon^{-1}),\infty\big)\Big),
	\end{equation}
	let $g_n\colon \R\to[0,1]$, $n\in\N$, be the functions which satisfy for all $n\in\N$, $x\in\R$ that
%	\begin{equation}\label{NNsquare:NNsgDef}
%	g_s(x)=\begin{cases}
%	2x & \colon n=1,x<\tfrac{1}{2}\\
%	2-2x & \colon n=1,x\geq\tfrac{1}{2}\\
%	g_1(g_{n-1}(t)) & \colon s\geq 1
%	\end{cases},
%	\end{equation}
		\begin{equation}\label{NNsquare:NNsgDef}
		g_1(x)=\begin{cases}
		2x & \colon x\in [0,\tfrac{1}{2})\\
		2-2x & \colon x\in [\tfrac{1}{2},1]\\
		0 & \colon x\in \R \backslash[0,1]\\
		\end{cases}
		\end{equation}
		and $g_{n+1}(x)=g_1(g_{n}(x))$,
	 let $f_n\colon [0,1]\to[0,1]$, $n\in\N_0$, 
	be the functions which satisfy for all 
	$n\in\N_0$, $k\in\zeroto{2^n-1}$, $x\in\big[\tfrac{k}{2^n},\tfrac{k+1}{2^n}\big)$ 
	that $f_n(1)=1$ and
	\begin{equation}\label{NNsquare:NNsfmDef}
	f_n(x)=\brac{\tfrac{2k+1}{2^n}}x-\tfrac{(k^2+k)}{2^{2n}},
	\end{equation}
	let $(A_k,b_k)\in\R^{4\times 4}\times\R^4$, $k\in\N\cap[2,\infty)$, 
	 satisfy for all $k\in\N\cap[2,\infty)$ that
	\begin{equation}\label{NNsquare:MatricesDef}
	A_k=\begin{pmatrix}
	2 & -4 & 2 & 0\\
	2 & -4 & 2 & 0\\
	2 & -4 & 2 & 0\\
	(-2)^{3-2k} & 2^{4-2k} & (-2)^{3-2k} & 1
	\end{pmatrix}
	\quad\mathrm{and}\quad
	b_k=\begin{pmatrix}
	0 \\ -\frac{1}{2} \\ -1 \\0
	\end{pmatrix},
	\end{equation}    
	let $\affineMap_k\in \R^{1\times 4}\times \R$, $k\in\N\cap [2,\infty)$, satisfy for all $k\in\N\cap [2,\infty)$ that 
	\begin{equation}\label{productLemma:DefinitionLastLayer}
		\affineMap_k=\pa{\begin{pmatrix}(-2)^{3-2k} & 2^{4-2k} & (-2)^{3-2k} & 1\end{pmatrix},0},
	\end{equation}
	let $\phi_k\in\ANNs$, $k\in\N\cap [2,\infty)$, satisfy  that 
		\begin{equation}\label{productLemma:DefinitionANNOne}
		\phi_2=\pa{\pa{\begin{pmatrix}1 \\ 1 \\ 1 \\ 1 \end{pmatrix},\begin{pmatrix}0 \\ -\frac{1}{2} \\ -1 \\0\end{pmatrix}},
			\affineMap_2}
		\end{equation}
	and
	\begin{equation}\label{productLemma:DefinitionANNTwo}
	\forallDist k\in\N\cap [3,\infty)\colon
	\phi_k=\pa{\pa{\begin{pmatrix}1 \\ 1 \\ 1 \\ 1 \end{pmatrix},\begin{pmatrix}0 \\ -\frac{1}{2} \\ -1 \\0\end{pmatrix}},(A_2,b_2),\dots,(A_{k-1},b_{k-1}),
		\affineMap_k}\!,
	\end{equation}
	and let $r_k=(r_{k,1},r_{k,2},r_{k,3},r_{k,4})\colon\R\to\R^4$, $k\in\N$, be the functions which satisfy for all $x\in\R$, $k\in\N$ that 
	\begin{equation}\label{NNsquare:StartingValueRfunctions}
	r_1(x)=(r_{1,1}(x),r_{1,2}(x),r_{1,3}(x),r_{1,4}(x))=\mathfrak{M}_{a, 4}\big(x,x-\tfrac{1}{2},x-1,x\big)
	\end{equation}  
	and 
	\begin{equation}\label{NNsquare:RecursiveValueRfunctions}
	r_{k+1}(x)=(r_{k+1,1}(x),r_{k+1,2}(x),r_{k+1,3}(x),r_{k+1,4}(x))=\mathfrak{M}_{a, 4}\big(A_{k+1} r_{k}(x)+b_{k+1}\big)
	\end{equation}
	(cf.\ Definition~\ref{Def:multidim_version}).
	 Note that \eqref{NNsquare:StartingValueRfunctions}, \eqref{multidim_version:Equation}, \eqref{NNsquare:NNsgDef}, and the hypothesis that for all $x\in\R$ it holds that $\relu(x)=\max\{x,0\}$  show that for all $x\in\R$ it holds that
	 \begin{equation}\label{NNsT01}
	 \begin{split}
	 &2r_{1,1}(x)-4r_{1,2}(x)+2r_{1,3}(x)
	 =2\relu(x)-4\relu(x-\tfrac{1}{2})+2\relu(x-1)
	 \\&=2\max\{x,0\}-4\max\{x-\tfrac{1}{2},0\}+2\max\{x-1,0\}
	 =g_1(x).
	 \end{split}
	 \end{equation}
	 Furthermore, observe that \eqref{NNsquare:StartingValueRfunctions}, \eqref{multidim_version:Equation}, the hypothesis that for all $x\in\R$ it holds that $\relu(x)=\max\{x,0\}$, and the fact that for all $x\in[0,1]$ it holds that $f_0(x)=x=\max\{x,0\}$ imply that for all $x\in\R$ it holds that 
	 \begin{equation}\label{NNsT01b}
	 r_{1,4}(x)=\max\{x,0\}
	 =\begin{cases}
	 f_{0}(x) & \colon x\in [0,1]\\
	 \max\{x,0\} & \colon x\in \R \backslash[0,1]\\
	 \end{cases}.		
	 \end{equation}
	Next we claim that for all $k\in\N$ it holds that
		\begin{equation}\label{NNsNIa}
	\big(\forallDist x\in\R\colon\,	2r_{k,1}(x)-4r_{k,2}(x)+2r_{k,3}(x)=g_k(x)\big)
		\end{equation}
%		\begin{equation}\label{NNsNIb}
%		\andq r_{k,4}(x)
%		=\begin{cases}
%		x-\Big[\smallsum\limits_{j=1}^{k-1} 2^{-2j}g_j(x)\Big] & \colon x\in [0,1]\\
%		\max\{x,0\} & \colon x\in \R \backslash[0,1]\\
%		\end{cases}.		
%		\end{equation}
and
		\begin{equation}\label{NNsNIb}
		\left(\forallDist x\in\R\colon\, r_{k,4}(x)
		=\begin{cases}
		f_{k-1}(x) & \colon x\in [0,1]\\
		\max\{x,0\} & \colon x\in \R \backslash[0,1]\\
		\end{cases}\right)\!.		
		\end{equation}
	We now prove \eqref{NNsNIa}--\eqref{NNsNIb} by induction on $k\in\N$. 	
	Note that \eqref{NNsT01} and \eqref{NNsT01b} prove \eqref{NNsNIa}--\eqref{NNsNIb} in the base case $k=1$.
	For the induction step $\N\ni k\to k+1\in\N\cap [2,\infty)$ assume that 
	there exists $k\in\N$ such that for all $x\in\R$ it holds that
	\begin{equation}\label{NNsNIa:InductionHypothesis}
	2r_{k,1}(x)-4r_{k,2}(x)+2r_{k,3}(x)=g_k(x)
	\end{equation}
	\begin{equation}\label{NNsNIb:InductionHypothesis}
		\andq r_{k,4}(x)
		=\begin{cases}
		f_{k-1}(x) & \colon x\in [0,1]\\
		\max\{x,0\} & \colon x\in \R \backslash[0,1]\\
		\end{cases}.	
	\end{equation}
	Observe that \eqref{NNsNIa:InductionHypothesis}, \eqref{NNsT01}, \eqref{NNsquare:MatricesDef}, \eqref{multidim_version:Equation}, and \eqref{NNsquare:RecursiveValueRfunctions} ensure that for all $x\in\R$ it holds that
	\begin{equation}\label{NNsNIa:InductionCalculation}
		\begin{split}
		g_{k+1}(x)&=g_1(g_{k}(x))
		=g_1(2r_{k,1}(x)-4r_{k,2}(x)+2r_{k,3}(x))
		\\&=2\,\relu\big(2r_{k,1}(x)-4r_{k,2}(x)+2r_{k,3}(x)\big)
		\\&\quad-4\,\relu\big(2r_{k,1}(x)-4r_{k,2}(x)+2r_{k,3}(x)-\tfrac{1}{2}\big)
		\\&\quad
		+2\,\relu\big(2r_{k,1}(x)-4r_{k,2}(x)+2r_{k,3}(x)-1\big)
		\\&=2r_{k+1,1}(x)-4r_{k+1,2}(x)+2r_{k+1,3}(x).
				\end{split}
	\end{equation}
%	In addition, observe that \eqref{NNsquare:MatricesDef}, \eqref{NNsquare:RecursiveValueRfunctions}, \eqref{NNsNIa:InductionHypothesis}, \eqref{NNsNIb:InductionHypothesis}, and the fact that for all $x\in [0,1]$ it holds that $f_k(x)\ge 0$ demonstrate that for all $x\in [0,1]$ it holds that 
%	\begin{equation}
%		\begin{split}
%		&r_{k+1,4}(x)
%		=\relu\big(-2^{-2k+1}r_{k,1}(x)+2^{-2k+2}r_{k,2}(x)-2^{-2k+1}r_{k,3}(x)+r_{k,4}(x)\big)
%		\\&=\relu\Big(-2^{-2k}\big[2r_{k,1}(x)-4r_{k,2}(x)+2r_{k,3}(x)\big]+x-\Big[\smallsum\limits_{j=1}^{k-1} 2^{-2j}g_j(x)\Big]\Big)
%		\\&=\relu\Big(-2^{-2k}g_k(x)+x-\Big[\smallsum\limits_{j=1}^{k-1} 2^{-2j}g_j(x)\Big]\Big)
%		\\&=\relu\Big(y-\Big[\smallsum\limits_{j=1}^{k} 2^{-2j}g_j(x)\Big]\Big)
%		=\relu(f_k(x))=f_k(x).
%		\end{split}
%	\end{equation}
	In addition, observe that \eqref{multidim_version:Equation}, \eqref{NNsquare:MatricesDef}, \eqref{NNsquare:RecursiveValueRfunctions}, and \eqref{NNsNIa:InductionHypothesis}  demonstrate that for all $x\in \R$ it holds that 
	\begin{equation}\label{NNsquare:CalculationFourthComponentOne}
	\begin{split}
	&r_{k+1,4}(x)
	\\&=\relu\big((-2)^{3-2(k+1)}r_{k,1}(x)+2^{4-2(k+1)}r_{k,2}(x)+(-2)^{3-2(k+1)}r_{k,3}(x)+r_{k,4}(x)\big)
	\\&=\relu\big((-2)^{1-2k}r_{k,1}(x)+2^{2-2k}r_{k,2}(x)+(-2)^{1-2k}r_{k,3}(x)+r_{k,4}(x)\big)
		\\&=\relu\big(2^{-2k}\big[-2r_{k,1}(x)+2^2r_{k,2}(x)-2r_{k,3}(x)\big]+r_{k,4}(x)\big)
	\\&=\relu\big(-\big[2^{-2k}\big]\big[2r_{k,1}(x)-4r_{k,2}(x)+2r_{k,3}(x)\big]+r_{k,4}(x)\big)
	\\&=\relu\big(-\big[2^{-2k}\big]g_k(x)+r_{k,4}(x)\big).
%	\\&=\relu\Big(y-\Big[\smallsum\limits_{j=1}^{k} 2^{-2j}g_j(x)\Big]\Big)
%	=\relu(f_k(x))=f_k(x).
	\end{split}
	\end{equation}
	Combining this with \eqref{NNsNIb:InductionHypothesis}, Lemma~\ref{Lemma:productAuxiliaryTwo}, the hypothesis that for all $x\in\R$ it holds that $\relu(x)=\max\{x,0\}$, and the fact that for all $x\in [0,1]$ it holds that $f_k(x)\ge 0$ shows that for all $x\in[0,1]$ it holds that
		\begin{equation}\label{NNsquare:CalculationFourthComponentTwo}
		\begin{split}
		r_{k+1,4}(x)
		&=\relu\big(-\big[2^{-2k}g_k(x)\big]+f_{k-1}(x)\big)
\\&=\relu\Big(-\big(2^{-2k}g_k(x)\big)+x-\Big[\smallsum\limits_{j=1}^{k-1} \big(2^{-2j}g_j(x)\big)\Big]\Big)
		\\&=\relu\Big(x-\Big[\smallsum\limits_{j=1}^{k} 2^{-2j}g_j(x)\Big]\Big)
		=\relu(f_k(x))=f_k(x).
		\end{split}
		\end{equation}
	Next note that \eqref{NNsquare:CalculationFourthComponentOne}, \eqref{NNsNIb:InductionHypothesis}, 
	item~\eqref{productAuxiliaryOne:exterior} in Lemma~\ref{Lemma:productAuxiliaryOne}, and the hypothesis that for all $x\in\R$ it holds that $\relu(x)=\max\{x,0\}$ prove that for all $x\in\R\backslash [0,1]$ it holds that 
	\begin{equation}\label{NNsquare:CalculationFourthComponentThree}
		r_{k+1,4}(x)
		=\relu\Big(-\big(2^{-2k}g_k(x)\big)+r_{k,4}(x)\Big)
		=\relu(\max\{x,0\})=\max\{x,0\}.
	\end{equation}
	Combining  \eqref{NNsNIa:InductionCalculation} and \eqref{NNsquare:CalculationFourthComponentTwo} hence proves \eqref{NNsNIa}--\eqref{NNsNIb} in the case $k+1$. Induction thus establishes \eqref{NNsNIa}--\eqref{NNsNIb}.
%	This, \eqref{NNsquare:CalculationFourthComponentTwo}, \eqref{NNsNIa:InductionCalculation}, and induction prove \eqref{NNsNIa} and \eqref{NNsNIb}.	
%	Next note that \eqref{NNsNIa}, \eqref{NNsNIb}, and Lemma~\ref{Lemma:productAuxiliaryTwo}
%	assure that for all $m\in\N$, $x\in\R$ it holds that
%	\begin{equation}\begin{split}
%	&[\functionANN(\phi_m)](y)
%	=-2^{-2m+3}r_{m-1,1}(y)+2^{-2m+4}r_{m-1,2}(y)-2^{-2m+3}r_{m-1.3}(y)+r_{m-1,4}(y)\\
%	&=-2^{-2(m-1)}\big[2r_{m-1,1}(y)-4r_{m-1,2}(y)+2r_{m-1,3}(y)\big]+y-\Big[\smallsum\limits_{j=1}^{m-2} 2^{-2j}g_j(y)\Big]\\
%	&=-2^{-2(m-1)}g_{m-1}(y)+y-\Big[\smallsum\limits_{j=1}^{m-2} 2^{-2j}g_j(y)\Big]
%\\&=y-\Big[\smallsum\limits_{j=1}^{m-1} 2^{-2j}g_j(y)\Big]
%=f_{m-1}(y).
%	\end{split}\end{equation} 
Next note that  \eqref{setting_NN:ass2}, \eqref{NNsquare:MatricesDef}, \eqref{productLemma:DefinitionLastLayer}, \eqref{NNsNIa}, \eqref{productLemma:DefinitionANNOne}, \eqref{productLemma:DefinitionANNTwo}, \eqref{NNsquare:StartingValueRfunctions}, and \eqref{NNsquare:RecursiveValueRfunctions}
assure that for all $m\in\N\cap [2,\infty)$, $x\in\R$ it holds that $\functionANN(\phi_m)\in C(\R,\R)$ and
\begin{equation}\label{NNsquare:CalculationANNOne}
\begin{split}
&(\functionANN(\phi_m))(x)
\\&=(-2)^{3-2m}r_{m-1,1}(x)+2^{4-2m}r_{m-1,2}(x)+(-2)^{3-2m}r_{m-1,3}(x)+r_{m-1,4}(x)\\
&=(-2)^{4-2m} \Big(\Big[\tfrac{r_{m-1,1}(x)+r_{m-1,3}(x)}{(-2)}\Big]+r_{m-1,2}(x)\Big)+r_{m-1,4}(x)\\
&=2^{4-2m} \Big(\Big[\tfrac{r_{m-1,1}(x)+r_{m-1,3}(x)}{(-2)}\Big]+r_{m-1,2}(x)\Big)+r_{m-1,4}(x)\\
&=2^{2-2m} \big(4r_{m-1,2}(x)-2r_{m-1,1}(x)-2r_{m-1,3}(x)\big)+r_{m-1,4}(x)\\
&=-\big[2^{-2(m-1)}\big]\big[2r_{m-1,1}(x)-4r_{m-1,2}(x)+2r_{m-1,3}(x)\big]+r_{m-1,4}(x)\\
&=-\big[2^{-2(m-1)}\big]g_{m-1}(x)+r_{m-1,4}(x).
%\\&=y-\Big[\smallsum\limits_{j=1}^{m-1} 2^{-2j}g_j(x)\Big]
%=f_{m-1}(x).
\end{split}\end{equation} 
	Combining this with \eqref{NNsNIb} and Lemma~\ref{Lemma:productAuxiliaryTwo} shows that for all $m\in\N\cap [2,\infty)$, $x\in[0,1]$ it holds that
	\begin{equation}\label{NNsquare:CalculationANNTwo}
	\begin{split}
	(\functionANN(\phi_m))(x)
	&=-\big(2^{-2(m-1)}g_{m-1}(x)\big)+f_{m-2}(x)
	\\&=-\big(2^{-2(m-1)}g_{m-1}(x)\big)+x-\Big[\smallsum\limits_{j=1}^{m-2} 2^{-2j}g_j(x)\Big]
	\\&=x-\Big[\smallsum\limits_{j=1}^{m-1} 2^{-2j}g_j(x)\Big]
=f_{m-1}(x).
	\end{split}
	\end{equation}
	Lemma~\ref{Lemma:productAuxiliaryTwo} therefore implies that for all $m\in\N\cap [2,\infty)$, $x\in [0,1]$ it holds that 
	\begin{equation}\label{NNsphimEst}
	\abs{x^2-(\functionANN(\phi_m))(x)}\leq 2^{-2(m-1)-2}= 2^{-2m}.
	\end{equation}
	Next note that \eqref{NNsquare:M_epsDef} assures that 
	\begin{equation}
		\begin{split}
		M&=\min\!\Big(\N\cap  \big[\max\!\big\{2,\tfrac{1}{2}\LogBin(\varepsilon^{-1})\big\},\infty\big)\Big)
		\\&\ge \min\!\Big(\big[\max\!\big\{2,\tfrac{1}{2}\LogBin(\varepsilon^{-1})\big\},\infty\big)\Big)
		\\&=\max\!\big\{2,\tfrac{1}{2}\LogBin(\varepsilon^{-1})\big\}\ge \tfrac{1}{2}\LogBin(\varepsilon^{-1}).
		\end{split}
	\end{equation}
	This and \eqref{NNsphimEst} demonstrate that for all $x\in [0,1]$ it holds that 
	\begin{equation}\label{NNsphimEst:Varepsilon}
	\abs{x^2-(\functionANN(\phi_M))(x)}\leq 2^{-2M}\le 2^{- \LogBin(\varepsilon^{-1})}= \varepsilon.
	\end{equation}	
	Moreover, observe that item~\eqref{productAuxiliaryOne:exterior} in Lemma~\ref{Lemma:productAuxiliaryOne}, \eqref{NNsNIb}, and \eqref{NNsquare:CalculationANNOne} ensure that for all $m\in\N\cap [2,\infty)$, $x\in\R\backslash[0,1]$ it holds that
	\begin{equation}\label{NNsquare:ANNoutside}
	\begin{split}
			(\functionANN(\phi_m))(x)&=-2^{-2(m-1)}g_{m-1}(x)+r_{m-1,4}(x)
			\\&=r_{m-1,4}(x)
			=\max\{x,0\}=a(x).
	\end{split}
	\end{equation}
	Furthermore, observe that \eqref{NNsquare:M_epsDef}, \eqref{productLemma:DefinitionANNOne}, and \eqref{productLemma:DefinitionANNTwo} assure that
	\begin{equation}\label{NNsLMphiLength}
		\L(\phi_M)=M\le \max\{\tfrac{1}{2} \LogBin(\varepsilon^{-1})+1,2\}.
	\end{equation}
This, \eqref{NNsquare:M_epsDef}, \eqref{productLemma:DefinitionANNOne}, and \eqref{productLemma:DefinitionANNTwo} show that
	\begin{equation}\label{NNsLMphi}
	\begin{split}
		\paramANN(\phi_M)
		&=4(1+1)+\left[\smallsum\limits_{j=2}^{M-1}4(4+1)\right]+(4+1)
		\\&=8+20(M-2)+5\le 20\max\{\tfrac{1}{2} \LogBin(\varepsilon^{-1})-1,0\}+13
		\\&= \max\{10 \LogBin(\varepsilon^{-1})-20,0\}+13
		=\max\{10 \LogBin(\varepsilon^{-1})-7,13\}.
	\end{split}
	\end{equation}
	Combining  \eqref{NNsphimEst:Varepsilon}, \eqref{NNsLMphiLength}, \eqref{NNsquare:ANNoutside}, and the fact that $\functionANN(\phi_M)\in C(\R,\R)$ hence establishes items~\eqref{NNsquare:RealizationCont}--\eqref{NNsquare:Length}.
	The proof of Proposition~\ref{Prop:NNsquare} is thus completed.    
\end{proof}

%% file: ANNApproxSquareOnR.tex
%\subsection{DNN approximation of the quadratic function on the real line}

%\begin{prop}\label{Lemma:ApproxOfSquare}
%Let $a\colon\R\to\R$ satisfy for all $x\in\R$ that $a(x)=\max\{0,x\}$.
%	Then there exists $(\Phi_\varepsilon)_{\varepsilon\in (0,1]}\subseteq \ANNs$ such that for all $x\in \R$, $\varepsilon\in (0,1]$ it holds that $\paramANN(\Phi_\varepsilon)\le 72 \!\left(\tfrac{1}{2} \LogBin(\tfrac{1}{\varepsilon})+1\right)$, $\paramNotZeroANN(\Phi_\varepsilon)\le 30 \!\left(\tfrac{1}{2} \LogBin(\tfrac{1}{\varepsilon})+1\right)$, $\lengthANN(\Phi_\varepsilon)\le \tfrac{1}{2} \LogBin(\tfrac{1}{\varepsilon})+1$, and
%	\begin{equation}
%	0\le (\functionANN (\Phi_\varepsilon))(x)\le\vert x\vert,\qquad \sup_{y\in[-1,1]}\vert y^2-(\functionANN (\Phi_\varepsilon))(y)\vert\le \varepsilon
%	\end{equation}
%	(cf.\ Definition~\ref{Def:ANN} and Definition~\ref{Definition:ANNrealization}).
%%	\begin{equation}
%%	\paramANN(\Phi_\varepsilon)\le 72 \!\left(\tfrac{1}{2} \LogBin(\tfrac{1}{\varepsilon})+1\right),
%%	\end{equation}
%%	\begin{equation}
%%	\paramNotZeroANN(\Phi_\varepsilon)\le , 
%%	\end{equation}
%%	$\lengthANN(\Phi_\varepsilon)\le $
%%	\begin{equation}
%%	\lengthANN(\Phi_\varepsilon)\le .
%%	\end{equation}
%\end{prop}

\begin{prop}\label{Lemma:ApproxOfSquare}
	Let $\varepsilon\in (0,1]$, $q\in (2,\infty)$, $a\in C(\R,\R)$ satisfy for all $x\in \R$ that $\activation(x)=\max\{x,0\}$.
	Then there exists $\Phi\in \ANNs$  
	such that
	\begin{enumerate}[(i)]
		\item \label{ApproxOfSquare:Realization}
		 it holds  that  $\functionANN (\Phi)\in C(\R,\R)$,
		 \item \label{ApproxOfSquare:Zero}
		 it holds  that $(\functionANN (\Phi))(0)\allowbreak=0$,
		\item \label{ApproxOfSquare:Estimates}
		it holds for all  $x\in\R$ that $0\le (\functionANN (\Phi))(x)\le \varepsilon+ \vert x\vert^2$,
		\item \label{ApproxOfSquare:EstimatesDifference}
		it holds for all  $x\in\R$ that 
				$\vert x^2-(\functionANN (\Phi))(x)\vert\le \varepsilon \max\{1,\vert x\vert^q\}$,
%		\begin{equation}
%		\vert x^2-(\functionANN (\Phi))(x)\vert\le \varepsilon \max\{1,\vert x\vert^q\},
%		\end{equation}
		\item \label{ApproxOfSquare:Params}
		it holds  that $\paramANN(\Phi)\le  \max\!\big\{\big[\tfrac{40q}{(q-2)}\big]\LogBin({\varepsilon}^{-1})+\tfrac{80}{(q-2)}-28,52\big\}$, and 
		\item \label{ApproxOfSquare:Length}
		it holds  that $\lengthANN(\Phi)\le \max\! \big\{\tfrac{q}{2(q-2)}\LogBin({\varepsilon}^{-1})+\tfrac{1}{(q-2)}+1,2\big\}$
%		\item \label{ApproxOfSquare:Params}
%		it holds  that 
%		\begin{equation}
%		\paramANN(\Phi)\le  \tfrac{40q}{(q-2)}\LogBin({\varepsilon}^{-1})+\tfrac{80}{(q-2)}+52
%		\end{equation}	
%		\begin{equation}
%		\andq \lengthANN(\Phi)\le \max \big\{\tfrac{q}{2(q-2)}\LogBin({\varepsilon}^{-1})+\tfrac{1}{(q-2)}+1,2\big\}
%		\end{equation}
	\end{enumerate}
 (cf.\ Definition~\ref{Def:ANN} and Definition~\ref{Definition:ANNrealization}).
\end{prop}

\begin{proof}[Proof of Proposition~\ref{Lemma:ApproxOfSquare}]
	Throughout this proof let  $\delta\in (0,1]$ satisfy that $\delta=2^{-\nicefrac{2}{(q-2)}} \varepsilon^{\nicefrac{q}{(q-2)}}$, let $\affineMap_1\in (\R^{2\times 1}\times \R^2)\subseteq \ANNs$, $\affineMap_2\in(\R^{1\times 2}\times \R)\subseteq \ANNs$ satisfy that
	\begin{equation}
	\affineMap_1=\left(\begin{pmatrix}
	(\tfrac{\varepsilon}{2})^{\nicefrac{1}{(q-2)}} \\[1ex]-(\tfrac{\varepsilon}{2})^{\nicefrac{1}{(q-2)}} 
	\end{pmatrix},\begin{pmatrix}
	0 \\0 
	\end{pmatrix}\right)\qandqShort
		\affineMap_2=\left(\begin{pmatrix}
		(\tfrac{\varepsilon}{2})^{-\nicefrac{2}{(q-2)}} &(\tfrac{\varepsilon}{2})^{-\nicefrac{2}{(q-2)}} 
		\end{pmatrix},0\right)\!,
	\end{equation}	
	let
	$\Psi\in \ANNs$ 
	satisfy that 
		\begin{enumerate}[(I)] 
			\item\label{ApproxOfSquare:BasisItemOne} it holds that  $\functionANN(\Psi)\in C(\R,\R)$,
			\item\label{ApproxOfSquare:BasisItemTwo} it holds for all  $x\in \R\backslash [0,1]$ that $(\functionANN(\Psi))(x) = a(x)$,
			\item \label{ApproxOfSquare:BasisItemThree}it holds for all $x\in [0,1]$ that 
%			$\sup_{y\in[0,1]}\abs{y^2-\brac{\functionANN(\Psi)}\!(y)}\leq \delta$, and 
$\abs{x^2-(\functionANN(\Psi))(x)}\leq \delta$, 
			\item \label{ApproxOfSquare:BasisItemFour}it holds that $\paramANN(\Psi)
			\le \max\{10\LogBin(\delta^{-1})-7,13\}$, and 
			\item \label{ApproxOfSquare:BasisItemFive} it holds that $\L(\Psi)\leq\max\{\tfrac{1}{2}\LogBin(\delta^{-1})+1,2\}$
	\end{enumerate}
%	
%	
%	
%	
%	
%	
%	
%	
%	for all $x\in \R\backslash [0,1]$  that
%	\begin{equation}\label{ApproxOfSquare:BasicParams}
%	\paramANN(\Psi)
%	\le 10\LogBin(\delta^{-1})+13,\qquad \L(\Psi)\leq\max\{\tfrac{1}{2}\LogBin(\delta^{-1})+1,2\},
%	\end{equation}
%	\begin{equation}\label{ApproxOfSquare:BasicEstimate}
%	\sup_{y\in[0,1]}\abs{y^2-\brac{\functionANN(\Psi)}\!(y)}\leq \delta,\qandq [\functionANN(\Psi)](x) = a(x)
%	\end{equation}     
	(cf.\ Proposition~\ref{Prop:NNsquare}), and let $\Phi\in \ANNs$ satisfy  that  
	\begin{equation}\label{ApproxOfSquare:ConstructionNetwork}
		\Phi=\compANN{\affineMap_2}{\compANN{\big[\parallelizationSpecial_2(\Psi,\Psi)\big]}{\affineMap_1}}
	\end{equation}
	(cf.\ Definition~\ref{Definition:ANNcomposition}, Definition~\ref{Definition:simpleParallelization}, and Lemma~\ref{Lemma:CompositionAssociative}).
	Note that  Proposition~\ref{Lemma:PropertiesOfParallelizationEqualLength} and item~\eqref{PropertiesOfCompositions:Realization} in  Proposition~\ref{Lemma:PropertiesOfCompositions} ensure that for all $x\in \R$ it holds that
	\begin{equation}
	\begin{split}
		&\big(\functionANN\big(\compANN{\big(\parallelizationSpecial_2(\Psi,\Psi)\big)}{\affineMap_1}\big)\big)(x)
%	=\big(\big[\functionANN\big(\parallelizationSpecial_2(\Psi,\Psi)\big)\big]\circ\big[\functionANN(\affineMap_1)\big]\big)(x)
%	\\&
	=\big(\functionANN\big(\parallelizationSpecial_2(\Psi,\Psi)\big)\big)\big(\big(\functionANN(\affineMap_1)\big)(x)\big)
		\\&=\big(\functionANN\big(\parallelizationSpecial_2(\Psi,\Psi)\big)\big)\big((\tfrac{\varepsilon}{2})^{\nicefrac{1}{(q-2)}} x,-(\tfrac{\varepsilon}{2})^{\nicefrac{1}{(q-2)}}x\big)
%	\\&=\big([\functionANN(\Psi)](x),[\functionANN(\Psi)](-x)\big).
\\&=\begin{pmatrix}
(\functionANN(\Psi))\big((\tfrac{\varepsilon}{2})^{\nicefrac{1}{(q-2)}}x\big)\\[1ex] (\functionANN(\Psi))\big(-(\tfrac{\varepsilon}{2})^{\nicefrac{1}{(q-2)}}x\big)
\end{pmatrix}.
	\end{split}
	\end{equation}
Item~\eqref{PropertiesOfCompositions:Realization} in	Proposition~\ref{Lemma:PropertiesOfCompositions} and \eqref{ApproxOfSquare:ConstructionNetwork} therefore demonstrate that for all $x\in \R$ it holds that
	\begin{equation}\label{ApproxOfSquare:ComputationRealization}
	\begin{split}
	(\functionANN (\Phi))(x)
%	= \functionANN\big(\compANN{\affineMap_2}{\big(\compANN{\big[\parallelizationSpecial_2(\Psi,\Psi)\big]}{\affineMap_1}\big)}\big)
		&= (\functionANN({\affineMap_2}))\big(\functionANN\big(\compANN{\big[\parallelizationSpecial_2(\Psi,\Psi)\big]}{\affineMap_1}\big)(x)\big)
	\\&=\begin{pmatrix}
	(\tfrac{\varepsilon}{2})^{-\nicefrac{2}{(q-2)}} &(\tfrac{\varepsilon}{2})^{-\nicefrac{2}{(q-2)}} 
	\end{pmatrix} \begin{pmatrix}
	[\functionANN(\Psi)]\big((\tfrac{\varepsilon}{2})^{\nicefrac{1}{(q-2)}}x\big)\\[1ex] [\functionANN(\Psi)]\big(-(\tfrac{\varepsilon}{2})^{\nicefrac{1}{(q-2)}}x\big)
	\end{pmatrix}
	\\&= (\tfrac{\varepsilon}{2})^{-\nicefrac{2}{(q-2)}} \Big(	[\functionANN(\Psi)]\big((\tfrac{\varepsilon}{2})^{\nicefrac{1}{(q-2)}}x\big)+ [\functionANN(\Psi)]\big(-(\tfrac{\varepsilon}{2})^{\nicefrac{1}{(q-2)}}x\big)\Big).
	\end{split}
	\end{equation}
%	\begin{equation}
%		\functionANN (\Phi)(x)= \functionANN (\Psi)(x)+\functionANN (\Psi)(-x).
%	\end{equation}
This, \eqref{ApproxOfSquare:BasisItemOne}, \eqref{ApproxOfSquare:BasisItemTwo}, and the hypothesis that for all $x\in\R$ it holds that $a(x)=\max\{x,0\}$ imply that 
	\begin{equation}\label{eq:ApproxOfSquare:Zero}
	\begin{split}
	(\functionANN (\Phi))(0)
	&= (\tfrac{\varepsilon}{2})^{-\nicefrac{2}{(q-2)}} \big(	[\functionANN(\Psi)](0)+ [\functionANN(\Psi)](0)\big)
		\\&= (\tfrac{\varepsilon}{2})^{-\nicefrac{2}{(q-2)}}  \big(a(0)+a(0)\big)=0.
	\end{split}
	\end{equation}
	Moreover, observe that \eqref{ApproxOfSquare:BasisItemOne} and \eqref{ApproxOfSquare:BasisItemTwo} ensure that for all $x\in \R\backslash [-1,1]$ it holds that 
		\begin{equation}\label{ApproxOfSquare:EstimateOutsideUnitCube}
		\begin{split}
%[\functionANN(\Psi)](x)+ [\functionANN(\Psi)](-x)
%		&=[\activation(x)+0]\indicator{[0,\infty)}(x)+[0+\activation(-x)]\indicator{(-\infty,0)}(x)
%		\\&=x\indicator{[0,\infty)}(x)-x\indicator{(-\infty,0)}(x)
%		=\vert x\vert
[\functionANN(\Psi)](x)+ [\functionANN(\Psi)](-x)
&=a(x)+a(-x)=\max\{x,0\}+\max\{-x,0\}
\\&=\max\{x,0\}-\min\{x,0\}=\vert x\vert.
		\end{split}
		\end{equation}
	Furthermore, note that \eqref{ApproxOfSquare:BasisItemTwo} and \eqref{ApproxOfSquare:BasisItemThree} show that 
	\begin{equation}\label{ApproxOfSquare:EstimateUnitCube}
	\begin{split}
	&\sup_{x\in [-1,1]}\left\vert x^2-\big([\functionANN(\Psi)](x)+ [\functionANN(\Psi)](-x)\big)\right\vert
	\\&=\max\!\Big\{\sup_{x\in [-1,0]}\left\vert x^2-\big(\activation(x)+[\functionANN (\Psi)](-x)\big)\right\vert,\sup_{x\in [0,1]}\left\vert x^2-\big([\functionANN (\Psi)](x)+\activation(-x)\big)\right\vert\Big\}
	\\&=\max\!\Big\{\sup_{x\in [-1,0]}\left\vert (-x)^2-(\functionANN (\Psi))(-x)\right\vert,\sup_{x\in [0,1]}\left\vert x^2-(\functionANN (\Psi))(x)\right\vert\Big\}	
	\\&=\sup_{x\in [0,1]}\left\vert x^2-(\functionANN (\Psi))(x)\right\vert\le \delta.
	\end{split}
	\end{equation}	
	Next observe that \eqref{ApproxOfSquare:ComputationRealization} and \eqref{ApproxOfSquare:EstimateOutsideUnitCube} prove that for all $x\in \R\backslash [-({\varepsilon}/{2})^{-\nicefrac{1}{(q-2)}},({\varepsilon}/{2})^{-\nicefrac{1}{(q-2)}}]$ it holds that
	\begin{equation}\label{ApproxOfSquare:BasicCalculation}
	\begin{split}
			0&\le [\functionANN (\Phi)](x)
			\\&=(\tfrac{\varepsilon}{2})^{-\nicefrac{2}{(q-2)}} \Big(	[\functionANN(\Psi)]\big((\tfrac{\varepsilon}{2})^{\nicefrac{1}{(q-2)}}x\big)+ [\functionANN(\Psi)]\big(-(\tfrac{\varepsilon}{2})^{\nicefrac{1}{(q-2)}}x\big)\Big)
			\\&=(\tfrac{\varepsilon}{2})^{-\nicefrac{2}{(q-2)}} \left\vert (\tfrac{\varepsilon}{2})^{\nicefrac{1}{(q-2)}}x\right\vert
			=(\tfrac{\varepsilon}{2})^{\nicefrac{-1}{(q-2)}}\vert x\vert\le \vert x\vert^2.
	\end{split}
	\end{equation}
The triangle inequality therefore ensures that for all $x\in \R\backslash[-({\varepsilon}/{2})^{-\nicefrac{1}{(q-2)}},({\varepsilon}/{2})^{-\nicefrac{1}{(q-2)}}]$ it holds that
	\begin{equation}\label{ApproxOfSquare:LargeValues}
			\begin{split}
			\left\vert x^2-(\functionANN(\Phi))(x)\right\vert
			&=\left\vert x^2-(\tfrac{\varepsilon}{2})^{\nicefrac{-1}{(q-2)}}\vert x\vert\right\vert
%			\\&=(\tfrac{\varepsilon}{2})^{-\nicefrac{2}{(q-2)}} \left\vert \vert(\tfrac{\varepsilon}{2})^{\nicefrac{1}{(q-2)}}x\vert^2- \vert (\tfrac{\varepsilon}{2})^{\nicefrac{1}{(q-2)}}x\vert \right\vert
%			\\&\le (\tfrac{\varepsilon}{2})^{-\nicefrac{2}{(q-2)}} \left(\vert(\tfrac{\varepsilon}{2})^{\nicefrac{1}{(q-2)}}x\vert^2+ \vert (\tfrac{\varepsilon}{2})^{\nicefrac{1}{(q-2)}}x\vert\right)  
			\le  \left(\vert x\vert^2+ (\tfrac{\varepsilon}{2})^{-\nicefrac{1}{(q-2)}}\vert x\vert\right)  
			\\&=  \left( \vert x\vert^q\vert x\vert^{-(q-2)}+ (\tfrac{\varepsilon}{2})^{-\nicefrac{1}{(q-2)}}\vert x\vert^q\vert x\vert^{-(q-1)}\right)
			\\&\le  \left( \vert x\vert^q (\tfrac{\varepsilon}{2})^{\nicefrac{(q-2)}{(q-2)}}+ (\tfrac{\varepsilon}{2})^{-\nicefrac{1}{(q-2)}}\vert x\vert^q(\tfrac{\varepsilon}{2})^{\nicefrac{(q-1)}{(q-2)}}\right)  
			\\&=(\tfrac{\varepsilon}{2}+\tfrac{\varepsilon}{2}) \vert x\vert^q=\varepsilon \vert x\vert^q \le \varepsilon \max\!\big\{1,\vert x\vert^q\big\}.
			\end{split}
	\end{equation}	
%		Moreover, observe that \eqref{ApproxOfSquare:BasicEstimate}, \eqref{ApproxOfSquare:BasicCalculation}, and the triangle inequality ensure that for all $\varepsilon\in (0,1]$, $x\in \R\backslash[-(\tfrac{\varepsilon}{2})^{-\nicefrac{1}{(q-2)}},(\tfrac{\varepsilon}{2})^{-\nicefrac{1}{(q-2)}}]$ it holds that
%		\begin{equation}\label{ApproxOfSquare:LargeValues}
%		\begin{split}
%		&\left\vert x^2-(\functionANN (\Phi))(x)\right\vert
%		\\&=(\tfrac{\varepsilon}{2})^{-\nicefrac{2}{(q-2)}} \left\vert ((\tfrac{\varepsilon}{2})^{\nicefrac{1}{(q-2)}}x)^2- (\functionANN (\Psi))((\tfrac{\varepsilon}{2})^{\nicefrac{1}{(q-2)}}x)\right\vert
%		\\&=(\tfrac{\varepsilon}{2})^{-\nicefrac{2}{(q-2)}} \left\vert \vert(\tfrac{\varepsilon}{2})^{\nicefrac{1}{(q-2)}}x\vert^2- \vert (\tfrac{\varepsilon}{2})^{\nicefrac{1}{(q-2)}}x\vert \right\vert
%		\\&\le (\tfrac{\varepsilon}{2})^{-\nicefrac{2}{(q-2)}} \left(\vert(\tfrac{\varepsilon}{2})^{\nicefrac{1}{(q-2)}}x\vert^2+ \vert (\tfrac{\varepsilon}{2})^{\nicefrac{1}{(q-2)}}x\vert\right)  
%		=  \left(\vert x\vert^2+ (\tfrac{\varepsilon}{2})^{-\nicefrac{1}{(q-2)}}\vert x\vert\right)  
%		\\&\le  \left( \vert x\vert^q\vert x\vert^{-(q-2)}+ (\tfrac{\varepsilon}{2})^{-\nicefrac{1}{(q-2)}}\vert x\vert^q\vert x\vert^{-(q-1)}\right)
%		\\&\le  \left( \vert x\vert^q (\tfrac{\varepsilon}{2})^{\nicefrac{(q-2)}{(q-2)}} \vert x\vert^{-(q-2)}+ (\tfrac{\varepsilon}{2})^{-\nicefrac{1}{(q-2)}}\vert x\vert^q(\tfrac{\varepsilon}{2})^{\nicefrac{(q-1)}{(q-2)}}\right)  
%		\\&=\varepsilon \vert x\vert^q\le \varepsilon \big(1+\vert x\vert^q\big). 
%		\end{split}
%		\end{equation}	
	Next note that \eqref{ApproxOfSquare:ComputationRealization}, \eqref{ApproxOfSquare:EstimateUnitCube}, and 
%	shows that for all $x\in\R$ it holds that
%	\begin{equation}
%	\begin{split}
%	&\left\vert x^2-[\functionANN(\Phi)](x)\right\vert
%	\\&= (\tfrac{\varepsilon}{2})^{-\nicefrac{2}{(q-2)}} \left\vert \big((\tfrac{\varepsilon}{2})^{\nicefrac{1}{(q-2)}}x\big)^2- \Big(	[\functionANN(\Psi)]\big((\tfrac{\varepsilon}{2})^{\nicefrac{1}{(q-2)}}x\big)+ [\functionANN(\Psi)]\big(-(\tfrac{\varepsilon}{2})^{\nicefrac{1}{(q-2)}}x\big)\Big)\right\vert.
%	\end{split}
%	\end{equation}
	the fact that $\delta=2^{-\nicefrac{2}{(q-2)}} \varepsilon^{\nicefrac{q}{(q-2)}}$ demonstrate that  for all $x\in  [-({\varepsilon}/{2})^{-\nicefrac{1}{(q-2)}},({\varepsilon}/{2})^{-\nicefrac{1}{(q-2)}}]$ it holds that
	\begin{equation}\label{ApproxOfSquare:ApproximationSmallValues}
	\begin{split}
	&\left\vert x^2-(\functionANN(\Phi))(x)\right\vert
	\\&= (\tfrac{\varepsilon}{2})^{-\nicefrac{2}{(q-2)}} \left\vert \big((\tfrac{\varepsilon}{2})^{\nicefrac{1}{(q-2)}}x\big)^2- \Big(	[\functionANN(\Psi)]\big((\tfrac{\varepsilon}{2})^{\nicefrac{1}{(q-2)}}x\big)+ [\functionANN(\Psi)]\big(-(\tfrac{\varepsilon}{2})^{\nicefrac{1}{(q-2)}}x\big)\Big)\right\vert
	\\&\le (\tfrac{\varepsilon}{2})^{-\nicefrac{2}{(q-2)}}\bigg[ \sup_{y\in [-1,1]}\left\vert y^2-\big([\functionANN(\Psi)](y)+ [\functionANN(\Psi)](-y)\big)\right\vert\bigg] 
	\\&\le  (\tfrac{\varepsilon}{2})^{-\nicefrac{2}{(q-2)}} \delta
	=(\tfrac{\varepsilon}{2})^{-\nicefrac{2}{(q-2)}} 2^{-\nicefrac{2}{(q-2)}} \varepsilon^{\nicefrac{q}{(q-2)}}=\varepsilon\le \varepsilon \max\!\big\{1,\vert x\vert^q\big\}.
	\end{split}
	\end{equation}
%		This, the fact that for all $\varepsilon\in (0,1]$, $x\in [-(\tfrac{\varepsilon}{2})^{-\nicefrac{1}{(q-2)}},(\tfrac{\varepsilon}{2})^{-\nicefrac{1}{(q-2)}}]$ it holds that $\vert(\tfrac{\varepsilon}{2})^{\nicefrac{1}{(q-2)}}x\vert\le 1$, and \eqref{ApproxOfSquare:BasicEstimate} prove that for all $\varepsilon\in (0,1]$, $x\in [-(\tfrac{\varepsilon}{2})^{-\nicefrac{1}{(q-2)}},(\tfrac{\varepsilon}{2})^{-\nicefrac{1}{(q-2)}}]$ it holds that 	
%		\begin{equation}\label{ApproxOfSquare:SmallValues}
%		\begin{split}
%		\left\vert x^2-(\functionANN (\Phi))(x)\right\vert
%		%	=\sup_{x\in [-(\tfrac{\varepsilon}{2})^{-\nicefrac{1}{(q-2)}},(\tfrac{\varepsilon}{2})^{-\nicefrac{1}{(q-2)}}]}\left\vert x^2-(\tfrac{\varepsilon}{2})^{-\nicefrac{2}{(q-2)}} (\functionANN (\Psi))((\tfrac{\varepsilon}{2})^{\nicefrac{1}{(q-2)}}x)\right\vert
%		%	\\&=(\tfrac{\varepsilon}{2})^{-\nicefrac{2}{(q-2)}} \sup_{x\in [-(\tfrac{\varepsilon}{2})^{-\nicefrac{1}{(q-2)}},(\tfrac{\varepsilon}{2})^{-\nicefrac{1}{(q-2)}}]}\left\vert ((\tfrac{\varepsilon}{2})^{\nicefrac{1}{(q-2)}}x)^2- (\functionANN (\Psi))((\tfrac{\varepsilon}{2})^{\nicefrac{1}{(q-2)}}x)\right\vert 
%		&\le(\tfrac{\varepsilon}{2})^{-\nicefrac{2}{(q-2)}} \sup_{y\in [-1,1]}\left\vert y^2- (\functionANN (\Psi))(y)\right\vert 
%		\\&\le (\tfrac{\varepsilon}{2})^{-\nicefrac{2}{(q-2)}}
%		2^{-\nicefrac{2}{(q-2)}} \varepsilon^{\nicefrac{q}{(q-2)}} =
%		\varepsilon\le \varepsilon \big(1+\vert x\vert^q\big).
%		\end{split}
%		\end{equation}	
		Combining this and \eqref{ApproxOfSquare:LargeValues} implies that for all  $x\in \R$ it holds that
		\begin{equation}\label{ApproxOfSquare:ApproximationAllValues}
		\begin{split}
		&\left\vert x^2-(\functionANN (\Phi))(x)\right\vert
		\le \varepsilon \max\!\big\{1,\vert x\vert^q\big\}
		\le \varepsilon \big(1+\vert x\vert^q\big). 
		\end{split}
		\end{equation}
				In addition, note that \eqref{ApproxOfSquare:ApproximationSmallValues} ensures that for all  $x\in [-({\varepsilon}/{2})^{-\nicefrac{1}{(q-2)}},({\varepsilon}/{2})^{-\nicefrac{1}{(q-2)}}]$ it holds that 	
				\begin{equation}\label{ApproxOfSquare:GrowthSmallValues}
				\left\vert (\functionANN (\Phi))(x)\right\vert\le \left\vert x^2-(\functionANN (\Phi))(x)\right\vert+	\vert x\vert^2\le \varepsilon+\vert x\vert^2.
				\end{equation}
				This and \eqref{ApproxOfSquare:BasicCalculation} show for all $x\in\R$ that 
				\begin{equation}\label{ApproxOfSquare:AllValues}
				\left\vert (\functionANN (\Phi))(x)\right\vert\le \varepsilon+\vert x\vert^2.
				\end{equation}
%Counting Params and Length
Furthermore, observe that the fact that $\delta=2^{-\nicefrac{2}{(q-2)}} \varepsilon^{\nicefrac{q}{(q-2)}}$ ensures that 
\begin{equation}\label{ApproxOfSquare:LogCalculation}
	\LogBin(\delta^{-1})=\LogBin(2^{\nicefrac{2}{(q-2)}} \varepsilon^{-\nicefrac{q}{(q-2)}})	
	=\tfrac{2}{(q-2)}+\Big[\big[\tfrac{q}{(q-2)}\big]\LogBin({\varepsilon}^{-1})\Big].
\end{equation}
Next note	 that  Corollary~\ref{Lemma:ParallelizationImprovedBoundsOne} implies that
	$\paramANN\big(\parallelizationSpecial_2(\Psi,\Psi)\big)\le 4\, \paramANN(\Psi)$.
	%	\begin{equation}
	%	\paramANN\big(\parallelizationSpecial_2(\Psi,\Psi)\big)=4\, \paramANN(\Psi).
	%	\end{equation}
	Corollary~\ref{Lemma:PropertiesOfCompositionsWithAffineMaps}, \eqref{ApproxOfSquare:ConstructionNetwork}, \eqref{ApproxOfSquare:BasisItemFour}, and \eqref{ApproxOfSquare:LogCalculation} hence ensure that
	\begin{equation}\label{ApproxOfSquare:ParamEstimate}
	\begin{split}
	\paramANN(\Phi)&\le \left[\max\!\left\{1,\tfrac{\outDimANN(\affineMap_2)}{\outDimANN(\parallelizationSpecial_2(\Psi,\Psi))}\right\}\right] 
	\left[\max\!\left\{1,\tfrac{\inDimANN(\affineMap_1)+1}{\inDimANN(\parallelizationSpecial_2(\Psi,\Psi))+1}\right\}\right]
		\paramANN\big(\parallelizationSpecial_2(\Psi,\Psi)\big)
	\\&= \left[\max\{1,\tfrac{1}{2}\}\right] \left[\max\{1,\tfrac{2}{3}\}\right] \paramANN\big(\parallelizationSpecial_2(\Psi,\Psi)\big)
	\\&=\paramANN\big(\parallelizationSpecial_2(\Psi,\Psi)\big)
	\le 4\, \paramANN(\Psi)
		\le 4\max\{10\LogBin(\delta^{-1})-7,13\}
		\\&= \max\!\big\{40\big[\tfrac{2}{(q-2)}\big]+40 \big[\tfrac{q}{(q-2)}\big]\LogBin({\varepsilon}^{-1})-28,52\big\}
		\\&= \max\!\big\{\big[\tfrac{40q}{(q-2)}\big]\LogBin({\varepsilon}^{-1})+\tfrac{80}{(q-2)}-28,52\big\}.
	\end{split}
	\end{equation}	
	In addition, observe that item~\eqref{PropertiesOfCompositions:Length} in Proposition~\ref{Lemma:PropertiesOfCompositions}, \eqref{ApproxOfSquare:ConstructionNetwork}, \eqref{ApproxOfSquare:BasisItemFive}, and \eqref{ApproxOfSquare:LogCalculation} demonstrate that 
	\begin{equation}
	\begin{split}
		\lengthANN(\Phi)&=\lengthANN\big(\parallelizationSpecial_2(\Psi,\Psi)\big)=\lengthANN(\Psi)
		\le \max\left\{\tfrac{1}{2}\LogBin(\delta^{-1})+1,2\right\}
		\\&= \max \!\big\{\big[\tfrac{q}{2(q-2)}\big]\LogBin({\varepsilon}^{-1})+\tfrac{1}{(q-2)}+1,2\big\}.
	\end{split}
	\end{equation}
Combining this with \eqref{eq:ApproxOfSquare:Zero}, \eqref{ApproxOfSquare:BasicCalculation}, \eqref{ApproxOfSquare:AllValues}, \eqref{ApproxOfSquare:ApproximationAllValues}, \eqref{ApproxOfSquare:ParamEstimate} establishes items~\eqref{ApproxOfSquare:Realization}--\eqref{ApproxOfSquare:Length}.
	The proof of Proposition~\ref{Lemma:ApproxOfSquare} is thus completed.
\end{proof}

%% file: ANNApproxOneDimProducts.tex
\begin{prop}\label{Lemma:ApproxOfProduct}
	Let  $\varepsilon\in (0,1]$, $q\in (2,\infty)$,  $a\in C(\R,\R)$ satisfy for all $x\in \R$ that $\activation(x)=\max\{x,0\}$.
	Then there exists $\Phi\in \ANNs$
	such that
	\begin{enumerate}[(i)]
		\item \label{ApproxOfProduct:RealizationProp}
		it holds that  $\functionANN (\Phi)\in C(\R^2,\R)$,
		\item \label{ApproxOfProduct:ZeroProp}
		it holds for all $x\in\R$ that 	$(\functionANN (\Phi))(x,0)=(\functionANN (\Phi))(0,x)=0$,
%		\begin{equation}
%		(\functionANN (\Phi))(x,0)=(\functionANN (\Phi))(0,y)=0,
%		\end{equation}
		\item \label{ApproxOfProduct:EstimatesProp}
		it holds for all  $x,y\in\R$ that 
		\begin{equation}
		\vert xy-(\functionANN (\Phi))(x,y)\vert\le \varepsilon \max\!\big\{1,\vert x\vert^q,\vert y\vert^q\big\},
		\end{equation}
		\item \label{ApproxOfProduct:EstimatesGrowth}
		it holds for all  $x,y\in\R$ that 
		\begin{equation}
		\vert(\functionANN (\Phi))(x,y)\vert
		\le 
		\tfrac{3}{2}\big( \tfrac{\varepsilon}{3}+ x^2+  y^2\big)\le 1+2 x^2+ 2 y^2,
%				\tfrac{3}{2}\big( \tfrac{\varepsilon}{3}+\vert x\vert^2+ \vert y\vert^2\big)\le 1+2\vert x\vert^2+ 2\vert y\vert^2,
		\end{equation}
				\item \label{ApproxOfProduct:ParamsProp}
				it holds  that 
				\begin{equation}
				\begin{split}
				\paramANN(\Phi)&\le\tfrac{360q}{(q-2)}\big[\LogBin(\eps^{-1})+\LogBin(2^{q-1}+1)\big]+\tfrac{1}{(q-2)}-252
				%	=\tfrac{360q}{(q-2)}\LogBin({\delta}^{-1})+\tfrac{720}{(q-2)}+468.
				\\&\le \tfrac{360q}{(q-2)}\big[\LogBin(\eps^{-1})+q+1\big]-252,
				\end{split}
				\end{equation}	
				and
		\item \label{ApproxOfProduct:LengthProp}
		it holds  that 
				\begin{equation}
				\begin{split}
				\lengthANN(\Phi)&\le
				\tfrac{q}{2(q-2)}\big[\LogBin(\eps^{-1})+\LogBin(2^{q-1}+1)\big]+\tfrac{(q-1)}{(q-2)}
				 \\&\le\tfrac{q}{(q-2)}\big[\LogBin(\eps^{-1})+q\big]
				\end{split}
				\end{equation}
%		\item \label{ApproxOfProduct:ParamsProp}
%		it holds for all $\varepsilon\in (0,1]$, $q\in (2,\infty)$ that 
%	\begin{equation}
%	\begin{split}
%	\paramANN(\Phi)&\le\tfrac{360q}{(q-2)}\LogBin(\eps^{-1})+\tfrac{360q}{(q-2)}\LogBin(2^{q-1}+1)+108
%	%	=\tfrac{360q}{(q-2)}\LogBin({\delta}^{-1})+\tfrac{720}{(q-2)}+468.
%	\end{split}
%	\end{equation}	
%and
%	\begin{equation}\label{ApproxOfProduct:LengthEstimate}
%	\begin{split}
%	\lengthANN(\Phi)&\le \max \big\{\tfrac{q}{2(q-2)}\LogBin(\eps^{-1})+\tfrac{q}{2(q-2)}\LogBin(2^{q-1}+1)+\tfrac{1}{2},2\big\}
%	\end{split}
%	\end{equation}
	\end{enumerate}
	(cf.\ Definition~\ref{Def:ANN} and Definition~\ref{Definition:ANNrealization}).
\end{prop}

\begin{proof}[Proof of Proposition~\ref{Lemma:ApproxOfProduct}]	
		Throughout this proof let $\delta\in (0,1]$ satisfy that $\delta=\varepsilon (2^{q-1}+1)^{-1}$, 
		let $\affineMap_1\in (\R^{3\times 2}\times \R^3)\subseteq \ANNs$, $\affineMap_2\in(\R^{1\times 3}\times \R)\subseteq \ANNs$ satisfy that
	\begin{equation}
	\affineMap_1=\left(\begin{pmatrix}
	1 &1\\1 &0\\0&1 
	\end{pmatrix},\begin{pmatrix}
	0 \\0 \\0
	\end{pmatrix}\right)\qandq
	\affineMap_2=\left(\begin{pmatrix}
	\tfrac{1}{2} &-\tfrac{1}{2}&-\tfrac{1}{2} 
	\end{pmatrix},0\right)\!,
	\end{equation}	
	 let $\Psi\in \ANNs$ satisfy that 
	\begin{enumerate}[(I)]
		\item \label{ApproxOfProduct:RealizationSquare}
		it holds that  $\functionANN (\Psi)\in C(\R,\R)$, 
		\item \label{ApproxOfProduct:Zero} it holds that $[\functionANN (\Psi)](0)=0$,
		\item \label{ApproxOfProduct:EstimatesSquareGrowth}
		it holds for all $x\in\R$  that $0\le [\functionANN (\Psi)](x)\le \delta+ \vert x\vert^2$, 
		\item\label{ApproxOfProduct:EstimatesSquare} it holds for all $x\in\R$  that $\vert x^2-[\functionANN (\Psi)](x)\vert\le \delta \max\!\big\{1,\vert x\vert^q\big\}$,
%		\begin{equation}
%		\vert x^2-[\functionANN (\Psi)](x)\vert\le \delta \max\!\big\{1,\vert x\vert^q\big\},
%		\end{equation}
		\item \label{ApproxOfProduct:ParamsSquare}
		it holds that   $\paramANN(\Psi)\le  \max\!\big\{\big[\tfrac{40q}{(q-2)}\big]\LogBin({\delta}^{-1})+\tfrac{80}{(q-2)}-28,52\big\}$, and 
%		\begin{equation}
%		\paramANN(\Psi)\le  \max\!\big\{\big[\tfrac{40q}{(q-2)}\big]\LogBin({\varepsilon}^{-1})+\tfrac{80}{(q-2)}-28,52\big\}
%		\end{equation}	
		\item \label{ApproxOfProduct:LengthSquare}
		it holds that $\lengthANN(\Psi)\le \max\! \big\{\big[\tfrac{q}{2(q-2)}\big]\LogBin({\delta}^{-1})+\tfrac{1}{(q-2)}+1,2\big\}$
%		\begin{equation}
%		\andq \lengthANN(\Psi)\le \max\! \big\{\big[\tfrac{q}{2(q-2)}\big]\LogBin({\varepsilon}^{-1})+\tfrac{1}{(q-2)}+1,2\big\}
%		\end{equation}
	\end{enumerate}
%			\item \label{ApproxOfSquare:Params}
%			it holds  that $\paramANN(\Phi)\le  \max\!\big\{\big[\tfrac{40q}{(q-2)}\big]\LogBin({\varepsilon}^{-1})+\tfrac{80}{(q-2)}-28,52\big\}$, and 
%			\item \label{ApproxOfSquare:Length}
%			it holds  that $\lengthANN(\Phi)\le \max\! \big\{\tfrac{q}{2(q-2)}\LogBin({\varepsilon}^{-1})+\tfrac{1}{(q-2)}+1,2\big\}$
%	 let $\Psi\in \ANNs$ such that it holds that $\functionANN (\Psi)\in C(\R,\R)$ and for all $x\in\R$ it holds that 
%	\begin{equation}\label{ApproxOfProduct:squareEstimate}
%	0\le (\functionANN (\Psi))(x)\le (\tfrac{\delta}{2})^{-\nicefrac{1}{(q-2)}} \vert x\vert,\qquad\vert x^2-(\functionANN (\Psi))(x)\vert\le \delta \big(1+\vert x\vert^q\big),
%	\end{equation}
%	\begin{equation}
%	\paramANN(\Psi)\le  \tfrac{72}{(q-2)}+\tfrac{36q}{(q-2)}\LogBin(\tfrac{1}{\delta})+72
%	,\qquad  \paramNotZeroANN(\Psi)\le  \tfrac{30}{(q-2)}+\tfrac{15q}{(q-2)}\LogBin(\tfrac{1}{\delta})+30,
%	\end{equation}
%	\begin{equation}
%	\lengthANN(\Psi)\le  \tfrac{1}{(q-2)}+\tfrac{q}{2(q-2)}\LogBin(\tfrac{1}{\delta})+1,
%	\end{equation}
%	and $(\functionANN (\Psi))(0)=0$
	(cf.\ Proposition~\ref{Lemma:ApproxOfSquare}),
 and let $\Phi\in \ANNs$ satisfy that  
	\begin{equation}\label{ApproxOfProduct:ConstructionNetwork}
	\Phi=\compANN{\affineMap_2}{\compANN{\big[\parallelizationSpecial_3(\Psi,\Psi,\Psi)\big]}{\affineMap_1}} 
	\end{equation}
	(cf.\ Definition~\ref{Definition:ANNcomposition}, Definition~\ref{Definition:simpleParallelization}, and Lemma~\ref{Lemma:CompositionAssociative}).
	Note that  item~\eqref{PropertiesOfCompositions:Realization} in	Proposition~\ref{Lemma:PropertiesOfCompositions} and Proposition~\ref{Lemma:PropertiesOfParallelizationEqualLength} ensure that for all $x,y\in \R$ it holds that $\functionANN\big(\compANN{\big[\parallelizationSpecial_3(\Psi,\Psi,\Psi)\big]}{\affineMap_1}\big)\allowbreak\in C(\R^2,\R^3)$ and
	\begin{equation}
	\begin{split}
	&\big[\functionANN\big(\compANN{\big[\parallelizationSpecial_3(\Psi,\Psi,\Psi)\big]}{\affineMap_1}\big)\big](x,y)
%	=\big(\big[\functionANN\big(\parallelizationSpecial_3(\Psi,\Psi,\Psi)\big)\big]\circ\big[\functionANN(\affineMap_1)\big]\big)(x,y)
	=\big[\functionANN\big(\parallelizationSpecial_3(\Psi,\Psi,\Psi)\big)\big]\big(\big[\functionANN(\affineMap_1)\big](x,y)\big)
	\\&=\big[\functionANN\big(\parallelizationSpecial_3(\Psi,\Psi,\Psi)\big)\big](x+y,x,y)
	%	\\&=\big([\functionANN(\Psi)](x,y),[\functionANN(\Psi)](-x)\big).
	=\begin{pmatrix}
	[\functionANN(\Psi)](x+y)\\ [\functionANN(\Psi)](x)\\ [\functionANN(\Psi)](y)
	\end{pmatrix}.
	\end{split}
	\end{equation}
	Item~\eqref{PropertiesOfCompositions:Realization} in Proposition~\ref{Lemma:PropertiesOfCompositions} and \eqref{ApproxOfProduct:ConstructionNetwork} therefore demonstrate that for all $x,y\in \R$ it holds that $\functionANN (\Phi)\in C(\R^2,\R)$ and
	\begin{equation}\label{ApproxOfProduct:TransformedANN}
	\begin{split}
	[\functionANN (\Phi)](x,y)
	&= \big(\functionANN\big(\compANN{\affineMap_2}{\compANN{\big[\parallelizationSpecial_3(\Psi,\Psi,\Psi)\big]}{\affineMap_1}}\big)\big)(x,y)
	\\&= [\functionANN({\affineMap_2})]\big(\functionANN\big(\compANN{\big[\parallelizationSpecial_3(\Psi,\Psi,\Psi)\big]}{\affineMap_1}\big)(x,y)\big)
	\\&=\begin{pmatrix}
	\tfrac{1}{2} &-\tfrac{1}{2}&-\tfrac{1}{2} 
	\end{pmatrix} \begin{pmatrix}
	[\functionANN(\Psi)](x+y)\\ [\functionANN(\Psi)](x)\\ [\functionANN(\Psi)](y)
	\end{pmatrix}
	\\&=	\tfrac{1}{2} [\functionANN(\Psi)](x+y)-\tfrac{1}{2} [\functionANN(\Psi)](x)-\tfrac{1}{2}[\functionANN(\Psi)](y).
	\end{split}
	\end{equation}
%	
%	
%	
%	Throughout this proof  and let $\Phi\in \ANNs$ such that for all $x,y\in\R$ it holds that 
%	\begin{equation}\label{ApproxOfProduct:TransformedANN}
%	(\functionANN (\Phi))(x,y)=\tfrac{1}{2} (\functionANN (\Psi))(x+y)-\tfrac{1}{2} (\functionANN (\Psi))(x)-\tfrac{1}{2} (\functionANN (\Psi))(y).
%	\end{equation}	
	The fact that for all $\alpha,\beta\in\R$ it holds that $\alpha\beta=\tfrac{1}{2}\vert \alpha+\beta\vert^2-\tfrac{1}{2}\vert \alpha\vert^2-\tfrac{1}{2}\vert \beta\vert^2$, the triangle inequality,  and \eqref{ApproxOfProduct:EstimatesSquare} hence ensure that for all $x,y\in\R$ it holds that 
	\begin{equation}
	\begin{split}
	&\left\vert [\functionANN (\Phi)](x,y)-xy\right\vert
	\\&=\left\vert \tfrac{1}{2} \big[[\functionANN (\Psi)](x+y)-\vert x+y\vert^2\big]-\tfrac{1}{2} \big[[\functionANN (\Psi)](x)-\vert x\vert^2\big]-\tfrac{1}{2} \big[[\functionANN (\Psi)](y)-\vert y\vert^2\big]\right\vert
	\\&\le\tfrac{1}{2}\left\vert  [\functionANN (\Psi)](x+y)-\vert x+y\vert^2\right\vert+\tfrac{1}{2} \left\vert [\functionANN (\Psi)](x)-\vert x\vert^2\right\vert+\tfrac{1}{2} \left\vert [\functionANN (\Psi)](y)-\vert y\vert^2\right\vert
	\\&\le \tfrac{\delta}{2} \big[\max\!\big\{1,\vert x+y\vert^q\big\}+ \max\!\big\{1,\vert x\vert^q\big\}+ \max\!\big\{1,\vert y\vert^q\big\}\big].
%		\\&\le \tfrac{1}{2}\delta \big(1+ 2^{q -1}[|x|^{q} + |y|^{q}]\big)+\tfrac{1}{2}\delta \big(1+\vert x\vert^q\big)+\tfrac{1}{2}\delta \big(1+\vert y\vert^q\big)
%		\\&= \tfrac{1}{2}\delta \big(3+ [2^{q -1}+1][|x|^{q} + |y|^{q}]\big)
%		\le \tfrac{1}{2}\delta  [2^{q -1}+1] \big(1+ |x|^{q} + |y|^{q}\big)
%		\\&= \varepsilon \big(1+ |x|^{q} + |y|^{q}\big).
	\end{split}
	\end{equation}
	This, the fact that for all $\alpha,\beta\in\R$, $p \in [1, \infty)$ it holds that   $|\alpha + \beta|^{p}  \leq 2^{p -1}(|\alpha|^{p} + |\beta|^{p})$, and the fact that $\delta=\varepsilon (2^{q-1}+1)^{-1}$ establish that for all $x,y\in\R$ it holds that 
	\begin{equation}\label{ApproxOfProduct:AccuracyEstimate}
	\begin{split}
	&\left\vert [\functionANN (\Phi)](x,y)-xy\right\vert
	\\&\le \tfrac{\delta}{2} \big[ \max\!\big\{1, 2^{q -1}|x|^{q} + 2^{q -1} |y|^{q}\big\}+ \max\!\big\{1,\vert x\vert^q\big\}+ \max\!\big\{1,\vert y\vert^q\big\}\big]
		\\&\le \tfrac{\delta}{2} \big[\max\!\big\{1, 2^{q -1}|x|^{q} \big\}+ 2^{q -1} |y|^{q}+ \max\!\big\{1,\vert x\vert^q\big\}+ \max\!\big\{1,\vert y\vert^q\big\}\big]
		\\&\le \tfrac{\delta}{2} \big[2^q+2\big]\max\!\big\{1, |x|^{q}, |y|^{q}\big\}= \varepsilon \max\!\big\{1, |x|^{q}, |y|^{q}\big\}.
	\end{split}
	\end{equation}
	Moreover, observe that \eqref{ApproxOfProduct:EstimatesSquareGrowth}, \eqref{ApproxOfProduct:TransformedANN}, the triangle inequality,
	the fact that for all $\alpha,\beta\in\R$ it holds that $\vert \alpha+\beta\vert^2\le 2(\vert \alpha\vert^2+\vert \beta\vert^2)$,
	and the fact that $\delta=\varepsilon (2^{q-1}+1)^{-1}$ prove that  for all $x,y\in\R$ it holds  that
%	\begin{equation}
%		\begin{split}
%			\vert[\functionANN (\Phi)](x,y)\vert
%		&\le\tfrac{1}{2} \vert [\functionANN(\Psi)](x+y)\vert +\tfrac{1}{2} \vert[\functionANN(\Psi)](x)\vert+\tfrac{1}{2}\vert[\functionANN(\Psi)](y)\vert
%		\\&\le  \tfrac{1}{2} (\tfrac{\delta}{2})^{-\nicefrac{1}{(q-2)}} \vert x+y\vert +\tfrac{1}{2} (\tfrac{\delta}{2})^{-\nicefrac{1}{(q-2)}} \vert x\vert+\tfrac{1}{2} (\tfrac{\delta}{2})^{-\nicefrac{1}{(q-2)}} \vert y\vert
%		\\&\le (\tfrac{\delta}{2})^{-\nicefrac{1}{(q-2)}} \big(\vert x\vert+\vert y\vert\big)
%		\le \varepsilon^{-\nicefrac{1}{(q-2)}} (2^{q-1}+1)^{\nicefrac{1}{(q-2)}} \big(\vert x\vert+\vert y\vert\big).
%		\end{split}
%	\end{equation}
	\begin{equation}\label{ApproxOfProduct:GrowthEstimate}
\begin{split}
\vert[\functionANN (\Phi)](x,y)\vert
&\le\tfrac{1}{2} \vert [\functionANN(\Psi)](x+y)\vert +\tfrac{1}{2} \vert[\functionANN(\Psi)](x)\vert+\tfrac{1}{2}\vert[\functionANN(\Psi)](y)\vert
\\&\le  \tfrac{1}{2}\big(\delta+ \vert x+y\vert^2\big) +\tfrac{1}{2} \big(\delta+ \vert x\vert^2\big)+\tfrac{1}{2} \big(\delta+ \vert y\vert^2\big)
\\&\le \tfrac{3\delta}{2} + \tfrac{3}{2} \big( \vert x\vert^2+ \vert y\vert^2\big)
= \big[\tfrac{3\varepsilon}{2}\big] [2^{q-1}+1]^{-1}+ \tfrac{3}{2} \big( \vert x\vert^2+ \vert y\vert^2\big)
\\&= \tfrac{3}{2}\big[\tfrac{\varepsilon}{(2^{q-1}+1)}+ \vert x\vert^2+ \vert y\vert^2\big]
\le \tfrac{3}{2}\big[\tfrac{\varepsilon}{3}+ \vert x\vert^2+ \vert y\vert^2\big].
\end{split}
\end{equation}
%Next note that \eqref{ApproxOfProduct:RealizationSquare} and \eqref{ApproxOfProduct:TransformedANN} prove that for all $x,y\in\R$ it holds that 
%\begin{equation}\label{ApproxOfProduct:ZeroAnnihilation}
%	[\functionANN (\Phi)](x,0)=[\functionANN (\Phi)](0,y)=0.
%\end{equation}
Next note that \eqref{ApproxOfProduct:RealizationSquare} and \eqref{ApproxOfProduct:TransformedANN} prove that for all $x,y\in\R$ it holds that 
\begin{equation}\label{ApproxOfProduct:ZeroAnnihilation}
\begin{split}
[\functionANN (\Phi)](x,0)
&=\tfrac{1}{2} [\functionANN(\Psi)](x)-\tfrac{1}{2} [\functionANN(\Psi)](x)-\tfrac{1}{2}[\functionANN(\Psi)](0)
=0
\\&=\tfrac{1}{2} [\functionANN(\Psi)](y)-\tfrac{1}{2} [\functionANN(\Psi)](0)-\tfrac{1}{2}[\functionANN(\Psi)](y)
=[\functionANN (\Phi)](0,y).
\end{split}
\end{equation}
Furthermore, observe that the fact that $\delta=\varepsilon (2^{q-1}+1)^{-1}$ shows that
	\begin{equation}\label{ApproxOfProduct:AuxiliaryCalculation}
	\begin{split}
	&\big[\tfrac{q}{2(q-2)}\big] \LogBin(\delta^{-1})+\tfrac{1}{(q-2)}=\big[\tfrac{q}{2(q-2)}\big] \LogBin\!\big(\varepsilon^{-1} (2^{q-1}+1)\big)+\tfrac{1}{(q-2)}
	\\&= \tfrac{q}{2(q-2)}\big[\LogBin(\eps^{-1})+\LogBin(2^{q-1}+1)\big]+\tfrac{1}{(q-2)}
	\\&= \big[\tfrac{q}{2(q-2)}\big]\LogBin(\eps^{-1})+\big[\tfrac{q}{2(q-2)}\big]\LogBin(2^{q-1}+1)+\tfrac{1}{(q-2)}.
	\end{split}
	\end{equation}
	Moreover, observe that Corollary~\ref{Lemma:ParallelizationImprovedBoundsOne} implies that
	$\paramANN\big(\parallelizationSpecial_3(\Psi,\Psi,\Psi)\big)\allowbreak\le 9\, \paramANN(\Psi)$.
	%	\begin{equation}
	%		\paramANN\big(\parallelizationSpecial_3(\Psi,\Psi,\Psi)\big)=9\, \paramANN(\Psi).
	%	\end{equation}
	Items~\eqref{PropertiesOfCompositionsWithAffineMaps:ItemFront}--\eqref{PropertiesOfCompositionsWithAffineMaps:ItemBehind} in Corollary~\ref{Lemma:PropertiesOfCompositionsWithAffineMaps}, \eqref{ApproxOfProduct:ParamsSquare}, \eqref{ApproxOfProduct:ConstructionNetwork}, and \eqref{ApproxOfProduct:AuxiliaryCalculation} hence ensure that
	\begin{equation}\label{ApproxOfProduct:ParamEstimate}
	\begin{split}
	&\paramANN(\Phi)\le \left[\max\!\left\{1,\tfrac{\outDimANN(\affineMap_2)}{\outDimANN(\parallelizationSpecial_3(\Psi,\Psi,\Psi))}\right\} \right]
	\left[\max\!\left\{1,\tfrac{\inDimANN(\affineMap_1)+1}{\inDimANN(\parallelizationSpecial_3(\Psi,\Psi,\Psi))+1}\right\}\right]
		\paramANN\big(\parallelizationSpecial_3(\Psi,\Psi,\Psi)\big)
	\\&= \left[\max\{1,\tfrac{1}{3}\}\right]	
	\left[\max\{1,\tfrac{3}{4}\}\right] \paramANN\big(\parallelizationSpecial_3(\Psi,\Psi,\Psi)\big)=\paramANN\big(\parallelizationSpecial_3(\Psi,\Psi,\Psi)\big)
	\\&\le 9\, \paramANN(\Psi)
	\le 9 \max\!\big\{\big[\tfrac{40q}{(q-2)}\big]\LogBin({\delta}^{-1})+\tfrac{80}{(q-2)}-28,52\big\}
	\\&=  \max\!\big\{720\big(\big[\tfrac{q}{2(q-2)}\big]\LogBin({\delta}^{-1})+\tfrac{1}{(q-2)}\big)-252,468\big\}
	\\&=  \max\!\big\{720\big(\big[\tfrac{q}{2(q-2)}\big]\LogBin(\eps^{-1})+\big[\tfrac{q}{2(q-2)}\big]\LogBin(2^{q-1}+1)+\tfrac{1}{(q-2)}\big)-252,468\big\}		
	\\&=  \max\!\big\{\tfrac{360q}{(q-2)}\big(\LogBin(\eps^{-1})+\LogBin(2^{q-1}+1)\big)+\tfrac{720}{(q-2)}-252,468\big\}.
	\end{split}
	\end{equation}	
	Next note that the fact that for all $r\in (-\infty,4]$ it holds that $r\ge 2r-4=2(r-2)$ ensures that for all $r\in (2,4]$ it holds that $\tfrac{r(r-1)}{(r-2)}\ge \tfrac{r}{(r-2)}\ge 2$.
	This and the fact that for all $r\in [3,\infty)$ it holds that $\tfrac{r(r-1)}{(r-2)}\ge r-1\ge 2$ imply that for all $r\in (2,\infty)$ it holds that $\tfrac{r(r-1)}{(r-2)}\ge  2$.
	Hence, we obtain that for all $r\in (2,\infty)$ it holds that 
	\begin{equation}
		\begin{split}
		\big[\tfrac{360r}{(r-2)}\big]\LogBin(2^{r-1}+1)-252
		&\ge \big[\tfrac{360r}{(r-2)}\big]\LogBin(2^{r-1})-252
		\\&=\tfrac{360r(r-1)}{(r-2)}-252\ge 720-252=468.
		\end{split}
	\end{equation}
	Combining this with \eqref{ApproxOfProduct:ParamEstimate} shows that 
		\begin{equation}\label{ApproxOfProduct:ParamEstimateTilde}
		\begin{split}
		\paramANN(\Phi)&\le 		
		\tfrac{360q}{(q-2)}\big(\LogBin(\eps^{-1})+\LogBin(2^{q-1}+1)\big)+\tfrac{720}{(q-2)}-252.
%		\\&\le \tfrac{360q}{(q-2)}\big(\LogBin(\eps^{-1})+q\big)+\tfrac{720}{(q-2)}-252.
		\end{split}
		\end{equation}
	The fact that 
%	\begin{equation}
%		\begin{split}
%		q-\LogBin(2^{q-1}+1)&=\LogBin(2^q)-\LogBin(2^{q-1}+1)=\LogBin\!\big(\tfrac{2^q}{2^{q-1}+1}\big)
%		\\&=\LogBin\!\big(\tfrac{1}{2^{-1}+2^{-q}}\big)\ge \LogBin\!\big(\tfrac{1}{2^{-1}+2^{-2}}\big)
%%		=\LogBin\!\big(\tfrac{4}{3}\big)
%		\\&=\LogBin\!\big(\tfrac{2^2}{2+1}\big)=2-\LogBin(3)
%		\end{split}
%	\end{equation}
	\begin{equation}
	\begin{split}
	\LogBin(2^{q-1}+1)&=\LogBin(2^{q-1}+1)-\LogBin(2^q)+q
	=\LogBin\!\big(\tfrac{2^{q-1}+1}{2^q}\big)+q
	\\&=\LogBin\!\big({2^{-1}+2^{-q}}\big)+q
	\le \LogBin\!\big({2^{-1}+2^{-2}}\big)+q
	\\&=\LogBin\!\big(\tfrac{3}{4}\big)+q
	=\LogBin(3)-2+q
	\end{split}
	\end{equation}
	hence proves that 
			\begin{equation}\label{ApproxOfProduct:ParamEstimateTildeTilde}
			\begin{split}
			\paramANN(\Phi)&\le 		
			\tfrac{360q}{(q-2)}\big(\LogBin(\eps^{-1})+\LogBin(2^{q-1}+1)\big)+\tfrac{720}{(q-2)}-252
			\\&\le\tfrac{360q}{(q-2)}\big(\LogBin(\eps^{-1})+q+\LogBin(3)-2\big)+\tfrac{720}{(q-2)}-252
			\\&= \tfrac{360q}{(q-2)}\big(\LogBin(\eps^{-1})+q+\LogBin(3)-2+\tfrac{2}{q}\big)-252
			\\&\le \tfrac{360q}{(q-2)}\big(\LogBin(\eps^{-1})+q+\LogBin(3)-1\big)-252.
			\end{split}
			\end{equation}
			In addition, observe that item~\eqref{PropertiesOfCompositions:Length} in Proposition~\ref{Lemma:PropertiesOfCompositions}, \eqref{ApproxOfProduct:ConstructionNetwork},  \eqref{ApproxOfProduct:LengthSquare}, the fact that $\delta=\varepsilon (2^{q-1}+1)^{-1}$, and \eqref{ApproxOfProduct:AuxiliaryCalculation} demonstrate that 
			\begin{equation}\label{ApproxOfProduct:LengthEstimate}
			\begin{split}
			\lengthANN(\Phi)&=\lengthANN\big(\parallelizationSpecial_3(\Psi,\Psi,\Psi)\big)=\lengthANN(\Psi)
			\\&\le \max\! \big\{\big[\tfrac{q}{2(q-2)}\big]\LogBin({\delta}^{-1})+\tfrac{1}{(q-2)}+1,2\big\}
			\\&\le \max \!\big\{\tfrac{q}{2(q-2)}\big[\LogBin(\eps^{-1})+\LogBin(2^{q-1}+1)\big]+\tfrac{(q-1)}{(q-2)},2\big\}.
			\end{split}
			\end{equation}
	Furthermore, note that the fact for all $r\in (2,\infty)$ it holds that $\tfrac{r(r-1)}{(r-2)}\ge  2$ assures that 
	\begin{equation}
		\begin{split}
		&\tfrac{q}{2(q-2)}\big[\LogBin(\eps^{-1})+\LogBin(2^{q-1}+1)\big]+\tfrac{(q-1)}{(q-2)}
		\\&\ge \big[\tfrac{q}{2(q-2)}\big]\LogBin(2^{q-1})+1
		= \tfrac{q(q-1)}{2(q-2)}+1\ge 2.
		\end{split}
	\end{equation}	
	Combining this with \eqref{ApproxOfProduct:LengthEstimate} proves that
			\begin{equation}
			\begin{split}
			\lengthANN(\Phi)
			&\le \tfrac{q}{2(q-2)}\big[\LogBin(\eps^{-1})+\LogBin(2^{q-1}+1)\big]+\tfrac{(q-1)}{(q-2)}
			\\&\le \tfrac{q}{2(q-2)}\big[\LogBin(\eps^{-1})+\LogBin(2^{q-1}+2^{q-1})\big]+\tfrac{q}{(q-2)}
			\\&= \tfrac{q}{(q-2)}\big[\tfrac{\LogBin(\eps^{-1})}{2}+\tfrac{q}{2}+1\big]
			\le  \tfrac{q}{(q-2)}\big[\LogBin(\eps^{-1})+\tfrac{q}{2}+\tfrac{q}{2}\big]
			\\&= \tfrac{q}{(q-2)}\big[\LogBin(\eps^{-1})+q\big].
			\end{split}
			\end{equation}	
	This, the fact that $\functionANN (\Phi)\in C(\R^2,\R)$, \eqref{ApproxOfProduct:AccuracyEstimate}, \eqref{ApproxOfProduct:GrowthEstimate}, \eqref{ApproxOfProduct:ZeroAnnihilation}, 
	and \eqref{ApproxOfProduct:ParamEstimateTildeTilde} establish items~\eqref{ApproxOfProduct:RealizationProp}--\eqref{ApproxOfProduct:LengthProp}.
	The proof of Proposition~\ref{Lemma:ApproxOfProduct} is thus completed.
\end{proof}

%% file: ANNApproxMultiDimProducts.tex
\begin{definition}[The Euclidean norm]
	\label{Def:euclideanNorm}
	We denote by $\norm{\cdot} \colon (\cup_{d\in\N} \R^d) \to [0,\infty)$ the function which satisfies for all $d\in\N$, $ x = ( x_1, \dots, x_{d} ) \in \R^{d} $ that
	\begin{equation}\label{euclideanNorm:Equation}
	\norm{x}=\big[\smallsum\nolimits_{j=1}^d \vert x_j\vert^2\big]^{\nicefrac{1}{2}}.
	\end{equation}
\end{definition}

\begin{prop}\label{Lemma:ApproxOfScalarVectorProduct}
	Let $\varepsilon\in (0,1]$, $q\in (2,\infty)$, $d\in\N$, $a\in C(\R,\R)$ satisfy for all $x\in \R$ that $\activation(x)=\max\{x,0\}$.
% 	let 
%	$
%	\left\| \cdot \right\| \colon \R^d \to [0,\infty)
%	$
%	be the  Euclidean norm on $\R^d$.
	Then there exists $\Phi\in \ANNs$ such that
	\begin{enumerate}[(i)]
		\item \label{ApproxOfScalarVectorProduct:RealizationProp}
		it holds that  $\functionANN (\Phi)\in C(\R^{d+1},\R^d)$,
		\item \label{ApproxOfScalarVectorProduct:ZeroProp}
		it holds for all  $t\in\R$, $x\in\R^d$ that $(\functionANN (\Phi))(t,0)=(\functionANN (\Phi))(0,x)
			=0$,
%	\begin{equation}
%	\begin{split}
%	&(\functionANN (\Phi))(0,x)
%	=(\functionANN (\Phi))(t,0)
%	=0,
%	\end{split}
%	\end{equation}	
		\item \label{ApproxOfScalarVectorProduct:EstimatesProp}
		it holds for all  $t\in\R$, $x\in\R^d$ that 
		\begin{equation}
		\left\|  t x-(\functionANN (\Phi))(t,x)\right\|
		\le \varepsilon \big(\sqrt{d}\left[\max\!\big\{1,\vert t\vert^q\big\}\right]+\| x\|^q\big),
%		\le \varepsilon \sqrt{d}\big(1+\vert t\vert^q+\left\| x\right\|^q\!\big),
		\end{equation}
			\item \label{ApproxOfScalarVectorProduct:EstimatesGrowth}
			it holds for all  $t\in\R$, $x\in\R^d$ that 
			\begin{equation}
			\left\| (\functionANN (\Phi))(t,x)\right\|
			\le \sqrt{d}\big(1+2 t^2\big)+2\| x\|^2,
%			\le 3 \sqrt{d}\big(1+\vert t\vert^2+\left\| x\right\|^2\!\big),
			\end{equation}
		\item \label{ApproxOfScalarVectorProduct:ParamsProp}
		it holds that $\paramANN(\Phi)\le d^2 \big[\tfrac{360q}{(q-2)}\big]\big[\LogBin(\eps^{-1})+q+1\big]-252d^2$,
%		\begin{equation}
%		\begin{split}
%		\paramANN(\Phi)&\le d^2 \big[\tfrac{360q}{(q-2)}\big]\big[\LogBin(\eps^{-1})+q+\LogBin(3)-1\big]-252d^2,
%		\end{split}
%		\end{equation}	
		and
		\item \label{ApproxOfScalarVectorProduct:LengthProp}
		it holds that $\lengthANN(\Phi)\le \tfrac{q}{(q-2)}[\LogBin(\eps^{-1})+q]$
%			\begin{equation}
%			\begin{split}
%			\lengthANN(\Phi)\le \tfrac{q}{(q-2)}[\LogBin(\eps^{-1})+q]
%			\end{split}
%			\end{equation}
	\end{enumerate}	
	(cf.\ Definition~\ref{Def:ANN}, Definition~\ref{Definition:ANNrealization}, and Definition~\ref{Def:euclideanNorm}).
\end{prop}

%					that $\paramANN(\Psi)\le \tfrac{360q}{(q-2)}\big[\LogBin(\eps^{-1})+q+\LogBin(3)-1\big]-252$, and 
%					\item \label{ApproxOfScalarVectorProduct:LengthLoad}
%					that $\lengthANN(\Psi)\le \tfrac{q}{(q-2)}[\LogBin(\eps^{-1})+q]$

\begin{proof}[Proof of Proposition~\ref{Lemma:ApproxOfScalarVectorProduct}]	
	Throughout this proof let 
	 $ v,w \in \R^{2 \times 1}$, $b\in \R^{2d}$, $A \in \R^{(2d) \times (d+1)}$  satisfy that
	\begin{equation}\label{ApproxOfScalarVectorProduct:Vectors}
	v = \begin{pmatrix}
	0\\
	1
	\end{pmatrix},\qquad 
	w = \begin{pmatrix}
	1\\
	0
	\end{pmatrix},
	 \qquad b =0,
	\end{equation}
	and
	\begin{equation}\label{ApproxOfScalarVectorProduct:MatrixA}
		A =
		\begin{pmatrix}
		w& v& 0& 0& \cdots& 0\\
		w& 0& v& 0&\cdots& 0\\
		w& 0& 0& v&\cdots& 0\\
		\vdots& \vdots&\vdots& \vdots& \ddots& \vdots\\
		w& 0& 0& 0&\cdots& v
		\end{pmatrix}, 
%		\in \R^{(2d) \times (d+1)},
	\end{equation}
%		\begin{equation}
%	w_1 = \begin{pmatrix}
%	1\\
%	-1
%	\end{pmatrix}
%	\in \R^{2 \times 1}, 
%	\quad
%	W_1 =
%	\begin{pmatrix}
%	w_1& 0& 0& \cdots& 0\\
%	0& w_1& 0& \cdots& 0\\
%	0& 0& w_1& \cdots& 0\\
%	\vdots& \vdots& \vdots& \ddots& \vdots\\
%	0& 0& 0& \cdots& w_1
%	\end{pmatrix} 
%	\in \R^{(2d) \times d}, \quad B_1 =0 \in \R^{2d},
%	\end{equation}
%	 let $\Psi\in \ANNs$ such that it holds that $\functionANN (\Psi)\in C(\R^2,\R)$ and for all $y,z\in\R$ it holds that 
%	 	\begin{equation}\label{ApproxOfScalarVectorProduct:scalarAnnihilation}
%	 (\functionANN (\Psi))(y,0)=(\functionANN (\Psi))(0,z)=0,
%	 \end{equation}
%	\begin{equation}\label{ApproxOfScalarVectorProduct:scalarEstimate}
%	\vert (\functionANN (\Psi))(y,z)\vert\le 3\big(1+\vert y\vert^2+\vert z\vert^2\big),
%	\qandqShort
%	\vert yz-(\functionANN (\Psi))(y,z)\vert\le \varepsilon \big(1+\vert y\vert^q+\vert z\vert^q\big),
%	\end{equation}
let $\Psi\in \ANNs$ satisfy that
	\begin{enumerate}[(I)]
		\item \label{ApproxOfScalarVectorProduct:RealizationLoad}
		it holds that   $\functionANN (\Psi)\in C(\R^2,\R)$,
		\item \label{ApproxOfScalarVectorProduct:ZeroLoad}
		it holds for all $x\in\R$ that  
		$[\functionANN (\Psi)](x,0)=[\functionANN (\Psi)](0,x)=0$,
%		\begin{equation}
%		[\functionANN (\Psi)](y,0)=[\functionANN (\Psi)](0,z)=0,
%		\end{equation}
		\item \label{ApproxOfScalarVectorProduct:EstimatesLoad}
		it holds for all  $x,y\in\R$ that 
				$\vert xy-[\functionANN (\Psi)](x,y)\vert\le \varepsilon \max\!\big\{1,\vert x\vert^q,\vert y\vert^q\big\}$,
		\item \label{ApproxOfScalarVectorProduct:EstimatesGrowthLoad}
		it holds for all  $x,y\in\R$ that $\vert[\functionANN (\Psi)](x,y)\vert
		\le 
		1+2x^2+2 y^2$,
		\item \label{ApproxOfScalarVectorProduct:ParamsLoad}
		it holds that $\paramANN(\Psi)\le \tfrac{360q}{(q-2)}\big[\LogBin(\eps^{-1})+q+1\big]-252$, and 
\item \label{ApproxOfScalarVectorProduct:LengthLoad}
		it holds that $\lengthANN(\Psi)\le \tfrac{q}{(q-2)}[\LogBin(\eps^{-1})+q]$
	\end{enumerate}	
	(cf.\ Proposition~\ref{Lemma:ApproxOfProduct}), and let $\affineMap\in (\R^{2d\times (d+1)}\times \R^{2d})\subseteq \ANNs$, $\Phi\in \ANNs$ 
	satisfy that 
	\begin{equation}\label{ApproxOfScalarVectorProduct:ConstructionNetwork}
		\affineMap=(A,b)\qandq \Phi=\compANN{\big[\parallelizationSpecial_d(\Psi,\Psi,\dots,\Psi)\big]}{\affineMap}
	\end{equation}
	(cf.\ Definition~\ref{Definition:ANNcomposition} and Definition~\ref{Definition:simpleParallelization}).	
	Observe that \eqref{ApproxOfScalarVectorProduct:Vectors} and \eqref{ApproxOfScalarVectorProduct:MatrixA} ensure that for all $y=(y_1,y_2,\dots, y_{d+1})\in\R^{d+1}$ it holds that 
	\begin{equation}
%		Ay=\Big(y_1 w+y_2 v,y_1 w+y_3 v,\dots, y_1 w+y_{d+1} v\Big)
Ay=\begin{pmatrix}
y_1 w+y_2 v\\ y_1 w+y_3 v \\\vdots\\ y_1 w+y_{d+1} v
\end{pmatrix}
=\begin{pmatrix}
y_1 \\ y_2 \\ y_1\\y_3  \\\vdots\\ y_1 \\y_{d+1}
\end{pmatrix}.
	\end{equation}
Combining this with   \eqref{ApproxOfScalarVectorProduct:ConstructionNetwork} proves that for all $t\in\R$, $x=(x_1,x_2,\dots,x_d)\in\R^d$ it holds that 
%	$\functionANN (\affineMap)\in C(\R^{d+1},\R^{2d})$ and 
	\begin{equation}
	\functionANN (\affineMap)\in C(\R^{d+1},\R^{2d})\qandq
			(\functionANN (\affineMap))(t,x)
			=(t,x_1,t,x_2,\dots,t,x_d)
			.
	\end{equation}
%	$(\functionANN (\affineMap))(t,x)
%	=(t,x_1,t,x_2,\dots,t,x_d)$.
	Proposition~\ref{Lemma:PropertiesOfParallelizationEqualLength}, \eqref{ApproxOfScalarVectorProduct:ConstructionNetwork}, and item~\eqref{PropertiesOfCompositions:Realization} in Proposition~\ref{Lemma:PropertiesOfCompositions} hence demonstrate that for all $t\in\R$, $x=(x_1,x_2,\dots,x_d)\in\R^d$ it holds that 
	 $\functionANN (\Phi)\in C(\R^{d+1},\R^d)$ and 
	 	\begin{equation}\label{ApproxOfScalarVectorProduct:TransformedANN}
	 	\begin{split}
	 	(\functionANN (\Phi))(t,x)
	 	&=\big(\big[\functionANN\big(\parallelizationSpecial_d(\Psi,\Psi,\dots,\Psi)\big)\big]\circ\big[\functionANN(\affineMap)\big]\big)(t,x)
	 	\\&=\big[\functionANN\big(\parallelizationSpecial_d(\Psi,\Psi,\dots,\Psi)\big)\big](t,x_1,t,x_2,\dots,t,x_d)	
	 	\\&=\big((\functionANN (\Psi))(t,x_1),(\functionANN (\Psi))(t,x_2),\dots,(\functionANN (\Psi))(t,x_d)\big).
	 	\end{split}
	 	\end{equation}	
%	\begin{equation}\label{ApproxOfScalarVectorProduct:TransformedANN}
%	\begin{split}
%		&(\functionANN (\Phi))(t,x_1,x_2,\dots, x_d)
%	\\&=\big(\big[\functionANN\big(\parallelizationSpecial_d(\Psi,\Psi,\dots,\Psi)\big)\big]\circ\big[\functionANN(\affineMap)\big]\big)(t,x_1,x_2,\dots, x_d)
%	\\&=\big[\functionANN\big(\parallelizationSpecial_d(\Psi,\Psi,\dots,\Psi)\big)\big](t,x_1,t,x_2,\dots,t,x_d)	
%	\\&=\big((\functionANN (\Psi))(t,x_1),(\functionANN (\Psi))(t,x_2),\dots,(\functionANN (\Psi))(t,x_d)\big).
%		\end{split}
%	\end{equation}	
	Combining this with \eqref{ApproxOfScalarVectorProduct:ZeroLoad} proves that for all $t\in\R$ it holds that
	\begin{equation}\label{ApproxOfScalarVectorProduct:ZeroAnnihilationOne}
\begin{split}
(\functionANN (\Phi))(t,0,0,\dots, 0)&=\big((\functionANN (\Psi))(t,0),(\functionANN (\Psi))(t,0),\dots,(\functionANN (\Psi))(t,0)\big)
\\&=(0,0,\dots, 0)=0.
\end{split}
\end{equation}	
Next note that \eqref{ApproxOfScalarVectorProduct:ZeroLoad} and \eqref{ApproxOfScalarVectorProduct:TransformedANN} imply that for all  $x=(x_1,x_2,\dots,x_d)\in\R^d$ it holds that
	\begin{equation}\label{ApproxOfScalarVectorProduct:ZeroAnnihilationTwo}
	\begin{split}
	(\functionANN (\Phi))(0,x)&=\big((\functionANN (\Psi))(0,x_1),(\functionANN (\Psi))(0,x_2),\dots,(\functionANN (\Psi))(0,x_d)\big)
	\\&=(0,0,\dots, 0)=0.
	\end{split}
	\end{equation}		
%In addition, observe that the triangle inequality, and the fact that for all $p\in[1,\infty)$, $r\in [p,\infty)$, $a_1,a_2,\dots, a_d\in\R$ it holds that 
%\begin{equation}
%	\Big[\smallsum\nolimits_{j=1}^d \vert a_j\vert^{r}\Big]^{1/r}\le \Big[\smallsum\nolimits_{j=1}^d \vert a_j\vert^{p}\Big]^{1/p}
%\end{equation}
%prove that for all $t\in\R$, $x=(x_1,x_2,\dots,x_d)\in\R^d$, $r\in [1,\infty)$ it holds that 
%	\begin{equation}\label{ApproxOfScalarVectorProduct:AuxiliarySums}
%		\begin{split}
%		&\Big[\smallsum\nolimits_{j=1}^d \big(1+\vert t\vert^r+\vert x_j\vert^r\big)^2\Big]^{\nicefrac{1}{2}}
%		\le \Big[\smallsum\nolimits_{j=1}^d \big(1+\vert t\vert^r\big)^2\Big]^{\nicefrac{1}{2}}
%		+\Big[\smallsum\nolimits_{j=1}^d \vert x_j\vert^{2r}\Big]^{\nicefrac{1}{2}}
%		\\&\le  \sqrt{d}\big(1+\vert t\vert^r\big)
%		+\Big[\smallsum\nolimits_{j=1}^d \vert x_j\vert^{2}\Big]^{\nicefrac{r}{2}}
%		=  \sqrt{d}\big(1+\vert t\vert^r\big)+\left\| x\right\|^r.
%		\end{split}
%	\end{equation}	
In addition, observe that the triangle inequality and the fact that for all $r\in [1,\infty)$,  $(x_1,x_2,\dots,x_d)\in\R^d$ it holds that 
\begin{equation}
\Big[\smallsum\nolimits_{j=1}^d \vert x_j\vert^{2r}\Big]^{\nicefrac{1}{2}}\le \Big[\smallsum\nolimits_{j=1}^d \vert x_j\vert^{2}\Big]^{\nicefrac{r}{2}}
\end{equation}
prove that for all $b\in\R$, $x=(x_1,x_2,\dots,x_d)\in\R^d$, $r\in [1,\infty)$ it holds that 
\begin{equation}\label{ApproxOfScalarVectorProduct:AuxiliarySums}
\begin{split}
\Big[\smallsum\nolimits_{j=1}^d \big(\vert b\vert+\vert x_j\vert^r\big)^2\Big]^{\nicefrac{1}{2}}
&\le \Big[\smallsum\nolimits_{j=1}^d b^2\Big]^{\nicefrac{1}{2}}
+\Big[\smallsum\nolimits_{j=1}^d \vert x_j\vert^{2r}\Big]^{\nicefrac{1}{2}}
\\&\le  \vert b\vert\sqrt{d}
+\Big[\smallsum\nolimits_{j=1}^d \vert x_j\vert^{2}\Big]^{\nicefrac{r}{2}}
=  \vert b\vert\sqrt{d}+\left\| x\right\|^r.
\end{split}
\end{equation}
This, \eqref{ApproxOfScalarVectorProduct:EstimatesLoad}, and \eqref{ApproxOfScalarVectorProduct:TransformedANN} assure that for all $t\in\R$, $x=(x_1,x_2,\dots,x_d)\in\R^d$ it holds that 
	\begin{equation}\label{ApproxOfScalarVectorProduct:AccuracyEstimate}
	\begin{split}
	&\left\|  t x-(\functionANN (\Phi))(t,x)\right\|
	=\Big[\smallsum\nolimits_{j=1}^d \vert tx_j-(\functionANN (\Psi))(t,x_j)\vert^2\Big]^{\nicefrac{1}{2}}
	\\&\le \Big[\smallsum\nolimits_{j=1}^d \big[\varepsilon\max\!\big\{1,\vert t\vert^q,\vert x_j\vert^q\big\}\big]^2\Big]^{\nicefrac{1}{2}}
	\le \varepsilon\Big[\smallsum\nolimits_{j=1}^d \big[\max\!\big\{1,\vert t\vert^q\big\}+\vert x_j\vert^q\big]^2\Big]^{\nicefrac{1}{2}}
		\\&\le \varepsilon \big(\sqrt{d}\left[\max\!\big\{1,\vert t\vert^q\big\}\right]+\left\| x\right\|^q\big).
	\end{split}
	\end{equation}	
	Furthermore, observe that \eqref{ApproxOfScalarVectorProduct:EstimatesGrowthLoad}, \eqref{ApproxOfScalarVectorProduct:TransformedANN}, and \eqref{ApproxOfScalarVectorProduct:AuxiliarySums} show that for all $t\in\R$, $x=(x_1,x_2,\dots,x_d)\in\R^d$ it holds that 
	\begin{equation}\label{ApproxOfScalarVectorProduct:GrowthEstimate}
	\begin{split}
	\left\| (\functionANN (\Phi))(t,x)\right\|
	&=\Big[\smallsum\nolimits_{j=1}^d \vert (\functionANN (\Psi))(t,x_j)\vert^2\Big]^{\nicefrac{1}{2}}
	\\&\le \Big[\smallsum\nolimits_{j=1}^d \big(1+2\vert t\vert^2+2\vert  x_j\vert^2\big)^2\Big]^{\nicefrac{1}{2}}
	\\&=\Big[\smallsum\nolimits_{j=1}^d \big(1+2\vert t\vert^2+\vert \sqrt{2} x_j\vert^2\big)^2\Big]^{\nicefrac{1}{2}}
	\\&\le  \sqrt{d}\big(1+2\vert t\vert^2\big)+\big\|  \sqrt{2} x\big\|^2
	=  \sqrt{d}\big(1+2\vert t\vert^2\big)+2\|   x\|^2.
	\end{split}
	\end{equation}
	In addition, note that Corollary~\ref{Lemma:ParallelizationImprovedBoundsOne} implies that
		\begin{equation}
		\paramANN\big(\parallelizationSpecial_d(\Psi,\Psi,\dots, \Psi)\big)\le d^2\, \paramANN(\Psi).
		\end{equation}
		Item~\eqref{PropertiesOfCompositionsWithAffineMaps:ItemBehind} in Corollary~\ref{Lemma:PropertiesOfCompositionsWithAffineMaps}, \eqref{ApproxOfScalarVectorProduct:ParamsLoad}, and \eqref{ApproxOfScalarVectorProduct:ConstructionNetwork} hence ensure that
		\begin{equation}\label{ApproxOfScalarVectorProduct:ParamEstimate}
		\begin{split}
		\paramANN(\Phi)&\le 
		\left[\max\!\left\{1,\tfrac{\inDimANN(\affineMap)+1}{\inDimANN(\parallelizationSpecial_d(\Psi,\Psi,\dots,\Psi))+1}\right\}\right]
		\paramANN\big(\parallelizationSpecial_d(\Psi,\Psi,\dots,\Psi)\big)
		\\&=
		\left[\max\{1,\tfrac{d+2}{2d+1}\}\right] \paramANN\big(\parallelizationSpecial_d(\Psi,\Psi,\dots,\Psi)\big)
		=\paramANN\big(\parallelizationSpecial_d(\Psi,\Psi,\dots,\Psi)\big)
		\\&\le  d^2\, \paramANN(\Psi)
		\le d^2 \big[\tfrac{360q}{(q-2)}\big]\big[\LogBin(\eps^{-1})+q+1\big]-252d^2.
		\end{split}
		\end{equation}
			Next note that item~\eqref{PropertiesOfCompositions:Length} in Proposition~\ref{Lemma:PropertiesOfCompositions}, \eqref{ApproxOfScalarVectorProduct:LengthLoad}, 
			and \eqref{ApproxOfScalarVectorProduct:ConstructionNetwork} demonstrate that 
			\begin{equation}\label{ApproxOfScalarVectorProduct:LengthEstimate}
			\begin{split}
			\lengthANN(\Phi)&=\lengthANN\big(\parallelizationSpecial_d(\Psi,\Psi,\dots, \Psi)\big)=\lengthANN(\Psi)
			\le \tfrac{q}{(q-2)}[\LogBin(\eps^{-1})+q].
			\end{split}
			\end{equation}
	This, the fact that $\functionANN (\Phi)\in C(\R^{d+1},\R^d)$, \eqref{ApproxOfScalarVectorProduct:ZeroAnnihilationOne}, \eqref{ApproxOfScalarVectorProduct:ZeroAnnihilationTwo}, \eqref{ApproxOfScalarVectorProduct:AccuracyEstimate}, \eqref{ApproxOfScalarVectorProduct:GrowthEstimate}, and \eqref{ApproxOfScalarVectorProduct:ParamEstimate} establish items~\eqref{ApproxOfScalarVectorProduct:RealizationProp}--\eqref{ApproxOfScalarVectorProduct:LengthProp}.
	The proof of Proposition~\ref{Lemma:ApproxOfScalarVectorProduct} is thus completed.
\end{proof}

%% file: SpaceTimeRepEuler.tex
\begin{lemma}\label{Lem:EulerCont}
	Let $N, d \in \N$,
	$\mu\in C(\R^d,\R^d)$,
	 $T\in (0,\infty)$, $(t_n)_{n\in\{-1,0,1,\dots,N+1\}}\allowbreak\subseteq\R$
	satisfy 
	 that $t_{-1}<0=t_0<t_1<\ldots <t_N=T<t_{N+1}$,
	 let $f_n\colon\R\to \R$, $n\in \{0,1,\dots, N\}$, be the functions which satisfy for all $n\in \{0,1,\dots, N\}$, $t\in\R$ that 
%		let $(f_n)_{n\in \{0,1,\dots, N\}}\subseteq C(\R,\R)$ satisfy for all $n\in \{0,1,\dots, N\}$, $t\in\R$ that 
		\begin{equation}\label{EulerCont:TimeNetworkMiddle}
		f_n(t)= \left[\tfrac{(t-t_{n-1})}{(t_n-t_{n-1})} \right]\indicator{(t_{n-1},t_n]}(t)+\left[\tfrac{(t_{n+1}-t)}{(t_{n+1}-t_n)} \right]\indicator{(t_{n},t_{n+1})}(t),
		\end{equation}
	and let 
	$Y= (Y^{x,y }_t)_{(t,x,y)\in [0,T]\times \R^d\times(\R^d)^N} \colon\allowbreak [0,T]\times \R^d\times (\R^d)^N \to \R^d $
	be the 
	function	
	which satisfies for all 
	$n\in\{0,1,\dots,N-1\}$,
	$ t \in [t_{n},t_{n+1}]$,
	$ x \in \R^d $, $y=(y_1,y_2,\dots, y_N)\in (\R^d)^N$
	that $\affineProcess^{x,y }_0=x$ and
	\begin{equation}
	\label{EulerCont:Y_processes}
	\begin{split}
	&\affineProcess^{x,y }_t 
	=
	\affineProcess^{x,y }_{t_n}+ \tfrac{(t-t_n)}{(t_{n+1}-t_n)}\,\left[(t_{n+1}-t_n)\,\mu \big( 
	\affineProcess^{x,y }_{t_n} 
	\big)
	+
	y_{n+1}\right]
	\end{split}
	\end{equation}
		(cf.\ Definition~\ref{Def:ANN} and Definition~\ref{Definition:ANNrealization}). Then  
		\begin{enumerate}[(i)]
			\item\label{ItemOne:EulerCont} it holds  that 
%			$([0,T]\times \R^d\ni (t,x)\mapsto\affineProcess^{x,y }_t\in \R^d)\in C([0,T]\times \R^d,\R^d)$ 
\begin{multline}
	\big([0,T]\times \R^d\times (\R^d)^N \ni (t,x,y)\mapsto\affineProcess^{x,y }_t\in \R^d\big)
	\\\in C([0,T]\times \R^d\times (\R^d)^N,\R^d)
\end{multline}
%\begin{equation}
%	([0,T]\times \R^d\times (\R^d)^N \ni (t,x,y)\mapsto\affineProcess^{x,y }_t\in \R^d)\in C([0,T]\times \R^d\times (\R^d)^N,\R^d) 
%\end{equation}
			and 
			\item\label{ItemTwo:EulerCont}  it holds for all $ t \in [0,T]$,
			$ x \in \R^d $, $y\in (\R^d)^N$ that
				\begin{equation}
			\begin{split}
			\affineProcess_t^{x,y }
			=\smallsum\limits_{n=0}^N f_n(t)\, Y_{t_n}^{x,y }.
			\end{split}
			\end{equation}
		\end{enumerate}
\end{lemma}

\begin{proof}[Proof of Lemma~\ref{Lem:EulerCont}]
	Observe that \eqref{EulerCont:Y_processes} 
	ensures that for all $n\in \{0,1,\dots, N-1\}$, $t\in [t_n,t_{n+1}]$, $x\in\R^d$, $y=(y_1,y_2,\dots, y_N)\in (\R^d)^N$ it holds that
	\begin{equation}
	\begin{split}
	&\affineProcess^{x,y }_{t_n}\big( \tfrac{t_{n+1}-t}{t_{n+1}-t_n}\big)+
	\affineProcess^{x,y }_{t_{n+1}} \big(\tfrac{t-t_n}{t_{n+1}-t_n}\big)
	\\&=\affineProcess^{x,y }_{t_n}\big(1- \tfrac{t-t_n}{t_{n+1}-t_n}\big)+
	\affineProcess^{x,y }_{t_{n+1}} \big(\tfrac{t-t_n}{t_{n+1}-t_n}\big)
	\\&= \affineProcess^{x,y }_{t_n}\big(1- \tfrac{t-t_n}{t_{n+1}-t_n}\big)
	\\&\quad+ \left(\affineProcess^{x,y }_{t_n}+ \tfrac{t_{n+1}-t_n}{t_{n+1}-t_n}\, \big[(t_{n+1}-t_n)\,\mu \big( 
	\affineProcess^{x,y }_{t_n} 
	\big)
	+
	y_{n+1}\big]\right)
	\big(\tfrac{t-t_n}{t_{n+1}-t_n}\big)
	\\&= \affineProcess^{x,y }_{t_n}+  \big[(t_{n+1}-t_n)\,\mu \big( 
	\affineProcess^{x,y }_{t_n} 
	\big)
	+
	y_{n+1}\big]
	\big(\tfrac{t-t_n}{t_{n+1}-t_n}\big)
	=Y_t^{x,y }.
	\end{split}
	\end{equation}
	Hence, we obtain that for all  $t\in [0,T]$, $x\in\R^d$, $y\in (\R^d)^N$ it holds that
	\begin{equation}
	\begin{split}
	&\affineProcess_t^{x,y }=Y_{t_0}^{x,y }\,\indicator{\{t_0\}}(t)
	+\smallsum\limits_{n=0}^{N-1}\big(\affineProcess_t^{x,y } \,\indicator{(t_n,t_{n+1}]}(t)\big)
	\\&=Y_{t_0}^{x,y }\,\indicator{\{t_0\}}(t)
	+\smallsum\limits_{n=0}^{N-1} \left(\left[\affineProcess^{x,y }_{t_n}\big( \tfrac{t_{n+1}-t}{t_{n+1}-t_n}\big)+
	\affineProcess^{x,y }_{t_{n+1}}\, \big(\tfrac{t-t_n}{t_{n+1}-t_n}\big) \right] \indicator{(t_n,t_{n+1}]}(t)\right)
	\\&=Y_{t_0}^{x,y }\,\indicator{\{t_0\}}(t)
	+\left[\smallsum\limits_{n=0}^{N-1}\affineProcess^{x,y }_{t_n}\big( \tfrac{t-t_{n+1}}{t_n-t_{n+1}}\big) \,\indicator{(t_n,t_{n+1}]}(t)\right]
	\\&\quad+\left[\smallsum\limits_{n=1}^{N}\affineProcess^{x,y }_{t_{n}} \,\big(\tfrac{t-t_{n-1}}{t_{n}-t_{n-1}}\big)  \,\indicator{(t_{n-1},t_{n}]}(t)\right]
	.
	\end{split}
	\end{equation}
	Combining this with \eqref{EulerCont:TimeNetworkMiddle}  implies that for all  $t\in [0,T]$, $x\in\R^d$, $y\in (\R^d)^N$ it holds that
	\begin{equation}\label{EulerCont:SumRepresentation}
	\begin{split}
	\affineProcess_t^{x,y }
	&=Y_{t_0}^{x,y }\,\indicator{\{t_0\}}(t) +
	\affineProcess^{x,y }_{t_0}\big( \tfrac{t-t_{1}}{t_0-t_{1}}\big) \,\indicator{(t_0,t_{1}]}(t)+ \affineProcess^{x,y }_{t_{N}} \big(\tfrac{t-t_{N-1}}{t_{N}-t_{N-1}}\big)  \,\indicator{(t_{N-1},t_{N}]}(t)
	\\&\quad+\smallsum\limits_{n=1}^{N-1}\affineProcess^{x,y }_{t_n}\left[\big( \tfrac{t-t_{n+1}}{t_n-t_{n+1}}\big) \,\indicator{(t_n,t_{n+1}]}(t)+\big(\tfrac{t-t_{n-1}}{t_{n}-t_{n-1}} \big) \,\indicator{(t_{n-1},t_{n}]}(t)\right]
	\\&=
	\affineProcess^{x,y }_{t_0}\big( \tfrac{t_{1}-t}{t_{1}-t_0}\big)\, \indicator{[t_0,t_{1}]}(t)+ \affineProcess^{x,y }_{t_{N}}\, f_N(t)
	+\smallsum\limits_{n=1}^{N-1}f_n(t)\, \affineProcess^{x,y }_{t_n} 
	\\&=\smallsum\limits_{n=0}^N f_n(t)\, Y_{t_n}^{x,y }.
	\end{split}
	\end{equation}
	Next	we claim that 	for all $n\in \{0,1,\dots, N\}$ it holds that 
		\begin{equation}\label{EulerCont:InductionHypo}
		 \big( \R^d\times (\R^d)^N \ni (x,y)\mapsto\affineProcess^{x,y }_{t_n}\in \R^d\big)\in C( \R^d\times (\R^d)^N,\R^d). 
		\end{equation}
		We now prove \eqref{EulerCont:InductionHypo} by induction on $n\in \{0,1,\dots, N\}$. 
		Note that the fact that for all $x\in\R^d$, $y\in (\R^d)^N$ it holds that $Y_{t_0}^{x,y}=Y_{0}^{x,y}=x$
		proves \eqref{EulerCont:InductionHypo} in the base case $n=0$.
		For the induction step assume there exists $n\in \{0,1,\dots, N-1\}$ which satisfies that
		\begin{equation}\label{EulerCont:InductionAssumption}
		\big( \R^d\times (\R^d)^N \ni (x,y)\mapsto\affineProcess^{x,y }_{t_n}\in \R^d\big)\in C( \R^d\times (\R^d)^N,\R^d). 
		\end{equation}
		Observe that \eqref{EulerCont:Y_processes} ensures that for all $x\in\R^d$, $y\in (\R^d)^N$ it holds that
		\begin{equation}
		\label{EulerCont:Y_processesInduction}
		\begin{split}
		&\affineProcess^{x,y }_{t_{n+1}} 
		=
		\affineProcess^{x,y }_{t_n}+ (t_{n+1}-t_n)\,\mu \big( 
		\affineProcess^{x,y }_{t_n} 
		\big)
		+
		y_{n+1}.
		\end{split}
		\end{equation}
		Combining this with \eqref{EulerCont:InductionAssumption} and the hypothesis that $\mu\in C(\R^d,\R^d)$ demonstrates that 
		\begin{equation}
		\big( \R^d\times (\R^d)^N \ni (x,y)\mapsto\affineProcess^{x,y }_{t_{n+1}}\in \R^d\big)\in C( \R^d\times (\R^d)^N,\R^d). 
		\end{equation}
		Induction thus proves \eqref{EulerCont:InductionHypo}.
%	Moreover, observe that item~\eqref{itemRepresentation:DNNrepresentationEulerSpaceFunction} in Proposition~\ref{Lemma:DNNrepresentationEulerSpace} ensures that 
%	for all $y\in (\R^d)^N$, $n\in \{0,1,\dots,N\}$ it holds that $( \R^d\ni x\mapsto\affineProcess^{x,y }_{t_n}\in \R^d)\in C(\R^d,\R^d)$.
	 Next observe that \eqref{EulerCont:SumRepresentation},  \eqref{EulerCont:InductionHypo}, and the fact that for all $n\in\{0,1,\dots, N\}$ it holds that $f_n\in C(\R,\R)$ show that
	 \begin{multline}
	 \big([0,T]\times \R^d\times (\R^d)^N \ni (t,x,y)\mapsto\affineProcess^{x,y }_t\in \R^d\big)
	 \\\in C([0,T]\times \R^d\times (\R^d)^N,\R^d).
	 \end{multline}
%	  for all $y\in (\R^d)^N$ it holds that $([0,T]\times \R^d\ni (t,x)\mapsto\affineProcess^{x,y }_t\in \R^d)\in C([0,T]\times \R^d,\R^d)$. 
	  Combining this with \eqref{EulerCont:SumRepresentation} establishes items~\eqref{ItemOne:EulerCont}--\eqref{ItemTwo:EulerCont}. 
	The proof of
	Lemma~\ref{Lem:EulerCont} is thus completed.
\end{proof}

%% file: ANNrepHatFct.tex
\begin{lemma}\label{Lemma:DNNhatFunction}
	Let $a\in C(\R,\R)$ satisfy for all $x\in\R$ that $a(x)=\max\{x,0\}$,
	let $\alpha,\beta,\gamma,h\in\R$ satisfy that $\alpha<\beta<\gamma$,
	let $W_1\in\R^{4\times 1}$, $B_1\in \R^4$, $W_2\in\R^{1\times 4}$, $B_2\in \R$ satisfy that 
	\begin{equation}
	W_1=\begin{pmatrix}
	\tfrac{1}{(\beta-\alpha)} \\[1ex] \tfrac{1}{(\beta-\alpha)}\\[1ex]\tfrac{1}{(\gamma-\beta)}\\[1ex]\tfrac{1}{(\gamma-\beta)}
	\end{pmatrix}
	,\qquad		B_1=\begin{pmatrix}
	-\tfrac{\alpha}{(\beta-\alpha)} \\[1ex] -\tfrac{\beta}{(\beta-\alpha)}\\[1ex]-\tfrac{\beta}{(\gamma-\beta)}\\[1ex]-\tfrac{\gamma}{(\gamma-\beta)}
	\end{pmatrix},
	\end{equation}
	\begin{equation}
	W_2=\begin{pmatrix}
	h& -h& -h& h
	\end{pmatrix}, \qquad
	B_2=0,
	\end{equation}
	and let $\Phi\in(\R^{4\times 1}\times \R^4)\times (\R^{1\times 4}\times \R)\subseteq\ANNs$ satisfy that $\Phi=((W_1,B_1),\allowbreak(W_2,B_2))$ (cf.\ Definition~\ref{Def:ANN}).
%	Then it holds  for all $t\in\R$ that $\functionANN(\Phi)\in C(\R,\R)$ and 
%	\begin{equation}\label{DNNrepresentationHatFunction}
%	(\functionANN(\Phi))(t)= \tfrac{(t-\alpha)h}{(\beta-\alpha)} \indicator{(\alpha,\beta]}(t)+\tfrac{(t-\gamma)h}{\beta-\gamma} \indicator{(\beta,\gamma)}(t)
%	\end{equation}
	Then
	\begin{enumerate}[(i)]
		\item it holds that $\functionANN(\Phi)\in C(\R,\R)$ and
		\item  it holds  for all $t\in\R$ that
%		\begin{equation}\label{DNNrepresentationHatFunction}
%		(\functionANN(\Phi))(t)= \tfrac{(t-\alpha)h}{(\beta-\alpha)} \indicator{(\alpha,\beta]}(t)+\tfrac{(t-\gamma)h}{\beta-\gamma} \indicator{(\beta,\gamma)}(t)
%		\end{equation}
				\begin{equation}\label{DNNrepresentationHatFunction}
				\begin{split}
				(\functionANN(\Phi))(t)&= \left[\tfrac{(t-\alpha)h}{(\beta-\alpha)}\right] \indicator{(\alpha,\beta]}(t)+\left[\tfrac{(\gamma-t)h}{(\gamma-\beta)} \right]\indicator{(\beta,\gamma)}(t)\\
				&=\begin{cases}
				0 &: t\in (-\infty,\alpha]\cup [\gamma,\infty)\\[1ex]
				\tfrac{(t-\alpha)h}{(\beta-\alpha)}&: t\in (\alpha,\beta]\\[1ex]
				\tfrac{(\gamma-t)h}{(\gamma-\beta)}&: t\in (\beta,\gamma)
				\end{cases}
				\end{split}
				\end{equation}
	\end{enumerate}	
	(cf.\ Definition~\ref{Definition:ANNrealization}).
\end{lemma}

\begin{proof}[Proof of Lemma~\ref{Lemma:DNNhatFunction}]	
	Observe that for all $t\in\R$ it holds that $\functionANN(\Phi)\in C(\R,\R)$ and 
	\begin{equation}
	\begin{split}
	&(\functionANN(\Phi))(t)= W_2 \big(\activationDim{4}(W_1t+B_1)\big)+B_2
	%\\&= \begin{pmatrix}
	%1&-1&-1& 1
	%\end{pmatrix} 
	%\begin{pmatrix}
	%\max\!\left\{\tfrac{t-a}{b-a},0\right\}\\
	%\max\!\left\{\tfrac{t-t_{n}}{b-a},0\right\}\\
	%\max\!\left\{\tfrac{t-t_{n}}{c-b},0\right\}\\
	%\max\!\left\{\tfrac{t-c}{c-b},0\right\}
	%\end{pmatrix}	
	\\&=h\max\!\left\{\tfrac{(t-\alpha)}{(\beta-\alpha)},0\right\}- h\max\!\left\{\tfrac{(t-\beta)}{(\beta-\alpha)},0\right\}
	-h\max\!\left\{\tfrac{(t-\beta)}{(\gamma-\beta)},0\right\}
	+h\max\!\left\{\tfrac{(t-\gamma)}{(\gamma-\beta)},0\right\}
	\\&=h[0-0- 0+0]\indicator{(-\infty,\alpha]}(t)
	+h\left[\tfrac{(t-\alpha)}{(\beta-\alpha)}-0-0+ 0\right] \indicator{(\alpha,\beta]}(t)
	\\&\quad+h\left[\tfrac{(t-\alpha)}{(\beta-\alpha)}-\tfrac{(t-\beta)}{(\beta-\alpha)}-\tfrac{(t-\beta)}{(\gamma-\beta)}+ 0\right] \indicator{(\beta,\gamma)}(t)
	\\&\quad+h\left[\tfrac{(t-\alpha)}{(\beta-\alpha)}-\tfrac{(t-\beta)}{(\beta-\alpha)}-\tfrac{(t-\beta)}{(\gamma-\beta)}+ \tfrac{(t-\gamma)}{(\gamma-\beta)}\right] \indicator{[\gamma,\infty)}(t)
	\\&=h\left[\tfrac{(t-\alpha)}{(\beta-\alpha)}\right] \indicator{(\alpha,\beta]}(t)
	+h\left[1-\tfrac{(t-\beta)}{(\gamma-\beta)}\right] \indicator{(\beta,\gamma)}(t)
	\\&=\left[\tfrac{(t-\alpha)h}{(\beta-\alpha)}\right] \indicator{(\alpha,\beta]}(t)+\left[\tfrac{(\gamma-t)h}{(\gamma-\beta)} \right]\indicator{(\beta,\gamma)}(t)
	\end{split}
	\end{equation}
	(cf.\ Definition~\ref{Def:multidim_version}).
	The proof of Lemma~\ref{Lemma:DNNhatFunction} is thus completed.
\end{proof}

%% file: PosterioriEstimates.tex
\begin{prop}\label{Cor:ApproxOfEuler}
	Let $N, d \in \N$,
	$a\in C(\R,\R)$ satisfy for all  $x\in \R$ that $\activation(x)=\max\{x,0\}$,
	let $T\in (0,\infty)$, $(t_n)_{n\in\{0,1,\dots,N\}}\allowbreak\subseteq\R$
	satisfy 
	for all $n\in\{0,1,\dots,N\}$ that $t_n=\tfrac{nT}{N}$,
%	let $p\in\R$ satisfy  that 
%	\begin{equation}\label{ApproxOfEuler:paramBound}
%	p=2 \Big(\tfrac{360q}{(q-2)}\big[\LogBin(\eps^{-1})+\LogBin(2^{q-1}+1)\big]+108\Big),
%	\end{equation}
		let $\mathfrak{D}\in[1,\infty)$, $\varepsilon\in (0,1]$, $q\in(2,\infty)$ satisfy that 
%		\begin{equation}\label{ApproxOfEuler:paramBound}
%		\mathfrak{D}=2 \Big(\tfrac{360q}{(q-2)}\big[\LogBin(\eps^{-1})+\LogBin(2^{q-1}+1)\big]+108\Big),
%		\end{equation}
				\begin{equation}\label{ApproxOfEuler:paramBound}
				\mathfrak{D}=\big[\tfrac{720q}{(q-2)}\big]\big[\LogBin(\eps^{-1})+q+1\big]-504,
				\end{equation}
%		let 
%		$
%		\Vert \cdot \Vert \colon \R^d \to [0,\infty)
%		$
%		be the Euclidean norm on $\R^d$,
		let 
		$\Phi\in \ANNs $
		satisfy 
		that $\inDimANN(\Phi)=\outDimANN(\Phi)=d$,
	%	let $b_n\in \R^d$, $n\in\{0,1,\dots,N-1\}$, 
%	and let 
%	$Y^n= \big(Y^{n,x,y }_t\big)_{t\in [t_n,t_{n+1}], x\in\R^d, y\in (\R^d)^{n+1}} \colon\allowbreak [t_n,t_{n+1}]\times \R^d\times (\R^d)^{n+1} \to \R^d $, $n\in\{0,1,\dots,N-1\}$, 
%	be the 
%	functions 	
%	which satisfy for all 
%	$n\in\{0,1,\dots,N-1\}$,
%		$ t \in [t_{n},t_{n+1}]$
%	$ x \in \R^d $, $y=(y_0,y_1,\dots, y_n)\in (\R^d)^{n+1}$
%	that $\affineProcess^{n,x,y }_0=x$ and
%	\begin{equation}
%	\label{ApproxOfEuler:Y_processes}
%	\begin{split}
%	&\affineProcess^{n,x,y }_t 
%	=
%	\affineProcess^{n,x,y }_{t_n}+ \left(\tfrac{tN}{T}-n\right)\left[\tfrac{T}{N}(\functionANN(\Phi)) \big( 
%	\affineProcess^{n,x,y }_{t_n} 
%	\big)
%	+
%	y_n\right]
%	\end{split}
%	\end{equation}
	and let 
	$Y= (Y^{x,y }_t)_{(t,x,y)\in [0,T]\times \R^d\times(\R^d)^N} \colon\allowbreak [0,T]\times \R^d\times (\R^d)^N \to \R^d $ 
	be the 
	function	
	which satisfies for all 
	$n\in\{0,1,\dots,N-1\}$,
	$ t \in [t_{n},t_{n+1}]$,
	$ x \in \R^d $, $y=(y_1,y_2,\dots, y_N)\in (\R^d)^N$
	that $\affineProcess^{x,y }_0=x$ and
	\begin{equation}
	\label{ApproxOfEuler:Y_processes}
	\begin{split}
	&\affineProcess^{x,y }_t 
	=
	\affineProcess^{x,y }_{t_n}+ \left(\tfrac{tN}{T}-n\right)\left[\tfrac{T}{N}(\functionANN(\Phi)) \big( 
	\affineProcess^{x,y }_{t_n} 
	\big)
	+
	y_{n+1}\right]
	\end{split}
	\end{equation}
	(cf.\ Definition~\ref{Def:ANN} and Definition~\ref{Definition:ANNrealization}).
%	Then there exist continuous functions $f_{\varepsilon,x,n,t}\colon  (\R^d)^{n+1}\to \R^d$, $ t \in [t_{n},t_{n+1}]$, $n\in\{0,1,\dots, N-1\}$, $x\in\R^d$, $\varepsilon\in (0,1]$,  and artificial neural networks $\Psi_{\varepsilon,y}\in \ANNs$, $\varepsilon\in(0,1]$, $y\in (\R^d)^N$, which satisfy
%	%	$(\Psi_{\varepsilon,y})_{\varepsilon\in(0,1]}\subseteq \ANNs$
	Then there exist $\Psi_{y}\in \ANNs$,  $y\in (\R^d)^N$, such that
	%	$(\Psi_{\varepsilon,y})_{\varepsilon\in(0,1]}\subseteq \ANNs$
	\begin{enumerate}[(i)]
		\item \label{ApproxOfEuler:Function} it holds for all  $y\in (\R^d)^N$ that $\functionANN (\Psi_{y})\in C(\R^{d+1},\R^d)$,
		\item it holds for all
		$n\in\{0,1,\dots,N-1\}$,
		$ t \in [t_{n},t_{n+1}]$,
		$x\in\R^d$,
		$y\in (\R^d)^N$ that 
		\begin{equation}
		\Vert \affineProcess^{x,y }_t -(\functionANN (\Psi_{y}))(t,x)\Vert
%		\le 4\varepsilon \sqrt{d} \left(1+\Vert \affineProcess^{x,y }_{t_n}\Vert^q+\Vert \affineProcess^{x,y }_{t_{n+1}}\Vert^q\right),
\le\varepsilon \big(2\sqrt{d}+\|  Y_{t_{n}}^{x,y }\|^q+\|  Y_{t_{n+1}}^{x,y }\|^q\big),
		\end{equation}
		\item it holds for all 
		$n\in\{0,1,\dots,N-1\}$,
		$ t \in [t_{n},t_{n+1}]$,
		 $x\in\R^d$, 
		 $y\in (\R^d)^N$
		that 
%		\begin{equation}
%		\Vert (\functionANN (\Psi_{y}))(t,x)\Vert
%		\le  12 \sqrt{d} \left(1+\Vert \affineProcess^{x,y }_{t_n}\Vert^2+\Vert \affineProcess^{x,y }_{t_{n+1}}\Vert^2\right),
%		\end{equation}
		\begin{equation}
		\Vert (\functionANN (\Psi_{y}))(t,x)\Vert
%		\le  12 \sqrt{d} \left(1+\Vert \affineProcess^{x,y }_{t_n}\Vert^2+\Vert \affineProcess^{x,y }_{t_{n+1}}\Vert^2\right),
\le 6\sqrt{d}+2\big(\Vert Y_{t_{n}}^{x,y }\Vert^2+\Vert Y_{t_{n+1}}^{x,y }\Vert^2\big),
		\end{equation}
		\item \label{ApproxOfEuler:ItemParams}
		it holds for all  $y\in (\R^d)^N$ that 
%		\begin{multline}
%		\paramANN(\Psi_{\varepsilon,y})
%		\\ \le\tfrac{1}{2} \bigg((2N+1)
%		\Big[d^2 \mathfrak{D}_\varepsilon+\big(
%		23
%		+6N(\lengthANN(\Phi)-1)
%		+  4d^2+3d+N \big[2(d^2+d)+\paramANN(\Phi)\big]^{\!2}\big)^{\!2}\Big]
%		\\-Nd(2d+1)+(\lengthANN(\Phi)-1) d(2d+1)N(N+1)
%		\bigg)^{\!2},
%		\end{multline}
				\begin{equation}
				\begin{split}
				&\paramANN(\Psi_{y})
				 \le\tfrac{1}{2} \bigg[  6d^2N^2 \hiddenLength(\Phi) 
				 \\&+3N
				\Big[d^2 \mathfrak{D}+\big(
				23
				+6N\hiddenLength(\Phi)
				+  7d^2+N \big[4d^2+\paramANN(\Phi)\big]^{\!2}\big)^{\!2}\Big]
				\bigg]^{\!2},
				\end{split}
				\end{equation}
%		\item \label{ApproxOfEuler:ItemMeasurability}
%		that for all $\varepsilon\in (0,1]$, $y=(y_1,y_2,\dots, y_N)\in (\R^d)^N$,  $n\in\{0,1,\dots, N-1\}$, $t\in [t_{n},t_{n+1}]$, $x\in\R^d$ it holds that 
%		\begin{equation}
%		(\functionANN (\Psi_{\varepsilon,y}))(t,x)=f_{\varepsilon,x,n,t}(y_0,y_1,\dots, y_n).
%		\end{equation}
		\item \label{ApproxOfEuler:Continuity}
		it holds  for all  $t\in [0,T]$, $x\in\R^d$ that 
		\begin{equation}
		\big[(\R^d)^N\ni y\mapsto (\functionANN (\Psi_{y}))(t,x)\in\R^d\big]\in C\big((\R^d)^N,\R^d\big),
		\end{equation}
		and
		\item \label{ApproxOfEuler:Adaptedness}
		it holds  for all    $n\in\{0,1,\dots, N\}$, $t\in [0,t_n]$, $x\in\R^d$, $y=(y_1,y_2,\dots, y_N),\allowbreak z=(z_1,z_2,\dots, z_N)\in (\R^d)^N$  with $\forall\, k\in \N\cap [0,n]\colon y_k=z_k$ that 
		\begin{equation}
			(\functionANN (\Psi_{y}))(t,x)=(\functionANN (\Psi_{z}))(t,x)
		\end{equation}
	\end{enumerate}
	(cf.\  Definition~\ref{Def:euclideanNorm}).
\end{prop}

\begin{proof}[Proof of Proposition~\ref{Cor:ApproxOfEuler}]	
	Throughout this proof
	let $t_n\in\R$, $n\in\{-1,N+1\}$,
	satisfy 
	for all $n\in\{-1,N+1\}$ that $t_n=\tfrac{nT}{N}$,
	let $(\mathbb{I}_{\mathfrak{d}})_{\mathfrak{d}\in\N}\subseteq \ANNs$ satisfy for all $\mathfrak{d}\in\N$, $x\in\R^\mathfrak{d}$ that $\functionANN(\mathbb{I}_\mathfrak{d})\in C(\R^\mathfrak{d},\R^\mathfrak{d})$, $\dims(\mathbb{I}_\mathfrak{d})=(\mathfrak{d},2\mathfrak{d},\mathfrak{d})$, and 
	\begin{equation}
		(\functionANN(\mathbb{I}_\mathfrak{d}))(x)=x
	\end{equation}
%	$(\functionANN(\mathbb{I}_\mathfrak{d}))(x)=x$ 
	(cf., e.g., \cite[Lemma~5.4]{Salimova2018}),
	let $(\Pi_n)_{n\in \{0,1,\dots, N\}}\subseteq \ANNs$ satisfy for all $n\in \{0,1,\dots, N\}$, $t\in\R$ that $\inDimANN(\Pi_n)=\outDimANN(\Pi_n)=1$, $\hiddenLength(\Pi_n)=1$, $\paramANN(\Pi_n)=13$, and
	%	\begin{equation}\label{ApproxOfEuler:TimeNetworkStart}
	%			(\functionANN(\rho_{1,0}))(t)=  \indicator{(-\infty,t_0]}(t)+\left(1-\tfrac{t}{t_{1}-t_{0}}\right) \indicator{(0,t_{1}]}(t),
	%	\end{equation}
	%		\begin{equation}\label{ApproxOfEuler:TimeNetworkEnd}
	%	(\functionANN(\rho_{1,N}))(t)= \tfrac{t-t_{N-1}}{t_N-t_{N-1}} \indicator{(t_{N-1},t_N]}(t)+ \indicator{(t_N,\infty)}(t),
	%	\end{equation}
	%	 and
	\begin{equation}\label{ApproxOfEuler:TimeNetworkMiddle}
	(\functionANN(\Pi_n))(t)= \left[\tfrac{(t-t_{n-1})}{(t_n-t_{n-1})}\right] \indicator{(t_{n-1},t_n]}(t)+\left[\tfrac{(t_{n+1}-t)}{(t_{n+1}-t_n)}\right] \indicator{(t_{n},t_{n+1})}(t)
	\end{equation}
	(cf.\ Lemma~\ref{Lemma:DNNhatFunction}), 
%	let $(\Xi_{n,y})_{ n\in \{0,1,\dots, N\}}\subseteq \ANNs$ 
%	satisfy 
%	\begin{enumerate}[(A)]
%		\item \label{ApproxOfEuler:ParamBoundsSpace} that for all  $n\in \{0,1,\dots, N\}$ it holds that 
%		\begin{equation}\label{ApproxOfEuler:PropSpatialApproxLength}
%		\lengthANN(\Xi_{n,y})=2+n(\lengthANN(\Phi)-1)
%		\end{equation}
%		\begin{equation}
%		\andq	\paramANN(\Xi_{n,y})\le \paramANN(\mathbb{I})+n \big[\tfrac{1}{2}\paramANN(\mathbb{I})+\paramANN(\Phi)\big]^{\!2}
%		\end{equation}
%		and
%		\item \label{ApproxOfEuler:RealizationSpatialNetwork}
%		 that 	for all  $n\in \{0,1,\dots, N\}$, $x\in\R^d$ that 	$\functionANN(\Xi_{n,y})\in C(\R^d,\R^d)$ and
%		$(\functionANN(\Xi_{n,y}))(x)=Y_{t_n}^x$
%	\end{enumerate}
%	(cf.\ Proposition~\ref{Lemma:DNNrepresentationEulerSpace} (with $a=a$, $d=d$, $\hiddenDimId=2d$, $\mathbb{I}=\mathbb{I}$, $\Phi=\Phi$,  $(\Omega, \mathcal{F}, \P )=( \Omega, \mathcal{F}, \P )$ $X_n=X_n \left(\tfrac{tN}{T}-n\right) \indicator{\{0,1,\dots, N-1\}}(n)$, $A_n=\left(\tfrac{tN}{T}-n\right) \tfrac{T}{N} \idMatrix_d \indicator{\{0,1,\dots, N-1\}}(n)$, $Y_n^x=Y_{t_n}^x \indicator{\{0,1,\dots, N-1\}}(n)+Y_{t_N}^x \indicator{\N\cap [N,\infty)}(n)$ for $x\in\R^d$,  $n\in\N_0$  in the notation of Proposition~\ref{Lemma:DNNrepresentationEulerSpace}),	
	let $(\Xi_{n,y})_{ (n,y)\in \{0,1,\dots, N\}\times(\R^d)^N}\subseteq \ANNs$ 
	satisfy that
		\begin{enumerate}[(I)]
			\item \label{ApproxOfEuler:SpatialFunction}
			it holds for all $n\in\{0,1,\dots,N\}$, $y\in(\R^d)^N$ that
			$\functionANN(\Xi_{n,y}) \in C(\R^d, \R^{d})$, 
			\item \label{ApproxOfEuler:SpatialFunctionValue} 	it holds for all $n\in\{0,1,\dots,N\}$, $y\in(\R^d)^N$, $x \in \R^d$  that $(\functionANN(\Xi_{n,y}))(x)\allowbreak=Y_{t_n}^{x,y}$,
			\item \label{ApproxOfEuler:SpatialHiddenLength}
			it holds for all $n\in\{0,1,\dots,N\}$, $y\in(\R^d)^N$ that $\hiddenLength(\Xi_{n,y})=1+n\hiddenLength(\Phi)$,
			\item \label{ApproxOfEuler:SpatialGeometry}	it holds for all $n\in\{0,1,\dots,N\}$, $y\in(\R^d)^N$ that
			\begin{equation}
				\paramANN(\Xi_{n,y})\le \paramANN(\mathbb{I}_d)+n \big[\tfrac{1}{2}\paramANN(\mathbb{I}_d)+\paramANN(\Phi)\big]^{\!2},
			\end{equation}
			\item \label{ApproxOfEuler:SpatialContinuity}
			it holds for all $n\in\{0,1,\dots,N\}$, $x\in\R^d$  that 
			\begin{equation}
			\big[(\R^d)^N\ni y\mapsto (\functionANN (\Xi_{n,y}))(x)\in\R^d\big]\in C\big((\R^d)^N,\R^d\big),
			\end{equation}
			and
			\item \label{ApproxOfEuler:SpatialAdaptedness}
			it holds
			for all $n\in\{0,1,\dots, N\}$, $m\in\N_0\cap [0,n]$,  $x\in\R^d$, $y=(y_1,y_2,\dots, y_N),\allowbreak z=(z_1,z_2,\dots, z_N)\in (\R^d)^N$  with $\forall\, k\in \N\cap [0,n]\colon y_k=z_k$ that 
			\begin{equation}
			(\functionANN(\Xi_{m,y}))(x)=(\functionANN(\Xi_{m,z}))(x)
			\end{equation}
		\end{enumerate}
	(cf.\ Proposition~\ref{Lemma:DNNrepresentationEulerSpace}), 
%	(with $a=a$, $N=N$, $d=d$, $\hiddenDimId=2d$, $\mathbb{I}=\mathbb{I}_d$, $\Phi=\Phi$,  $A_{n+1}=\left(\tfrac{tN}{T}-n\right) \tfrac{T}{N} \idMatrix_d$ for $n\in\{0,1,\dots,N-1\}$ (cf.\ Definition~\ref{Definition:identityMatrix}), $Y_n^x=Y_{t_n}^x$ for $x\in\R^d$,  $n\in\{0,1,\dots,N\}$  in the notation of Proposition~\ref{Lemma:DNNrepresentationEulerSpace}),		
	let $\Gamma\in \ANNs$ satisfy that
%	\begin{enumerate}[(a)]
%				\item \label{ApproxOfEuler:productLemmaParams} it holds	that
%				\begin{equation}
%				\begin{split}
%				\paramANN(\Gamma)&\le d^2\Big(\tfrac{360q}{(q-2)}\LogBin(\eps^{-1})+\tfrac{360q}{(q-2)}\LogBin(2^{q-1}+1)+108\Big),
%				\end{split}
%				\end{equation}
%		\item \label{ApproxOfEuler:annihilation} it holds for all $t\in\R$,	$x\in\R^{d}$ that $\functionANN (\Gamma)\in C(\R^{{d}+1},\R^{d})$ and
%		\begin{equation}\label{ApproxOfEulerAbstractProductNetworkCancellation}
%		(\functionANN (\Gamma))(t,0)=(\functionANN (\Gamma))(0,x)=0,
%		\end{equation}
%		and
%		\item it holds for all $t\in\R$,	$x\in\R^{d}$  that
%			\begin{equation}\label{ApproxOfEulerAbstractProductNetwork}
%			\Vert  t x-(\functionANN (\Gamma))(t,x)\Vert\le \varepsilon \sqrt{d} \big(1+\vert t\vert^q+\Vert x\Vert^q\big)
%			\end{equation}
%				\begin{equation}\label{ApproxOfEulerAbstractProductNetworkGrowth}
%				\andq\Vert(\functionANN (\Gamma))(t,x)\Vert\le 3 \sqrt{d} \big(1+\vert t\vert^2+\Vert x\Vert^2\big)
%				\end{equation}
%	\end{enumerate}
		\begin{enumerate}[(a)]
			\item \label{ApproxOfEuler:VectorProductRealizationProp}
			it holds that  $\functionANN (\Gamma)\in C(\R^{d+1},\R^d)$,
			\item \label{ApproxOfEuler:VectorProductZeroProp}
			it holds for all  $t\in\R$, $x\in\R^d$ that $(\functionANN (\Gamma))(t,0)=(\functionANN (\Gamma))(0,x)
			=0$,	
			\item \label{ApproxOfEuler:VectorProductEstimatesDifference}
			it holds for all  $t\in\R$, $x\in\R^d$ that 
			\begin{equation} 
			\left\|  t x-(\functionANN (\Gamma))(t,x)\right\|
			\le \varepsilon \big(\sqrt{d}\left[\max\!\big\{1,\vert t\vert^q\big\}\right]+\| x\|^q\big),
			\end{equation}
			\item \label{ApproxOfEuler:VectorProductEstimatesGrowth}
			it holds for all  $t\in\R$, $x\in\R^d$ that 
			\begin{equation}
			\left\| (\functionANN (\Gamma))(t,x)\right\|
			\le \sqrt{d}\big(1+2 t^2\big)+2\| x\|^2,
			%			\le 3 \sqrt{d}\big(1+\vert t\vert^2+\left\| x\right\|^2\!\big),
			\end{equation}
			\item \label{ApproxOfEuler:VectorProductParamsProp}
			it holds that $\paramANN(\Gamma)\le d^2 \big[\tfrac{360q}{(q-2)}\big]\big[\LogBin(\eps^{-1})+q+1\big]-252d^2$,
			and
			\item \label{ApproxOfEuler:VectorProductLengthProp}
			it holds that $\lengthANN(\Gamma)\le \tfrac{q}{(q-2)}[\LogBin(\eps^{-1})+q]$
		\end{enumerate}	
	(cf.\ Proposition~\ref{Lemma:ApproxOfScalarVectorProduct}),	
	let 
	$(\Psi_{n,y})_{(n,y)\in\{0,1,\dots,N\}\times(\R^d)^N}\subseteq \ANNs$ satisfy 
	%	that $\functionANN (\Psi_{n})\in C(\R^{d+1},\R^{d})$ and satisfy 
	%	for all $t\in\R$,
	%	%	 $x=(x_1,x_2,\dots,x_d)\in\R^d$ 
	%	$x\in\R^d$, 
	for all  $n\in\{0,1,\dots,N\}$, $y\in(\R^d)^N$
	 that $\inDimANN(\Psi_{n,y})=d+1$, $\outDimANN(\Psi_{n,y})=d$, and
	\begin{equation}\label{ApproxOfEulerCompProductNetwork}
	\Psi_{n,y}={\Gamma}\odot_{\mathbb{I}_{d+1}}\big[{\operatorname{P}_{2,(\mathbb{I}_1, \mathbb{I}_d)}(\Pi_n,\Xi_{n,y})}\big]
	\end{equation}
	(cf.\ Definition~\ref{Definition:ANNconcatenation}, Definition~\ref{Definition:generalParallelization}, Proposition~\ref{Lemma:PropertiesOfConcatenations}, and Corollary~\ref{Lemma:PropertiesOfParallelizationRealization}),
	let $L_y\in\N$, $y\in (\R^d)^N$, satisfy for all $y\in (\R^d)^N$ that $L_y=\max_{n\in \{0,1,\dots,N\}}\lengthANN(\Psi_{n,y})$, and let $(\Phi_y)_{y\in (\R^d)^N}\subseteq\ANNs$ satisfy that
%	that for all  $y\in\R^{d+1}$ it holds that
%	$\functionANN (\Phi_y)\in C(\R^{d+1},\R^{d})$ and
%	\begin{equation}\label{ApproxOfEulerSumNetwork}
%	\begin{split}
%	(\functionANN (\Phi_y))(y)=\smallsum\limits_{n=0}^N (\functionANN (\Psi_{n}))(y)
%	\end{split}
%	\end{equation}
%	(cf.\ Proposition~\ref{Lemma:SumsOfANNS}).
%	
	\begin{enumerate}[(A)]
		\item \label{ApproxOfEuler:SumContinuity} it holds for all $y\in (\R^d)^N$ that $\functionANN (\Phi_y)\in C(\R^{d+1},\R^{d})$,
		\item\label{ApproxOfEuler:SumRealization} it holds for all $y\in (\R^d)^N$,
		$z\in\R^{d+1}$ that 
		\begin{equation}\label{ApproxOfEuler:SumContinuityEquation}
		\begin{split}
		(\functionANN (\Phi_y))(z)=\smallsum\limits_{n=0}^N (\functionANN (\Psi_{n,y}))(z),
		\end{split}
		\end{equation}
%		\item  it holds for all $y\in (\R^d)^N$ that $L=\max_{n\in \{0,1,\dots,N\}}\lengthANN(\Psi_{n,y})$,
		and
		\item\label{ApproxOfEuler:SumParams}  it holds for all $y\in (\R^d)^N$ that 
%		\begin{equation}
%		\begin{split}
%							\paramANN(\Phi_y)
%							&\le \tfrac{1}{2} \bigg(\left[\smallsum\nolimits_{n=0}^N
%							2\,\paramANN(\Psi_{n,y})
%							\indicator{(\lengthANN(\Psi_{n,y}),\infty)}(L)\right]
%							\\&\qquad+\left[\smallsum\nolimits_{n=0}^N\big(
%							(L-\lengthANN(\Psi_{n,y})-1) \,2d(2d+1)
%							+d(2d+1)\big)
%							\indicator{(\lengthANN(\Psi_{n,y}),\infty)}(L)\right]
%							\\&\qquad+\left[\smallsum\nolimits_{n=0}^N  \paramANN(\Psi_{n,y})\indicator{\{\lengthANN(\Psi_{n,y})\}}(L)\right]\bigg)^{\!2}.
%		\end{split}
%		\end{equation}
		\begin{multline}\label{ApproxOfEuler:SumParamsTwo}
										\paramANN(\Phi_y)
										\le \tfrac{1}{2} \bigg[\left[\smallsum\nolimits_{n=0}^N
										2\,\paramANN(\Psi_{n,y})
										\indicator{(\lengthANN(\Psi_{n,y}),\infty)}(L_y)\right]
										\\+\left[\smallsum\nolimits_{n=0}^N\big(
										(L_y-\lengthANN(\Psi_{n,y})-1) \,2d(2d+1)
										+d(2d+1)\big)
										\indicator{(\lengthANN(\Psi_{n,y}),\infty)}(L_y)\right]
										\\+\left[\smallsum\nolimits_{n=0}^N  \paramANN(\Psi_{n,y})\indicator{\{\lengthANN(\Psi_{n,y})\}}(L_y)\right]\bigg]^{\!2}
		\end{multline}
	\end{enumerate}		
	(cf.\ Proposition~\ref{Lemma:SumsOfANNS}).	
	Note that \eqref{ApproxOfEuler:SpatialHiddenLength} and the fact that for all $n\in \{0,1,\dots,N\}$ it holds that $\hiddenLength(\Pi_n)=1$ ensure that for all $n\in \{0,1,\dots,N\}$, $y\in (\R^d)^N$ it holds that 
	$\mathcal{L}(\Xi_{n,y})=2+n\hiddenLength(\Phi)\geq 2$, $\mathcal{L}(\Pi_n)=2$, and 
	\begin{equation}
	\label{eq:R}
		\max\{\mathcal{L}(\Pi_n),\mathcal{L}(\Xi_{n,y})\}=\max\{2,2+n\hiddenLength(\Phi)\}=2+n\hiddenLength(\phi)=\mathcal{L}(\Xi_{n,y}).
	\end{equation}
	Corollary~\ref{Lemma:PropertiesOfParallelization} (with $a=a$, $n=2$, $L=\max\{\mathcal{L}(\Pi_n),\mathcal{L}(\Xi_{n,y})\}$, $\hiddenDimId_1=2$, $\hiddenDimId_2=2d$, $\Psi=(\mathbb{I}_1,\mathbb{I}_d)$, $\Phi=(\Pi_n,\Xi_{n,y})$ for $n\in\{0,1,\dots,N\}$, 
	$y\in(\R^d)^N$ in the notation of Corollary~\ref{Lemma:PropertiesOfParallelization}),
	 \eqref{ApproxOfEuler:SpatialGeometry}, and the fact that for all  $n\in \{0,1,\dots, N\}$ it holds that $\paramANN(\Pi_n)=13$ hence prove that for all  $n\in \{0,1,\dots, N\}$, $y\in (\R^d)^N$ it holds that 
	\begin{equation}\label{ApproxOfEuler:ParamsParallelization}
	\begin{split}
	&\paramANN\big(\!\operatorname{P}_{2,(\mathbb{I}_1, \mathbb{I}_d)}(\Pi_n,\Xi_{n,y})\big)
	\le \tfrac{1}{2} \big(
	2\paramANN(\Pi_n)
	+6(\lengthANN(\Xi_{n,y})-3) 
	+3
	+  \paramANN(\Xi_{n,y})\big)^{\!2}
	\\&= \tfrac{1}{2} \big(
	11
	+6\lengthANN(\Xi_{n,y}) 
	+  \paramANN(\Xi_{n,y})\big)^{\!2}
	\\&\le 
	\tfrac{1}{2} \big(
	11
	+6\big(2+n\hiddenLength(\Phi)\big)
	+  \paramANN(\mathbb{I}_d)+n \big[\tfrac{1}{2}\paramANN(\mathbb{I}_d)+\paramANN(\Phi)\big]^{\!2}\big)^{\!2}
	\\&=	\tfrac{1}{2} \big(
	23
	+6n\hiddenLength(\Phi)
	+  \paramANN(\mathbb{I}_d)+n \big[\tfrac{1}{2}\paramANN(\mathbb{I}_d)+\paramANN(\Phi)\big]^{\!2}\big)^{\!2}.
	\end{split}
	\end{equation}	
		Moreover, observe that \eqref{ApproxOfEuler:paramBound} and \eqref{ApproxOfEuler:VectorProductParamsProp} imply that $2\paramANN(\Gamma)\le d^2 \mathfrak{D}.$ 
		Combining this with Proposition~\ref{Lemma:PropertiesOfConcatenations}, \eqref{ApproxOfEuler:ParamsParallelization}, and the fact that $\paramANN(\mathbb{I}_d)=4d^2+3d\le 4(d^2+d)$ ensures that for all  $n\in \{0,1,\dots, N\}$, $y\in (\R^d)^N$ it holds that 
	\begin{equation}\label{ApproxOfEuler:ParamsEachSummand}
	\begin{split}
	\paramANN(\Psi_{n,y})&=\paramANN\big({\Gamma}\odot_{\mathbb{I}_{d+1}}[{\operatorname{P}_{2,(\mathbb{I}_1, \mathbb{I}_d)}(\Pi_n,\Xi_{n,y})}]\big)
	\\&\le \max\!\left\{1,\tfrac{2(d+1)}{{(d+1)}}\right\} \big(\paramANN({\Gamma})+\paramANN(\operatorname{P}_{2,(\mathbb{I}_1, \mathbb{I}_d)}(\Pi_n,\Xi_{n,y}))\big)
	\\&\le d^2 \mathfrak{D}+\big(
	23
	+6n\hiddenLength(\Phi)
	+  \paramANN(\mathbb{I}_d)+n \big[\tfrac{1}{2}\paramANN(\mathbb{I}_d)+\paramANN(\Phi)\big]^{\!2}\big)^{\!2}
		\\&\le d^2 \mathfrak{D}+\big(
		23
		+6n\hiddenLength(\Phi)
		+  4d^2+3d+n \big[2(d^2+d)+\paramANN(\Phi)\big]^{\!2}\big)^{\!2}.
	\end{split}
	\end{equation}
Next note that
	\eqref{ApproxOfEuler:SpatialHiddenLength}, \eqref{eq:R}, \eqref{ANNenlargement:Equation}, \eqref{generalParallelization:DefinitionFormula}, \eqref{parallelisationSameLengthDef}, item \eqref{PropertiesOfANNenlargementGeometry:ItemLonger} in Proposition~\ref{Lemma:PropertiesOfConcatenations}, and item \eqref{PropertiesOfANNenlargementGeometry:BulletPower} in Lemma~\ref{Lemma:PropertiesOfANNenlargementGeometry} demonstrate that for all  $n\in \{0,1,\dots, N\}$, $y\in (\R^d)^N$ it holds that 
	\begin{equation}\label{ApproxOfEuler:LengthEachSummand}
	\begin{split}
	\lengthANN(\Psi_{n,y})
	&=\lengthANN(\Gamma)+\lengthANN\big(\!\operatorname{P}_{2,(\mathbb{I}_1, \mathbb{I}_d)}(\Pi_n,\Xi_{n,y})\big)
	\\&=\lengthANN(\Gamma)+\lengthANN\big(\parallelizationSpecial_{2}(\longerANN{\max\{\lengthANN(\Pi_n),\lengthANN(\Xi_{n,y})\},\mathbb{I}_1}({\Pi_n}),\longerANN{\max\{\lengthANN(\Pi_n),\lengthANN(\Xi_{n,y})\},\mathbb{I}_d}({\Xi_{n,y}}))
\big)
	\\&=\lengthANN(\Gamma)+\lengthANN\big(\parallelizationSpecial_{2}(\longerANN{\lengthANN(\Xi_{n,y}),\mathbb{I}_1}({\Pi_n}),\longerANN{\lengthANN(\Xi_{n,y}),\mathbb{I}_d}({\Xi_{n,y}}))
\big)
	\\&=\lengthANN(\Gamma)+\lengthANN\big(\longerANN{\lengthANN(\Xi_{n,y}),\mathbb{I}_d}({\Xi_{n,y}})\big)
	\\&=\lengthANN(\Gamma)+\lengthANN\big(\compANN{\big((\mathbb{I}_d)^{\bullet 0}\big)}{\Xi_{n,y}}\big)
	\\&=\lengthANN(\Gamma)+\lengthANN\big((\mathbb{I}_d)^{\bullet 0}\big) +\lengthANN\big(\Xi_{n,y}\big)-1
	\\&=\lengthANN(\Gamma)+\lengthANN(\Xi_{n,y})
	= \lengthANN(\Gamma)+\hiddenLength(\Xi_{n,y})+1
	\\&= \lengthANN(\Gamma)+2+n\hiddenLength(\Phi).
%	=\lengthANN(\Gamma)+2+n(\lengthANN(\Phi)-1).
	\end{split}
	\end{equation}
Therefore, we obtain that
for all  $n\in \{0,1,\dots, N\}$, $y\in (\R^d)^N$ it holds that 
\begin{equation}\label{ApproxOfEuler:LengthDifferenceEachSummand}
\begin{split}
\lengthANN(\Psi_{N,y})-\lengthANN(\Psi_{n,y})-1
&=(\lengthANN(\Gamma)+2+N\hiddenLength(\Phi))-(\lengthANN(\Gamma)+2+n\hiddenLength(\Phi))-1
\\&=(N-n)\hiddenLength(\Phi)-1.
\end{split}
\end{equation}
In addition, note that \eqref{ApproxOfEuler:LengthEachSummand} proves that for all $y\in (\R^d)^N$ it holds that  $L_y=\lengthANN(\Psi_{N,y})=\lengthANN(\Psi_{N,0})=L_0$.
The fact that 
	$\smallsum\nolimits_{n=0}^{N-1}
	(N-n)=\smallsum\nolimits_{m=1}^{N}
	m=\tfrac{1}{2}N(N+1)$, \eqref{ApproxOfEuler:SumParamsTwo}, and \eqref{ApproxOfEuler:LengthDifferenceEachSummand} hence assure that for all $y\in(\R^d)^N$ it holds that
	
	\begin{equation}
	\begin{split}
	&\paramANN(\Phi_y)\\
	&\le \tfrac{1}{2} \bigg[\Big[\smallsum\nolimits_{n=0}^{N-1}\big(
	2\paramANN(\Psi_{n,y})\\
	&
	+\max\!\big\{\big(\lengthANN(\Psi_{N,y})-\lengthANN(\Psi_{n,y})-1\big) 2d(2d+1)
	+d(2d+1),0\big\}\big)
	\Big]
	+ \paramANN(\Psi_{N,y})\bigg]^{\!2}
	\\&=\tfrac{1}{2} \bigg[\left[\smallsum\nolimits_{n=0}^{N-1}\big(
	2\paramANN(\Psi_{n,y})
	+\max\!\big\{(N-n)\hiddenLength(\Phi) 2d(2d+1)
	-d(2d+1),0\big\}\big)
	\right]
	\\&+ \paramANN(\Psi_{N,y})\bigg]^{\!2}
	\\&\le\tfrac{1}{2} \bigg[(2N+1)
	\paramANN(\Psi_{N,y})\\
	&+\max\!\big\{\hiddenLength(\Phi) 2d(2d+1)\Big[\smallsum\nolimits_{n=0}^{N-1}
	(N-n)
	\Big]-Nd(2d+1),0\big\}
	\bigg]^{\!2}
	\\&=\tfrac{1}{2} \bigg[(2N+1)
	\paramANN(\Psi_{N,y})+\max\!\big\{\hiddenLength(\Phi) d(2d+1)N(N+1)-Nd(2d+1),0\big\}
	\bigg]^{\!2}.
	\end{split}
	\end{equation}
This and \eqref{ApproxOfEuler:ParamsEachSummand} imply that for all $y\in(\R^d)^N$ it holds that
	\begin{multline}
			\paramANN(\Phi_y)\le\tfrac{1}{2} \bigg[(2N+1)
			\Big[d^2 \mathfrak{D} +\big(
			23
			+6N\hiddenLength(\Phi)
			+  4d^2+3d+N \big[2(d^2+d)+\paramANN(\Phi)\big]^{\!2}\big)^{\!2}\Big]
			\\+\max\!\big\{\hiddenLength(\Phi) d(2d+1)N(N+1)-Nd(2d+1),0\big\}
			\bigg]^{\!2}.
	\end{multline}
	Therefore, we obtain that for all $y\in(\R^d)^N$ it holds that
		\begin{equation}\label{ApproxOfEuler:ParamProofFinal}
		\begin{split}
		&\paramANN(\Phi_y)
		\\&\le\tfrac{1}{2} \bigg[(2N+1)
		\Big[d^2 \mathfrak{D} +\big(
		23
		+6N\hiddenLength(\Phi)
		+  7d^2+N \big[4d^2+\paramANN(\Phi)\big]^{\!2}\big)^{\!2}\Big]
		\\&\qquad+\max\!\big\{\hiddenLength(\Phi) d(2d+1)N(N+1)-Nd(2d+1),0\big\}
		\bigg]^{\!2}
%		\\&\le\tfrac{1}{2} \bigg[3N
%		\Big[d^2 \mathfrak{D} +\big(
%		23
%		+6N\lengthANN(\Phi)
%		+  7d^2+N \big[4d^2+\paramANN(\Phi)\big]^{\!2}\big)^{\!2}\Big]
%		+6d^2 N^2\lengthANN(\Phi) 
%		\bigg]^{\!2}.
		\\&\le\tfrac{1}{2} \bigg[3N
		\Big[d^2 \mathfrak{D} +\big(
		23
		+6N\hiddenLength(\Phi)
		+  7d^2+N \big[4d^2+\paramANN(\Phi)\big]^{\!2}\big)^{\!2}\Big]
		+6d^2 N^2\hiddenLength(\Phi) 
		\bigg]^{\!2}.
		\end{split}
		\end{equation}
In addition, note that \eqref{ApproxOfEuler:TimeNetworkMiddle}, \eqref{ApproxOfEuler:SpatialFunctionValue}, and Lemma~\ref{Lem:EulerCont} (with $N=N$, $d=d$, $\mu=\functionANN(\Phi)$, $T=T$, $(\{-1,0,1,\dots,N+1\}\ni n\mapsto t_n\in\R)=(\{-1,0,1,\dots,N+1\}\ni n\mapsto t_n\in\R)$, $(\{0,1,\dots,N\}\ni n\mapsto f_n\in C(\R,\R))=(\{0,1,\dots,N\}\ni n\mapsto \functionANN(\Pi_n)\in C(\R,\R))$, $Y=Y$ in the notation of Lemma~\ref{Lem:EulerCont}) ensure that for all $t\in [0,T]$, $x\in\R^d$, $y\in (\R^d)^N$ it holds that 
		\begin{equation}\label{ApproxOfEuler:SumRepresentation}
		\begin{split}
		&\affineProcess_t^{x,y }
		=\smallsum\limits_{n=0}^N [(\functionANN(\Pi_n))(t)] \,Y_{t_n}^{x,y }
		=\smallsum\limits_{n=0}^N [(\functionANN(\Pi_n))(t)] \,[(\functionANN(\Xi_{n,y}))(x)].
		\end{split}
		\end{equation}
	Moreover, observe that 
	\eqref{ApproxOfEulerCompProductNetwork}, \eqref{ApproxOfEuler:SumContinuityEquation}, item \eqref{PropertiesOfConcatenations:Realization} in Proposition~\ref{Lemma:PropertiesOfConcatenations} (with $\Psi=\mathbb{I}_{d+1}$, 
	$\Phi_1=\Gamma$, $\Phi_2=\operatorname{P}_{2,(\mathbb{I}_1, \mathbb{I}_d)}(\Pi_n,\Xi_{n,y})$, $\mathfrak{i}=2(d+1)$ for $n\in \{0,1,\dots,N\}$, $y\in (\R^d)^N$ in the notation of Proposition~\ref{Lemma:PropertiesOfConcatenations}), and
	Corollary~\ref{Lemma:PropertiesOfParallelizationRealization} (with $a=a$, $n=2$, $\mathbb{I}=(\mathbb{I}_1,\mathbb{I}_d)$, $\Phi=(\Pi_n,\Xi_{n,y})$ for $n\in\{0,1,\dots,N\}$, $y\in(\R^d)^N$ in the notation of 
	Corollary~\ref{Lemma:PropertiesOfParallelizationRealization})
	demonstrate that for all $t\in [0,T]$, $x\in\R^d$, $y\in (\R^d)^N$ it holds that
	\begin{equation}\label{ApproxOfEuler:SumRepresentationPhi}
	(\functionANN (\Phi_y))(t,x)=\smallsum\limits_{n=0}^N (\functionANN (\Gamma))\big((
	\functionANN (\Pi_n))(t), (\functionANN (\Xi_{n,y}))(x)\big).
	\end{equation}		
	Next note that 
	%\eqref{ApproxOfEuler:TimeNetworkStart}, \eqref{ApproxOfEuler:TimeNetworkEnd}, 
	\eqref{ApproxOfEuler:TimeNetworkMiddle} shows that for all $k\in\{0,1,\dots,N\}$, $t\in \R\backslash (t_{k-1},t_{k+1})$ it holds that 
	\begin{equation}
	\label{eq:Z}
		(\functionANN(\Pi_{k}))(t)=0.
	\end{equation}
	Combining this, \eqref{ApproxOfEuler:SumRepresentation}, and \eqref{ApproxOfEuler:SumRepresentationPhi} with \eqref{ApproxOfEuler:VectorProductZeroProp} 	
	proves that for all $k\in\{0,1,\dots, N-1\}$, $t\in [t_{k},t_{k+1}]$, $x\in\R^d$, $y\in (\R^d)^N$ it holds that 
	\begin{equation}
	\begin{split}
	\affineProcess_t^{x,y }
	&=\smallsum\limits_{n=0}^N [(\functionANN(\Pi_n))(t)] \,[(\functionANN(\Xi_{n,y}))(x)]
	\\&=[(\functionANN(\Pi_{k}))(t)] \,[(\functionANN(\Xi_{k,y}))(x)]
	+[(\functionANN(\Pi_{k+1}))(t)]\, [(\functionANN(\Xi_{k+1,y}))(x)]
	\end{split}
	\end{equation}
	and 
	\begin{equation}\label{ApproxOfEuler:SumRepresentationReduced}
	\begin{split}
	(\functionANN (\Phi_y))(t,x)&=\smallsum\limits_{n=0}^N (\functionANN (\Gamma))\big((
	\functionANN (\Pi_n))(t), (\functionANN (\Xi_{n,y}))(x)\big)
	\\&=(\functionANN (\Gamma))\big((
	\functionANN (\Pi_{k}))(t), (\functionANN (\Xi_{k,y}))(x)\big)
	\\&\quad+ (\functionANN (\Gamma))\big((
	\functionANN (\Pi_{k+1}))(t), (\functionANN (\Xi_{k+1,y}))(x)\big).
	\end{split}
	\end{equation}
%	
%				\begin{equation} 
%				\left\|  t x-(\functionANN (\Gamma))(t,x)\right\|
%				\le \varepsilon \big(\sqrt{d}\left[\max\!\big\{1,\vert t\vert^q\big\}\right]+\| x\|^q\big),
%				\end{equation}
%				\item \label{ApproxOfEuler:VectorProductEstimatesGrowth}
%				it holds for all  $t\in\R$, $x\in\R^d$ that 
%				\begin{equation}
%				\left\| (\functionANN (\Gamma))(t,x)\right\|
%				\le \sqrt{d}\big(1+2 t^2\big)+2\| x\|^2,
%				%			\le 3 \sqrt{d}\big(1+\vert t\vert^2+\left\| x\right\|^2\!\big),
%				\end{equation}
%	
%	
	The triangle inequality, \eqref{ApproxOfEuler:VectorProductEstimatesDifference}, and \eqref{ApproxOfEuler:VectorProductEstimatesGrowth}
%	\eqref{ApproxOfEulerAbstractProductNetwork}, and \eqref{ApproxOfEulerAbstractProductNetworkGrowth}
	hence establish that for all $k\in\{0,1,\dots, N-1\}$, $t\in [t_{k},t_{k+1}]$, $x\in\R^d$, $y\in (\R^d)^N$ it holds that 
	\begin{equation}\label{ApproxOfEuler:PreFinalEstimateOne}
	\begin{split}
	&\Vert Y_t^{x,y } -(\functionANN (\Phi_y))(t,x)\Vert
	%	=
	%	\Vert \left[\smallsum\limits_{n=0}^N  \big[(\functionANN (\Pi_n)(t)\big]\, \big[(\functionANN (\Xi_{n,y}))(x)\big]\right] -(\functionANN (\Phi_y))(t,x)\Vert
	\\&\le \smallsum\limits_{n=k}^{k+1} \big\Vert  [(\functionANN (\Pi_n))(t)] \,[(\functionANN (\Xi_{n,y}))(x)] -(\functionANN(\Gamma))\big((\functionANN(\Pi_n))(t),
	(\functionANN(\Xi_{n,y}))(x)\big)\big\Vert
%	\\&\le \smallsum\limits_{n=k}^{k+1}  \varepsilon \sqrt{d} \big(1+\left\vert (\functionANN(\Pi_n))(t)\right\vert^q+\Vert (\functionANN(\Xi_{n,y}))(x)\Vert^q\big)
	\\&\le \smallsum\limits_{n=k}^{k+1}\varepsilon \big(\sqrt{d}\left[\max\!\big\{1,\vert (\functionANN(\Pi_n))(t)\vert^q\big\}\right]+\| (\functionANN(\Xi_{n,y}))(x)\|^q\big)
	\end{split}
	\end{equation}	
	and 
	\begin{equation}\label{ApproxOfEuler:PreFinalEstimateTwo}
	\begin{split}
	\Vert (\functionANN (\Phi_y))(t,x)\Vert
	&\le \smallsum\limits_{n=k}^{k+1} \big\Vert (\functionANN(\Gamma))\big((\functionANN(\Pi_n))(t),
	(\functionANN(\Xi_{n,y}))(x)\big)\big\Vert
%	\\&\le 3 \sqrt{d}
%	\left[\smallsum\limits_{n=k}^{k+1}  \big(1+\left\vert (\functionANN(\Pi_n))(t)\right\vert^2+\Vert (\functionANN(\Xi_{n,y}))(x)\Vert^2\big)\right]\!
%	.
	\\&\le \smallsum\limits_{n=k}^{k+1} \big(\sqrt{d}\big(1+2 \vert (\functionANN(\Pi_n))(t)\vert^2\big)+2\Vert (\functionANN(\Xi_{n,y}))(x)\Vert^2\big).
	\end{split}
	\end{equation}	
	Next note that \eqref{ApproxOfEuler:TimeNetworkMiddle} ensures that for all $n\in\{0,1,\dots, N\}$, $t\in\R$ it holds that  $0\le (\functionANN(\Pi_n))(t)\le 1$. Combining this with  \eqref{ApproxOfEuler:PreFinalEstimateOne}, \eqref{ApproxOfEuler:PreFinalEstimateTwo}, and 
%	the fact that for all $n\in \{0,1,\dots, N\}$, $x\in\R^d$, $y\in (\R^d)^N$ it holds that 
%	$(\functionANN(\Xi_{n,y}))(x)=Y_{t_n}^{x,y }$ 
	\eqref{ApproxOfEuler:SpatialFunctionValue}
	demonstrates that for all $k\in\{0,1,\dots, N-1\}$, $t\in [t_{k},t_{k+1}]$, $x\in\R^d$, $y\in (\R^d)^N$ it holds that 
%	\begin{equation}\label{ApproxOfEuler:FinalEstimateOne}
%	\begin{split}
%	\Vert Y_t^{x,y } -(\functionANN (\Phi_y))(t,x)\Vert
%	%		&\le \varepsilon \sqrt{d}
%	%		\smallsum_{n=0}^N   \big(1+\left\vert (\functionANN(\Pi_n))(t)\right\vert^q+\Vert (\functionANN(\Xi_{n,y}))(x)\Vert^q\big)
%	%		\\
%	&\le \varepsilon \sqrt{d}
%	\left[\smallsum\limits_{n=k}^{k+1}   \big(2+\Vert Y_{t_n}^{x,y }\Vert^q\big)\right]
%	\\&\le 4\varepsilon \sqrt{d} \left(1+\Vert Y_{t_k}^{x,y }\Vert^q+\Vert Y_{t_{k+1}}^{x,y }\Vert^q\right)
%	\end{split}
%	\end{equation}
%	and
%	\begin{equation}\label{ApproxOfEuler:FinalEstimateTwo}
%	\begin{split}
%	\Vert (\functionANN (\Phi_y))(t,x)\Vert
%	%	&\le 3 \sqrt{d}
%	%	\smallsum\limits_{n=0}^N   \big(1+\left\vert (\functionANN(\Pi_n))(t)\right\vert^2+\Vert (\functionANN(\Xi_{n,y}))(x)\Vert^2\big)
%	%	\\
%	&\le 3 \sqrt{d}
%	\left[\smallsum\limits_{n=k}^{k+1}   \big(2+\Vert Y_{t_n}^{x,y }\Vert^2\big)\right]\!
%	\\&\le  12 \sqrt{d} \left(1+\Vert Y_{t_k}^{x,y }\Vert^2+\Vert Y_{t_{k+1}}^{x,y }\Vert^2\right).
%	\end{split}
%	\end{equation}
		\begin{equation}\label{ApproxOfEuler:FinalEstimateOne}
		\begin{split}
		\Vert Y_t^{x,y } -(\functionANN (\Phi_y))(t,x)\Vert
		&\le \smallsum\limits_{n=k}^{k+1}\varepsilon \big(\sqrt{d}+\| (\functionANN(\Xi_{n,y}))(x)\|^q\big)
		\\&=\varepsilon \big(\sqrt{d}+\|  Y_{t_{k}}^{x,y }\|^q\big)+\varepsilon \big(\sqrt{d}+\|  Y_{t_{k+1}}^{x,y }\|^q\big)
		\\&=\varepsilon \big(2\sqrt{d}+\|  Y_{t_{k}}^{x,y }\|^q+\|  Y_{t_{k+1}}^{x,y }\|^q\big)
		\end{split}
		\end{equation}	
		and 
		\begin{equation}\label{ApproxOfEuler:FinalEstimateTwo}
		\begin{split}
		\Vert (\functionANN (\Phi_y))(t,x)\Vert
		&\le \smallsum\limits_{n=k}^{k+1} \big(3\sqrt{d}+2\Vert (\functionANN(\Xi_{n,y}))(x)\Vert^2\big)
		\\&= 3\sqrt{d}+2\Vert Y_{t_{k}}^{x,y }\Vert^2+3\sqrt{d}+2\Vert Y_{t_{k+1}}^{x,y }\Vert^2
		\\&= 6\sqrt{d}+2\big(\Vert Y_{t_{k}}^{x,y }\Vert^2+\Vert Y_{t_{k+1}}^{x,y }\Vert^2\big).
		\end{split}
		\end{equation}
	Furthermore, observe that \eqref{ApproxOfEuler:SumRepresentationPhi}, \eqref{ApproxOfEuler:SpatialContinuity}, and \eqref{ApproxOfEuler:VectorProductRealizationProp} ensure 
	that for all $t\in [0,T]$,  $x\in\R^d$ it holds that 
	\begin{equation}\label{ApproxOfEuler:ContinuityProof}
	\big[(\R^d)^N\ni y\mapsto (\functionANN (\Phi_{y}))(t,x)\in\R^d\big]\in C\big((\R^d)^N,\R^d\big).
	\end{equation}
	In addition, observe that \eqref{ApproxOfEuler:VectorProductZeroProp}, \eqref{ApproxOfEuler:SumRepresentationPhi} and \eqref{eq:Z} demonstrate	
	that for all $n\in\{0,1,\dots, N\}$, $t\in [0,t_{n}]$, $x\in\R^d$, $y\in (\R^d)^N$ it holds that
		\begin{equation}
		\begin{split}
		(\functionANN (\Phi_y))(t,x)&=\smallsum\limits_{k=0}^n (\functionANN (\Gamma))\big((
		\functionANN (\Pi_k))(t), (\functionANN (\Xi_{k,y}))(x)\big).
		\end{split}
		\end{equation}
	This and 
	 \eqref{ApproxOfEuler:SpatialAdaptedness} show that 
					for all $n\in\{0,1,\dots, N\}$, $t\in [0,t_{n}]$,  $x\in\R^d$, $y=(y_1,y_2,\dots, y_N)$, $z=(z_1,z_2,\dots, z_N)\in (\R^d)^N$  with $\forall\, k\in \N\cap [0,n]\colon y_k=z_k$ it holds that 
					\begin{equation}
					\begin{split}
					(\functionANN(\Phi_{y}))(t,x)&=\smallsum\limits_{m=0}^n(\functionANN(\Gamma))((\functionANN(\Pi_m))(t),(\functionANN(\Xi_{m,y}))(x))\\
					&=\smallsum\limits_{m=0}^n(\functionANN(\Gamma))((\functionANN(\Pi_m))(t),(\functionANN(\Xi_{m,z}))(x))\\
					&=(\functionANN(\Phi_{z}))(t,x).
					\end{split}
					\end{equation}
	Combining this with \eqref{ApproxOfEuler:SumContinuity}, \eqref{ApproxOfEuler:ParamProofFinal}, \eqref{ApproxOfEuler:FinalEstimateOne},
	\eqref{ApproxOfEuler:FinalEstimateTwo}, and \eqref{ApproxOfEuler:ContinuityProof} establishes
	items~\eqref{ApproxOfEuler:Function}--\eqref{ApproxOfEuler:Adaptedness}.
	The proof of
	Proposition~\ref{Cor:ApproxOfEuler} is thus completed.
\end{proof}

%% file: PrioriEstimatesEuler.tex
\begin{lemma}\label{Lemma:Gronwall}
	Let $N, d \in \N$, $c,C\in [0,\infty)$, $A_1,A_2,\dots,A_N\in\R^{d\times d}$,
	let 
	$
	\NORm{\cdot} \colon \R^d \to [0,\infty)
	$
	be a norm on $\R^d$, let $\tripleNorm{\cdot}\colon \R^{d\times d}\to [0,\infty)$ 
	be the function which satisfies for all $A\in\R^{d\times d}$ that $\tripleNorm{A}=\sup_{\{x\in\R^d\colon\NOrm{x} \le 1\}}\NOrm{Ax}$,
	let $\mu\colon \R^d\to \R^d$ be a function which satisfies for all $x\in\R^d$ that 
	%	$\Vert \mu(x)\Vert_{\R^d} \le C(1+\Vert x\Vert_{\R^d} )$,
	\begin{equation}\label{Gronwall:linearGrowth}
	\NOrm{\mu(x)} \le C+c\NOrm{x},
	\end{equation}
	and 
	%	let 
	%	$Y= (Y^{x,y }_t)_{(t,x,y)\in [0,T]\times \R^d\times(\R^d)^N} \colon\allowbreak [0,T]\times \R^d\times (\R^d)^N \to \R^d $, 
	%	be the 
	%	function	
	%	which satisfies for all 
	%	$n\in\{0,1,\dots,N-1\}$,
	%	$ t \in [t_{n},t_{n+1}]$,
	%	$ x \in \R^d $, $y=(y_1,y_2,\dots, y_N)\in (\R^d)^N$
	%	that $\affineProcess^{x,y }_0=x$ and
	%	\begin{equation}
	%	\label{ApproxOfEuler:Y_processes}
	%	\begin{split}
	%	&\affineProcess^{x,y }_t 
	%	=
	%	\affineProcess^{x,y }_{t_n}+ \left(\tfrac{tN}{T}-n\right)\left[\tfrac{T}{N}\,\mu\big( 
	%	\affineProcess^{x,y }_{t_n} 
	%	\big)
	%	+
	%	y_{n+1}\right].
	%	\end{split}
	%	\end{equation}
	let 
	$Y_n= (Y^{x,y }_n)_{ (x,y)\in\R^d\times (\R^d)^N} \colon\allowbreak  \R^d\times(\R^d)^N \to \R^d $, $n\in\{0,1,\dots, N\}$,
	be the
	functions 	
	which satisfy for all 
	$n\in\{0,1,\dots, N-1\}$,
	$ x \in \R^d $,
	$y=(y_1,y_2,\dots, y_N)\in (\R^d)^N$  
	that $\affineProcess^{x,y} _0=x$ and
	\begin{equation}\label{Gronwall:Recursion}
	\begin{split}
	&\affineProcess^{x,y} _{n+1} 
	=
	\affineProcess^{x,y} _{n}+ A_{n+1}\,\mu\big( 
	\affineProcess^{x,y} _{n} 
	\big)
	+
	y_{n+1}.
	\end{split}
	\end{equation}
	Then 
	\begin{enumerate}[(i)]
		\item \label{Gronwall:ItemOne}it holds for all 
		$n\in\{0,1,\dots, N\}$,
		$ x \in \R^d $, $y=(y_1,y_2,\dots,, y_N)\in (\R^d)^N$ that
		\begin{equation}
		\begin{split}
		&\affineProcess^{x,y} _{n} 
		=
		x+\smallsum\limits_{k=0}^{n-1} \left[A_{k+1}\,\mu\big( 
		\affineProcess^{x,y} _{k} 
		\big)
		+
		y_{k+1}\right]
		\end{split}
		\end{equation}
		and
		\item\label{Gronwall:ItemTwo}
		it holds for all $n\in\{0,1,\dots, N\}$,
		$ x \in \R^d $, $y=(y_1,y_2,\dots, y_N)\in (\R^d)^N$ that 
		%		\begin{equation}
		%					\begin{split}
		%					&\Vert \affineProcess^{x,y} _{n}\Vert 
		%					\le 
		%					\Big(\Vert x\Vert+C \Big[\smallsum\limits_{k=1}^{n}\tripleNorm{{A_k}}\Big]+\Big\Vert\Big[\smallsum\limits_{k=1}^{n}y_k\Big]\Big\Vert\Big) \exp\!\Big(C \Big[\smallsum\limits_{k=1}^{n}\tripleNorm{{A_k}}\Big]\Big).
		%					\end{split}
		%		\end{equation}
		\begin{equation}\label{Gronwall:EstimateItemTwo}
		\begin{split}
		&\NOrm{\affineProcess^{x,y} _{n}}  \\
		& \le\left(\NOrm{x}  +C \left[\smallsum\limits_{k=1}^{n} \tripleNorm{{A_{k}}}\right]+
		\max_{m\in\{0,1,\dots,n\}}\NOrmmmm{\smallsum\limits_{k=1}^{m} y_{k}} \right)\exp\!\left(\!c \!\left[\smallsum\limits_{k=1}^{n} \tripleNorm{{A_{k}}}\right]\!\right)
		\!.
		\end{split}
		\end{equation}
	\end{enumerate}	 
\end{lemma}

\begin{proof}[Proof of Lemma~\ref{Lemma:Gronwall}]	
	We claim that 
	for all 
	$n\in\{0,1,\dots, N\}$,
	$ x \in \R^d $, $y=(y_1,y_2,\allowbreak\dots, y_N)\in (\R^d)^N$ it holds that
	\begin{equation}\label{Gronwall:SumInduction}
	\begin{split}
	&\affineProcess^{x,y} _{n} 
	=
	x+\smallsum\limits_{k=0}^{n-1} \left[A_{k+1}\,\mu\big( 
	\affineProcess^{x,y} _{k} 
	\big)
	+
	y_{k+1}\right].
	\end{split}
	\end{equation}
	We now prove \eqref{Gronwall:SumInduction} by induction on $n\in\{0,1,\dots,N\}$. Observe that the hypothesis that for all $ x \in \R^d $, $y\in (\R^d)^N$ it holds that  $\affineProcess^{x,y} _0=x$ proves \eqref{Gronwall:SumInduction} in the base case $n=0$. For the induction step note that \eqref{Gronwall:Recursion} implies that  for all $n\in\{0,1,\dots,N-1\}$,
	$ x \in \R^d $, $y=(y_1,y_2,\allowbreak\dots, y_N)\in (\R^d)^N$ with
	\begin{equation}\label{Gronwall:SumInductionHypo}
	\begin{split}
	&\affineProcess^{x,y} _{n} 
	=
	x+\smallsum\limits_{k=0}^{n-1} \left[A_{k+1}\,\mu\big( 
	\affineProcess^{x,y} _{k} 
	\big)
	+
	y_{k+1}\right]
	\end{split}
	\end{equation}
	it holds that 
	\begin{equation}\label{Gronwall:SumInductionStep}
	\begin{split}
	\affineProcess^{x,y} _{n+1} 
	&=
	\affineProcess^{x,y} _{n}+ A_{n+1}\,\mu\big( 
	\affineProcess^{x,y} _{n} 
	\big)
	+
	y_{n+1}
	\\&=
	x+\left[\smallsum\limits_{k=0}^{n-1} \left(A_{k+1}\,\mu\big( 
	\affineProcess^{x,y} _{k} 
	\big)
	+
	y_{k+1}\right)\right]
	+ \left(A_{n+1}\,\mu\big( 
	\affineProcess^{x,y} _{n} 
	\big)
	+
	y_{n+1}\right)
	\\&=
	x+\left[\smallsum\limits_{k=0}^{n} \left(A_{k+1}\,\mu\big( 
	\affineProcess^{x,y} _{k} 
	\big)
	+
	y_{k+1}\right)\right].
	\end{split}
	\end{equation}
	Induction thus proves \eqref{Gronwall:SumInduction}. Observe that \eqref{Gronwall:SumInduction} establishes item~\eqref{Gronwall:ItemOne}.
	In addition, note that \eqref{Gronwall:SumInduction}, the triangle inequality, 
%	and the hypothesis that for all $A\in\R^{d\times d}$ that $\tripleNorm{A}=\sup_{x\in\R^d\colon\norm{x}\le 1}\norm{Ax}$ 
	and the fact that for all $A\in\R^{d\times d}$, $x\in\R^d$ it holds that $\NOrm{Ax}\le \tripleNorm{A}\,\NOrm{x} $ 
	demonstrate that for all 
	$m\in\{0,1,\dots,N\}$,
	$ x \in \R^d $, $y=(y_1,y_2,\allowbreak\dots, y_N)\in (\R^d)^N$ it holds that
	\begin{equation}
	\begin{split}
	&\NOrm{\affineProcess^{x,y} _{m}}
	\le
	\NOrm{x}  +\left[\smallsum\limits_{k=0}^{m-1} \tripleNorm{{A_{k+1}}}\, \NOrmm{\mu\big( 
	\affineProcess^{x,y} _{k} 
	\big)} \right]
	+
	\NOrmmmm{ \smallsum\limits_{k=0}^{m-1} y_{k+1}}.
	\end{split}
	\end{equation}
	Combining this with \eqref{Gronwall:linearGrowth} ensures that 
	for all 
	$n\in\{0,1,\dots,N\}$, $m\in\{0,1,\dots,n\}$,
	$ x \in \R^d $, $y=(y_1,y_2,\allowbreak\dots, y_N)\in (\R^d)^N$ it holds that
	\begin{equation}\label{Gronwall:TriangleSum}
	\begin{split}
	&\NOrm{\affineProcess^{x,y} _{m}} 
	\le
	\NOrm{x}  + \left[\smallsum\limits_{k=0}^{m-1} \tripleNorm{{A_{k+1}}}\,\big(C+c\NOrmm{
	\affineProcess^{x,y} _{k}} \big)\right]
	+
	\NOrmmmm{ \smallsum\limits_{k=1}^{m} y_{k}} 
	\\&=\NOrm{x}  +C \left[\smallsum\limits_{k=1}^{m} \tripleNorm{{A_{k}}}\right]+
	\NOrmmmm{ \smallsum\limits_{k=1}^{m} y_{k}} 
	+c \left[\smallsum\limits_{k=0}^{m-1} \tripleNorm{{A_{k+1}}}\,\NOrmm{ 
	\affineProcess^{x,y} _{k} } \right]
	\\&\le\NOrm{x}  +C \left[\smallsum\limits_{k=1}^{n} \tripleNorm{{A_{k}}}\right]+
	\bigg[\max_{m\in\{0,1,\dots,n\}}\NOrmmmm{\smallsum\limits_{k=1}^{m} y_{k}} \bigg]
	+c \left[\smallsum\limits_{k=0}^{m-1} \tripleNorm{{A_{k+1}}}\,\NOrmm{
	\affineProcess^{x,y} _{k} } \right]
	\!.
	\end{split}
	\end{equation}
	The time-discrete Gronwall inequality 
	(cf., e.g., Hutzenthaler et al.~\!\cite[Lemma 2.1]{PDEapproximation}
	(with $N=n$, 
	$\alpha=\big(\NOrm{x}  +C \big[\sum_{k=1}^{n} \tripleNorm{{A_{k}}}\big]+
	\max_{m\in\{0,1,\dots,n\}}\NOrm{ \sum_{k=1}^{m} y_{k}} \big)$,
	$\beta_0=c\tripleNorm{{A_1}},\beta_1=c\tripleNorm{{A_2}},\dots,\beta_{n-1}=c\tripleNorm{{A_n}}, \epsilon_0=\NOrm{\affineProcess^{x,y} _0}, \epsilon_1=\NOrm{\affineProcess^{x,y} _1},\dots,\epsilon_n=\NOrm{\affineProcess^{x,y} _n}$ 
	for $n\in \{1,2,\dots,N\}$ in the notation of Hutzenthaler et al.~\cite[Lemma 2.1]{PDEapproximation})) 
	hence implies that for all 
	$n\in\{1,2,\dots,N\}$,
	$ x \in \R^d $, $y=(y_1,y_2,\allowbreak\dots, y_N)\in (\R^d)^N$ it holds that
	\begin{equation}
	\begin{split}
	\NOrm{\affineProcess^{x,y} _{n}} 
	&\le \left(\NOrm{x}  +C \left[\smallsum\limits_{k=1}^{n} \tripleNorm{{A_{k}}}\right]+
	\max_{m\in\{0,1,\dots,n\}} \NOrmmm{ \smallsum\limits_{k=1}^{m} y_{k}} \right)
	\exp\!\left(\!c \!\left[\smallsum\limits_{k=0}^{n-1} \tripleNorm{{A_{k+1}}}\right]\!\right)
	\!.
	\end{split}
	\end{equation}
	The hypothesis that for all $ x \in \R^d $,
	$y\in (\R^d)^N$ it holds that 
	$\affineProcess^{x,y} _0=x$ therefore assures that for all 
	$n\in\{0,1,\dots,N\}$,
	$ x \in \R^d $, $y=(y_1,y_2,\allowbreak\dots, y_N)\in (\R^d)^N$ it holds that
	\begin{equation}\label{Gronwall:EndEstimate}
	\begin{split}
	\NOrm{\affineProcess^{x,y} _{n} } 
	&\le \left(\NOrm{ x}  +C \left[\smallsum\limits_{k=1}^{n} \tripleNorm{{A_{k}}}\right]+
	\max_{m\in\{0,1,\dots,n\}}\NOrmmm{ \smallsum\limits_{k=1}^{m} y_{k}} \right)
	\exp\!\left(\!c \!\left[\smallsum\limits_{k=1}^{n} \tripleNorm{{A_{k}}}\right]\!\right)
	\!.
	\end{split}
	\end{equation}
	This establishes item~\eqref{Gronwall:ItemTwo}.
	The proof of Lemma~\ref{Lemma:Gronwall} is thus completed.
\end{proof}

%% file: PrioriEstimatesANN.tex
\begin{theorem}\label{Thm:ApproxOfEulerWithGronwall}
	Let $N, d \in \N$, $\mathfrak{C}\in [0,\infty)$, 
	$a\in C(\R,\R)$ satisfy for all  $x\in \R$ that $\activation(x)=\max\{x,0\}$,
	let $T\in (0,\infty)$, $(t_n)_{n\in\{0,1,\dots,N\}}\allowbreak\subseteq\R$
	satisfy 
	for all $n\in\{0,1,\dots,N\}$ that $t_n=\tfrac{nT}{N}$,
	let $\mathfrak{D}\in[1,\infty)$, $\varepsilon\in (0,1]$, $q\in(2,\infty)$ satisfy that 
	\begin{equation}\label{ApproxOfEulerWithGronwall:paramBound}
%	\mathfrak{D}=2 \Big(\tfrac{360q}{(q-2)}\big[\LogBin(\eps^{-1})+\LogBin(2^{q-1}+1)\big]+108\Big),
\mathfrak{D}=\big[\tfrac{720q}{(q-2)}\big]\big[\LogBin(\eps^{-1})+q+1\big]-504,
	\end{equation}
%		let 
%		$
%		\Vert \cdot \Vert \colon \R^d \to [0,\infty)
%		$
%		be the Euclidean norm on $\R^d$,
		let 
		$\Phi\in \ANNs $
		satisfy 
		for all $x\in \R^d$ that $\inDimANN(\Phi)=\outDimANN(\Phi)=d$ and $\Vert (\functionANN(\Phi))(x)\Vert\le \mathfrak{C}\big(1+\Vert x\Vert\big)$,
	%	let $b_n\in \R^d$, $n\in\{0,1,\dots,N-1\}$, 
	%	and let 
	%	$Y^n= \big(Y^{n,x,y }_t\big)_{t\in [t_n,t_{n+1}], x\in\R^d, y\in (\R^d)^{n+1}} \colon\allowbreak [t_n,t_{n+1}]\times \R^d\times (\R^d)^{n+1} \to \R^d $, $n\in\{0,1,\dots,N-1\}$, 
	%	be the 
	%	functions 	
	%	which satisfy for all 
	%	$n\in\{0,1,\dots,N-1\}$,
	%		$ t \in [t_{n},t_{n+1}]$
	%	$ x \in \R^d $, $y=(y_0,y_1,\dots, y_n)\in (\R^d)^{n+1}$
	%	that $\affineProcess^{n,x,y }_0=x$ and
	%	\begin{equation}
	%	\label{ApproxOfEulerWithGronwall:Y_processes}
	%	\begin{split}
	%	&\affineProcess^{n,x,y }_t 
	%	=
	%	\affineProcess^{n,x,y }_{t_n}+ \left(\tfrac{tN}{T}-n\right)\left[\tfrac{T}{N}(\functionANN(\Phi)) \big( 
	%	\affineProcess^{n,x,y }_{t_n} 
	%	\big)
	%	+
	%	y_n\right]
	%	\end{split}
	%	\end{equation}
	 let 
	$Y= (Y^{x,y }_t)_{(t,x,y)\in [0,T]\times \R^d\times(\R^d)^N} \colon\allowbreak [0,T]\times \R^d\times (\R^d)^N \to \R^d $ 
	be the 
	function	
	which satisfies for all 
	$n\in\{0,1,\dots,N-1\}$,
	$ t \in [t_{n},t_{n+1}]$,
	$ x \in \R^d $, $y=(y_1,y_2,\dots, y_N)\in (\R^d)^N$
	that $\affineProcess^{x,y }_0=x$ and
	\begin{equation}
	\label{ApproxOfEulerWithGronwall:Y_processes}
	\begin{split}
	&\affineProcess^{x,y }_t 
	=
	\affineProcess^{x,y }_{t_n}+ \left(\tfrac{tN}{T}-n\right)\left[\tfrac{T}{N}(\functionANN(\Phi)) \big( 
	\affineProcess^{x,y }_{t_n} 
	\big)
	+
	y_{n+1}\right]\!,
	\end{split}
	\end{equation}
	and let $g_n\colon \R^d\times (\R^d)^N\to [0,\infty)$, $n\in\{0,1,\dots,N\}$, be the functions which satisfy for all $n\in\{0,1,\dots,N\}$,
	$ x \in \R^d $, $y=(y_1,y_2,\dots, y_N)\in (\R^d)^N$ that 
	\begin{equation}
		g_{n}(x,y)=\bigg(\Vert x\Vert + \mathfrak{C}t_n+
		\max_{m\in\{0,1,\dots,n\}}\Big\Vert \smallsum\limits_{k=1}^{m} y_{k}\Big\Vert\bigg)
		\exp(\mathfrak{C} t_n)
	\end{equation}
	(cf.\ Definition~\ref{Def:ANN}, Definition~\ref{Definition:ANNrealization}, and Definition~\ref{Def:euclideanNorm}).
	%	Then there exist continuous functions $f_{\varepsilon,x,n,t}\colon  (\R^d)^{n+1}\to \R^d$, $ t \in [t_{n},t_{n+1}]$, $n\in\{0,1,\dots, N-1\}$, $x\in\R^d$, $\varepsilon\in (0,1]$,  and artificial neural networks $\Psi_{\varepsilon,y}\in \ANNs$, $\varepsilon\in(0,1]$, $y\in (\R^d)^N$, which satisfy
	%	%	$(\Psi_{\varepsilon,y})_{\varepsilon\in(0,1]}\subseteq \ANNs$
	Then there exist $\Psi_{y}\in \ANNs$,  $y\in (\R^d)^N$, such that
	%	$(\Psi_{y})_{\varepsilon\in(0,1]}\subseteq \ANNs$
	\begin{enumerate}[(i)]
		\item \label{ApproxOfEulerWithGronwall:Function} it holds for all  $y\in (\R^d)^N$ that $\functionANN (\Psi_{y})\in C(\R^{d+1},\R^d)$,
		\item it holds for all
		$n\in\{0,1,\dots,N-1\}$,
		$ t \in [t_{n},t_{n+1}]$,
		$x\in\R^d$, 
		$y\in (\R^d)^N$ that 
		\begin{equation}
		\Vert \affineProcess^{x,y }_t -(\functionANN (\Psi_{y}))(t,x)\Vert
		\le\varepsilon \big(2\sqrt{d}+(g_{n}(x,y))^q+(g_{n+1}(x,y))^q\big),
		\end{equation}
		\item it holds for all 
		$n\in\{0,1,\dots,N-1\}$,
		$ t \in [t_{n},t_{n+1}]$,
		$x\in\R^d$, 
		$y\in (\R^d)^N$
		that 
		%		\begin{equation}
		%		\Vert (\functionANN (\Psi_{y}))(t,x)\Vert
		%		\le  12 \sqrt{d} \left(1+\Vert \affineProcess^{x,y }_{t_n}\Vert^2+\Vert \affineProcess^{x,y }_{t_{n+1}}\Vert^2\right),
		%		\end{equation}
		\begin{equation}
		\Vert (\functionANN (\Psi_{y}))(t,x)\Vert
		\le 6\sqrt{d}+2\big((g_{n}(x,y))^2+(g_{n+1}(x,y))^2\big),
		\end{equation}
		\item \label{ApproxOfEulerWithGronwall:ItemParams}
		it holds for all  $y\in (\R^d)^N$ that 
		%		\begin{multline}
		%		\paramANN(\Psi_{y})
		%		\\ \le\tfrac{1}{2} \bigg((2N+1)
		%		\Big[d^2 \mathfrak{D}_\varepsilon+\big(
		%		23
		%		+6N\hiddenLength(\Phi)
		%		+  4d^2+3d+N \big[2(d^2+d)+\paramANN(\Phi)\big]^{\!2}\big)^{\!2}\Big]
		%		\\-Nd(2d+1)+\hiddenLength(\Phi) d(2d+1)N(N+1)
		%		\bigg)^{\!2},
		%		\end{multline}
%		\begin{multline}
%		\paramANN(\Psi_{y})
%		\le\tfrac{1}{2} \bigg( \hiddenLength(\Phi) d(2d+1)N(N+1)-Nd(2d+1)
%		\\
%		+(2N+1)
%		\Big[d^2 \mathfrak{D}+\big(
%		23
%		+6N\hiddenLength(\Phi)
%		+  7d^2+N \big[4d^2+\paramANN(\Phi)\big]^{\!2}\big)^{\!2}\Big]
%		\bigg)^{\!2},
%		\end{multline}
			\begin{equation}
			\begin{split}
			\paramANN(\Psi_{y})
			&\le\tfrac{9}{2}\, N^6 d^{16} \Big[  2 \hiddenLength(\Phi)
			+
			 \mathfrak{D}+\big(
			30
			+6\hiddenLength(\Phi)
			+   \big[4+\paramANN(\Phi)\big]^{\!2}\big)^{\!2}
			\Big]^{\!2},
			\end{split}
			\end{equation}
		%		\item \label{ApproxOfEulerWithGronwall:ItemMeasurability}
		%		that for all $\varepsilon\in (0,1]$, $y=(y_1,y_2,\dots, y_N)\in (\R^d)^N$,  $n\in\{0,1,\dots, N-1\}$, $t\in [t_{n},t_{n+1}]$, $x\in\R^d$ it holds that 
		%		\begin{equation}
		%		(\functionANN (\Psi_{y}))(t,x)=f_{\varepsilon,x,n,t}(y_0,y_1,\dots, y_n).
		%		\end{equation}
		\item \label{ApproxOfEulerWithGronwall:Continuity}
		it holds  for all $t\in [0,T]$, $x\in\R^d$ that 
		\begin{equation}
		\big[(\R^d)^N\ni y\mapsto (\functionANN (\Psi_{y}))(t,x)\in\R^d\big]\in C\big((\R^d)^N,\R^d\big),
		\end{equation}
		and
		\item \label{ApproxOfEulerWithGronwall:Adaptedness}
		it holds  for all    $n\in\{0,1,\dots, N\}$, $t\in [0,t_n]$, $x\in\R^d$, $y=(y_1,y_2,\dots, y_N)$, $z=(z_1,z_2,\dots, z_N)\in (\R^d)^N$  with $\forall\, k\in \N\cap [0,n]\colon y_k=z_k$ that 
		\begin{equation}
		(\functionANN (\Psi_{y}))(t,x)=(\functionANN (\Psi_{z}))(t,x).
		\end{equation}
	\end{enumerate}
\end{theorem}

\begin{proof}[Proof of Theorem~\ref{Thm:ApproxOfEulerWithGronwall}]	
	Throughout this proof let 
	$\Psi_{y}\in \ANNs$,  $y\in (\R^d)^N$, satisfy  that
	%	$(\Psi_{y})_{\varepsilon\in(0,1]}\subseteq \ANNs$
	\begin{enumerate}[(I)]
		\item \label{ProofApproxOfEulerWithGronwall:Function} it holds for all  $y\in (\R^d)^N$  that $\functionANN (\Psi_{y})\in C(\R^{d+1},\R^d)$,
		\item\label{ProofApproxOfEulerWithGronwall:Error} it holds for all  
		$n\in\{0,1,\dots,N-1\}$,
		$ t \in [t_{n},t_{n+1}]$,
		$x\in\R^d$,
		$y\in (\R^d)^N$  that 
		\begin{equation}
		\Vert \affineProcess^{x,y }_t -(\functionANN (\Psi_{y}))(t,x)\Vert
%		\le 4\varepsilon \sqrt{d} \left(1+\Vert \affineProcess^{x,y }_{t_n}\Vert^q+\Vert \affineProcess^{x,y }_{t_{n+1}}\Vert^q\right),
		\le \varepsilon \big(2\sqrt{d}+\|  Y_{t_{n}}^{x,y }\|^q+\|  Y_{t_{n+1}}^{x,y }\|^q\big),
		\end{equation}
		\item \label{ProofApproxOfEulerWithGronwall:Growth} it holds for all 
		$n\in\{0,1,\dots,N-1\}$,
		$ t \in [t_{n},t_{n+1}]$,
		$x\in\R^d$,
		$y\in (\R^d)^N$  that 
		%		\begin{equation}
		%		\Vert (\functionANN (\Psi_{y}))(t,x)\Vert
		%		\le  12 \sqrt{d} \left(1+\Vert \affineProcess^{x,y }_{t_n}\Vert^2+\Vert \affineProcess^{x,y }_{t_{n+1}}\Vert^2\right),
		%		\end{equation}
		\begin{equation}
		\Vert (\functionANN (\Psi_{y}))(t,x)\Vert
%		\le  12 \sqrt{d} \left(1+\Vert \affineProcess^{x,y }_{t_n}\Vert^2+\Vert \affineProcess^{x,y }_{t_{n+1}}\Vert^2\right),
		\le 6\sqrt{d}+2\big(\Vert Y_{t_{n}}^{x,y }\Vert^2+\Vert Y_{t_{n+1}}^{x,y }\Vert^2\big),
		\end{equation}
		\item \label{ProofApproxOfEulerWithGronwall:ItemParams}
		it holds for all  $y\in (\R^d)^N$ that 
		%		\begin{multline}
		%		\paramANN(\Psi_{y})
		%		\\ \le\tfrac{1}{2} \bigg((2N+1)
		%		\Big[d^2 \mathfrak{D}+\big(
		%		23
		%		+6N\hiddenLength(\Phi)
		%		+  4d^2+3d+N \big[2(d^2+d)+\paramANN(\Phi)\big]^{\!2}\big)^{\!2}\Big]
		%		\\-Nd(2d+1)+\hiddenLength(\Phi) d(2d+1)N(N+1)
		%		\bigg)^{\!2},
		%		\end{multline}
						\begin{multline}
						\paramANN(\Psi_{y})
						\le\tfrac{1}{2} \bigg[  6d^2N^2 \hiddenLength(\Phi)
						\\
						+3N
						\Big[d^2 \mathfrak{D}+\big(
						23
						+6N\hiddenLength(\Phi)
						+  7d^2+N \big[4d^2+\paramANN(\Phi)\big]^{\!2}\big)^{\!2}\Big]
						\bigg]^{\!2},
						\end{multline}
		%		\item \label{ProofApproxOfEulerWithGronwall:ItemMeasurability}
		%		that for all $\varepsilon\in (0,1]$, $y=(y_1,y_2,\dots, y_N)\in (\R^d)^N$,  $n\in\{0,1,\dots, N-1\}$, $t\in [t_{n},t_{n+1}]$, $x\in\R^d$ it holds that 
		%		\begin{equation}
		%		(\functionANN (\Psi_{y}))(t,x)=f_{\varepsilon,x,n,t}(y_0,y_1,\dots, y_n).
		%		\end{equation}
		\item \label{ProofApproxOfEulerWithGronwall:Continuity}
		it holds  for all  $t\in [0,T]$, $x\in\R^d$ that 
		\begin{equation}
		\big[(\R^d)^N\ni y\mapsto (\functionANN (\Psi_{y}))(t,x)\in\R^d\big]\in C\big((\R^d)^N,\R^d\big),
		\end{equation}
		and
		\item \label{ProofApproxOfEulerWithGronwall:Adaptedness}
		it holds  for all    $n\in\{0,1,\dots, N\}$, $t\in [0,t_n]$, $x\in\R^d$, $y=(y_1,y_2,\dots, y_N),\allowbreak z=(z_1,z_2,\dots, z_N)\in (\R^d)^N$  with $\forall\, k\in \N\cap [0,n]\colon y_k=z_k$  that 
		\begin{equation}
		(\functionANN (\Psi_{y}))(t,x)=(\functionANN (\Psi_{z}))(t,x)
		\end{equation}
	\end{enumerate}
	(cf.\ Proposition~\ref{Cor:ApproxOfEuler}).
		Note that 	\eqref{ProofApproxOfEulerWithGronwall:ItemParams} ensures  for all  $y\in (\R^d)^N$  that 
		\begin{equation}
		\begin{split}
		&\paramANN(\Psi_{y})
		\\&\le\tfrac{1}{2} \bigg[  6d^2N^2 \hiddenLength(\Phi)
		+3N
		\Big[d^2 \mathfrak{D}+\big(
		23
		+6N\hiddenLength(\Phi)
		+  7d^2+N d^4 \big[4+\paramANN(\Phi)\big]^{\!2}\big)^{\!2}\Big]
		\bigg]^{\!2}
				\\&\le\tfrac{1}{2} \bigg[  6d^2N^2 \hiddenLength(\Phi)
				+3N
				\Big[d^2 \mathfrak{D}+N^2 d^8\big(
				30
				+6\hiddenLength(\Phi)
				+   \big[4+\paramANN(\Phi)\big]^{\!2}\big)^{\!2}\Big]
				\bigg]^{\!2}.
				\end{split}
		\end{equation}
	Hence, we obtain that   for all  $y\in (\R^d)^N$ it holds that 
	\begin{equation}\label{ProofApproxOfEulerWithGronwall:EstimateParams}
	\begin{split}
	\paramANN(\Psi_{y})
	&\le\tfrac{1}{2} \bigg[  6d^2N^2 \hiddenLength(\Phi)
	+3N^3 d^8
	\Big[ \mathfrak{D}+\big(
	30
	+6\hiddenLength(\Phi)
	+   \big[4+\paramANN(\Phi)\big]^{\!2}\big)^{\!2}\Big]
	\bigg]^{\!2}
		\\&\le\tfrac{9}{2}\, N^6 d^{16} \Big[  2 \hiddenLength(\Phi)
		+
		 \mathfrak{D}+\big(
		30
		+6\hiddenLength(\Phi)
		+   \big[4+\paramANN(\Phi)\big]^{\!2}\big)^{\!2}
		\Big]^{\!2}.
	\end{split}
	\end{equation}
	In addition, observe that Lemma~\ref{Lemma:Gronwall} and the hypothesis that for all $n\in\{0,1,\dots,N\}$ it holds that $t_n=\tfrac{nT}{N}$ demonstrate that 	 for all $n\in\{0,1,\dots,N\}$,
	$ x \in \R^d $, $y=(y_1,y_2,\dots, y_N)\in (\R^d)^N$ it holds that 
	\begin{equation}\label{ApproxOfEuler:GronwallEstimate}
	\begin{split}
	\Vert \affineProcess^{x,y} _{t_n} \Vert
	&\le \bigg[\Vert x\Vert + \tfrac{\mathfrak{C}nT}{N}+
	\max_{m\in\{0,1,\dots,n\}}\Big\Vert \smallsum\limits_{k=1}^{m} y_{k}\Big\Vert\bigg]
	\exp\!\left(\tfrac{\mathfrak{C}nT}{N}\right)
	\\&=\bigg[\Vert x\Vert + \mathfrak{C}t_n+
	\max_{m\in\{0,1,\dots,n\}}\Big\Vert \smallsum\limits_{k=1}^{m} y_{k}\Big\Vert\bigg]
	\exp(\mathfrak{C} t_n)
	=g_n(x,y).
	\end{split}
	\end{equation}
	Combining this with  \eqref{ProofApproxOfEulerWithGronwall:Error} and \eqref{ProofApproxOfEulerWithGronwall:Growth} ensures that for all
	$n\in\{0,1,\dots,N-1\}$,
	$ t \in [t_{n},t_{n+1}]$,
	$x\in\R^d$, $y\in (\R^d)^N$ it holds that 
	\begin{equation}
	\begin{split}
		\Vert \affineProcess^{x,y }_t -(\functionANN (\Psi_{y}))(t,x)\Vert
		&\le\varepsilon \big(2\sqrt{d}+\|  Y_{t_{n}}^{x,y }\|^q+\|  Y_{t_{n+1}}^{x,y }\|^q\big)
				\\&\le\varepsilon \big(2\sqrt{d}+(g_{n}(x,y))^q+(g_{n+1}(x,y))^q\big)
%		\\&\le 4\varepsilon \sqrt{d} \left(1+\Vert \affineProcess^{x,y }_{t_n}\Vert^q+\Vert \affineProcess^{x,y }_{t_{n+1}}\Vert^q\right)
%		\\&\le 4\varepsilon \sqrt{d} \Big[1+\big(g_{n}(x,y)\big)^q+\big(g_{n+1}(x,y)\big)^q\Big]
	\end{split}
	\end{equation}
	and
			\begin{equation}
			\begin{split}
						\Vert (\functionANN (\Psi_{y}))(t,x)\Vert
						&\le 6\sqrt{d}+2\big(\Vert Y_{t_{n}}^{x,y }\Vert^2+\Vert Y_{t_{n+1}}^{x,y }\Vert^2\big)
						\\&\le 6\sqrt{d}+2\big((g_{n}(x,y))^2+(g_{n+1}(x,y))^2\big).
%						\\&\le  12 \sqrt{d} \left(1+\big(g_{n}(x,y)\big)^2+\big(g_{n+1}(x,y)\big)^2\right)
%						\\&\le  12 \sqrt{d} \Big[1+\big(g_{n}(x,y)\big)^2+\big(g_{n+1}(x,y)\big)^2\Big].
			\end{split}
			\end{equation}
	This, \eqref{ProofApproxOfEulerWithGronwall:Function},  \eqref{ProofApproxOfEulerWithGronwall:Continuity}, \eqref{ProofApproxOfEulerWithGronwall:Adaptedness}, and \eqref{ProofApproxOfEulerWithGronwall:EstimateParams} establish items~\eqref{ApproxOfEulerWithGronwall:Function}--\eqref{ApproxOfEulerWithGronwall:Adaptedness}.
	The proof of Theorem~\ref{Thm:ApproxOfEulerWithGronwall} is thus completed.
\end{proof}

\begin{cor}\label{Cor:ApproxOfEulerWithGronwall}
	Let  $\mathfrak{C}, T,\mathfrak{d} \in (0,\infty)$, $a\in C(\R,\R)$ satisfy for all  $x\in \R$ that $\activation(x)=\max\{x,0\}$, 
	let $\Phi_d\in \ANNs$, $d\in\N$, satisfy for all $d\in\N$, $x\in \R^d$ that $\inDimANN(\Phi_d)=\outDimANN(\Phi_d)=d$, $\Vert (\functionANN(\Phi_d))(x)\Vert\le \mathfrak{C}(1+\Vert x\Vert)$, and 
	$\paramANN(\Phi_d)\le \mathfrak{C}d^{\mathfrak{d}}$,
	let $Y^{d,N}= (Y^{d,N}_{t,x,y})_{(t,x,y)\in [0,T]\times \R^d\times(\R^d)^N} \colon\allowbreak [0,T]\times \R^d\times (\R^d)^N \to \R^d $, $N, d\in \N$, 
	be the functions which satisfy for all $d, N\in\N$, $n\in\{0,1,\dots,N-1\}$, $ t \in \big[\frac{nT}{N},\frac{(n+1)T}{N}\big]$, $ x \in \R^d $, $y=(y_1,y_2,\dots, y_N)\in (\R^d)^N$ that $\affineProcess^{d,N }_{0,x,y}=x$ and
	\begin{equation}
	\label{ApproxOfEulerWithGronwall:Y_processes}
	\begin{split}
	&\affineProcess^{d,N}_{t,x,y} 
	=
	\affineProcess^{d,N }_{\frac{nT}{N},x,y}+ \left(\tfrac{tN}{T}-n\right)\big[\tfrac{T}{N}(\functionANN(\Phi_d)) ( 
	\affineProcess^{d,N }_{\frac{nT}{N},x,y} )
	+
	y_{n+1}\big]\!.
	\end{split}
	\end{equation}
	(cf.\ Definition~\ref{Def:ANN}, Definition~\ref{Definition:ANNrealization}, and Definition~\ref{Def:euclideanNorm}).
	Then there exist $C\in\R$ and $\Psi_{\varepsilon,d,N,y}\in \ANNs$, $y\in (\R^d)^N$, $N, d\in \N$, $\varepsilon\in (0,1]$, such that 
	\begin{enumerate}[(i)]
		\item \label{cor1} it holds for all  $\varepsilon\in (0,1]$, $d,N\in\N$, $y\in (\R^d)^N$ that $\functionANN (\Psi_{\varepsilon,d,N,y})\in C(\R^{d+1},\R^d)$,
		\item  it holds for all
		$\varepsilon\in (0,1]$, $d,N\in\N$,
		$ t \in [0,T]$,
		$x\in\R^d$, 
		$y\in (\R^d)^N$ that 
		\begin{equation}
			\Vert \affineProcess^{d,N}_{t,x,y} -(\functionANN (\Psi_{\varepsilon,d,N,y}))(t,x)\Vert\le Cd^{\nicefrac{1}{2}} N^{\nicefrac{3}{2}}\varepsilon(1+\Vert x\Vert^3+\Vert y\Vert^3),
		\end{equation}
		\item it holds for all 
		$\varepsilon\in (0,1]$, $d,N\in\N$,
		$ t \in [0,T]$,
		$x\in\R^d$, 
		$y\in (\R^d)^N$
		that 
		\begin{equation}
			\Vert (\functionANN (\Psi_{\varepsilon,d,N,y}))(t,x)\Vert \le Cd^{\nicefrac{1}{2}}N(1+\Vert x\Vert^2+\Vert y\Vert^2),
		\end{equation}
		\item it holds for all  $\varepsilon\in (0,1]$, $d,N\in\N$, $y\in (\R^d)^N$ that 
			\begin{equation}
				\paramANN(\Psi_{\varepsilon,d,N,y})\le Cd^{16+8\mathfrak{d}}N^6\big[1+|\!\ln(\varepsilon)|^2\big],
			\end{equation}
		\item it holds  for all $\varepsilon\in (0,1]$, $d,N\in\N$, $t\in [0,T]$, $x\in\R^d$ that 
		\begin{equation}
		\big[(\R^d)^N\ni y\mapsto (\functionANN (\Psi_{\varepsilon,d,N,y}))(t,x)\in\R^d\big]\in C\big((\R^d)^N,\R^d\big),
		\end{equation}
		and
		\item \label{cor6}
		it holds for all $\varepsilon\in (0,1]$, $d,N\in\N$, $n\in\{0,1,\dots, N\}$, $t\in [0,\frac{nT}{N}]$, $x\in\R^d$, $y=(y_1,y_2,\dots, y_N)$, $z=(z_1,z_2,\dots, z_N)\in (\R^d)^N$  with $\forall\, k\in \N\cap [0,n]\colon y_k=z_k$ that 
		\begin{equation}
		(\functionANN (\Psi_{\varepsilon,d,N,y}))(t,x)=(\functionANN (\Psi_{\varepsilon,d,N,z}))(t,x).
		\end{equation}

	\end{enumerate}
\end{cor}

\begin{proof}[Proof of Corollary~\ref{Cor:ApproxOfEulerWithGronwall}]
	Throughout this proof let $\mathfrak{D}_{\varepsilon,q}\in [1,\infty)$, $q\in (2,\infty)$, $\varepsilon\in(0,1]$, satisfy for all $\varepsilon\in (0,1]$, $q\in(2,\infty)$ that
	\begin{equation}
	\label{eq: ConstD}
		\mathfrak{D}_{\varepsilon,q}=\left[\frac{720 q}{(q-2)}\right]\left[\log_2(\varepsilon^{-1})+q+1\right]-504,
	\end{equation}
	 let $c=\max\{\exp(\mathfrak{C}T),\mathfrak{D}_{1,3}, 62+6\mathfrak{C}(\mathfrak{C}+1)\}$, and let $g^{d,N}_n \colon\R^d\times (\R^d)^N\to [0,\infty)$, $n\in\{0,1,\dots,N\}$, $N, d\in\N$, be the functions which satisfy for all $d, N\in\N$, $n\in\{0,1,\dots,N\}$, $x\in \R^d$, 
	$y=(y_1,y_2,\dots,y_N)\in(\R^d)^N$ that 
	\begin{equation}
	\label{eq: familiyG}
		g^{d,N}_n(x,y)=\bigg(\Vert x\Vert + \frac{\mathfrak{C} nT}{N}+
		\max_{m\in\{0,1,\dots,n\}}\Big\Vert \smallsum\limits_{k=1}^{m} y_{k}\Big\Vert\bigg)
		\exp(\frac{\mathfrak{C}nT}{N}).
	\end{equation}
	Note that Theorem~\ref{Thm:ApproxOfEulerWithGronwall} (with $N=N$, $d=d$, $\mathfrak{C}=\mathfrak{C}$, $a=a$, $T=T$, $t_n=\frac{nT}{N}$, $\mathfrak{D}=\mathfrak{D}_{\varepsilon,3}$, 
	$\varepsilon=\varepsilon$, $q=3$, $\Phi=\Phi_d$,
	$Y=Y^{d,N}$, $g_n=g^{d,N}_n$ for $N,d\in\N$, $n\in \{0,1,\dots,N\}$, $\varepsilon\in (0,1]$ in the notation of Theorem~\ref{Thm:ApproxOfEulerWithGronwall})
	 implies that there exist 
	$\Psi_{\varepsilon,d,N,y}\in \ANNs$, $y\in (\R^d)^N$, $N, d\in\N$, $\varepsilon\in(0,1]$, which satisfy that
	\begin{enumerate}[(I)]
		\item \label{CorApproxOfEulerWithGronwall:Function} it holds for all  $\varepsilon\in (0,1]$, $d,N\in\N$, $y\in (\R^d)^N$ that $\functionANN (\Psi_{\varepsilon,d,N,y})\in C(\R^{d+1},\R^d)$,
		\item \label{CorApproxOfEulerWithGronwall:SecondClaim} it holds for all
		$\varepsilon\in (0,1]$, $d,N\in\N$, $n\in\{0,1,\dots,N-1\}$,
		$ t \in [\frac{nT}{N},\frac{(n+1)T}{N}]$,
		$x\in\R^d$, 
		$y\in (\R^d)^N$ that 
		\begin{equation}
		\Vert \affineProcess^{d,N}_{t,x,y} -(\functionANN (\Psi_{\varepsilon,d,N,y}))(t,x)\Vert
		\le\varepsilon \big(2\sqrt{d}+(g^{d,N}_n(x,y))^3+(g^{d,N}_{n+1}(x,y))^3\big),
		\end{equation}
		\item \label{CorApproxOfEulerWithGronwall:ThirdClaim} it holds for all 
		$\varepsilon\in (0,1]$, $d,N\in\N$, $n\in\{0,1,\dots,N-1\}$,
		$ t \in [\frac{nT}{N},\frac{(n+1)T}{N}]$,
		$x\in\R^d$, 
		$y\in (\R^d)^N$
		that 
		\begin{equation}
		\Vert (\functionANN (\Psi_{\varepsilon,d,N,y}))(t,x)\Vert
		\le 6\sqrt{d}+2\big((g^{d,N}_n(x,y))^2+(g^{d,N}_{n+1}(x,y))^2\big),
		\end{equation}
		\item \label{CorApproxOfEulerWithGronwall:ItemParams}
		it holds for all  $\varepsilon\in (0,1]$, $d,N\in\N$, $y\in (\R^d)^N$ that 
			\begin{equation}
			\begin{split}
			&\paramANN(\Psi_{\varepsilon,d,N,y})\\
			&\le\tfrac{9}{2}\, N^6 d^{16} \Big[  2 \hiddenLength(\Phi_d)
			+
			 \mathfrak{D}_{\varepsilon,3}+\big(
			30
			+6\hiddenLength(\Phi_d)
			+   \big[4+\paramANN(\Phi_d)\big]^{\!2}\big)^{\!2}
			\Big]^{\!2},
			\end{split}
			\end{equation}
		\item \label{CorApproxOfEulerWithGronwall:Continuity}
		it holds  for all $\varepsilon\in (0,1]$, $d,N\in\N$, $t\in [0,T]$, $x\in\R^d$ that 
		\begin{equation}
		\big[(\R^d)^N\ni y\mapsto (\functionANN (\Psi_{\varepsilon,d,N,y}))(t,x)\in\R^d\big]\in C\big((\R^d)^N,\R^d\big),
		\end{equation}
		and
		\item \label{CorApproxOfEulerWithGronwall:Adaptedness}
		it holds for all $\varepsilon\in (0,1]$, $d,N\in\N$, $n\in\{0,1,\dots, N\}$, $t\in [0,\frac{nT}{N}]$, $x\in\R^d$, $y=(y_1,y_2,\dots, y_N)$, $z=(z_1,z_2,\dots, z_N)\in (\R^d)^N$  with $\forall\, k\in \N\cap [0,n]\colon y_k=z_k$ it holds that 
		\begin{equation}
		(\functionANN (\Psi_{\varepsilon,d,N,y}))(t,x)=(\functionANN (\Psi_{\varepsilon,d,N,z}))(t,x).
		\end{equation}
	\end{enumerate}
	Observe that Jensen's inequality implies that for all $n\in\N$, $p\in [1,\infty)$, $(x_1,x_2,\dots,x_n)\in\R^n$ it holds that
	\begin{equation}
	\label{eq:discreteJensen}
		|x_1+\dots+x_n|^p\leq n^{p-1}(|x_1|^p+\dots+|x_n|^p).
	\end{equation}
	Moreover, note that H\"older's inequality shows that for all $N\in\N$, $y=(y_1,y_2,\dots,y_N)\in (\R^d)^N$ it holds that
	\begin{equation}
	\label{eq:YoungInequality}
		\smallsum\limits_{k=1}^N\Vert y_k\Vert =\smallsum\limits_{k=1}^N(1\Vert y_k\Vert)\le N^{\nicefrac{1}{2}}\Big(\smallsum\limits_{k=1}^N\Vert y_k\Vert^2\Big)^{\nicefrac{1}{2}}=N^{\nicefrac{1}{2}}\Vert y\Vert.
	\end{equation}
	Combining \eqref{eq:discreteJensen}, \eqref{CorApproxOfEulerWithGronwall:SecondClaim}, and \eqref{eq: familiyG} therefore ensures that for all $\varepsilon\in (0,1]$, $d,N\in\N$, $n\in\{0,1,\dots,N-1\}$, 
	$ t \in [\frac{nT}{N},\frac{(n+1)T}{N}]$, $x\in\R^d$, $y=(y_1,y_2,\dots,y_N)\in (\R^d)^N$ it holds that 
		\begin{equation}
		\label{eq:estimate 411}
		\begin{split}
			&\Vert \affineProcess^{d,N}_{t,x,y} -(\functionANN (\Psi_{\varepsilon,d,N,y}))(t,x)\Vert\le 2d^{\nicefrac{1}{2}}\varepsilon (1+(g^{d,N}_N(x,y))^3)\\
			&= 2d^{\nicefrac{1}{2}}\varepsilon \bigg(1+\bigg(\Vert x\Vert + \mathfrak{C}T+\max_{m\in\{0,1,\dots,N\}}\Big\Vert \smallsum\limits_{k=1}^{m} y_{k}\Big\Vert\bigg)^3\exp(3\mathfrak{C}T)\bigg)\\
			&\le 2d^{\nicefrac{1}{2}}\varepsilon \bigg(1+9\bigg(\Vert x\Vert^3+c^3+\Big(\smallsum\limits_{k=1}^{N} \Vert y_{k}\Vert\Big)^3\bigg)c^3\bigg)\\
			&\le 2d^{\nicefrac{1}{2}}\varepsilon \big(1+9\big(\Vert x\Vert^3+c^3+N^{\nicefrac{3}{2}}\Vert y\Vert ^3\big)c^3\big)\\
			&\le 2c^6d^{\nicefrac{1}{2}}N^{\nicefrac{3}{2}}\varepsilon \big(1+9\big(\Vert x\Vert^3+1+\Vert y\Vert ^3\big)\big)\\
			&= 2c^6d^{\nicefrac{1}{2}}N^{\nicefrac{3}{2}}\varepsilon \big(10+9\Vert x\Vert^3+9\Vert y\Vert ^3\big)\\
			&\le 20c^6d^{\nicefrac{1}{2}}  N^{\nicefrac{3}{2}} \varepsilon(1+\Vert x\Vert^3 +\Vert y\Vert^3).
		\end{split}
		\end{equation}
%	Note that \eqref{CorApproxOfEulerWithGronwall:SecondClaim}, \eqref{eq: familiyG}, and midpoint convexity of the mapping $\R\ni t\mapsto |t|^{\alpha}\in\R$, $\alpha\in[1,\infty)$, ensure for all $d,N\in\N$, $n\in\{0,1,\dots,N-1\}$, 
%	$q\in (0,\infty)$, $ t \in [\frac{nT}{N},\frac{(n+1)T}{N}]$, $x\in\R^d$, $y\in (\R^d)^N$ that 
%		\begin{equation}
%		\label{eq:estimate 411}
%		\begin{split}
%			&\Vert \affineProcess^{d,N}_{t,x,y} -(\functionANN (\Psi_{\varepsilon,d,N,y}))(t,x)\Vert\le 2\sqrt{d}\varepsilon (1+(g^{d,N}_N(x,y))^q)\\
%			&\le 2\sqrt{d}\varepsilon \bigg(1+\bigg(\Vert x\Vert + \mathfrak{C}T+\max_{m\in\{0,1,\dots,N\}}\Big\Vert \smallsum\limits_{k=1}^{m} y_{k}\Big\Vert\bigg)^q\exp(q\mathfrak{C}T)\bigg)\\
%			&\le 2\sqrt{d}\varepsilon \bigg(1+2^{q-1}\bigg(\Vert x\Vert^q+2^{q-1}\bigg(C^q+\max_{m\in\{0,1,\dots,N\}}\Big\Vert \smallsum\limits_{k=1}^{m} y_{k}\Big\Vert^q\bigg)\bigg)\exp(qC)\bigg)\\
%			&\le C\sqrt{d}\varepsilon4^q\bigg(2+\Vert x\Vert^q+\max_{m\in\{0,1,\dots,N\}}\max_{k\in \{0,1,\dots,m\}}(m^q\Vert y_k\Vert^q)\bigg)\\
%			&\le C\sqrt{d}\varepsilon 4^q(2+\Vert x\Vert^q +N^q\Vert y\Vert^q)\le C\sqrt{d}\varepsilon 4^q N^q (1+\Vert x\Vert^q +\Vert y\Vert^q).
%		\end{split}
%		\end{equation}
%	Hence, with $q=3$ we obtain that for all $d,N\in\N$, $n\in\{0,1,\dots,N-1\}$, $ t \in [\frac{nT}{N},\frac{(n+1)T}{N}]$, $x\in\R^d$, $y\in (\R^d)^N$ that 
%	\begin{equation}
%		\Vert \affineProcess^{d,N}_{t,x,y} -(\functionANN (\Psi_{\varepsilon,d,N,y}))(t,x)\Vert\le C\sqrt{d}\varepsilon N^3(1+\Vert x\Vert^3+\Vert y\Vert^3).
%	\end{equation}
	Next note that \eqref{CorApproxOfEulerWithGronwall:ThirdClaim}, \eqref{eq: familiyG}, \eqref{eq:discreteJensen}, and \eqref{eq:YoungInequality} imply that for all $\varepsilon\in (0,1]$, $d,N\in\N$, $n\in\{0,1,\dots,N-1\}$, 
	$ t \in [\frac{nT}{N},\frac{(n+1)T}{N}]$, $x\in\R^d$, $y=(y_1,y_2,\dots,y_N)\in (\R^d)^N$ it holds that
	\begin{equation}
	\label{eq:estimate 412}
	\begin{split}
		&\Vert (\functionANN (\Psi_{\varepsilon,d,N,y}))(t,x)\Vert\le 6\sqrt{d}+4(g^{d,N}_N(x,y))^2\\
		&= 6\sqrt{d}+4\bigg(\Vert x\Vert+\mathfrak{C}T+\max_{m\in\{0,1,\dots,N\}}\Big(\Vert \smallsum\limits_{k=1}^{m} y_{k}\Big\Vert\bigg)^2\exp(2\mathfrak{C}T)\\
		&\le 6\sqrt{d}+12\bigg(\Vert x\Vert^2+c^2+\Big( \smallsum\limits_{k=1}^{N} \Vert y_{k}\Vert\Big)^2\bigg)c^2\\
		&\le 6\sqrt{d}+12\big(\Vert x\Vert^2+c^2+N\Vert y\Vert^2\big)c^2\\
		&\le 18c^4\sqrt{d}N(1+\Vert x\Vert^2+\Vert y\Vert^2).
	\end{split}
	\end{equation}
%	Finally, observe that \eqref{CorApproxOfEulerWithGronwall:Adaptedness}, the hypothesis that all $\Phi_d\in\ANNs$, $d\in\N$, satisfy $\paramANN(\Phi_d)\le \mathfrak{C}d^{\mathfrak{d}}$, and Young's inequality
%	assure that for all $\varepsilon\in (0,1]$, $d,N\in\N$, $y\in (\R^d)^N$ it holds that
%	\begin{equation}
%	\begin{split}
%		&\paramANN(\Psi_{\varepsilon,d,N,y})\le \tfrac{9}{2}N^6d^{16}\Big[2\paramANN(\Phi_d)+\mathfrak{D}_{\varepsilon,q}+\big(30+6\paramANN(\Phi_d)+[4+\paramANN(\Phi_d)]^2\big)^2\Big]^2\\
%		&\le \tfrac{9}{2}N^6d^{16}\Big[2\paramANN(\Phi_d)+\mathfrak{D}_{\varepsilon,q}+\big(46+6\paramANN(\Phi_d)+\paramANN(\Phi_d)^2\big)^2\Big]^2\\
%		&\le \tfrac{9}{2}N^6d^{16}\Big[16(\paramANN(\Phi_d))^4+72(\paramANN(\Phi_d))^2+2\paramANN(\Phi_d)+\mathfrak{D}_{\varepsilon,q}+7688\Big]^2\\
%		&\le C\tfrac{N^6d^{16}q}{(q-2)}\Big[(\paramANN(\Phi_d))^8+(\paramANN(\Phi_d))^4+(\paramANN(\Phi_d))^2+\big(\log_2(\varepsilon^{-1})+q+1\big)^2+1\Big]\\
%		&\le C\tfrac{N^6d^{16}q}{(q-2)}\Big[(\paramANN(\Phi_d))^8+(\paramANN(\Phi_d))^4+(\paramANN(\Phi_d))^2+\Big(\tfrac{9}{2\varepsilon^2}+\tfrac{(2q-3)^2}{2}\Big)+1\Big]\\
%		&\le C\tfrac{N^6d^{16+8\mathfrak{d}}\varepsilon^{-2}q(2q-3)^2}{(q-2)}.
%	\end{split}
%	\end{equation}
%	Hence, with $q=3$ we obtain that for all $\varepsilon\in (0,1]$, $d,N\in\N$, $y\in (\R^d)^N$ it holds that 
%	\begin{equation}
%		\paramANN(\Psi_{\varepsilon,d,N,y})\le CN^6d^{16+8\mathfrak{d})}\varepsilon^{-2}.
%	\end{equation}
	Furthermore, observe that \eqref{eq: ConstD} shows that for all $\varepsilon\in (0,1]$ it holds that
	\begin{equation}
	\begin{split}
		(\mathfrak{D}_{\varepsilon,3})^2&=\big(\tfrac{2160}{\ln(2)}\!\ln(\varepsilon^{-1})+\mathfrak{D}_{1,3}\big)^2\le (\mathfrak{D}_{1,3})^2(\ln(\varepsilon^{-1})+1)^2\\
		&\le c^2|1-\!\ln(\varepsilon)|^2\le 2c^2(1+|\!\ln(\varepsilon)|^2).
	\end{split}
	\end{equation}
	This, \eqref{CorApproxOfEulerWithGronwall:ItemParams}, the hypothesis that for all $d\in\N$ it holds that $\paramANN(\Phi_d)\le \mathfrak{C}d^{\mathfrak{d}}$, and \eqref{eq:discreteJensen}
	assure that for all $\varepsilon\in (0,1]$, $d,N\in\N$, $y\in (\R^d)^N$ it holds that
	\begin{equation}
	\label{eq:estimate 413}
	\begin{split}
		&\paramANN(\Psi_{\varepsilon,d,N,y})\le \tfrac{9}{2}N^6d^{16}\Big[2\mathfrak{C}d^{\mathfrak{d}}+\mathfrak{D}_{\varepsilon,3}+\big(30+6\mathfrak{C}d^{\mathfrak{d}}+[4+\mathfrak{C}d^{\mathfrak{d}}]^2\big)^2\Big]^2\\
		&\le \tfrac{27}{2}N^6d^{16}\Big[4\mathfrak{C}^2d^{2\mathfrak{d}}+(\mathfrak{D}_{\varepsilon,3})^2+\big(30+6\mathfrak{C}d^{\mathfrak{d}}+2\big[16+\mathfrak{C}^2d^{2\mathfrak{d}}\big]\big)^4\Big]\\
		&\le \tfrac{27}{2}N^6d^{16}\Big[4\mathfrak{C}^2d^{2\mathfrak{d}}+2c^2\big(1+|\!\ln(\varepsilon)|^2\big)+\big(62+6\mathfrak{C}d^{\mathfrak{d}}+2\mathfrak{C}^2d^{2\mathfrak{d}}\big)^4\Big]\\
		&\le \tfrac{27}{2}N^6d^{16}\Big[4\mathfrak{C}^2d^{2\mathfrak{d}}+2c^2\big(1+|\!\ln(\varepsilon)|^2\big)+\big(62+6\mathfrak{C}\big(\mathfrak{C}+1\big)d^{2\mathfrak{d}}\big)^4\Big]\\
		&\le \tfrac{27}{2}N^6d^{16}\Big[c d^{2\mathfrak{d}}+2c^2\big(1+|\!\ln(\varepsilon)|^2\big)+\big(c d^{2\mathfrak{d}}\big)^4\Big]\\
		&\le 27N^6d^{16}\Big[c^2\big(1+|\!\ln(\varepsilon)|^2\big)+\big(c d^{2\mathfrak{d}}\big)^4\Big]\\
		&\le 54c^4N^6d^{16+8\mathfrak{d}}\Big[1+|\!\ln(\varepsilon)|^2\Big].\\
	\end{split}
	\end{equation}
	Combining \eqref{CorApproxOfEulerWithGronwall:Function}, \eqref{eq:estimate 411}, \eqref{eq:estimate 412}, \eqref{eq:estimate 413}, \eqref{CorApproxOfEulerWithGronwall:Continuity}, and \eqref{CorApproxOfEulerWithGronwall:Adaptedness} 
	establishes items \eqref{cor1}-\eqref{cor6}. The proof of Corollary~\ref{Cor:ApproxOfEulerWithGronwall} is thus completed.
\end{proof}